\documentclass{article}
\usepackage{graphicx,color}
\usepackage{amsmath,amsthm,amsfonts,amssymb}
\newtheorem{theorem}{\indent  Theorem}
\newtheorem{df}{\indent  Definition}
\newtheorem{prop}{\indent Proposition}

\newtheorem{rk}{\indent Remark}
\newtheorem{lm}{\indent Lemma}
\newtheorem{notation}[lm]{Notation}

\newcommand{\bp}{\begin{proof}}
\newcommand{\ep}{\end{proof}}

\newtheorem{cy}{\indent Corollary}
\newtheorem{Claim}{\indent Claim}
\newcommand{\beq}{\begin{equation}}
\newcommand{\eeq}{\end{equation}}
\newcommand{\be}{\begin{enumerate}}
\newcommand{\ee}{\end{enumerate}}
\newcommand{\bea}{\begin{eqnarray*}}
\newcommand{\eea}{\end{eqnarray*}}
\newcommand{\bi}{\begin{itemize}}
\newcommand{\ei}{\end{itemize}}
\usepackage{amssymb}
\usepackage{psfig}
\usepackage{amsmath}
\usepackage{amssymb}

\edef\storedcatcodeat{\the\catcode`\@} \catcode`\@=11
\def\ig#1#2#3{{\setlength\arraycolsep{0pt}\def\arraystretch{0.8}
\setbox250\vbox{\strut}
\begin{array}[t]{|lcr|}
\multicolumn3{|c|}{\: #1 \:} \\ \hline \scriptstyle \: #2 &
\scriptstyle\hspace{1em} & \scriptstyle #3 \:
\end{array}}}
\def\igl#1#2{{\setlength\arraycolsep{0pt}\def\arraystretch{0.8}
\setbox250\vbox{\strut}
\begin{array}[t]{|l}
\multicolumn1{|r}{\: #1} \\[-\ht250]
\rlap{\vbox{\hrule width0.5em}} \strut \\[-\dp250]
\scriptstyle \: #2 \end{array}}}
\def\igr#1#2{{\setlength\arraycolsep{0pt}\def\arraystretch{0.8}
\setbox250\vbox{\strut}
\begin{array}[t]{r|}
\multicolumn1{l|}{\: #1} \\[-\ht250]
\strut \llap{\vbox{\hrule width0.5em}} \\[-\dp250]
\scriptstyle #2 \: \end{array}}}
\catcode`\@=\storedcatcodeat\relax

\begin{document}

\begin{center}{\bf\Huge { Implicit function theorem over free
groups }}\end{center}

\bigskip
\begin{center}{\large {Olga  Kharlampovich }}\end{center}

\medskip
\begin{center}{Department of Mathematics and
Statistics \\ McGill University\\ Montreal, QC H3A 2K6, Canada \\
e- mail: olga@math.mcgill.ca} \end{center}

\bigskip
\begin{center}{\large { Alexei Myasnikov}}\end{center}

\medskip
\begin{center}{Department of Mathematics\\
City College of New York\\ New York, NY 10031, USA  \\ e-mail:
alexeim@att.net}\end{center}

\pagestyle{myheadings}
\begin{abstract}
 We introduce the notion of a
regular quadratic equation and a regular NTQ system over a free
group.  We prove the results that can be described as Implicit
function theorems for algebraic varieties corresponding to regular
quadratic and NTQ systems.  We will also show that the Implicit
function theorem is true only for these varieties. In algebraic
geometry such results would be described as lifting solutions of
equations into generic points. From the model theoretic view-point
we claim the existence of simple Skolem functions for particular
$\forall\exists$-formulas over free groups.  Proving these
theorems we describe in details a new version of the
Makanin-Razborov process for solving equations in free groups. We
also prove a weak version of the Implicit function theorem for NTQ
systems which is one of the key results in the solution of the
Tarski's problems about the elementary theory of a free group. We
call it the Parametrization theorem.
\end{abstract}
\vspace{5in}

\tableofcontents
\bigskip
\section{Introduction} \label{se:intro}

In this paper we prove  so-called Implicit function theorems for
regular quadratic and NTQ systems over free groups (Theorems
\ref{IMTQ}, \ref{1,2,3,4}, \ref{reg}). They can be viewed  as
analogs of the corresponding result from analysis, hence the name.
To show this we formulate a very basic version of the Implicit
function theorem.

Let
$$S(x_1,\ldots ,x_n,a_1,\ldots ,a_k)=1$$
be  a ``regular'' quadratic equation in variables $X = (x_1,\ldots
,x_n)$ with constants
 $a_1,\ldots,a_k$ in a free group $F$ (roughly speaking ``regular''
means that the radical of $S$  coincides with the normal closure
of $S$ and $S$ is not an equation of one of few very specific
types). Suppose now that for each solution of the equation $S(X) =
1$ some other equation
$$T(x_1,\ldots,x_n,y_1, \ldots y_m,a_1,\ldots ,a_k)=1$$
 has a solution in $F$, then $T(X,Y)=1$ has a solution
$Y = (y_1,\ldots ,y_m)$ in the coordinate group $G_{R(S)}$ of the
equation $S(X)= 1$.

 This implies, that locally (in terms of Zariski topology), i.e., in the
neighborhood defined by  the equation $S(X) = 1$ , the implicit
functions $y_1,\ldots ,y_m$ can be expressed as  explicit words in
variables $x_1, \ldots, x_n$ and constants from $F$, say $Y =
P(X)$.  This result  allows one to eliminate a quantifier from the
following formula

$$ \Phi = \forall X \exists Y (S(X) = 1 \  \rightarrow  \ T(X,Y) = 1).$$
Indeed, the sentence $\Phi$ is equivalent in $F$ to the following
one:

$$
\Psi = \forall X (S(X) = 1 \  \rightarrow  \ T(X, P(X)) = 1).
$$

   From model theoretic view-point the
theorems claim existence of very simple Skolem functions for
particular $\forall\exists$-formulas over free groups. While in
algebraic geometry such results would be described as lifting
solutions of equations into generic points. We discuss definitions
and general properties of liftings in Section 6. We also prove
Theorem \ref{tqe} which is a weak version of the Implicit function
theorem for NTQ systems. We call it the Parametrization theorem.
This weak version of the Implicit function theorems  forms an
important part of the solution of Tarski's problems in \cite{KM7}.
All Implicit function theorems will be proved in Section 7.

In Sections 4 and 5 we describe  a new version of the
Makanin-Razborov process for a system of equations with
parameters, describe a solution set of such a system (Theorems
\ref{th:5.3.1} and \ref{ge}) and introduce a new type of equations
over groups, so-called {\em cut equations} (see Definition
\ref{df:cut} and Theorem \ref{th:cut}).

We collect some preliminary results and basic notions of algebraic
geometry for free groups in Section 2. In Section 3 we discuss
first order formulas over a free group and reduce an arbitrary
sentence to a relatively simple form.

 This paper is an extended version of the paper \cite{KM5}, the basic  version of the Implicit function theorem  was
announced at the Model Theory conference  at MSRI in  1998 (see
\cite{MSRI} and \cite{KM4}).

We thank Igor Lysenok who carefully read the manuscript and whose
numerous  remarks and suggestions have substantially improved the
exposition.

\section{Preliminaries}

\subsection{Free monoids and free groups}
\label{se:2-1}

Let $A= \{a_1, \ldots, a_m\}$ be a set. By $F_{mon}(A)$ we denote
the free monoid generated by  $A$ which is defined as the set of
all words (including the empty word $1$)  over the alphabet $A$
with concatenation as multiplication. For a word $w = b_1 \ldots
b_n,$ where $ b_i \in A$,    by $|w|$ or $d(w)$ we denote the
length $n$ of $w$.

To each $a\in A$ we associate a symbol $a^{-1}$. Put $A^{-1} =
\{a^{-1} \mid a \in A\},$ and suppose that $A\cap
A^{-1}=\emptyset.$   We assume that $a^{1} = a$, $(a^{-1})^{-1}=
a$ and $A^1 = A$. Denote  $A^{\pm 1} = A \cup A^{-1}$. If $w =
b_1^{\varepsilon_1} \ldots b_n^{\varepsilon_n} \in F_{mon}(A^{\pm
1})$, where  $ (\varepsilon_i \in \{1, -1\}) $, then we put
$w^{-1} = b_n^{-\varepsilon_n} \ldots b_1^{-\varepsilon_1};$ we
see that $w^{-1}\in M(A^{\pm 1})$ and say that $w^{-1}$ is an
inverse of $w$. Furthermore, we put $1^{-1}=1$.

A word $w \in F_{mon}(A^{\pm 1})$ is called {\it reduced}  if it
does not contain subwords  $bb^{-1}$ for $b \in A^{\pm 1}$. If $w
= w_1bb^{-1}w_2,$ $w  \in F_{mon}(A^{\pm 1})$ then we say that
$w_1w_2$ is obtained from $w$ by an elementary reduction $ bb^{-1}
\rightarrow 1$. A reduction process for $w$ consists of finitely
many reductions which bring $w$ to a reduced word $\bar w$. This
$\bar w$ does not depend on a particular reduction process and is
called the {\em reduced form} of $w$.

Consider a congruence relation on $F_{mon}(A^{\pm 1})$, defined
the following way:  two words are congruent if a reduction process
brings them to the same reduced word. The set of congruence
classes with respect to this relation forms a free group $F(A)$
with basis $A$. If not said otherwise, we assume that $F(A)$ is
given as the set of all reduced words in $A^{\pm 1}$.
Multiplication in $F(A)$ of two words $u, w$ is given by the
reduced form of their concatenation, i.e.,  $u\dot v = {\bar
{uv}}$. A word $w \in F_{mon}(A^{\pm 1})$ determines the element
$\bar w \in F(A)$, in this event we sometimes say that $w$ is an
element of $F(A)$ (even though $w$ may not be reduced).

 Words $u, w \in F_{mon}(A^{\pm 1})$ are {\it graphically equal} if they are
equal in the monoid $F_{mon}(A^{\pm 1})$ (for example,
$a_1a_2a_2^{-1}$ is not graphically equal to $a_1$).

Let $X=\{x_1,\ldots ,x_n\}$be a finite set of elements disjoint
with $A$. Let $w(X)= w(x_1,\ldots,x_n)$ be a word in the alphabet
$(X\cup A)^{\pm 1}$ and $U = (u_1(A), \ldots, u_n(A))$ be a tuple
of words in the alphabet $A^{\pm 1}$. By $w(U)$ we denote the word
which is obtained from $w$ by replacing each $x_i$ by $u_i$.
Similarly, if $W = (w_1(X), \ldots, w_m(X))$ is an $m$-tuple of
words in variables $X$ then by $W(U)$ we denote the tuple
$(w_1(U),\ldots,w_m(U))$. For any set $S$ we denote by $S^n$ the
set of all $n$-tuples of elements from $S$. Every word $w(X)$
gives rise to a map $p_w : (F_{mon}(A^{\pm 1}))^n \rightarrow
F_{mon}(A^{\pm 1})$ defined by $p_w(U) = w(U)$ for $U \in
F_{mon}(A^{\pm 1})^n$.  We call $p_w$ the word map defined by
$w(X)$. If $W(X) = (w_1(X), \ldots, w_m(X))$ is an $m$-tuple of
words in variables $X$ then we define a word map $P_W:
(F_{mon}(A^{\pm 1}))^n \rightarrow F_{mon}(A^{\pm 1})^m$ by the
rule $P_W(U) = W(U)$.

\subsection{ On G-groups}
\label{se:2-2}

 For the purpose of algebraic geometry over a given fixed group
 $G$,
one has to consider the category of $G$-groups, i.e., groups which
contain the group $G$ as a distinguished subgroup. If $H$ and $K$
are $G$-groups then a homomorphism $\phi: H \rightarrow K$ is a
$G$- homomorphism if $g^\phi = g$ for every $g \in G$, in this
event we write $\phi: H \rightarrow_G K$. In this category
morphisms are $G$-homomorphisms; subgroups are $G$-subgroups, etc.
By $Hom_G(H,K)$ we denote the set of all $G$-homomorphisms from
$H$ into $K$. It is not hard to see that the free product $G \ast
F(X)$ is a free object in the category of $G$-groups. This group
is called  a free $G$-group with basis $X$,  and we denote it  by
$G[X]$.  A $G$-group $H$ is termed {\em finitely generated
$G$-group} if there exists a finite subset $A \subset H$ such that
the set $G \cup A$ generates $H$. We refer to \cite{BMR1} for a
general discussion on $G$-groups.

To deal with cancellation in the group $G[X]$ we need the
following notation. Let  $u = u_1 \ldots u_n \in G[X] = G \ast
F(X)$. We say that $u$  is {\em reduced} (as written) if $u_i \neq
1$, $u_i$ and $u_{i+1}$ are in
 different factors of the free product,  and if $u_i \in F(X)$ then it is reduced in the free group
 $F(X)$. By $red(u)$ we denote the reduced form of $u$.
  If  $red(u)  = u_1 \ldots u_n \in G[X]$, then we define  $|u| = n$, so $|u|$ is
  the syllable length  of $u$ in the free product $G[X]$.
 For reduced $u,v \in G[X]$, we write $u\circ v$ if the product $uv$ is reduced as written.
 If $u = u_1 \ldots u_n $ is  reduced  and $u_1, u_n$ are in
 different factors, then  we say that $u$ is {\em cyclically reduced}.

If $u = r \circ s$, $ v = s^{-1} \circ  t$, and $rt = r\circ t$
then we say that the word $s$ {\it cancels out in reducing}  $uv$,
or, simply, $s$ cancels out in $uv$. Therefore $s$ corresponds to
the {\it maximal} cancellation in $uv$.

\subsection{Formulas in the language $L_A$}
\label{se:2-3}

Let $G$ be a group generated by a set of generators $A$. The
standard first-order language of group theory, which we denote by
$L,$ consists of a symbol for multiplication $\cdot$, a symbol for
inversion $^{-1},$ and a symbol for the identity $1.$  To deal
with $G$-groups, we have to enlarge the language $L$ by all
non-trivial elements from $G$ as  constants.  In fact, we do not
need to add all the elements of $G$ as constants, it suffices to
add only new constants corresponding to the generating set $A$. By
$L_A$ we denote  the language $L$ with constants from $A$.

 A group word in variables
 $X$ and constants $A$ is a word $S(X,A)$ in the alphabet $(X\cup
A)^{\pm 1}$. One
 may consider the word $S(X,A)$ as a term in the language $L_A$.
Observe that
 every term in the language $L_A$ is  equivalent modulo the
axioms of group
 theory to a group word in variables $X$ and constants $A \cup\{1\}$.
 An {\it atomic formula}  in the language $L_A$ is a formula of the type
$S(X,A) =
 1$, where $S(X,A)$ is a group word in $X$ and $A$.  With a slight abuse
of language
 we  will consider atomic formulas in $L_A$  as equations over $G$, and
vice versa.
   A {\it Boolean combination}
 of atomic formulas in the language $L_A$ is a disjunction of
conjunctions of
  atomic  formulas or their negations. Thus every Boolean combination
$\Phi$  of atomic  formulas in $L_A$
 can be written in the form $\Phi =  \bigvee_{i=1}^n\Psi_i$, where each
$\Psi_i$ has one
 of the following forms:

 $$\bigwedge_{j = 1}^n(S_j(X,A) = 1 ),   \ \ \  \bigwedge_{j =
1}^n(T_j(X,A) \neq 1
 ),  \ \ \
\bigwedge_{j = 1}^n (S_j(X,A) = 1 ) \&  \bigwedge_{k = 1}^m (T_k
(X,A) \neq 1).
$$

Observe that if the group $G$ is not trivial, then every formula
$\Psi$, as above, can be written in the from
 $$\bigwedge_{j = 1}^n (S_j(X,A) = 1 \ \ \& \
\ T_j (X,A) \neq 1), $$ where (if necessary) we add into the
formula the trivial equality  $1 = 1$, or an inequality of  the
type $a \neq 1$ for a given fixed non-trivial $a \in A$.

 It follows from general results on
disjunctive normal forms in propositional logic that every
quantifier-free formula in the language $L_A$ is logically
equivalent (modulo the axioms of group theory) to a Boolean
combination of  atomic ones. Moreover, every formula $\Phi$ in
$L_A$ with variables $Z\{z_1,\ldots ,z_k\}$  is logically
equivalent to a formula of the type
 $$Q_1x_1 Q_2 x_2 \ldots Q_n x_n \Psi(X,Z,A),$$
   where  $Q_i \in \{\forall, \exists \}$,
and  $\Psi(X,Z,A)$ is a Boolean combination of atomic formulas in
variables from $X \cup Z$. Using vector notations $Q Y  = Q y_1
\ldots Q y_n$ for strings of similar quantifiers  we can rewrite
such formulas  in the form
 $$\Phi(Z)  =  Q_1 Z_1  \ldots Q_k Z_k \Psi(Z_1,  \ldots, Z_k, X).$$
Introducing  fictitious quantifiers, one can always rewrite the
formula $\Phi$ in the form
   $$\Phi(Z)  = \forall X_1 \exists Y_1 \ldots \forall X_k \exists Y_k
\Psi(X_1,
 Y_1, \ldots, X_k, Y_k, Z).$$

If $H$ is a $G$-group, then the set $Th_A(H)$ of all sentences in
$L_A$ which are valid in $H$  is called the {\it elementary
theory} of $H$ in the language $L_A$. Two $G$-groups $H$ and $K$
are {\it elementarily equivalent } in the language $L_A$ (or
$G$-elementarily equivalent) if $Th_A(H) = Th_A(K).$

Let $T$ be a set of sentences in the language $L_A$. For a formula
$\Phi(X)$ in the language $L_A$, we write  $T \vdash \Phi$ if
$\Phi$ is a logical consequence of the theory  $T$. If $K$ is  a
$G$-group, then we write $K \models T$ if every sentence from $T$
holds in $K$ (where we interpret  constants from $A$ by
corresponding elements in the subgroup $G$ of $K$).  Notice, that
$Th_A(H) \vdash \Phi$ holds if and only if $K \models \forall X
\Phi(X)$ for every $G$-group $K$ which is $G$-elementarily
equivalent to $H$. Two formulas $\Phi(X)$ and $\Psi(X)$ in the
language $L_A$ are said to be {\it equivalent modulo} $T$  (we
write $\Phi \sim_T \Psi$) if $T \vdash \forall X (\Phi(X)
\leftrightarrow \Psi(X))$. Sometimes, instead of $\Phi
\sim_{Th_A(G)} \Psi$ we write $\Phi \sim_G \Psi$ and say that
$\Phi$ is equivalent to $\Psi$ over $G$.

\subsection{Elements of algebraic geometry over groups}
\label{se:2-4}

Here  we  introduce some basic notions of algebraic geometry over
groups. We refer to \cite{BMR1} and \cite{KM9} for details.

Let $G$ be a group generated by a finite set $A$, $F(X)$ be a free
group with basis $X = \{x_1, x_2, \ldots  x_n\}$, $G[X] = G \ast
F(X)$ be a free product of $G$ and $F(X)$. If $S \subset G[X]$
then the expression $S = 1$ is called {\em a system of equations}
over $G$. As an element of the free product, the left side of
every equation in $S = 1$ can be written as a product of some
elements from $X \cup X^{-1}$ (which are called {\em variables})
and some elements from $A$ ({\em constants}). To emphasize this we
sometimes write $S(X,A) = 1$.

A {\em solution} of the system $ S(X) = 1$ over a group $G$ is a
tuple of elements $g_1, \ldots, g_n \in G$ such that after
replacement of each $x_i$ by $g_i$ the left hand side of every
equation in $S = 1$ turns into the trivial element of $G$.
Equivalently, a solution of the system $S = 1$ over $G$ can be
described as a $G$-homomorphism $\phi : G[X] \longrightarrow G$
such that $\phi(S) = 1$. Denote by $ncl(S)$ the normal closure of
$S$ in $G[X]$, and by $G_S$ the quotient group $G[X]/ncl(S)$. Then
every solution of $S(X) = 1$ in $G$ gives rise to a
$G$-homomorphism $G_S \rightarrow G$, and vice versa. By $V_G(S)$
we denote the set of all solutions in $G$ of the system $ S = 1$,
it is called the {\em algebraic set defined by} $S$. This
algebraic set $V_G(S)$ uniquely corresponds to the normal subgroup
$$ R(S) = \{ T(x) \in G[X] \ \mid \ \forall A\in G^n (S(A) = 1
\rightarrow T(A) = 1) \} $$ of the group $G[X]$. Notice that if
$V_G(S) = \emptyset$, then $R(S) = G[X]$. The subgroup $R(S)$
contains $S$, and it is called the {\it radical of $S$}. The
quotient group
$$G_{R(S)}=G[X]/R(S)$$ is the {\em coordinate group} of the
algebraic set  $V(S).$ Again, every solution of $S(X) = 1$ in $G$
can be described as a $G$-homomorphism $G_{R(S)} \rightarrow G$.

We recall from \cite{MyasExpo2} that a group $G$ is called a {\em
CSA group} if every maximal Abelian subgroup $M$ of $G$ is
malnormal, i.e., $M^g \cap M = 1$ for any $g \in G - M.$ The class
of CSA-groups is quite substantial. It includes all Abelian
groups, all torsion-free hyperbolic groups \cite{MyasExpo2}, all
groups acting freely on $\Lambda$-trees \cite{bass}, and  many
one-relator groups  \cite{gkm}.

We define a Zariski topology on $G^n$ by taking algebraic sets in
$G^n$ as a sub-basis for the closed sets of this topology. If $G$
is a  non-Abelian CSA group, in particular, a non-Abelian freely
discriminated group, then the union of two algebraic sets is again
algebraic (see Lemma \ref{le:gur}). Therefore the closed sets in
the Zariski topology over $G$ are precisely the algebraic sets.

A $G$-group $H$ is called {\it equationally Noetherian} if every
system $S(X) = 1$ with coefficients from $G$ is equivalent over
$G$ to a finite subsystem $S_0 = 1$, where $S_0 \subset S$, i.e.,
$V_G(S) = V_G(S_0)$. If $G$ is $G$-equationally Noetherian, then
we say that $G$ is equationally Noetherian. It is known that
linear groups (in particular, freely discriminated groups) are
equationally Noetherian (see \cite{Gub}, \cite{Br}, \cite{BMR1}).
If $G$ is equationally Noetherian then the Zariski topology over
$G^n$ is {\em Noetherian} for every $n$, i.e., every proper
descending chain of closed sets in $G^n$ is finite. This implies
that every algebraic set $V$ in $G^n$ is a finite union of
irreducible subsets (called {\it irreducible components} of $V$),
and such a decomposition of $V$ is unique. Recall that a closed
subset $V$ is {\it irreducible} if it is not a union of two proper
closed (in the induced topology) subsets.

Two algebraic sets $V_F(S_1)$ and $V_F(S_2)$ are {\em rationally
equivalent} if  there exists an isomorphism between their
coordinate groups which is identical on $F$.
\subsection{Discrimination and big powers}
\label{se:2-5}

Let $H$ and $K$ be $G$-groups.  We say that a family of
$G$-homomorphisms ${\mathcal F} \subset Hom_G(H,K)$ {\it
separates} [{\it discriminates}] $H$ into $K$ if for every
non-trivial element $h \in H$ [every finite set of non-trivial
elements $H_0 \subset H$] there exists $\phi \in {\mathcal F}$
such that $h^\phi \neq 1$ [$h^\phi \neq 1$ for every $h \in H_0$].
In this case we say that $H$ is $G$-{\it separated} ($G$-{\it
discriminated}) by $K$. Sometimes we do not mention $G$ and simply
say that $H$ is separated [discriminated] by $K$. In the event
when $K$ is a free group we say that $H$ is  {\it freely
separated} [{\it freely discriminated}].

Below we describe a  method of discrimination which is called a
{\em big powers} method.  We refer to \cite{MyasExpo2} and
\cite{KvM} for details about BP-groups.

Let $G$ be a group. We say that a tuple $u = (u_{1},...,u_{k})\in
G^{k}$ has {\it commutation} if $[u_{i}, u_{i+1}] = 1$ for some
$i=1,...k-1.$ Otherwise we call $u$ {\it commutation-free}.

\begin{df}  A group $G$ satisfies the {\em big powers condition} (BP) if
for any commutation-free tuple $ u = (u_{1},...,u_{k})$ of
elements from $G$ there exists an integer $n(u)$ (called {\em a
boundary of separation} for $u$) such that
$${u_{1}}^ {{\alpha}_{1}}...u_{k}^{{\alpha}_{k}} \neq 1$$
 for
any integers $\alpha_1, \ldots, \alpha_k \geq n(u)$. Such groups
are called {\em BP-groups}. \end{df}

The following  provides a host of examples of BP-groups.
Obviously, a subgroup of a BP-group is a BP-group; a group
discriminated by a BP-group is a BP-group (\cite{MyasExpo2});
every torsion-free hyperbolic group is a BP-group (\cite{Ol'sh1}).
From those facts it follows that every freely discriminated group
is a BP-group.

Let $G$ be a non-Abelian CSA group and $u \in G$  not be a proper
power. The following HNN-extension $$G(u,t) =  \langle G, t \mid
g^t = g (g \in C_G(u))\rangle$$  is called a {\em free extension}
of the centralizer $C_G(u)$ by a letter $t$. It is not hard to see
that for any integer $k$ the map $t \rightarrow u^k$ can be
extended uniquely to a $G$-homomorphism $\xi_k : G(u,t)
\rightarrow G$.

The result below is the essence of the big powers method of
discrimination.

\bigskip
{\bf Theorem}(\cite{MyasExpo2}) {\it Let $G$ be a non-Abelian CSA
BP-group. If $G(u,t)$ is a free extension of the centralizer of
the nonproper power $u$ by $t$,   Then  the family of
$G$-homomorphisms $\{\xi_k \mid k \ is \ an \ integer \}$
discriminates $G(u,v)$ into $G$. More precisely, for every $w \in
G(u,t)$ there exists an integer  $N_w$ such that for every $k \geq
N_w$ $w^{\xi_k} \neq 1$.}

\smallskip
If $G$ is  a non-Abelian CSA BP-group and $X$ is a finite set,
then   the group $G[X]$ is $G$-embeddable into $G(u,t)$ for any
nonproper power $u \in G$. It follows from the theorem above that
$G[X]$ is $G$-discriminated by $G$.

Unions of chains of extensions of centralizers play an important
part in this paper. Let   $G$ be a non-Abelian CSA BP- group and
 $$G = G_0 <  G_1 <  \ldots < G_n $$
  be a chain of extensions of centralizers $G_{i+1} = G_i(u_i,t_i)$. Then every
  $n$-tuple of integers $p = (p_1, \ldots, p_n)$ gives rise to a $G$-homomorphism
  $\xi_p: G_n \rightarrow G$ which is  composition of homomorphisms
  $\xi_{p_i}: G_i \rightarrow G_{i-1}$ described above. If a
  centralizer of $u_i$ is extended several times, we can suppose it
  is extended on the consecutive steps by letters $t_i,\ldots
  ,t_{i+j}$. Therefore $u_{i+1}=t_i,\ldots ,u_{i+j}=t_{i+j-1}.$

  A set $P$  of
  $n$-tuples  of integers is called {\it unbounded} if  for every integer $d$
  there exists a tuple $p  = (p_1, \ldots, p_n) \in P$ with $p_i \geq d$ for each $i$.
   The following result is a consequence of  the theorem above.

\bigskip
{\bf Corollary} {\it
   Let $G_n$ be as above. Then for every unbounded set of tuples $P$
   the set of $G$-homomorphisms $\Xi_P = \{\xi_p \mid p \in P\}$ $G$-discriminates
   $G_n$ into $G$.}


   Similar results hold for  infinite chains
   of extensions of centralizers (see \cite{MyasExpo2}) and \cite{BMR3}).
   For example,   Lyndon's  free
   $Z[x]$-group $F^{Z[x]}$ can be realized as  union of a  countable chain of extensions of
   centralizers which starts with the free group $F$ (see  \cite{MyasExpo2}),
   hence  there exists a family of $F$-homomorphisms which discriminates
   $F^{Z[x]}$ into $F$.

\subsection{Freely discriminated groups}
\label{se:2-6}

Here we formulate several results on freely discriminated groups
which are crucial for our considerations.

It is not hard to see that every freely discriminated group is a
torsion-free CSA group \cite{BMR1}.

Notice that every CSA group is commutation transitive
\cite{MyasExpo2}. A group $G$ is called {\em commutation
transitive} if commutation is transitive on the set of all
non-trivial elements of $G$, i.e., if $a,b,c \in G - \{1\}$ and
$[a,b] = 1, [b,c] = 1$, then $[a,c] = 1.$ Clearly, commutation
transitive groups are precisely the groups in which centralizers
of non-trivial elements are commutative.  It is easy to see that
every commutative transitive group $G$ which satisfies the
condition $[a,a^b] = 1 \rightarrow [a,b] = 1$ for all $a, b \in G$
is CSA.

\medskip
{\bf Theorem }(\cite{R1}).   {\it Let $F$ be a free non-abelian
group. Then a finitely generated $F$-group $G$ is freely
$F$-discriminated by $F$ if and only if $G$ is $F$-universally
equivalent to $F$ (i.e., $G$ and $F$ satisfy precisely the same
universal sentences in the  language $L_A$).}

\smallskip

\medskip
{\bf Theorem }(\cite{BMR1}, \cite{KM9}).  {\it Let $F$ be a free
non-abelian group. Then a finitely generated $F$-group $G$ is the
coordinate group of a non-empty irreducible algebraic set over $F$
if and only if $G$ is freely $F$-discriminated by $F$.}

\smallskip

\medskip
{\bf Theorem } (\cite{KMNull}). {\it Let $F$ be a non-abelian free
group. Then a finitely generated $F$-group is  the coordinate
group $F_{R(S)}$ of an irreducible non-empty algebraic set $V(S)$
over $F$ if and only if $G$ is $F$-embeddable into the  free
Lyndon's $Z[t]$-group $F^{Z[t]}$.  }

\smallskip

This theorem implies that  finitely generated freely discriminated
groups are finitely presented, also it allows one to present such
groups as fundamental groups of graphs of groups of a very
particular type (see \cite{KMNull} for details).

\subsection{Quadratic equations over freely discriminated groups}
\label{se:2-7}

In this section we collect some known results about quadratic
equations over fully residually free groups, which will be in use
throughout this paper.

\smallskip
 Let $S \subset G[X]$. Denote by $var(S)$ the set of variables that
occur in $S$.
\begin{df} A set $S \subset G[X]$ is called quadratic  if every variable
from
 $var(S)$ occurs in $S$ not more then twice. The set $S$ is strictly
quadratic if every letter from $var(S)$ occurs in $S$ exactly
twice.

A system $S = 1$ over $G$ is {\em quadratic [strictly quadratic]}
if the corresponding set $S$ is quadratic [strictly quadratic].
\end{df}

\begin{df}
A standard quadratic equation over a group $G$ is an equation of
the one of the following forms (below $d,c_i$ are nontrivial
elements from $G$):
\begin{equation}\label{eq:st1}
\prod_{i=1}^{n}[x_i,y_i] = 1, \ \ \ n > 0;
\end{equation}
\begin{equation}\label{eq:st2}
\prod_{i=1}^{n}[x_i,y_i] \prod_{i=1}^{m}z_i^{-1}c_iz_i d = 1,\ \ \
n,m \geq 0, m+n \geq 1 ;
\end{equation}
\begin{equation}\label{eq:st3}
\prod_{i=1}^{n}x_i^2 = 1, \ \ \ n > 0;
\end{equation}
\begin{equation}\label{eq:st4}
\prod_{i=1}^{n}x_i^2 \prod_{i=1}^{m}z_i^{-1}c_iz_i d = 1, \ \ \
n,m \geq 0, n+m \geq 1.
\end{equation}

Equations (\ref{eq:st1}), (\ref{eq:st2}) are called {\em
orientable} of genus $n$, equations (\ref{eq:st3}), (\ref{eq:st4})
are called {\em non-orientable} of genus $n$.
\end{df}

\begin{lm} \label{EC1}
Let $W$ be a strictly quadratic word over $G$. Then there is a
$G$- automorphism $f \in Aut_G(G[X])$ such that  $W^f$ is a
standard quadratic word over $G.$
\end{lm}
{\em Proof.}   See \cite{EC}.

\begin{df}
Strictly quadratic words  of the type $ [x,y], \ x^2, \ z^{-1}cz,
$ where $c \in G$, are called {\em atomic quadratic words} or
simply {\em atoms}.
\end {df}

By definition a standard quadratic equation $S = 1$  over $G$  has
the form $$ r_1 \ r_2 \ldots r_k d = 1,$$ where $r_i$ are atoms,
$d \in G$.  This number $k$ is called the {\it atomic rank } of
this equation,  we denote it by $r(S)$.  In Section \ref{se:2-4}
we defined the notion of the coordinate group $G_{R(S)}.$ Every
solution of the system $S=1$ is a homomorphism $\phi : G_{R(S)}
\rightarrow G$.
\begin{df}  Let $S = 1$ be a standard quadratic equation written in the
atomic form
 $r_1r_2\ldots r_kd = 1 $ with $k \geq 2$.  A solution $\phi : G_{R(S)}
\rightarrow G$
 of $S = 1$  is called:
 \begin{enumerate}
 \item degenerate, if $r_i^\phi = 1$ for some $i$, and
 non-degenerate otherwise;
 \item  commutative, if $[r_i^{\phi},r_{i+1}^{\phi}]=1$ for all
$i=1,\ldots ,k- 1,$  and  non-commutative otherwise;
 \item in a general position, if $[r_i^{\phi},r_{i+1}^{\phi}] \neq 1$ for all
$i=1,\ldots ,k-1,$.
 \end{enumerate}
 \end{df}
 Observe that if a standard quadratic equation $S(X) = 1$ has a
 degenerate non-commutative solution then it has a non-degenerate
 non-commutative solution {see \cite{KMNull}).
\begin{theorem}[\cite{KMNull}]  Let $G$ be a freely discriminated
group and $S = 1$ a standard quadratic equation over $G$ which has
a solution in $G$.
 In the following cases a standard quadratic equation  $S=1$
 always has a solution in a general position:
\begin{enumerate}
 \item $S=1$ is of the form (\ref{eq:st1}), $n>2$;
 \item $S=1$ is of the form (\ref{eq:st2}),  $n>0, \ n+m>1$;
 \item $S=1$ is of the form (\ref{eq:st3}), $n>3;$
 \item $S=1$ is of the form (\ref{eq:st4}), $n>2;$
 \item  $r(S) \geq 2$  and $S=1$ has a non-commutative solution.
\end{enumerate}
\end{theorem}
The following theorem describes the radical $R(S)$ of a standard
quadratic equation $S=1$ which has at least one  solution in a
freely discriminated group $G$.

\begin{theorem}[\cite{KMNull}]
\label{th:Nul} Let $G$ be a freely discriminated  group and let
$S=1$ be a standard quadratic equation over $G$ which has a
solution in $G$. Then
 \begin{enumerate}
\item If $S = [x,y]d$ or $S = [x_1,y_1][x_2,y_2]$, then $R(S) =
ncl(S)$; \item If $S = x^2d$, then $R(S) = ncl(xb)$ where $b^2 =
d$; \item If $S=c^zd$, then $R(S) = ncl([zb^{-1},c])$ where
$d^{-1} = c^b$; \item If $S = x_1^2x_2^2$, then $R(S) =
ncl([x_1,x_2])$; \item If $S = x_1^2x_2^2x_3^2$, then $R(S) =
ncl([x_1,x_2], [x_1,x_3], [x_2,x_3])$; \item If  $r(S) \geq 2$ and
$S=1$
 has a non-commutative solution,  then $R
(S)=ncl (S)$; \item If  $S = 1$ is of the type (\ref{eq:st4}) and
all solutions of $S=1$
 are commutative, then $R(S)$ is the normal closure of the following
system:
  $$  \{x_1\ldots x_n=s_1\ldots s_n, [x_k,x_l]=1, [a_i^{-1}z_i,x_k]=1,
[x_k,C] = 1, [a_i^{-1}z_i,C] = 1,  $$ $$ [a_i^{-1}z_i,a_j^{-1}z_j]
= 1 \ \ (k,l = 1,\ldots,n; i,j = 1, \ldots , m) \} , $$ where $x_k
\rightarrow s_k, z_i \rightarrow a_i$ is a solution of $S=1$ and
$C=C_G(c_1^{a_1},\ldots ,c_m^{a_m},s_1,\ldots ,s_n)$ is the
corresponding centralizer.  The group $G_{R(S)}$ is an extension
of the centralizer $C$.
  \end{enumerate}
  \end{theorem}

\begin{df}\label{regular}
 A standard quadratic equation $S = 1$ over $F$ is called {\em regular} if
either it is an equation of the type $[x,y] = d \ (d\neq 1)$, or
the equation $[x_1,y_1][x_2,y_2] = 1$, or $r(S) \geq 2$ and $S(X)
= 1$  has a non-commutative solution  and it is not an equation of
the type $c_1^{z_1}c_2^{z_2} = c_1c_2$, $x^2c^z = a^2c$,
$x_1^2x_2^2 = a_1^2a_2^2$.
\end{df}

Put
 $$\kappa(S) = |X| + \varepsilon(S),$$
  where $\varepsilon(S) = 1$  if
the coefficient $d$ occurs in $S$, and  $\varepsilon(S) = 0$
otherwise.

Equivalently, a standard quadratic equation $S(X) = 1$ is  {\em
regular} if $\kappa(S) \geq 4$ and there is a non-commutative
solution of $S(X) = 1$ in $G$, or it is an equation of the type
$[x,y]d = 1$.

Notice, that if $S(X) = 1$ has a solution in $G$,  $\kappa(S) \geq
4$, and $n
> 0$ in the orientable case ($n > 2$ in the non-orientable case),
then the equation $S = 1$ has a non-commutative solution, hence
regular.

\begin{cy}\begin{enumerate}
\item Every consistent orientable quadratic equation $S(X) = 1$ of
positive genus is regular, unless it is the equation $[x,y] = 1$;

\item Every consistent non-orientable equation of positive genus
is regular, unless  it is an equation of the type $x^2c^z = a^2c$,
$x_1^2x_2^2 = a_1^2a_2^2, x_1^2x_2^2x_3^2 = 1$, or $ S(X) =1$ can
be transformed to the form $[\bar z_i, \bar z_j]=[\bar z_i, a]=1,
\ i,j=1,\ldots ,m$ by changing variables.

\item Every standard quadratic equation $S(X) = 1$ of genus 0 is
regular unless either it is an equation of the type $c_1^{z_1} =
d, c_1^{z_1}c_2^{z_2} = c_1c_2$, or $ S(X) =1$ can be transformed
to the form $[\bar z_i, \bar z_j]=[\bar z_i, a]=1, \ i,j=1,\ldots
,m$ by changing variables.\end{enumerate}\end{cy}

\subsection{Formulation of the basic  Implicit function theorem}
\label{se:2-8}

In this section we formulate the implicit function theorem over
free groups in its basic simplest form. We refer to Sections
\ref{se:7.2}, \ref{se:7.3}  for the proofs and to Section
\ref{se:7.4} for  generalizations.

\begin{theorem} \label{IMTQ}
Let $S(X) = 1$ be a regular standard quadratic equation over a
non-abelian free group $F$ and let $T(X,Y) = 1$ be an equation
over $F$, $|X| = m, \ |Y| = n.$ Suppose that for any solution $U
\in V_F(S)$ there exists a tuple of elements $W \in F^n$ such that
$T(U,W) = 1.$ Then there exists a tuple of words $P = (p_1(X),
\ldots, p_n(X))$, with constants from $F$, such that $T(U,P(U)) =
1$ for any $U \in V_F(S)$. Moreover, one can fund a tuple $P$ as
above effectively.
\end{theorem}

From algebraic geometric view-point the implicit function theorem
tells one that (in the notations above)  $T(X,Y) = 1$ has a
solution at a generic point of the equation $S(X) = 1$.

\section{Formulas over freely discriminated groups}

In this section we collect some results (old and new) on how to
effectively rewrite formulas over a non-Abelian freely
discriminated group $G$ into more simple or more convenient
"normal" forms. Some of these results hold for many other groups
beyond the class of freely discriminated ones. We do not present
the most general formulations here, instead,  we limit our
considerations to  a class of groups ${\mathcal T}$ which will
just suffice for our purposes.

Let us fix a finite set of constants $A$ and the corresponding
group theory language $L_A$, let also $a, b$ be two  fixed
elements in $A$.
\begin{df}
A group $G$ satisfies Vaught's conjecture if the following
universal sentence holds in $G$
\begin{description}
\item [(V)] $ \ \ \forall x \forall y \forall z (x^2 y^2 z^2 = 1
\rightarrow [x,y] = 1\  \&\  [x,z] = 1 \ \& \ [y,z] = 1)$

\end{description}
\end{df}
Lyndon proved that every free group satisfies the condition (V)
(see \cite{Ly0}).

Denote by ${\mathcal T}$ the class of all groups $G$ such that:

1) $G$ is torsion-free;

2) $G$ satisfies  Vaught's conjecture;

3) $G$ is CSA;

4) $G$ has  two distinguished elements $a, b$ with $[a,b] \neq 1$.

It is easy to write down axioms for the class ${\mathcal T}$ in
the language $L_{\{a,b\}}$. Indeed, the following  universal
sentences describe the conditions 1)-4) above:
\begin{description}
\item [(TF)] $x^n = 1 \rightarrow x = 1 \ \ (n = 2,3, \ldots);$
\item [(V)]   $  \ \ \forall x \forall y \forall z (x^2 y^2 z^2 =
1 \rightarrow [x,y] = 1\  \&\  [x,z] = 1 \ \& \ [y,z] = 1);$ \item
[(CT)] $  \ \ \forall x \forall y \forall z(x \neq 1 \ \& \ y \neq
1 \ \& \ z\neq 1 \ \& \ [x,y] = 1 \ \& \ [x,z] = 1 \rightarrow
[y,z] = 1);$ \item [(WCSA)] $  \ \ \forall x \forall y ([x,x^y] =
1 \rightarrow [x,y] = 1);$ \item [(NA)] $[a,b] \neq 1.$
\end{description} Observe that the condition (WCSA) is a weak form of (CSA) but (WCSA) and (CT) together
 provide the CSA condition. Let $GROUPS$  be a set of
axioms of
 group theory. Denote by  $A_{\mathcal
T}$ the union of axioms (TF), (V), (CT), (WCSA),  (NA) and
$GROUPS$.  Notice that the axiom (V) is equivalent modulo $GROUPS$
to the following quasi-identity $$ \forall x \forall y \forall z
(x^2 y^2 z^2 = 1 \rightarrow [x,y] = 1) .$$ It follows that all
axioms in $A_{\mathcal T}$, with  exception of  (CT) and (NA),
are quasi-identities.

\begin{lm}
\label{le:ax} The class ${\mathcal T}$ contains all freely
discriminated non-abelian groups.
 \end{lm}
  {\it Proof. } We show here that every  freely discriminated group $G$
satisfies (V). Similar arguments work for the other conditions.
 If $ u^2v^2w^2 = 1$ for some $u, v, w \in G$ and, say, $[u,v] \neq
1$, then there exists a homomorphism $\phi: G \rightarrow F$ from
$G$  onto a free group $F$ such  that $[u^\phi,v^\phi] \neq 1$.
This shows that the elements $u^\phi, v^\phi, w^\phi$ in $F$ give
a counterexample to  Vaught's conjecture. This contradicts to the
Lyndon's result. Hence (V) holds in $G$. This proves the lemma.

Almost all results in  this section state that a formula $\Phi(X)$
in $L_A$ is equivalent modulo $A_{\mathcal T}$ to a formula
$\Psi(X)$ in $L_A$. We will use these results in the following
particular form. Namely, if $G$ is a group  generated by $A$ and
$H$ is a $G$-group from $\mathcal T$ then for any tuple of
elements $U \in H^n$ (here $n = |X|$) the formula  $\Phi(U)$ holds
in $H$ if and only if $\Psi(U)$ holds in $H$.
\subsection{Quantifier-free  formulas}

In this section by letters $X, Y, Z$ we denote finite tuples of
variables.

 The following result is due to A.Malcev \cite{Mal2}. He proved
it for free groups, but his argument is valid in a more general
context.
\begin{lm} \label{le:malcev} Let
$G \in {\mathcal T}$.  Then the equation \beq \label{eq:malcev}
x^2ax^2a^{- 1} = (ybyb^{-1})^2 \eeq has only the trivial solution
$x = 1$ and $y = 1$ in $G$.
\end{lm}

Proof. Let $G$ be as above and let $x,y$ be a solution in $G$ of
the equation (\ref{eq:malcev}) such that $x \neq 1.$ Then \beq
\label{eq:malcev2} (x^2a)^2 a^{-2} = ((yb)^2b^{-2})^2. \eeq In
view of the condition (V), we deduce from (\ref{eq:malcev2}) that
$[x^2a,a^{-1}]= 1$,  hence $[x^2,a^{-1}]= 1$. By transitivity of
commutation    $[x,a] = 1$ (here we use inequality $x \neq 1$).
Now,  we can rewrite (\ref{eq:malcev2}) in the form $$ x^2 x^2 =
((yb)^2b^{-2})^2, $$ which implies (according to (V)), that
$[x^2,(yb)^2b^{-2}] = 1$, and hence (since $G$ is torsion-free)
\beq \label{eq:malcev3} x^2 = (yb)^2b^{-2}. \eeq Again, it follows
from (V) that $[y,b] = 1$. Henceforth, $x^2 = y^2$ and, by the
argument above,  $x = y$. We proved that $[x,a] = 1$ and $[x,b] =
1$ therefore,  by  transitivity of commutation,   $[a,b] =1$,
which contradicts to the choice of $a,b$.  This contradiction
shows that $x = 1$. In this event,  the equation
(\ref{eq:malcev2}) transforms into $$((yb)^2b^{-2})^2 = 1,$$ which
implies  $ (yb)^2b^{-2} = 1$. Now from (V) we deduce that $[yb,b]
= 1$, and hence $[y,b] = 1.$  It follows  that $y^2 = 1$, so $y =
1$, as desired.

\begin{cy}
\label{co:1} Let $G \in {\mathcal T}$. Then for any finite system
of equations $S_1(X) = 1, \ldots, S_k(X) = 1$ over $G$ one can
effectively find a single equation $S(X) = 1$ over $G$ such that
$$V_G(S_1, \ldots, S_n) = V_G(S).$$
\end{cy}
{\it Proof}. By induction it suffices to prove the result for $k
 = 2$. In this case, by the lemma above, the following equation (after
bringing
 the right side to the left)
$$ S_1(X)^2aS_1(X)^2 a^{-1} = (S_2(X)bS_2(X)b^{-1})^2 $$ can be chosen as the
equation $S(X) = 1.$

 \begin{cy}
\label{co:2} For any finite system of equations $$S_1(X) = 1,
\ldots, S_k(X) = 1$$ in $L_A,$  one can effectively find a single
equation $S(X) = 1$ in $L_A$ such that $$ ( \bigwedge_{i = 1}^{k}
S_i(X) = 1)\  \sim_{ A_{\mathcal T}} \ S(X) = 1. $$
\end{cy}
\begin{rk}
In the proof of Lemma  \ref{le:malcev} and Corollaries \ref{co:1}
and \ref{co:2} we did not use the condition (WCSA) so the results
hold for an arbitrary non-abelian torsion-free commutation
transitive group satisfying  Vought's conjecture.
\end{rk}
 The next lemma shows how to rewrite  finite disjunctions of
equations into conjunctions of equations. In the case of free
groups this result was known for years (in \cite{Mak84}
 Makanin
attributes this to Y.Gurevich). We give here a different proof.

\begin{lm}
\label{le:gur} Let $G$ be a CSA group and let $a,b$ be arbitrary
non-commuting elements in $G$. Then for any solution $x,y \in G$
of the system \beq \label{eq:gur} [x,y^a] = 1, \ \  [x,y^b] = 1, \
\ [x,y^{ab}] = 1, \eeq either $x = 1$ or $y = 1$. The converse is
also true.
\end{lm}

{\it Proof}. Suppose $x,y$ are non-trivial elements from $G$, such
that $$ [x,y^a] = 1, \ \ [x,y^b] = 1, \ \ [x,y^{ab}] = 1. $$ Then
by the transitivity of commutation $[y^b,y^{ab}] = 1$ and
$[y^a,y^b] = 1$. The first relation implies that $[y,y^a]= 1$, and
since a maximal Abelian subgroup $M$ of $G$ containing $y$ is
malnormal in $G$, we have $[y,a] = 1$. Now from $[y^a,y^b] = 1$ it
follows that $[y,y^b] = 1$ and, consequently, $[y,b]= 1$. This
implies $[a,b] = 1$, a contradiction, which completes the proof.

Combining Lemmas \ref{le:gur} and \ref{le:malcev} yields an
algorithm to encode an arbitrary finite disjunction of equations
into a single equation.
\begin{cy}
Let $G \in {\mathcal T}$. Then for any finite set of equations
$S_1(X) = 1, \ldots, S_k(X) = 1$ over $G$ one can effectively find
a single equation $S(X) = 1$ over $G$ such that $$ V_G(S_1) \cup
\ldots \cup V_G(S_k) = V_G(S). $$
\end{cy}

Inspection of the proof above shows that the following corollary
holds.
\begin{cy}
\label{co:4} For any finite set of equations $S_1(X) = 1, \ldots,
S_k(X) = 1$ in $L_A,$ one can effectively find a single equation
$S(X) = 1$ in $L_A$ such that
 $$
(\bigvee_{i = 1}^{k} S_i(X) = 1) \ \sim_{A_{\mathcal T}} \ S(X) =
1.
$$
\end{cy}
\begin{cy}
\label{co:qfreepos} Every positive quantifier-free formula
$\Phi(X)$  in $L_A$ is equivalent modulo $A_{\mathcal T}$ to a
single equation $S(X) = 1$.
\end{cy}

 The next result shows that  one can effectively encode
finite conjunctions and finite disjunctions of {\em inequalities}
into a single inequality modulo $A_{\mathcal T}$.

\begin{lm}
\label{le:ineq}  For any finite set of inequalities
 $$S_1(X) \neq 1, \ldots,
S_k(X) \neq 1$$
 in $L_A$, one can effectively find an inequality
$R(X) \neq 1$ and an inequality $T(X) \neq 1$ in $L_A$ such that
 $$
(\bigwedge_{i = 1}^{k} S_i(X) \neq 1)\  \sim_{A_{\mathcal T}} \
R(X) \neq 1 $$ and
 $$
 (\bigvee_{i = 1}^{k}S_i(X) \neq 1) \ \sim_{A_{\mathcal T}} \  T(X) \neq 1.
$$
\end{lm}

{\em Proof.}   By Corollary \ref{co:4} there exists an equation
$R(X) =  1$ such that $$ \bigvee_{i = 1}^{k} (S_i(X) = 1) \
\sim_{A_{\mathcal T}} \  R(X) = 1. $$ Hence $$( \bigwedge_{i =
1}^{k} S_i(X) \neq 1 )  \  \sim_{A_{\mathcal T}}   \ \neg(\bigvee
_{i = 1}^{k} S_i(X) = 1)  \ \sim_{A_{\mathcal T}} \ \neg(R(X) = 1)
\ \sim_{A_{\mathcal T}} \ R(X) \neq 1. $$ This proves the first
part of the result. Similarly, by Corollary \ref{co:2} there
exists an equation $T(X) = 1$ such that
 $$
 ( \bigwedge_{i = 1}^{k} S_i(X) = 1) \  \sim_{A_{\mathcal T}} \ T(X)
= 1. $$ Hence

 $$( \bigvee_{i = 1}^{k}S_i(X) \neq 1) \   \sim_{A_{\mathcal T}} \
\neg (\bigwedge_{i = 1}^{k}S_i(X) = 1)$$ $$\  \sim_{A_{\mathcal
T}} \ \neg (T(X) = 1) \  \sim_{A_{\mathcal T}}  \ T(X) \neq 1. $$
This completes the proof.

\begin{cy}
\label{co:atomic}
 For every quantifier-free formula  $\Phi(X)$ in the language $L_A$, one
can effectively find a formula  $$\Psi(X) =  \bigvee_{i=1}^n
(S_i(X) = 1 \ \ \& \ \ T_i(X) \neq 1) $$ in $L_A$ which is
equivalent to $\Phi(X)$ modulo $A_{\mathcal T}$. In particular, if
$G \in {\mathcal T}$, then  every quantifier-free formula
$\Phi(X)$ in $L_G$ is equivalent over $G$ to a formula  $\Psi(X)$
as above.
\end{cy}

\subsection{Universal formulas over $F$}

In this section we discuss canonical forms of universal formulas
in the language $L_A$ modulo the theory $A_{\mathcal T}$ of the
class ${\mathcal T}$ of all torsion-free non-Abelian CSA groups
satisfying  Vaught's conjecture.  We show that every universal
formula in $L_A$ is equivalent modulo $A_{\mathcal T}$ to a
universal formula in canonical radical form. This implies that if
$G \in {\mathcal T}$ is generated by $A$, then  the universal
theory of $G$ in the language $L_A$ consists of the the axioms
describing the diagram of $G$ (multiplication table for $G$ with
all the equalities and inequalities between group words in $A$),
the set of axioms $A_{\mathcal T}$,  and   a set of axioms $A_R$
which describes the radicals of finite systems over $G$.

Also, we describe an effective quantifier elimination for
universal positive formulas in $L_A$ modulo $Th_A(G)$, where $G
\in {\mathcal T}$ and $G$ is a BP-group (in particular, a
non-Abelian freely discriminated group).  Notice, that in Section
\ref{se:mer} in the case when $G$ is a free group,  we describe an
effective quantifier elimination procedure (due to  Merzljakov and
Makanin) for arbitrary positive sentences modulo $Th_A(G)$.

Let $G \in {\mathcal T}$ and $A$ be a generating set for $G$.

 We say that a universal formula in $L_A$ is  in  {\it canonical
radical form} (is a {\it radical formula}) if it has the following
form
 \begin{equation}
 \label{eq:rad}
  \Phi_{S,T}(X) = \forall Y (S(X,Y) = 1 \rightarrow T(Y) = 1)
  \end{equation}
for some $S \in G[X \cup Y], T \in G[Y].$

For an arbitrary finite system  $S(X) = 1$ with coefficients from
$A$ denote by ${\tilde S}(X) = 1$  an equation which is equivalent
over $G$ to the system $S(X) = 1$ (such ${\tilde S}(X)$ exists by
Corollary  \ref{co:2}). Then for the radical $R(S)$ of the system
$S=1$ we have
 $$
 R(S) = \{T \in G[X] \mid G \models \Phi_{{\tilde S},T}\}.$$
 It follows that the set of radical sentences
 $$A_S = \{\Phi_{{\tilde S},T} \mid G \models \Phi_{{\tilde S},T}\}$$
  describes precisely the radical $R(S)$  of the system $S = 1$ over
$G$, hence the name.

\begin{lm}
\label{le:unrad}
 Every universal formula in $L_A$ is equivalent modulo $A_{\mathcal T}$ to a
radical formula.
\end{lm}
{\em Proof.}    By Corollary \ref{co:atomic}  every Boolean
combination of atomic formulas in the language $L_A$ is equivalent
modulo $A_{\mathcal T}$ to a formula of the type
 $$ \bigvee_{i=1}^n (S_i = 1 \ \ \& \ \ T_i \neq 1). $$
This implies that every existential formula in $L_A$ is equivalent
to a formula in the form
 $$ \exists Y( \bigvee_{i=1}^n (S_i(X,Y) = 1 \ \ \& \ \ T_i(X,Y) \neq
1)).$$ This formula is equivalent modulo $A_{\mathcal T}$ to the
formula
 $$\exists z_1 \ldots \exists z_n  \exists Y
 ((\bigwedge_{i = 1}^n z_i \neq 1)\ \  \& \ \ (\bigvee_{i=1}^n (S_i(X,Y)
= 1 \ \ \& \ \ T_i(X,Y)
  = z_i)) ) .$$
  By Corollaries \ref{co:4} and \ref{co:2} one can effectively find $S
\in
  G[X,Y,Z]$ and $T \in G[Z]$ (where $Z = (z_1, \ldots,z_n)$) such that
  $$ \bigvee_{i=1}^n (S_i(X,Y) = 1 \ \ \& \ \ T_i(X,Y) = z_i) \ \
\sim_{A_{\mathcal T}}
  \ \ S(X,Y,Z) =
  1$$
  and
  $$\bigwedge_{i = 1}^n (z_i \neq 1)\ \   \sim_{A_{\mathcal T}} \ \ T(Z)
\neq 1.$$
  It follows that every existential formula in $L_A$ is equivalent
modulo  $A_{\mathcal T}$  to a
  formula of the type
  $$ \exists Z \exists Y (S(X,Y,Z) = 1\ \  \& \ \ T(Z) \neq 1).$$
  Hence  every universal formula in $L_A$ is equivalent  modulo
$A_{\mathcal T}$ to a
  formula in the form
  $$\forall Z \forall Y (S(X,Y,Z) \neq 1\ \  \bigvee \ \ T(Z) = 1),$$
  which is equivalent to the radical  formula
  $$\forall Z \forall Y (S(X,Y,Z) = 1 \  \rightarrow   \  T(Z) = 1).$$
  This proves the lemma.

Now we consider universal positive formulas.
\begin{lm}
\label{le:3.2.5} Let $G$ be a BP-group from $\mathcal T$. Then $$
G \models \forall X (U(X) = 1) \Longleftrightarrow G[X] \models
U(X) = 1, $$ i.e., only the trivial equation has the whole set
$G^n$ as its  solution set.
\end{lm}
{\em Proof.}   The group $G[X]$ is  discriminated by $G$
\cite{BMR1}. Therefore, if the word $U(X)$ is a non-trivial
element of $G[X]$, then there exists a $G$-homomorphism $\phi:
G[X] \rightarrow G$ such that $U^\phi \neq 1.$ But then $U(X^\phi)
\neq 1$ in $G$ $--$ contradiction with conditions of the lemma. So
$U(X) = 1$ in $G[X].$

\begin{rk}
The proof above holds for every non-Abelian group $G$ for which
$G[X]$ is discriminated by $G$.
\end{rk}

 The next result shows how to eliminate quantifiers from positive
universal formulas over non-Abelian freely discriminated groups.

\begin{lm}
\label{le:elimin}  Let $G$ be a BP-group from $\mathcal T$. For a
given word $U(X,Y) \in G[X \cup Y]$, one can effectively find a
word $W(Y) \in G[Y]$ such that
 \beq \label{eq:elimin}  \forall X (U(X,Y) = 1)  \ \ \sim_G \ \ W(Y) =
1. \eeq
\end{lm}
{\it Proof}.  By Lemma \ref{le:3.2.5},  for any tuple of constants
$C$ from $G$, the following equivalence holds: $$ G \models
\forall X(U(X,C) = 1) \ \Longleftrightarrow \   G[X] \models
U(X,C) = 1. $$ Now it suffices to prove that for a given $U(X,Y)
\in G[X \cup Y]$ one can effectively find a word $W(Y) \in G[Y]$
such that for any tuple of constants $C$ over $F$ the following
equivalence holds $$ G[X] \models U(X,C) = 1 \Longleftrightarrow G
\models W(C) = 1. $$ We do this by induction on the syllable
length of $U(X,Y)$  which comes from the free product \newline
$G[X \cup Y] = G[Y] \ast F(X)$ (notice that $F(X)$ does not
contain constants from $G$, but $G[Y]$ does). If $U(X,Y)$ is of
the  syllable length  1, then either $U(X,Y) = U(X) \in F(X)$ or
$U(X,Y) = U(Y) \in G[Y]$. In the first event $F \models U(X) = 1$
means exactly that the reduced form of $U(X)$ is trivial, so we
can take $W(Y)$ trivial also.  In the event $U(X,Y) = U(Y)$ we can
take $W(Y) = U(Y)$.

Suppose now that $U(X,Y) \in G[Y] \ast F(X)$ and it has the
following reduced form: $$ U(X,Y) = g_1(Y)v_1(X)g_2(Y)v_2(X)
\ldots v_m(X)g_{m+1}(Y) $$ where $v_i$'s are reduced nontrivial
words in $F(X)$ and $g_i(Y)$'s are reduced words in $G[Y]$ which
are all nontrivial except, possibly, $g_1(Y)$ and $g_{m+1}(Y)$.

If for a tuple of constants $C$ over $G$ we have $ G[X] \models
U(X,C) = 1$ then at least one of the elements $g_2(C),
\ldots,g_m(C)$ must be trivial in $G$. This observation leads to
the following construction. For each $i = 2, \ldots, m$ delete the
subword $g_i(Y)$ from $U(X,Y)$ and reduce the new word to the
reduced form in the free product $F(X) \ast G[Y]$. Denote the
resulting word by $U_i(X,Y)$. Notice that the syllable length of
$U_i(X,Y)$ is less then the length of $U(X,Y)$. It follows from
the argument above that for any tuple of constants $C$ the
following equivalence holds: $$ G[X] \models U(X,C) = 1
\Longleftrightarrow G[X] \models \bigvee_{i=2}^{m} (g_i(C) = 1 \
\& \ U_i(X,C) = 1). $$ By induction one can effectively find words
$W_2(Y), \ldots, W_m(Y) \in G[Y]$ such that for any tuple of
constants $C$ we have $$ G[X] \models U_i(X,C) = 1
\Longleftrightarrow G \models W_i(C) = 1, $$ for each $i =
2,\ldots,m$. Combining the equivalences above we see that $$ G[X]
\models U(X,C) = 1 \Longleftrightarrow G \models \bigvee_{i=2}^{m}
(g_i(C) = 1\  \&\  W_i(C) = 1).
$$ By Corollaries \ref{co:2} and \ref{co:4} from the previous section we can
effectively rewrite the disjunction $$\bigvee_{i=2}^{m} (g_i(Y) =
1 \ \& \ W_i(Y) = 1)$$ as a single equation $W(Y) = 1$. That
finishes the proof.

\subsection{Positive and general formulas}

In this section we describe normal forms of general formulas and
positive formulas.  We show that every positive formula is
equivalent modulo $A_{\mathcal T}$ to a formula which consists of
an equation and a string of quantifiers in front of it; and for an
arbitrary formula $\Phi$ either $\Phi$ or $\neg\Phi$ is equivalent
modulo $A_{\mathcal T}$ to a formula in a general radical form (it
is a radical formula with a string of quantifiers in front of it).

\begin{lm}
\label{le:genpos} Every positive formula $\Phi(X)$ in $L_A$ is
equivalent modulo $A_{\mathcal T}$ to a formula of the type
 $$  Q_1X_1  \ldots Q_k X_k(S(X,X_1, \ldots , X_k) = 1),$$
 where  $Q_i \in \{\exists, \forall \} \ (i = 1, \ldots,k).$
 \end{lm}
 {\it Proof}. The result follows immediately from Corollaries \ref{co:2}
and
 \ref{co:4}.

 \begin{lm}
\label{le:impform} Let $\Phi(X)$ be a  formula in $L_A$ of the
form
$$\Phi(X) = Q_1X_1  \ldots Q_k X_k  \forall
   Y \Phi_0(X,X_1,  \ldots,X_k, Y),$$
where $Q_i \in \{\exists, \forall \}$ and  $\Phi_0$ is a
quantifier-free formula. Then one can effectively find a formula
$\Psi(X)$ of the form
 $$
\Psi(X)  = Q_1X_1  \ldots Q_k X_k \forall
   Y  \forall Z(S(X,
X_1,\ldots ,X_k,Y,Z) =1 \rightarrow T(Z)=1) $$ such that $\Phi(X)$
is equivalent to $\Psi(X)$ modulo $A_{\mathcal T}$.
 \end{lm}
  {\it Proof}.  Let
   $$\Phi(X) = Q_1X_1  \ldots Q_k X_k \forall  Y\Phi_0(X,X_1,
\ldots,X_k, Y),$$
  where $Q_i \in \{\exists, \forall \}$ and  $\Phi_0$ is a quantifier-
free formula.
   By  Lemma \ref{le:unrad} there exists  equations
   $S(X,X_1, \ldots ,X_k,Y,Z) = 1$ and $T(Z) = 1$ such that
  $$ \forall Y\Phi_0(X,X_1,  \ldots,X, Y) \sim_{A_{\mathcal T}} \forall Y
\forall Z
 (S(X_1,\ldots X_k,Y,Z) = 1 \
   \rightarrow   \  T(Z) = 1).$$
 It follows that
 $$\Phi(X) =  Q_1X_1  \ldots Q_k X_k \forall
   Y\Phi_0(X,X_1,  \ldots,X_k, Y) \sim_{A_{\mathcal T}}   Q_1X_1  \ldots Q_k
X_k \forall
   Y \forall Z(S(X, \ldots ,X_k,Y,Z) = 1
  \rightarrow T(Z) = 1,$$
 as desired.

\begin{lm}\label{fos}
For any formula $\Phi(X)$ in the language $L_A$, one can
effectively find a formula $\Psi(X)$ in the language $L_A$ in the
 form
 $$ \Psi(X) = \exists
X_1\forall Y_1\ldots\exists X_k \forall Y_k \forall
Z(S(X,X_1,Y_1,\ldots ,X_k,Y_k,Z) =1 \rightarrow T(Z)=1), $$ such
that $\Phi(X)$ or its negation $\neg \Phi(X)$ (and we can check
effectively which one of them) is equivalent to $ \Psi(X)$ modulo
$A_{\mathcal T}$.
\end{lm}

{\it Proof}. For any formula  $\Phi(X)$  in the language $L_A$ one
can effectively find a  disjunctive normal form $\Phi_1(X)$  of
$\Phi(X)$, as well as a disjunctive normal form $\Phi_2$ of the
negation $\neg \Phi(X)$ of $\Phi(X)$ (see, for example,
\cite{CK}).  We can assume that either  in $\Phi_1(X)$ or in
$\Phi_2(X)$ the quantifier prefix  ends with a universal
quantifier.  Moreover, adding (if necessary) an existential
quantifier $\exists v$ in front of the formula (where $v$ does not
occur in the formula) we may also assume that the formula begins
with an existential quantifier.  Now by Lemma \ref{le:impform} one
can effectively find a formula $\Psi$ with the required
conditions. $\Box$

\section{Generalized equations and positive theory of free groups}

Makanin \cite{Mak82} introduced the concept of a generalized
equation constructed for a finite system of equations in a free
group $F=F(A)$.   Geometrically a generalized equation consists of
three kinds of objects: bases, boundaries and items. Roughly it is
 a long interval with marked division points. The marked division
points are the boundaries. Subintervals between division points
are items (we assign a  variable to each item). Line segments
below certain subintervals, beginning at some boundary and ending
at some other boundary, are bases.  Each base either corresponds
to a letter from $A$ or has a double.

\begin{figure}[here]
\centering{\mbox{\psfig{figure=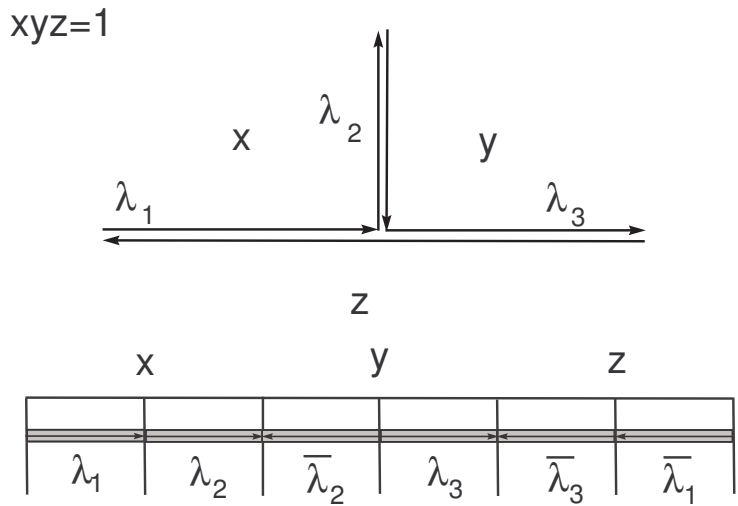,height=3in}}} \caption{From
the cancellation tree for the equation $xyz=1$ to the generalized
equation ($x=\lambda _1\circ\lambda _2,\ y=\lambda
_2^{-1}\circ\lambda _3,\ z=\lambda _3^{-1}\circ\lambda _1^{-1}).$}
\label{ET0}
\end{figure}
This concept becomes crucial to our subsequent work and is
difficult to understand. This is one of the main tools used to
describe solution sets of systems of equations. In subsequent
papers we will use it also to obtain effectively different
splittings of groups. Before we give a formal definition we will
try to motivate it with a simple example.

Suppose we have the simple equation $xyz = 1$ in a free group.
Suppose that we have a solution to this equation denoted by
$x^\phi,y^\phi,z^\phi$ where is $\phi$ is a given homomorphism
into a free group $F(A)$.  Since $x^\phi,y^\phi,z^\phi$ are
reduced words in the generators $A$ there must be complete
cancellation. If we take a concatenation of the geodesic subpaths
corresponding to $x^{\phi}, y^{\phi}$ and $z^{\phi}$ we obtain a
path in the Cayley graph corresponding to this complete
cancellation. This is called a cancellation tree (see Fig. 1). In
the simplest situation $x = \lambda_1\circ \lambda_2$,$y =
\lambda_2^{-1}\circ \lambda_3$ and $z = \lambda_3^{-1}\circ
\lambda_1^{-1}$.  The generalized equation would then be the
following interval.

 The
boundaries would be the division points, the bases are the
$\lambda's$ and the items in this simple case are also the
$\lambda's$.  In a more complicated equation where the variables
$X,Y,Z$ appear more than one time this basic interval would be
extended, Since the solution of any equation in a free group must
involve complete cancellation this drawing of the interval is
essentially the way one would solve such an equation. Our picture
above depended on one fixed solution $\phi$. However for any
equation there are only finitely many such cancellation trees and
hence only finitely many generalized equations.

\subsection{Generalized equations}
\label{se:4-1}

 Let $A= \{a_1, \ldots, a_m\}$ be a set of constants and $X = \{x_1,
\ldots, x_n\}$
  be a set of variables. Put   $G = F(A)$ and  $G[X] = G \ast F(X).$

\begin{df}
A combinatorial generalized equation $\Omega$ (with constants from
$A^{\pm 1}$) consists  of the following objects:

1. A finite set of {\bf bases} $BS = BS(\Omega)$.  Every base is
either a constant base or a variable base. Each constant base is
associated with exactly one letter from $A^{\pm 1}$. The set of
variable bases ${\mathcal M}$ consists of $2n$ elements ${\mathcal
M} = \{\mu_1, \ldots, \mu_{2n}\}$. The set ${\mathcal M}$ comes
equipped with two functions: a function $\varepsilon: {\mathcal M}
\rightarrow \{1,-1\}$ and an involution $\Delta: {\mathcal M}
\rightarrow {\mathcal M}$ (i.e., $\Delta$ is a bijection such that
$\Delta^2$ is an identity on  ${\mathcal M}$). Bases $\mu$ and
$\Delta(\mu)$ (or $\bar\mu$) are called {\it dual bases}.  We
denote variable bases by $\mu, \lambda, \ldots.$

2.  A set of {\bf boundaries} $BD = BD(\Omega)$. $BD$ is  a finite
initial segment of the set of positive integers  $BD = \{1, 2,
\ldots, \rho+1\}$. We use letters $i,j, \ldots$ for boundaries.

3. Two functions $\alpha : BS \rightarrow BD$ and $\beta : BS
\rightarrow BD$. We call $\alpha(\mu)$ and $\beta(\mu)$ the
initial and terminal boundaries of the base $\mu$ (or endpoints of
$\mu$). These functions satisfy the following conditions:
$\alpha(b) <  \beta(b)$  for every base $b \in BS$; if $b$ is a
constant base then $\beta(b) = \alpha(b) + 1$.

4. A finite set of {\bf boundary connections} $BC = BC(\Omega)$. A
boundary connection is a triple $(i,\mu,j)$ where $i, j \in BD$,
$\mu \in {\mathcal M}$ such that $\alpha(\mu) <  i < \beta(\mu)$
and $\alpha(\Delta(\mu)) <  j < \beta(\Delta(\mu))$. We will
assume for simplicity,  that if $(i,\mu,j) \in BC$ then
$(j,\Delta(\mu), i) \in BC$. This allows one to identify
connections $ (i,\mu,j)$ and  $(j,\Delta(\mu), i)$.
\end{df}

   For a combinatorial generalized equation $\Omega$, one can canonically
associate a system of equations in {\bf variables} $h_1, \ldots,
h_\rho$ over $F(A)$ (variables $h_i$ are sometimes  called {\it
items}). This system is called a {\bf generalized equation}, and
(slightly abusing the language) we  denote it by the same symbol
$\Omega$. The generalized equation  $\Omega$  consists of the
following three types of equations.

1. Each pair of dual variable bases $(\lambda, \Delta(\lambda))$
provides an equation
 $$[h_{\alpha
(\lambda )}h_{\alpha (\lambda )+1}\cdots h_{\beta (\lambda )-1}]^
{\varepsilon (\lambda)}= [h_{\alpha (\Delta (\lambda ))}h_{\alpha
(\Delta (\lambda ))+1} \cdots h_{\beta (\Delta (\lambda ))-1}]^
{\varepsilon (\Delta (\lambda))}.$$ These equations are called
{\bf basic equations}.

2. For each constant base $b$ we write down a {\bf coefficient
equation}
 $$ h_{\alpha(b)} = a,$$
 where $a \in A^{\pm 1}$ is the constant associated with $b$.

3.  Every boundary connection $(p,\lambda,q)$ gives rise to a {\bf
boundary equation}
 $$[h_{\alpha (\lambda )}h_{\alpha (\lambda
)+1}\cdots h_{p-1}]= [h_{\alpha (\Delta (\lambda ))}h_{\alpha
(\Delta (\lambda ))+1} \cdots h_{q-1}],$$ if $\varepsilon
(\lambda)= \varepsilon (\Delta(\lambda))$ and $$[h_{\alpha
(\lambda )}h_{\alpha (\lambda )+1}\cdots h_{p-1}]= [h_{q}h_{q+1}
\cdots h_{\beta (\Delta (\lambda ))-1}]^{-1} ,$$ if $\varepsilon
(\lambda)= -\varepsilon (\Delta(\lambda)).$

\begin{rk}
We  assume that every generalized equation comes associated with a
combinatorial one;
\end{rk}

{\bf Example.} Consider as an example the Malcev equation
$[x,y][b,a]=1$, where $a,b\in A.$ Consider the following solution
of this equation:
$$  x^{\phi}=(( b^{n_1}a)^{n_2}b)^{n_3}b^{n_1}a, \ \ y^{\phi}=(
b^{n_1}a)^{n_2}b.$$ Fig. 2 shows the cancellation tree and the
generalized equation for this solution.
\begin{figure}[here]
\centering{\mbox{\psfig{figure=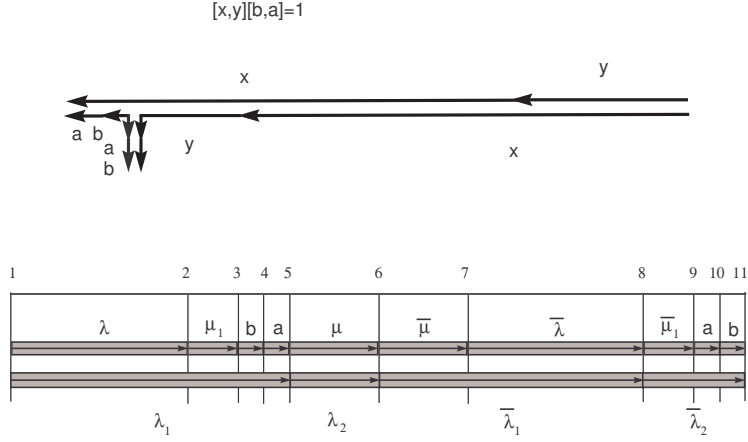,height=3.3in}}} \caption{A
cancellation tree and the generalized equation corresponding to
this tree for the equation $[x,y][b,a]=1.$} \label{ET00}
\end{figure}
This generalized equation has ten variables $h_1,\ldots ,h_{10}$
and eleven boundaries. The system of basic equations for this
generalized equation is the following
$$h_1=h_7, \ h_2=h_8,\ h_5=h_6,\ h_1h_2h_3h_4=h_6h_7,\
h_5=h_8h_9h_{10}.$$ The system of coefficient equations is
$$h_3=b,\ h_4=a,\ h_9=a,\ h_{10}=b.$$

\begin{df}
Let $\Omega(h) = \{L_1(h)=R_1(h), \ldots, L_s(h) = R_s(h)\}$ be a
generalized equation in variables $h = (h_1, \ldots,h_{\rho})$
with constants from $A^{\pm 1 }$. A sequence of reduced   nonempty
words $U = (U_1(A), \ldots, U_{\rho}(A))$  in the alphabet $A^{\pm
1}$ is a {\em solution} of $\Omega $ if:
\begin{enumerate}
\item [1)] all words $L_i(U), R_i(U)$  are reduced as written;

\item [2)] $L_i(U) =  R_i(U),  \ \ i = 1, \ldots s.$
\end{enumerate}\end{df} The notation $(\Omega, U)$ means that $U$ is a
solution of the generalized equation $\Omega $.

\begin{rk}
Notice that a solution $U$ of a generalized equation $\Omega$ can
be viewed as  a solution of $\Omega$ in the free monoid
$F_{mon}(A^{\pm 1})$ (i.e., the equalities $L_i(U) = R_i(U)$ are
graphical) which satisfies an additional  condition  $ U \in F(A)
\leq F_{mon}(A^{\pm 1})$.
\end{rk}

 Obviously, each solution  $U$ of $\Omega$ gives rise
to  a solution of $\Omega$ in the free group $F(A)$. The converse
does not hold in general, i.e., it might happen that  $U$ is a
solution of $\Omega$ in $F(A)$ but not in $F_{mon}(A^{\pm 1})$,
i.e., all equalities $L_i(U) = R_i(U)$ hold only after a free
reduction but not graphically. We introduce the following notation
which will allow us  to distinguish in which structure
($F_{mon}(A^{\pm 1})$ or $F(A)$) we are looking for solutions  for
$\Omega$.

 If
 $$S =  \{L_1(h)=R_1(h), \ldots, L_s(h) =
R_s(h)\}$$ is an arbitrary system of equations with constants from
$A^{\pm 1}$, then by $S^*$ we denote the system of equations $$S^*
= \{L_1(h)R_1(h)^{-1} = 1, \ldots, L_s(h)R_s(h)^{-1} = 1\}$$  over
the free group $F(A)$.

\begin{df}
A generalized equation $\Omega$ is called {\em formally
consistent} if it satisfies the following conditions. \be \item
[1)] If $\varepsilon (\mu)=-\varepsilon (\Delta (\mu))$, then the
bases $\mu$ and $\Delta (\mu )$ do not intersect, i.e. non of the
the items $h_{\alpha(\mu)}, h_{\beta(\mu)-1}$ is contained in
$\Delta (\mu )$. \item [2)] If two boundary equations have
respective parameters $(p,\lambda ,q)$ and $(p_1,\lambda ,q_1)$
with $p\leq p_1,$ then $q\leq q_1$ in the case when $\varepsilon
(\lambda)\varepsilon (\Delta (\lambda))=1,$ and $q\geq q_1$ in the
case $\varepsilon (\lambda)\varepsilon (\Delta (\lambda))=-1,$ in
particular, if $p=p_1$ then $q = q_1$.
 \item [3)] Let $\mu$ be
a base such that  $\alpha (\mu)=\alpha (\Delta (\mu))$ (in this
case we say that bases $\mu$ and $\Delta(\mu)$ form a matched pair
of dual bases). If $(p,\mu ,q)$ is  a boundary connection related
to $\mu$ then  $p=q.$

\item  [4)]A variable cannot occur in two distinct coefficient
equations, i.e., any two constant bases with the same left
end-point are labelled by the same letter from $A^{\pm 1}$.

\item [5)]If $h_i$ is a variable from some coefficient equation,
and if $(i,\mu,q_1), (i+1,\mu,q_2)$ are boundary connections, then
$|q_1- q_2|=1.$ \ee
\end{df}

\begin{lm}
\label{le:4.1}\begin{enumerate} \item If a generalized equation
$\Omega$ has a solution then $\Omega$ is formally consistent;

\item There is an algorithm which for every generalized equation
checks whether it is formally consistent or not.\end{enumerate}
\end{lm}
The proof is easy and we omit it.
\smallskip

\begin{rk} In the sequel we  consider only formally consistent
generalized equations.
\end{rk}

It is convenient to  visualize  a generalized equation $\Omega$ as
follows.

\begin {center}
\begin{picture}(200,100)(0,0)
\put(10,85){\line(0,-1){80}} \put(30,85){\line(0,-1){80}}
\put(50,85){\line(0,-1){80}} \put(70,85){\line(0,-1){80}}
\put(90,85){\line(0,-1){80}} \put(110,85){\line(0,-1){80}}
\put(130,85){\line(0,-1){80}} \put(150,85){\line(0,-1){80}}
\put(170,85){\line(0,-1){80}} \put(10,80){\line(1,0){160}} \put
(10,87){$1$} \put (30,87){$2$} \put (50,87){$3$} \put
(140,87){$\rho -1$} \put (170,87){$\rho$} \put
(10,65){\line(1,0){60}} \put (70,45){\line(1,0){20}} \put
(90,55){\line(1,0){40}} \put (110,30){\line(1,0){40}} \put
(20,67){$\lambda$} \put (92,57){$\Delta(\lambda )$} \put
(80,47){$\mu$} \put (112,37){$\Delta(\mu)$}
\end{picture}
\end{center}

\subsection{Reduction to generalized equations}

In this section, following Makanin \cite{Mak82},   we show how for
a given finite system of equations $S(X,A) = 1$ over a free group
$F(A)$ one can canonically associate a finite collection  of
generalized equations ${\mathcal GE}(S)$ with constants from
$A^{\pm 1}$, which to some extent describes all solutions of the
system $S(X,A) = 1$.

 Let $S(X,A) = 1$ be a finite system of equations $S_ 1 = 1, \ldots, S_m
= 1$
  over a free group $F(A)$.  We write $S(X,A) = 1$  in the form
\begin{equation}\label{*}\begin{array}{c}
 r_{11}r_{12}\ldots r_{1l_1}=1,\\
 r_{21}r_{22}\ldots r_{2l_2}=1,\\
 \ldots \\
 r_{m1}r_{m2}\ldots r_{ml_m}=1,\\ \end{array} \end{equation} where $r_{ij}$ are letters in
the alphabet $ X^{\pm 1}\cup A^{\pm 1}.$

A {\it partition table} $T$  of the system above  is  a set of
reduced words
$$T = \{V_{ij}(z_1, \ldots ,z_p)\} \ \ (1\leq i\leq m, 1\leq j\leq l_i)$$ from a free
group $F[Z] = F(A \cup Z)$, where $Z = \{z_1,\ldots ,z_p\}$,
which satisfies the following conditions:
\begin{enumerate}
\item [1)] The equality $V_{i1}V_{i2} \ldots V_{il_i}=1, 1\leq
i\leq m,$ holds in  $F[Z]$; \item [2)] $|V_{ij}|\leq l_i - 1$;
\item [3)] if $r_{ij} = a \in A^{\pm 1}$, then $V_{ij}=
a.$\end{enumerate}

Since $|V_{ij}|\leq l_i - 1$ then  at most $|S| = \sum_{i = 1}^m
(l_i - 1)l_i$ different  letters $z_i$ can occur in a partition
table of the equation $S(X,A) = 1$. Therefore we will always
assume  that $p \leq |S|$.

Each partition table encodes a particular type of cancellation
that happens when one substitutes  a particular solution $W(A) \in
F(A)$ into  $S(X,A) = 1$ and then freely reduces the words in
$S(W(A),A)$ into the empty word.

\begin{lm}
\label{le:4.2} Let $S(X,A) = 1$ be a finite system of equations
over $F(A)$. Then

1) the set $PT(S)$ of all partition tables  of $S(X,A) = 1$ is
finite, and its cardinality is bounded by a number which depends
only on  $|S(X,A)|$;

2) one can effectively enumerate the set $PT(S)$.
\end{lm}
{\em Proof.}   Since the words $V_{ij}$ have  bounded length,  one
can effectively enumerate the finite set  of all collections of
words $\{V_{ij}\}$ in $F[Z]$ which satisfy the conditions 2), 3)
above.  Now for each such collection $\{V_{ij}\}$,  one can
effectively check  whether the  equalities $V_{i1}V_{i2} \ldots
V_{il_i}=1, 1\leq i\leq m$ hold in the free group $F[Z]$ or not.
This allows one to list effectively all partition tables for
$S(X,A) = 1$. $\Box$

To each partition table $T=\{V_{ij}\}$ one can assign a
generalized equation $\Omega _T$ in the following way (below we
use the notation $\doteq $ for graphical equality). Consider the
following word $V$ in $M(A^{\pm 1} \cup Z^{\pm 1}):$
 $$V\doteq V_{11}V_{12}\ldots V_{1l_1}\ldots V_{m1}V_{m2} \ldots  V_{ml_m} = y_1
 \ldots y_\rho, $$
where $y_i \in A^{\pm 1} \cup Z^{\pm 1}$ and $\rho = l(V)$ is the
length of $V$. Then the generalized equation $\Omega_T =
\Omega_T(h)$ has $\rho + 1$ boundaries and $\rho$ variables
$h_1,\ldots ,h_{\rho}$ which are denoted by   $h = (h_1,\ldots
,h_{\rho})$.

Now we define bases of $\Omega_T$ and the functions $\alpha,
\beta, \varepsilon$.

 Let $z \in Z$.  For any two distinct occurrences of  $z$ in $V$ as
 $$y_i = z^{\varepsilon _i}, \ \ \ y_j = z^{\varepsilon _j} \ \ \
 (\varepsilon _i, \varepsilon _j \in \{1,-1\})$$
 we introduce a pair of dual variable bases $\mu_{z,i}, \mu_{z,j}$
such that $\Delta(\mu_{z,i}) = \mu_{z,j}$ (say, if $i < j$). Put
$$\alpha(\mu_{z,i}) = i, \ \ \ \beta(\mu_{z,i}) = i+1, \ \ \
\epsilon(\mu_{z,i}) = \varepsilon _i.$$
 The basic equation that
corresponds to this pair of dual bases is
 $h_{i}^{\varepsilon_i}=h_{j}^{\varepsilon _j} .$

Let $x \in X$.  For any two distinct occurrences of  $x$ in
$S(X,A) = 1$  as
 $$r_{i,j} = x^{\varepsilon_{ij}}, \ \ \ r_{s,t} =
x^{\varepsilon_{st}} \ \ \ (\varepsilon _{ij}, \varepsilon _{st}
\in \{1,-1\})$$
 we introduce a pair of dual bases $\mu_{x,i,j}$ and $\mu_{x,s,t}$
 such that $\Delta(\mu_{x,i,j}) = \mu_{x,s,t}$ (say, if $(i,j) < (s,t)$
 in the left lexicographic order).  Now let $V_{ij}$ occurs in the
 word $V$ as a subword
 $$V_{ij} = y_c \ldots y_d.$$
 Then we put
$$\alpha(\mu_{x,i,j}) = c, \ \ \ \beta(\mu_{x,i,j}) = d+1, \ \ \
\epsilon(\mu_{x,i,j}) = \varepsilon_{ij}.$$
 The basic equation which
corresponds to these dual bases can be written in the form
 $$[h_{\alpha(\mu_{x,i,j})}\ldots h_{\beta(\mu_{x,i,j})-1}]^{\varepsilon _{ij}}=
 [h_{\alpha(\mu_{x,s,t})}\ldots h_{\beta(\mu_{x,s,t})-
1}]^{\varepsilon _{st}}.$$

Let $r_{ij} = a \in A^{\pm 1}$.  In this case we introduce a
constant base $\mu_{ij}$ with the label $a$. If $V_{ij}$ occurs in
$V$ as $V_{ij} = y_c$, then we put
 $$ \alpha(\mu_{ij}) = c, \beta(\mu_{ij}) = c+1.$$
 The corresponding coefficient equation is written as  $h_{c}=a$.

The list of boundary connections here (and hence the boundary
equations) is empty. This defines the generalized equation
$\Omega_T$. Put
 $${\mathcal GE}(S) = \{\Omega_T \mid T \ is \  a \ partition \ table \
 for\ S(X,A)= 1 \}.$$
  Then ${\mathcal GE}(S)$ is a finite collection of generalized
  equations which can be effectively constructed for a given
  $S(X,A) = 1$.

  For a generalized equation $\Omega $ we can also  consider the same system of equations in a free group.
    We denote this system by $\Omega ^*$. By $F_{R(\Omega)}$ we denote the coordinate group of $\Omega ^*.$
    Now we explain relations between  the coordinate groups of
  $S(X,A) = 1$ and $\Omega_T^*$.

For a letter $x$ in $X$ we choose an arbitrary   occurrence of $x$
in $S(X,A) = 1$ as
 $$r_{ij} = x^{\varepsilon_{ij}}.$$
  Let $\mu = \mu_{x,i,j}$ be the base that corresponds to this occurrence
  of $x$.  Then $V_{ij}$ occurs  in $V$ as the subword
  $$V_{ij} = y_{\alpha(\mu)} \ldots y_{\beta(\mu) -1}.$$
 Define a word $P_x(h) \in F[h]$ (where $h = \{h_1, \ldots,h_\rho\}$) as
 $$ P_x(h,A) = h_{\alpha(\mu)} \ldots h_{\beta(\mu)-1}^{\varepsilon_{ij}},$$
 and put
 $$P(h) = (P_{x_1}, \ldots, P_{x_n}).$$
 The tuple of words $P(h)$ depends on a choice of occurrences of letters from
  $X$ in $V$. It follows from the construction above that  the map $ X
 \rightarrow F[h]$ defined by $x \rightarrow  P_x(h,A)$
 gives rise to an $F$-homomorphism
 $$\pi : F_{R(S)}\rightarrow F_{R(\Omega _T)}.$$
Observe that the image  $\pi (x)$ in $F_{R(\Omega _T)}$ does not
depend on a particular choice of the occurrence of $x$ in $S(X,A)$
(the basic equations of $\Omega_T$ make these images equal). Hence
$\pi$ depends only on $\Omega_T$.

  Now we relate solutions of $S(X,A) = 1$  with solutions of generalized
  equations  from ${\mathcal GE}(S)$.
  Let $W(A)$ be a solution of $S(X,A) = 1$ in $F(A)$. If in the system
  (\ref{*}) we make the substitution  $\sigma : X \rightarrow W(A)$,
   then
   $$(r_{i1}r_{i2}\ldots r_{il_i})^{\sigma} =
   r_{i1}^{\sigma}r_{i2}^{\sigma}\ldots r_{il_i}^{\sigma} = 1$$
    in $F(A)$ for  every $i = 1, \ldots, m$.
Hence every product $R_i = r_{i1}^{\sigma}r_{i2}^{\sigma}\ldots
r_{il_i}^{\sigma}$   can be reduced to the empty word by a
sequence of free   reductions. Let us fix a particular reduction
process for  each $R_i$.   Denote by ${\tilde z}_1, \ldots,
{\tilde z}_p$ all the (maximal) non-trivial subwords of
$r_{ij}^{\sigma}$  that cancel out in some $R_i$ ($i = 1, \ldots
,m$) during the chosen reduction process. Since every word
$r_{ij}^\sigma$ in this process  cancels out completely, ithat
mplies that
$$r_{ij}^\sigma = V_{ij}({\tilde z}_1, \ldots,  {\tilde z}_p)$$
 for some reduced words $V_{ij}(Z)$ in variables $Z = \{z_1, \ldots, z_p\}$.
  Moreover, the equality above is graphical. Observe also that if
  $r_{ij} = a \in A^{\pm 1}$ then $r_{ij}^{\sigma} = a$ and we have
 $V_{ij} = a$. Since every word $r_{ij}^\sigma$
 in $R_i$ has at most one cancellation with any other word
 $r_{ik}^\sigma$
 and does not have cancellation with itself, we have $l(V_{ij}) \leq l_i - 1$.
  This shows that
   the set $T = \{V_{ij}\}$ is  a partition
table for $S(X,A) = 1$. Obviously,
$$U(A) = ({\tilde z}_1, \ldots, {\tilde z}_p)$$
is the  solution of the generalized equation $\Omega_T$, which is
induced by  $W(A)$. From the construction of the map $P(H)$ we
deduce that $W(A) = P(U(A))$.

The reverse is also true: if $U(A)$ is an arbitrary  solution of
the generalized equation $\Omega_T$, then $P(U(A))$ is a solution
of $S(X,A) = 1$.

We  summarize the discussion above in the following lemma, which
is essentially due to Makanin \cite{Mak82}.

\begin{lm}
\label{le:R1} For a given  system of equations $S(X,A)=1$ over  a
free group $F = F(A)$,  one can  effectively  construct a finite
set
$${\mathcal GE}(S) = \{\Omega_T \mid T \ is \  a \ partition \ table \
 for\ S(X,A)= 1 \}$$
 of generalized equations  such  that
\begin{enumerate}
\item If the set ${\mathcal GE}(S)$ is empty, then  $S(X,A)= 1$
has no solutions in $F(A)$; \item   for each $\Omega (H) \in
{\mathcal GE}(S)$ and  for each $x \in X$  one can effectively
find  a word $P_x(H,A) \in F[H]$   of length at most $|H|$ such
that the map $x: \rightarrow P_x(H,A)$ ($x \in X$) gives rise to
an $F$-homomorphism $\pi_\Omega: F_{R(S)}\rightarrow
F_{R(\Omega)}$; \item for any solution $W(A)  \in F(A)^n$ of the
system $S(X,A)=1$ there exists $\Omega (H) \in {\mathcal GE}(S)$
and a solution $ U(A)$ of $\Omega(H)$ such that $W(A) = P(U(A))$,
where $P(H) = (P_{x_1}, \ldots, P_{x_n})$,
 and this equality is graphical;
\item for any $F$-group $\tilde F$, if a  generalized equation
$\Omega (H) \in {\mathcal GE}(S)$
 has a solution $\tilde U$ in $\tilde F$, then $P(\tilde U)$
is a solution of $S(X,A) = 1$ in $\tilde F$.
\end{enumerate}
\end{lm}
\begin{cy}
\label{co:R1} In the notations of Lemma \ref{le:R1}  for any
solution $W(A)  \in F(A)^n$ of the system $S(X,A)=1$ there exists
$\Omega (H) \in {\mathcal GE}(S)$  and a solution $ U(A)$ of
$\Omega(H)$ such that the following diagram commutes.
\end{cy}
\medskip

\begin{center}

\begin{picture}(100,100)(0,0)
\put(0,100){$F_{R(S)}$} \put(100,100){$F_{R(\Omega )}$}
\put(0,0){$F$} \put(15,103){\vector(1,0){80}}
\put(5,93){\vector(0,-1){78}} \put(95,95){\vector(-1,-1){80}}
\put(-10,50){$\pi_{W}$} \put(55,45){$\pi _{ U}$}
\put(50,108){$\pi$}
\end{picture}

\end{center}

\medskip

\subsection{Generalized equations with parameters}
\label{se:parametric}

In this section, following \cite{Razborov3} and \cite{KMIrc},  we
consider generalized equations with {\it parameters}. This kind of
equations appear naturally in Makanin's type rewriting processes
and provide a convenient tool to organize induction properly.

Let $\Omega$ be a generalized equation. An  item $h_i$ {\it
belongs} to a base $\mu$ (and, in this event, $\mu$ {\it contains}
$h_i$) if $\alpha (\mu)\leq i\leq \beta (\mu)-1.$ An item $h_i$ is
{\it constant} if it belongs to a constant base, $h_i$ is {\it
free } if it does not belong to any base. By $\gamma(h_i) =
\gamma_i$ we denote the number of bases which contain $h_i$. We
call $\gamma_i$ the {\it degree } of $h_i$.

A boundary $i$ {\it crosses} (or {\it intersects})  the base $\mu$
if $\alpha (\mu)<i<\beta (\mu).$ A boundary $i$ {\it touches} the
base $\mu$ (or $i$ is an end-point of $\mu$)  if $i=\alpha (\mu)$
or $i=\beta (\mu)$. A boundary is said to be {\em open} if it
crosses at least one base, otherwise it is called {\it closed}. We
say that a boundary $i$ is {\it tied} (or {\it bound}) by a base
$\mu$ (or {\it $\mu$-tied})  if there exists a boundary connection
$(p,\mu,q)$ such that $i = p$ or $i = q$. A boundary is {\em free}
if it does not touch any base and it is not tied by a boundary
connection.

 A set of consecutive  items $[i,j] =
\{h_i,\ldots, h_{i+j-1}\}$   is called a {\em section}.
 A section
is said to be {\em closed} if the boundaries $i$ and $i+j$ are
closed and all the boundaries between them are open. A base $\mu$
is {\it contained } in a base $\lambda$ if $\alpha(\lambda) \leq
\alpha(\mu) < \beta(\mu) \leq \beta(\lambda)$. If $\mu$ is a base
then by $\sigma(\mu)$ we denote the section
$[\alpha(\mu),\beta(\mu)]$ and by $h(\mu)$ we denote the product
of items $h_{\alpha(\mu)}\ldots h_{\beta(\mu)-1}$. In general for
a section $[i,j]$ by $h[i,j]$ we denote the product  $h_i \ldots
h_{j-1}$.

\begin{df}
\label{de:gepar} Let $\Omega$ be a generalized equation. If the
set $\Sigma = \Sigma_\Omega$ of all closed sections  of $\Omega$
is partitioned into a disjoint union of subsets
  \begin{equation}
  \label{eq:ge1-0}
\Sigma_\Omega = V\Sigma\cup  P\Sigma\cup C\Sigma ,
\end{equation}
then $\Omega $ is called a {\em generalized equation with
parameters} or a {\it parametric} generalized equation. Sections
from $V\Sigma, P\Sigma$, and $C\Sigma$ are called correspondingly,
{\em variable, parametric}, and {\em constant} sections. To
organize the branching process properly, we usually divide
variable sections into two disjoint parts:
 \begin{equation}
 \label{eq:geq1-1}
 V\Sigma = A\Sigma\cup NA\Sigma
\end{equation}
Sections from $A\Sigma$ are called {\em active}, and sections from
$NA\Sigma$ are {\em non-active}. In the case when partition
(\ref{eq:geq1-1}) is not specified we assume that $A\Sigma =
V\Sigma$. Thus, in general, we have a partition
  \begin{equation}
  \label{eq:ge1}
\Sigma_\Omega = A\Sigma\cup NA\Sigma\cup P\Sigma\cup C\Sigma
\end{equation}
If $\sigma \in \Sigma$, then every base or item from $\sigma$ is
called active, non-active, parametric, or constant, with respect
to the type of $\sigma$.
\end{df}

We will see later that every parametric generalized equation can
be written in a particular {\it standard} form.

\begin{df}
 We say that a parametric generalized equation $\Omega$ is {\em in a standard
 form}
 if the following conditions hold:
\begin{enumerate}
\item [1)]  all non-active sections from $NA\Sigma_\Omega$ are
located to the right of all active sections from $A\Sigma$, all
parametric sections from $P\Sigma_\Omega$ are located to
 the right  of all non-active sections, and all constant sections from $C\Sigma$
  are located to the right of all parametric sections; namely, there are
  numbers $1 \leq \rho_A  \leq  \rho_{NA} \leq \rho_P \leq \rho_C \leq \rho_ = \rho_
  \Omega$ such that $[1,\rho_A +1]$, $[\rho_A
  +1, \rho_{NA}+1]$, $[\rho_{NA}+1, \rho_P+1]$, and
  $[\rho_P +1, \rho_\Omega +1]$ are, correspondingly, unions of
  all active,  all non-active, all parametric, and
  all constant sections;
 \item [2)] for every letter $a \in A^{\pm 1}$ there is at most one
 constant base in $\Omega$ labelled by $a$, and  all such bases are located in the
 $C\Sigma$;
\item [3)] every free variable (item) $h_i$ of $\Omega$ is located
in $C\Sigma$.
 \end{enumerate}
 \end{df}

Now we describe a typical method for  constructing generalized
equations with parameters starting with a system of ordinary group
equations with constants from $A$.

 {\bf Parametric generalized equations
corresponding to group equations}

 Let
\begin{equation}
\label{eq:systS} S(X,Y_1,Y_2, \ldots, Y_k,A) = 1 \end{equation}
 be a finite system of equations  with constants from $A^{\pm 1}$ and with the
 set of variables partitioned into a  disjoint union
  \begin{equation}
  \label{eq:ge3}
   X \cup Y_1 \cup  \ldots \cup Y_k
   \end{equation}
 Denote by  ${\mathcal GE}(S)$ the  set
of generalized equations corresponding to $S = 1$ from Lemma
\ref{le:R1}. Put $Y = Y_1 \cup \ldots \cup Y_k$.   Let $\Omega \in
{\mathcal GE}(S)$.  Recall that every  base $\mu$ occurs in
$\Omega$
 either related to some occurrence of a variable  from $X \cup Y$ in the
 system $S(X,Y,A) = 1$, or
related to an occurrence of a letter $z \in Z$ in the word $V$
(see Lemma \ref{le:4.2}), or is a constant base. If $\mu$
corresponds to a variable $x \in X$ ($y \in Y_i$)  then we say
that $\mu$ is an {\it $X$-base} ({\it $Y_i$-base}). Sometimes we
refer to $Y_i$-bases as to $Y$ bases.  For a base $\mu$ of
$\Omega$ denote by $\sigma_\mu$ the section $ \sigma_\mu =
[\alpha(\mu),\beta(\mu)]$. Observe that the section $\sigma_\mu$
is closed in $\Omega$ for every $X$-base, or $Y$-base. If $\mu$ is
an $X$-base ($Y$-base or $Y_i$-base), then the section
$\sigma_\mu$ is called an {\it $X$-section} ({\it $Y$-section} or
{\it $Y_i$-section}). If $\mu$ is a constant base and the section
$\sigma_\mu$ is closed then we call $\sigma_\mu$ a {\it constant }
section. Using the derived transformation D2 we transport all
closed $Y_1$-sections to the right end of the generalized
equations behind all the sections of the equation (in an arbitrary
order), then we transport all $Y_2$-sections an put them behind
all $Y_1$-sections, and so on. Eventually, we transport  all
$Y$-sections to the very end of the interval and they appear there
with respect to the partition (\ref{eq:ge3}). After that we take
all the constant sections and put them behind all the parametric
sections. Now, let  $A\Sigma$ be  the set of all $X$-sections,
$NA\Sigma = \emptyset$, $P\Sigma$ be the set of all $Y$-sections,
and $C\Sigma$ be the set of all constant sections. This defines a
parametric generalized equation $\Omega = \Omega_Y$ with
parameters corresponding to the set of variables $Y$. If the
partition of variables (\ref{eq:ge3})
 is fixed we will omit $Y$ in the notation above and call $\Omega$ the
{\it parameterized} equation obtained from $\Omega$.
 Denote by
 $$ {\mathcal GE}_{par}(\Omega) = \{\Omega_{Y} \mid \Omega \in  {\mathcal
 GE}(\Omega)\}$$
the set of all parameterized equations of the system
(\ref{eq:systS}).

\subsection{Positive theory of free groups}
\label{se:mer} In this section we prove first  the Merzljakov's
result on elimination of quantifiers for positive sentences over
free group $F = F(A)$ \cite{Merz}.  This proof is based on the
notion of a generalized equation. Combining Merzljakov's theorem
with Makanin's result on decidability of equations over free
groups we obtain decidability of the positive theory of free
groups. This argument is due to Makanin \cite{Mak84}.

Recall that every positive formula $\Psi(Z)$ in the language $L_A$
is equivalent modulo $A_{\mathcal T}$ to a formula of the type $$
\forall x_1 \exists y_1 \ldots \forall x_k \exists y_k (S(X,Y,Z,A)
= 1),$$ where $S(X,Y,Z,A) = 1$ is an equation with constants from
$A^{\pm 1}$,  $X = (x_1, \ldots, x_k), Y = (y_1, \ldots,y_k), Z =
(z_1, \ldots, z_m)$.
 Indeed, one
can insert fictitious quantifiers to ensure the direct alteration
of quantifiers in the prefix.
 In particular, every positive sentence in $L_A$ is equivalent modulo
$A_{\mathcal T}$  to  a formula of the type
  $$\forall x_1 \exists y_1 \ldots \forall x_k \exists y_k (S(X,Y,A)= 1) .$$
 Now we prove the Merzlyakov's theorem from \cite{Merz}, though in a
 slightly different form.

\medskip
{\bf Merzljakov's Theorem.} {\it If
 $$F\models \forall x_1\exists y_1\ldots\forall x_k\exists y_k
(S(X,Y,A)=1),$$
 then there exist words (with constants from $F$)
$q_1(x_1),\ldots , q_k(x_1,\ldots ,x_k) \in F[X],$ such that
$$F[X]\models S(x_1, q_1(x_1),\ldots ,x_k, q_k(x_1,\ldots
,x_k,A))=1,$$
 i.e., the  equation
  $$S(x_1,y_1, \ldots, x_k,y_k,A) = 1$$
  (in variables  $Y$) has a solution in the free group $F[X]$.
}

\medskip
 {\em Proof.}   Let ${\mathcal GE}(u) = \{\Omega_1(Z_1),
\ldots, \Omega_r(Z_r)\}$ be generalized equations associated with
equation $S(X,Y,A) = 1$ in Lemma \ref{le:R1}. Denote by  $\rho_i =
|Z_i|$ the number of variables in $\Omega_i$.

Let $a, b \in A, [a,b] \neq 1$, and put
$$g_{1}=ba^{m_{11}}ba^{m_{12}}b\ldots a^{m_{1n_{1}}}b,$$ where
 $m_{11} < m_{12} < \ldots m_{1n_1}$ and
  $\max\{\rho_1, \ldots,\rho_r\}|S(X,A)|  < n_1 $.
Then there exists $h_1$ such that
 $$F \models \forall x_2 \exists y_2\ldots\forall x_k\exists y_k
 (S(g_1,h_1, x_2, y_2, \ldots, x_k,y_k)=1).$$
 Suppose now that elements $g_1, h_1, \ldots g_{i- 1},h_{i-1} \in
F$ are given. We define
\begin{equation}
\label{eq:two-stars}
 g_{i} = ba^{m_{i1}}ba^{m_{i2}}b\ldots
a^{m_{in_{i}}}b
\end{equation}

 such that:
\begin{enumerate}
\item [1)] $m_{i1} < m_{i2} < \ldots m_{in_i}$; \item [2)]
$\max\{\rho_1, \ldots,\rho_r\} |S(X,A)| < n_i $; \item [3)] no
subword of the type $ba^{m_{ij}}b$   occur in any of the words
$g_l, h_l$ for $l < i$.
\end{enumerate}
We call words \ref{eq:two-stars} {\em Merzljakov's words.} Then
there exists an element $h_i \in F$ such that
 $$F \models \forall x_{i+1}
\exists y_{i+1}\ldots\forall x_k\exists y_k
 (S(g_1,h_1, \ldots,  g_i, h_i,x_{i+1},y_{i+1},  \ldots, x_k,y_k) =
1).$$
  By induction we have constructed elements $g_1, h_1,
\ldots, g_k, h_k \in F$ such that
 $$S(g_1,h_1, \ldots,g_k,h_k) = 1$$
and each $g_i$ has the form (\ref{eq:two-stars}) and satisfies the
conditions 1), 2), 3).

By Lemma \ref{le:R1} there exists a generalized equation
$\Omega(Z) \in {\mathcal GE}(S)$,  words $P_{i}(Z,A), Q_{i}(Z,A)
\in F[Z]$ ($i = 1, \ldots,k$) of length not more then $\rho =
|Z|$, and a solution $U = (u_1, \ldots, u_\rho)$ of $\Omega(Z)$ in
$F$ such that the following words are graphically equal:
 $$g_i = P_{i}(U),  \ \ \  h_i = Q_{i}(U)   \ \ \ (i = 1,
\ldots,k).$$
 Since $n_i > \rho|S(X,A)|$ (by condition 2)) and $P_i(U) = y_1 \ldots y_q$ with
 $y_i \in U^{\pm 1}, q \leq \rho$,  the graphical  equalities
 \begin{equation}
 \label{eq:gp}
  g_i =  ba^{m_{i1}}ba^{m_{i2}}b\ldots
a^{m_{in_{i}}}b = P_{i}(U) \ \ (i = 1, \ldots,k)
  \end{equation}
  show that  there exists a subword $v_i = ba^{m_{ij}}b$ of $g_i$  such
that every occurrence of this subword in (\ref{eq:gp})  is an
occurrence inside some $u_j^{\pm 1}$.    For each $i$ fix  such a
subword  $v_i = ba^{m_{ij}}b$ in $g_i$.  In view of condition 3)
the word $v_i$ does not occur in any of the words $g_j$ ($j \neq
i$), $h_s$ ($s < i$), moreover, in $g_i$ it occurs precisely once.
Denote by $j(i)$ the unique index such that $v_i$ occurs inside
$u_{j(i)}^{\pm 1}$ in $P_i(U)$  from (\ref{eq:gp}) (and $v_i$
occurs in it precisely once).

 The argument above shows that the  variable $z_{j(i)}$ does not occur
in words $P_{t}(Z,A)$ ($t \neq i$), $Q_{s}(Z,A)$ ($s < i$).
Moreover, in $P_{i}(Z)$ it occurs precisely once.  It follows that
the variable $z_{j(i)}$ in the generalized equation $\Omega(Z)$
does not occur neither in
 coefficient equations nor in basic equations corresponding to the dual
bases related to $x_t$ ($t \neq i$), $y_s$ ($s < i$).

We "mark" (or select) the unique  occurrence of $v_i$ (as
$v_i^{\pm 1}$) in $u_{j(i)}$ $i = 1, \ldots, k$. Now we are going
to mark some other occurrences of $v_i$ in words $u_1, \ldots,
u_\rho$ as follows. Suppose some $u_d$ has a marked occurrence of
some $v_i$. If $\Omega$
  contains  an equation of the type
 $z_d^{\varepsilon} =  z_r^{\delta}$,  then $u_d^{\varepsilon} =
 u_r^{\delta}$ graphically. Hence $u_r$ has an  occurrence of subword
$v_i^{\pm 1}$ which correspond to the marked occurrence of
$v_i^{\pm 1}$ in $u_d$. We mark this occurrence of $v_i^{\pm 1}$
in $u_r$.

  Suppose  $\Omega$
  contains  an equation of the type
  $$[h_{\alpha _1}\ldots h_{\beta
_1-1}]^{\varepsilon _1}= [h_{\alpha _2}\ldots h_{\beta _2-
1}]^{\varepsilon _2}$$
 such that  $z_d$ occurs in it, say in the left. Then
 $$ [u_{\alpha _1}\ldots u_{\beta
_1-1}]^{\varepsilon _1}= [u_{\alpha _2}\ldots u_{\beta _2-
1}]^{\varepsilon _2}$$
 graphically. Since  $v_i^{\pm +1}$ is a subword of $u_d$, it occurs also
in the right-hand part
  of the equality above, say in some $u_r$. We marked this
  occurrence of $v_i^{\pm +1}$
  in  $u_r$. The marking process will be over in finitely many steps.
 Observe that one and
 the same $u_r$ can have several  marked occurrences
of some $v_i^{\pm 1}$.

Now in all words $u_1, \ldots, u_\rho$ we replace every marked
occurrence of $v_i = ba^{m_{ij}}b$ with a new  word
$ba^{m_{ij}}x_{i}b$ from the group $F[X]$. Denote the resulting
words from $F[X]$ by $\tilde u_1, \ldots, \tilde u_\rho$. It
follows from description of the marking process that the tuple
$\tilde U = (\tilde u_1, \ldots, \tilde u_{\rho})$  is a solution
of the generalized equation $\Omega$ in the free group $F[X]$.
Indeed, all the equations in $\Omega$ are graphically satisfied by
the substitution $z_i \rightarrow u_i$ hence the substitution $u_i
\rightarrow \tilde u_i$ still makes them graphically equal. Now by
Lemma \ref{le:R1} $X = P(\tilde U), Y = Q(\tilde U)$ is a solution
of the equation $S(X,A) = 1$ over  $F[X]$ as desired. $\Box$

\begin{cy}\cite{Mak84}
\label{co:merz} There is an algorithm which for a given positive
sentence
 $$ \forall x_1\exists y_1\ldots\forall x_k\exists y_k (S(X,Y,A)=1)$$
  in $L_A$ determines  whether or not this formula holds in $F$, and if
it does, the algorithm  finds  words
  $$q_1(x_1),\ldots , q_k(x_1,\ldots ,x_k) \in F[X]$$
   such that
$$F[X]\models u(x_1, q_1(x_1),\ldots ,x_k, q_k(x_1,\ldots ,x_k))=1.$$
\end{cy}
{\em Proof.}   The proof follows from Proposition \ref{merzl} and
decidability of equations over free groups with constraints $y_i
\in F[X_i]$, where $X_i = \{x_1, \ldots, x_i\}$ \cite{Mak82}.

\begin{df} Let $\phi$ be a sentence  in the language $L_A$
written in the standard form
 $$ \phi = \forall x_1\exists
y_1\ldots\forall x_k\exists y_k \phi_0(x_1,y_1, \ldots,
x_k,y_k),$$ where $\phi_0$ is a quantifier-free formula in $L_A$.
We say that $G$ {\em freely lifts} $\phi$ if there exist words
(with constants from $F$) $q_1(x_1),\ldots , q_k(x_1,\ldots ,x_k)
\in F[X],$ such that
$$F[X]\models \phi_0(x_1, q_1(x_1),\ldots ,x_k, q_k(x_1,\ldots
,x_k,A))=1.$$
\end{df}

\begin{theorem}
 \label{merzl} $F$ freely lifts every sentence in $L_A$ that
is true in $F$.
 \end{theorem}

\medskip
{\it Proof.} Suppose a sentence
  \begin{equation}
  \label{eq:sent}
   \phi = \forall x_1\exists
y_1\ldots\forall x_k\exists y_k (U(x_1,y_1, \ldots, x_k,y_k) =1
 \wedge V(x_1,y_1, \ldots, x_k,y_k) \neq 1),
  \end{equation}
 is true in $F$. We choose $x_1 = g_1, y_1 = h_1, \ldots, x_k= g_k, y_k=h_k $ precisely like in
the  Merzlyakov's Theorem. Then the formula
 $$U(g_1,h_1, \ldots, g_k,h_k) =1
 \wedge V(g_1,h_1, \ldots, g_k,h_k) \neq 1 $$
 holds in $F$. In particular, $U(g_1,h_1, \ldots, g_k,h_k) =1$ in
 $F$. It follows from the argument in Theorem \ref{merzl} that
 there are words $q_1(x_1) \in F[x_1], \ldots, q_k(x_1, \ldots, x_k) \in F[x_1,
 \ldots,x_k]$ such that
   $$F[X] \models U(x_1,q_1(x_1, \ldots,x_k), \ldots, x_k,q_k(x_1, \ldots,x_k))
   =1.$$
    Moreover, it follows from the construction that $h_1 = q_1(g_1), \ldots, h_k = q_k(g_1,
    \ldots,g_k)$. We claim that
     $$F[X] \models V(x_1,q_1(x_1, \ldots,x_k), \ldots, x_k,q_k(x_1, \ldots,x_k)) \neq
     1.$$
Indeed, if
 $$V(x_1,q_1(x_1, \ldots,x_k), \ldots, x_k,q_k(x_1, \ldots,x_k)) =
 1$$
in $F[X]$, then its image in $F$ under any specialization $X
\rightarrow F$ is also trivial, but this is not the case for
specialization  $x_1 \rightarrow g_1, \ldots, x_k \rightarrow g_k$
- contradiction. This proves the theorem for  sentences $\phi$ of
the form (\ref{eq:sent}). A similar argument works for formulas of
the type
 $$ \phi = \forall x_1\exists
y_1\ldots\forall x_k\exists y_k \bigvee_{i = 1}^{n}(U_i(x_1,y_1,
\ldots, x_k,y_k) =1
 \wedge V_i(x_1,y_1, \ldots, x_k,y_k) \neq 1),$$
 which is, actually,  the general case by Corollary
 \ref{co:atomic}. This finishes the proof.
  \hfill $\Box$

 \section{Makanin's process and  Cut equations}
\label{se:5}

\subsection{ Elementary transformations}
\label{se:5.1}

 In this section we describe {\em elementary transformations} of generalized
 equations which were introduced by Makanin in \cite{Mak82}.
 Recall that we consider only formally consistent equations. In
general, an elementary transformation $ET$ associates   to  a
generalized equation $\Omega$ a finite set of generalized
equations $ET(\Omega) = \{\Omega_1, \ldots, \Omega_r\}$  and a
collection of surjective homomorphisms $\theta _i:G_{R(\Omega
)}\rightarrow G_{R(\Omega _i)}$ such that for every pair $(\Omega
, U)$ there exists a unique pair of the type $(\Omega _i, U_i)$
for which the following diagram commutes.
\medskip

\begin{center}

\begin{picture}(100,100)(0,0)
\put(0,100){$F_{R(\Omega )}$} \put(100,100){$F_{R(\Omega _i)}$}
\put(0,0){$F(A)$} \put(15,103){\vector(1,0){80}}
\put(5,93){\vector(0,- 1){78}} \put(95,95){\vector(-1,-1){80}}
\put(-10,50){$\pi _{U_i}$} \put(55,45){$\pi _{U_i}$}
\put(50,108){$\theta _i$}
\end{picture}

\end{center}

\medskip
Here $\pi _{U}(X)=U.$  Since the pair $(\Omega_i,U_i)$ is defined
uniquely, we have a well-defined map  \newline $ET : (\Omega,U)
\rightarrow (\Omega_i,U_i).$

 ET1  ({\em Cutting a base}). Suppose  $\Omega$ contains a boundary
connection $<p,\lambda ,q>$. Then we replace (cut in $p$)  the
base $\lambda$ by two new bases  $\lambda _1$ and $\lambda _2$ and
also replace (cut in $q$) $\Delta(\lambda)$ by two new bases
$\Delta (\lambda _1)$ and $\Delta (\lambda _2)$ such that the
following conditions hold.

If $\varepsilon(\lambda) =  \varepsilon(\Delta(\lambda))$, then
$$\alpha(\lambda_1) = \alpha(\lambda), \ \  \beta(\lambda_1) = p, \ \ \ \
\alpha(\lambda_2) = p, \ \  \beta(\lambda_2) = \beta(\lambda); $$
 $$
\alpha(\Delta(\lambda_1)) = \alpha(\Delta(\lambda)), \ \
\beta(\Delta(\lambda_1)) = q, \ \ \ \ \alpha(\Delta(\lambda_2)) =
q, \ \ \beta(\Delta(\lambda_2)) = \beta(\Delta(\lambda)); $$

If  $\varepsilon(\lambda) = - \varepsilon(\Delta(\lambda))$, then
$$\alpha(\lambda_1) = \alpha(\lambda), \ \  \beta(\lambda_1) = p, \ \ \ \
\alpha(\lambda_2) = p, \ \  \beta(\lambda_2) = \beta(\lambda); $$
 $$
\alpha(\Delta(\lambda_1)) = q,  \ \  \beta(\Delta(\lambda_1)) =
\beta(\Delta(\lambda)), \ \ \ \ \alpha(\Delta(\lambda_2)) =
\alpha(\Delta(\lambda)),  \ \  \beta(\Delta(\lambda_2)) = q; $$

Put $\varepsilon(\lambda_i) = \varepsilon(\lambda), \ \
\varepsilon(\Delta(\lambda_i)) = \varepsilon(\Delta(\lambda)), \ i
= 1,2.$

Let $(p', \lambda, q')$ be  a boundary connection in $\Omega$.

If    $p' < p$,  then replace $(p', \lambda, q')$ by  $(p',
\lambda_1, q')$.

If    $p' > p$, then replace   $(p', \lambda, q')$ by  $(p',
\lambda_2, q')$.

Notice, since the equation $\Omega$ is formally consistent, then
the conditions above define boundary connections in the new
generalized equation. The resulting generalized equation
$\Omega^\prime$ is formally consistent.   Put $ET(\Omega) =
\{\Omega^\prime \}$.  Fig. 3 below explains the name of the
transformation ET1.

\begin{figure}[here]
\centering{\mbox{\psfig{figure=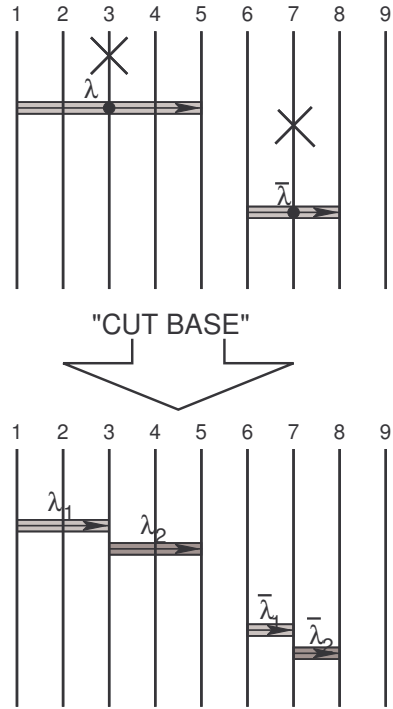,height=4in}}}
\caption{Elementary transformation ET1.} \label{ET1}
\end{figure}

ET2  ({\em Transfer of a base}).  Let    a base $\theta $ of a
generalized equation $\Omega$  be contained in the base $\mu$,
i.e., $\alpha (\mu)\leq \alpha (\theta)<\beta (\theta)\leq\beta
(\mu)).$ Suppose that the boundaries $\alpha(\theta)$ and
$\beta(\theta))$ are $\mu$-tied, i.e.,  there are boundary
connections of the type $<\alpha (\theta),\mu, \gamma_1>$ and
$<\beta (\theta),\mu,\gamma_2>$.  Suppose also that every
$\theta$-tied boundary is $\mu$-tied. Then we transfer $\theta$
from its location  on the base $\mu$ to the corresponding location
on the base $\Delta (\mu)$ and adjust all the basic and boundary
equations (see Fig. 4). More formally, we replace $\theta$ by a
new base $\theta^\prime$ such that $\alpha(\theta^\prime) =
\gamma_1, \beta(\theta^\prime) = \gamma_2$ and replace each
$\theta$-boundary connection $(p,\theta,q)$ with a new one
$(p^\prime,\theta^\prime,q)$ where $p$ and $p^\prime$ come from
the $\mu$-boundary connection $(p,\mu, p^\prime)$. The resulting
equation is denoted by $\Omega^\prime = ET2(\Omega)$.

\begin{figure}[here]\centering
{\mbox{\psfig{figure=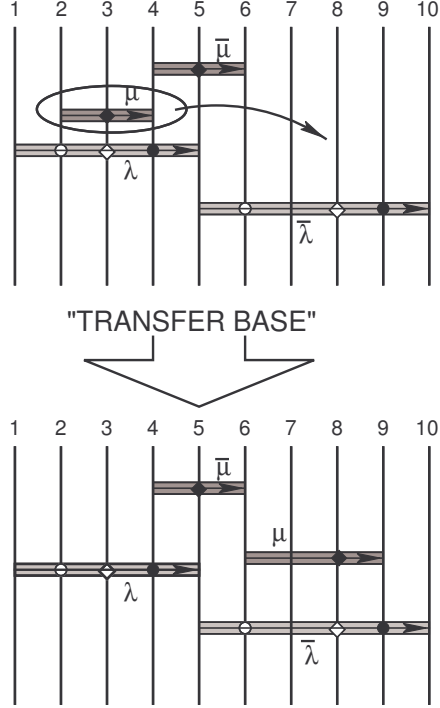,height=4in}}} \caption{Elementary
transformation ET2.} \label{ET2}
\end{figure}

 ET3 ({\em Removal of a pair of  matched bases} (see Fig. 5)). Let
$\mu$ and $\Delta(\mu)$ be a pair of matched bases in $\Omega$.
Since $\Omega$ is formally consistent one  has $\varepsilon(\mu) =
\varepsilon(\Delta(\mu))$, $\beta(\mu) = \beta(\Delta(\mu))$ and
every $\mu$-boundary connection is of the type $(p,\mu,p)$. Remove
the pair of bases  $\mu, \Delta(\mu)$ with all  boundary
connections related to $\mu$. Denote the new generalized equation
by $\Omega^\prime$.

\smallskip
 {\bf Remark.} Observe, that for $i = 1,2,3$  $ETi(\Omega)$ consists of
a single equation $\Omega^\prime$,  such that  $\Omega$ and
$\Omega^\prime$ have the same set of variables $H$,   and the
identity map $F[H] \rightarrow F[H]$ induces an $F$-isomorphism
$F_{R(\Omega )} \rightarrow F_{R(\Omega^{\prime })}$. Moreover,
 $U$ is a solution of $\Omega$ if and only if $U$ is  a solution of
$\Omega^\prime$.

\begin{figure}[here]
\centering{\mbox{\psfig{figure=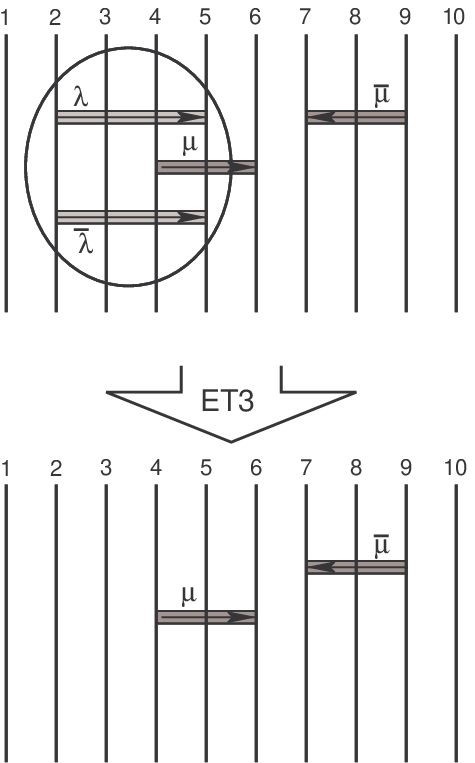}}} \caption{Elementary
transformation ET3.} \label{ET3}
\end{figure}

\begin{figure}[here]
\centering{\mbox{\psfig{figure=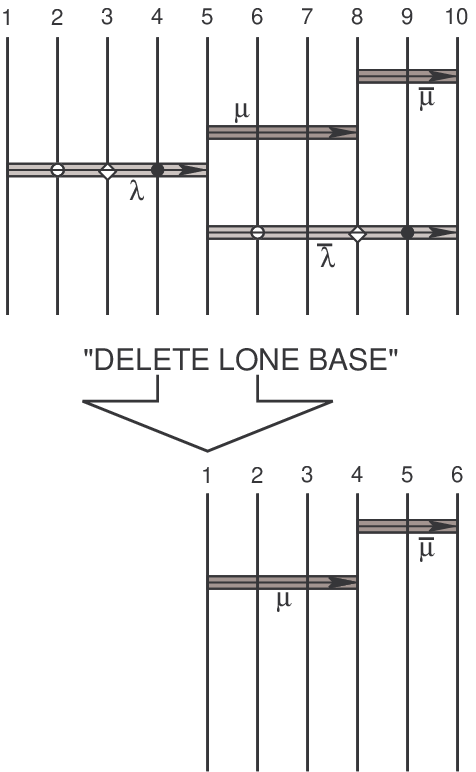}}} \caption{Elementary
transformation ET4.} \label{ET4}
\end{figure}

\begin{figure}[here]
\centering{\mbox{\psfig{figure=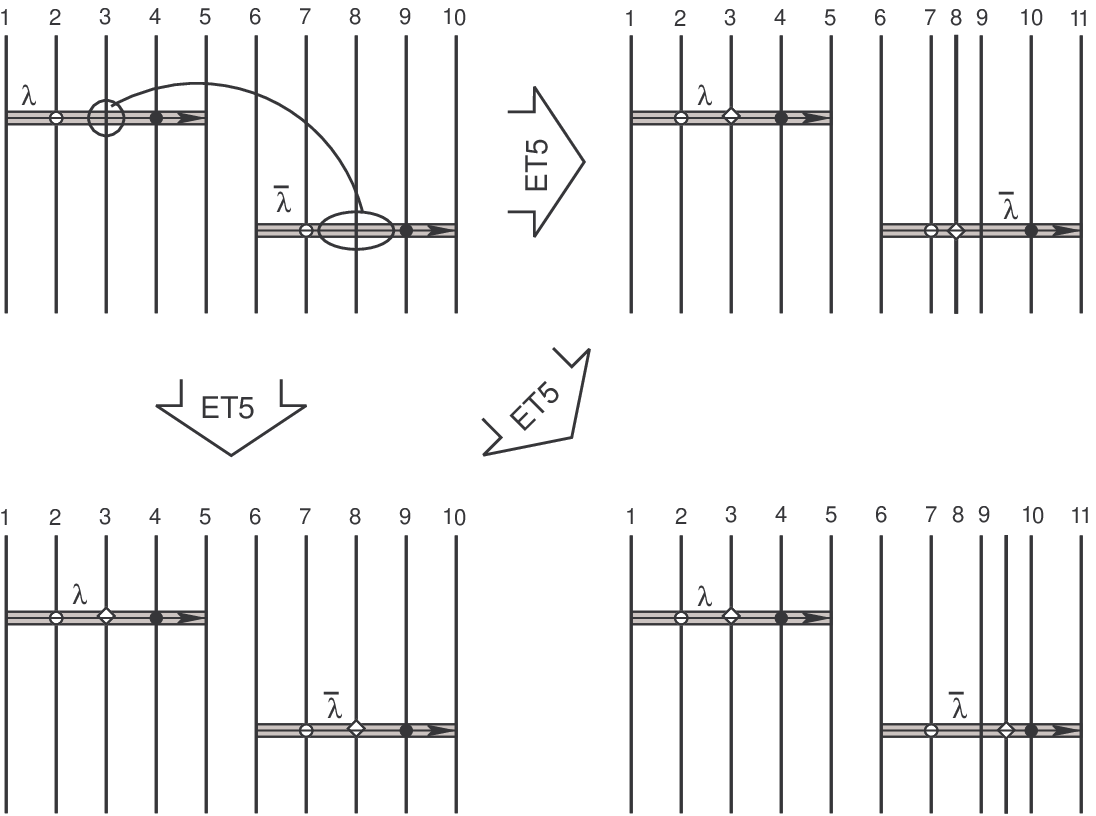}}} \caption{Elementary
transformation ET5.} \label{ET5}
\end{figure}

ET4 ({\em Removal of a lonely base} (see Fig. 6)).  Suppose in
$\Omega$ a variable base $\mu$ does not intersect any other
variable base, i.e., the items $h_{\alpha(\mu)}, \ldots,
h_{\beta(\mu)-1}$ are contained in only one variable base $\mu$.
Suppose also that all boundaries in $\mu$ are $\mu$-tied, i.e.,
for every $i$ ($\alpha(\mu)+1 \leq i\leq \beta -1$) there exists a
boundary $b(i)$ such that $(i,\mu ,b(i))$ is a boundary connection
in $\Omega$.  For convenience  we define: $b(\alpha(\mu)) =
\alpha(\Delta(\mu))$ and  $b(\beta(\mu)) = \beta(\Delta(\mu))$ if
$\varepsilon(\mu)\varepsilon(\Delta(\mu)) = 1$, and
$b(\alpha(\mu)) = \beta(\Delta(\mu))$ and $b(\beta(\mu)) =
\alpha(\Delta(\mu))$ if $\varepsilon(\mu)\varepsilon(\Delta(\mu))
= -1$.

\begin{figure}
\centering{\mbox{\psfig{figure=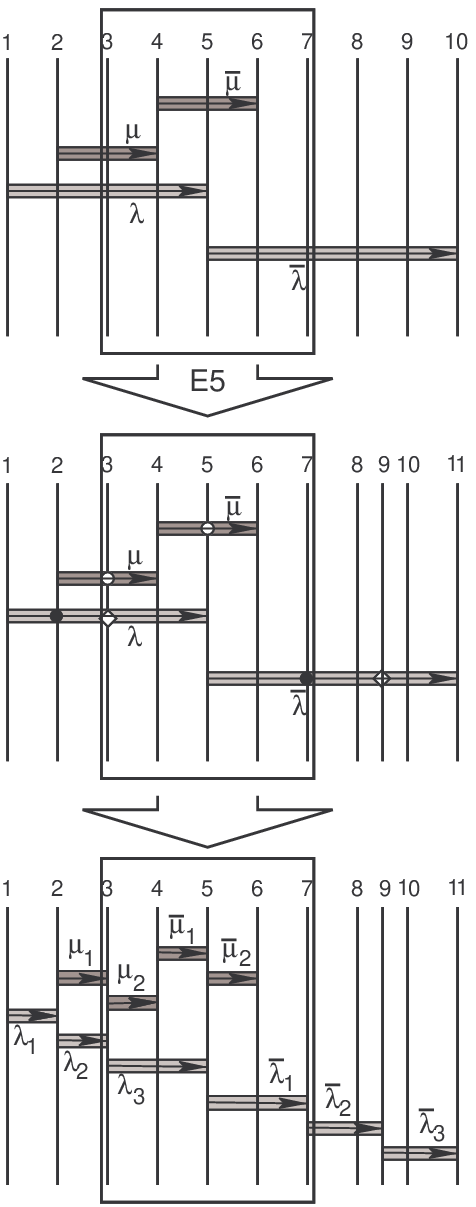}}} \caption{Derived
transformation D1.} \label{D1}
\end{figure}
\begin{figure}[here]
\centering{\mbox{\psfig{figure=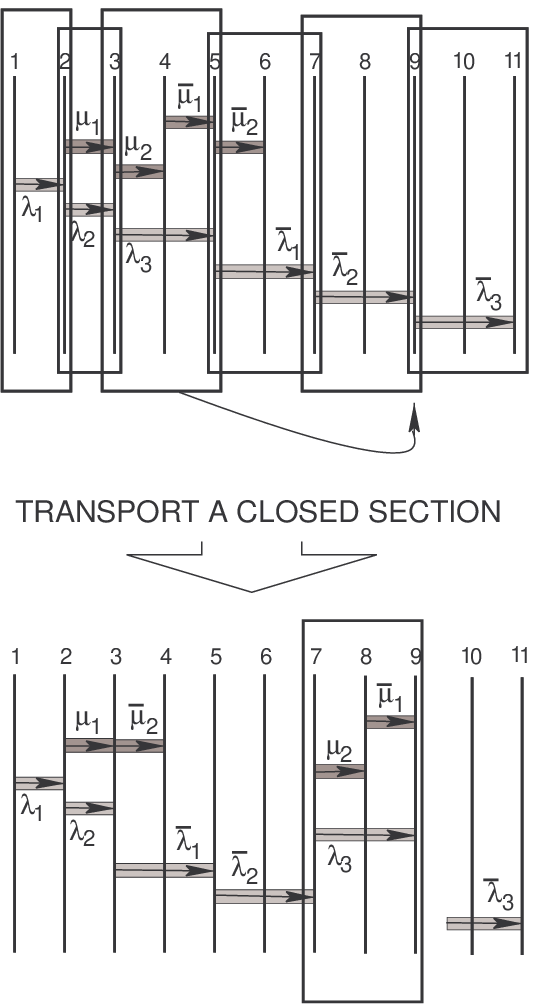}}} \caption{Derived
transformation D2.} \label{D2}
\end{figure}
 The transformation ET4 carries $\Omega $ into a  unique
generalized equation $\Omega _1$ which is obtained from $\Omega $
by deleting the pair of bases $\mu$ and $\Delta(\mu)$; deleting
all the boundaries $\alpha(\mu)+1, \ldots, \beta(\mu)-1$ ( and
renaming  the rest $\beta(\mu) - \alpha(\mu) - 1$ boundaries)
together with all $\mu$-boundary connections; replacing every
constant base $\lambda$ which is contained in  $\mu$ by a constant
base $\lambda^\prime$ with the same label as $\lambda$ and such
that $\alpha(\lambda^\prime) = b(\alpha(\lambda)),
\beta(\lambda^\prime) = b(\beta(\lambda))$.

 We define the homomorphism $\pi: F_{R(\Omega )} \rightarrow
F_{R(\Omega^{\prime })}$
 as follows:
$\pi(h_j)=h_j$ if $j<\alpha (\mu)$ or $j\geq \beta (\mu);$

$$\pi(h_{i})=\{\begin{array}{ll} h_{b(i)}\ldots h_{b(i)-1},& if \varepsilon
(\mu)=\varepsilon (\Delta\mu),\\ h_{b(i)}\ldots h_{b(i-1)-1},& if
\varepsilon (\mu)=-\varepsilon (\Delta\mu)
\end{array}$$
for $\alpha +1 \leq i\leq\beta (\mu)-1.$ It is not hard to see
that $\pi$ is  an $F$-isomorphism.

 ET5 ({\em Introduction of a boundary} (see Fig. 7)).  Suppose  a
point $p$ in a base $\mu$ is not
 $\mu$-tied. The transformation ET5 $\mu$-ties it in all possible ways,
producing
 finitely many different generalized equations. To this end, let $q$ be
a boundary
 on $\Delta (\mu)$. Then we perform one of the following two
transformations:

1. Introduce the boundary connection $<p,\mu ,q>$ if the resulting
equation $\Omega_q$ is formally consistent. In this case the
corresponding $F$-homomorphism
 $\pi_q: F_{R(\Omega )}$ into $F_{R(\Omega _q)}$ is induced by
the identity isomorphism on  $F[H]$. Observe that $\pi_q$  is not
necessary an isomorphism.

2. Introduce a new boundary $q^\prime$ between $q$ and $q+1$ (and
rename all the boundaries); introduce  a new boundary connection
$(p,\mu,q^\prime)$. Denote the resulting equation by
$\Omega_q^\prime$. In this case the corresponding $F$-homomorphism
 $\pi_{q^\prime}: F_{R(\Omega )}$ into $F_{R(\Omega _{q^{\prime}})}$
is induced by the map $\pi(h) = h,$ if $h\neq h_q$,  and $\pi(h_q)
= h_{q^\prime}h_{q^\prime +1}$.  Observe that  $\pi _{q^\prime}$
is an $F$-isomorphism.

Let $\Omega$ be a generalized equation and $E$ be an elementary
transformation. By $E(\Omega)$ we denote a generalized equation
obtained from $\Omega$ by elementary transformation $E$ (pehaps
several such equations) if $E$ is applicable to $\Omega$,
otherwise we put $E(\Omega) = \Omega$.   By $\phi_E :
F_{R(\Omega)} \rightarrow F_{R(E(\Omega))}$ we denote the
canonical homomorphism of the coordinate groups (which has been
described above in the case $E(\Omega) \neq \Omega)$, otherwise,
the identical  isomorphism.

\begin{lm}
\label{le:hom-check} There exists an algorithm which for every
generalized equation $\Omega$ and every elementary transformation
$E$ determines whether the canonical homomorphism $\phi_E:
F_{R(\Omega)} \rightarrow F_{R(E(\Omega))}$  is an isomorphism or
not.
\end{lm}
{\em Proof.}   The only non-trivial case is when $E = E5$ and no
new boundaries were introduced. In this case $E(\Omega)$ is
obtained from $\Omega$ by adding a new particular equation, say $s
= 1$, which is effectively determined by $\Omega$ and $E(\Omega)$.
In this event, the coordinate group
$$F_{R(E(\Omega))} = F_{R(\Omega \cup \{s\})}$$
is a quotient group of $F_{R(\Omega)}$. Now $\phi_E$ is an
isomorphism if and only if $R(\Omega) = R(\Omega \cup \{s\})$, or,
equivalently, $s \in R(\Omega)$. The latter condition holds if and
only if $s$ vanishes on all solutions of the system of
(group-theoretic) equations $\Omega = 1$ in $F$, i.e., if the
following formula holds in $F$:
 $$ \forall x_1 \ldots \forall x_\rho (\Omega(x_1, \ldots, x_\rho) = 1
 \rightarrow s(x_1, \ldots, x_\rho) = 1).$$
 This can be checked effectively, since the universal theory of a free
 group $F$ is decidable (\cite{Mak84}).

\subsection{Derived transformations and auxiliary transformations}
\label{se:5.2half}

 In this section we describe several useful transformations of
generalized equations. Some of them can be realized as  finite
sequences of elementary transformations, we call them {\it derived
} transformations. Other transformations result in equivalent
generalized equations but cannot be realized by finite sequences
of elementary moves.

 D1  ({\em Closing a section}).

Let $\sigma$ be a section of $\Omega $. The transformation D1
makes the section $\sigma$ closed. To perform D1 we  introduce
boundary connections (transformations ET5) through the end-points
of $\sigma$ until these end-points are tied by every base
containing them, and then cut through the end-points all the bases
containing them (transformations ET1)(see Fig. 8)

D2 ({\it Transporting a closed section}).

Let $\sigma$ be a closed section of a generalized equation
$\Omega$. We cut $\sigma$ out of the interval $[1,\rho_\Omega]$
together with all the bases and boundary connections on $\sigma$
and put $\sigma$ at the end of the interval or between any two
consecutive   closed sections of $\Omega$. After that we
correspondingly re-enumerate all the items and boundaries of the
latter equation to bring it to the proper form. Clearly, the
original equation $\Omega$ and the new one $\Omega^\prime$ have
the same solution  sets and their coordinate groups are isomorphic
(see Fig. 79)

D3 ({\it  Complete cut}).

 Let $\Omega $ be  a generalized
equation. For every boundary connection $(p,\mu,q)$ in $\Omega$ we
cut the base $\mu$ at $p$  applying ET1. The resulting generalized
equation $\widetilde\Omega$ is obtained from $\Omega$ by a
consequent application of all possible ET1 transformations.
Clearly, $\widetilde\Omega$ does not depend on a particular choice
of the sequence of transformations ET1. Since ET1 preserves
isomorphism between the coordinate groups, equations $\Omega$ and
$\widetilde\Omega$ have isomorphic coordinate groups, and the
isomorphism arises from the identity map $F[H] \rightarrow F[H]$.

D4 ({\it Kernel of a generalized equation}).

Suppose  that a generalized equation $\Omega$  does not contain
boundary connections. An active base $\mu  \in A\Sigma_\Omega$ is
called {\em eliminable} if at least one of the following holds:

 a) $\mu$ contains an  item $h_i$ with $\gamma(h_i)=1$;

b) at least one of the boundaries $\alpha (\mu),\beta (\mu)$ is
different from $1,\rho +1 $ and it does not touch any other base
(except  $\mu$).

An {\it elimination  process} for $\Omega$ consists of consequent
removals ({\em eliminations}) of  eliminable bases until no
eliminable bases left in the equation. The resulting generalized
equation is called a {\it kernel} of $\Omega$ and we denote it by
$Ker(\Omega )$. It is easy to see that $Ker(\Omega )$ does not
depend on a particular elimination process. Indeed, if $\Omega$
has two different eliminable bases $\mu_1$, $\mu_2$, and deletion
of  $\mu_i$ results in an equation $\Omega_i$ then  by induction
(on the number of eliminations) $Ker(\Omega_i)$ is uniquely
defined for $i = 1,2$. Obviously, $\mu_1$ is still eliminable in
$\Omega_2$, as well as $\mu_2$ is eliminable in $\Omega_1$. Now
eliminating $\mu_1$ and $\mu_2$ from $\Omega_2$ and $\Omega_1$ we
get one and  the same equation $\Omega_0$. By induction
$Ker(\Omega_1) = Ker(\Omega_0) = Ker(\Omega_2)$ hence the result.
 We say that a variable $h_i$ {\it belongs to
the kernel} ($h_i \in Ker(\Omega)$), if either $h_i$ belongs to at
least one base in the kernel, or it is parametric, or it is
constant.

Also, for an equation $\Omega$ by $\overline{\Omega}$ we denote
the equation which is obtained from $\Omega$ by deleting all free
variables. Obviously,
$$F_{R(\Omega)} = F_{R(\overline {\Omega})} \ast F(Y)$$
where $Y$ is the set of free variables in $\Omega$.

Let us consider what happens on the group level in the elimination
process.

 We start with the case  when just one base is eliminated.
 Let $\mu$ be an eliminable base in $\Omega = \Omega(h_1, \ldots,
h_\rho)$. Denote by $\Omega_1$ the equation resulting from
$\Omega$ by eliminating $\mu$.

1) Suppose  $h_i \in \mu$ and $\gamma(h_i) = 1$. Then the variable
$h_i$ occurs only once in $\Omega$ - precisely in the equation
$s_\mu = 1$ corresponding to the base $\mu$.  Therefore, in the
coordinate group $F_{R(\Omega )}$ the relation $s_\mu = 1$ can be
written as $h_i = w$, where $w$ does not contain $h_i$. Using
Tietze transformations we can rewrite the presentation of
$F_{R(\Omega )}$ as $F_{R(\Omega^\prime)}$, where $\Omega^\prime$
is obtained from $\Omega$ by deleting $s_\mu$ and the item $h_i$.
It follows immediately that
$$F_{R(\Omega_1)} \simeq  F_{R(\Omega^\prime)} \ast \langle h_i \rangle$$
and
 \begin{equation}
 \label{eq:ker1}
 F_{R(\Omega)} \simeq F_{R(\Omega^\prime)} \simeq F_{R(\overline
{\Omega_1})} \ast F(Z) \end{equation}
 for some free group $F(Z)$.
Notice that all the groups and equations which occur above can be
found effectively.

2) Suppose now that  $\mu$ satisfies case b) above with respect to
a boundary $i$. Then in the equation $s_\mu = 1$ the variable
$h_{i-1}$ either occurs only once or it occurs precisely twice and
in this event the second occurrence of $h_{i-1}$ (in
$\Delta(\mu)$) is a part of the subword $(h_{i-1}h_i)^{\pm 1}$. In
both cases it is easy to see that the tuple
$$(h_1, \ldots, h_{i-2}, s_\mu, h_{i-1}h_i, h_{i+1}, \ldots,h_\rho)$$
forms a basis of the ambient free group generated by $(h_1,
\ldots, h_\rho)$  and constants from $A$. Therefore, eliminating
the relation $s_\mu = 1$,  we can rewrite the presentation of
$F_{R(\Omega)}$ in generators $Y = (h_1, \ldots, h_{i-2},
h_{i-1}h_i, h_{i+1}, \ldots,h_\rho)$. Observe also  that  any
other equation $s_\lambda = 1$ ($\lambda \neq \mu$) of $\Omega$
either does not contain variables $h_{i-1}, h_i$ or it contains
them as parts of the subword $(h_{i-1}h_i)^{\pm 1}$, i.e., any
such a word  $s_\lambda$ can be expressed as a word $w_\lambda(Y)$
in terms of generators $Y$ and constants from $A$. This shows that
$$F_{R(\Omega)} \simeq F(Y \cup A)_{R(w_\lambda(Y) \mid \lambda \neq
\mu)} \simeq F_{R(\Omega^\prime)},$$ where $\Omega^\prime$ is a
generalized equation obtained from $\Omega_1$ by deleting the
boundary $i$. Denote by $\Omega^{\prime \prime}$ an equation
obtained from $\Omega^\prime$ by adding a free variable $z$ to the
right end of $\Omega^\prime$.
 It follows now that
 $$F_{R(\Omega_1)} \simeq  F_{R(\Omega^{\prime \prime})} \simeq
 F_{R(\Omega)} \ast \langle z \rangle$$
and
\begin{equation}
\label{eq:ker2} F_{R(\Omega)} \simeq F_{R(\overline
{\Omega^\prime})} \ast F(Z) \end{equation}
 for some free group $F(Z)$. Notice that all the
groups and equations which occur above can be found effectively.

By induction on the number of steps in elimination process we
obtain the following lemma.
\begin{lm}\label{7-10}
 $$F_{R(\Omega)} \simeq F_{R(\overline {Ker \Omega})} \ast F(Z)$$
 where $F(Z)$ is a free group on $Z$. Moreover, all the groups and equations
which occur above can be found effectively.\end{lm}
 {\em Proof.}   Let
 $$\Omega = \Omega_0 \rightarrow \Omega_1 \rightarrow \ldots
 \rightarrow \Omega_l = Ker\Omega$$
 be an elimination process for $\Omega$. It is easy to see (by induction on $l$)
 that for every $j = 0, \ldots,l-1$
 $$\overline{Ker\Omega_j} = \overline{Ker \overline{\Omega_j}}.$$
Moreover, if $\Omega_{j+1}$ is obtained from $\Omega_j$ as in the
case 2) above, then (in the notations above)
 $$\overline{Ker (\Omega_j)_1} = \overline{Ker \Omega_j^\prime} .$$
Now the statement of the lemma follows from the remarks above and
equalities (\ref{eq:ker1}) and (\ref{eq:ker2}).

D5 ({\it  Entire transformation}).

We need a few further definitions.
 A base $\mu $ of the equation $\Omega$ is called a {\it leading} base if
$\alpha(\mu)=1$. A leading base is said to be {\it maximal} (or a
{\it carrier}) if $\beta (\lambda)\leq \beta (\mu),$ for any other
leading base $\lambda $.  Let $\mu $ be a carrier base of
 $\Omega.$ Any active base $\lambda \neq \mu$ with $\beta(\lambda )\leq \beta (\mu )$
 is called a {\it transfer} base (with respect to $\mu$).

  Suppose now that $\Omega$ is a generalized equation with
$\gamma(h_i)\geq 2$ for each $h_i$ in the active part of $\Omega$.
An {\em entire transformation} is a sequence of elementary
transformations which are performed as follows. We fix a carrier
base $\mu$ of $\Omega$. For any transfer base $\lambda $ we
 $\mu$-tie (applying $ET5$) all boundaries in $\lambda $. Using ET2 we transfer
 all transfer
bases from $\mu$ onto $\Delta (\mu)$. Now, there exists some $i
<\beta (\mu)$ such that $h_1,\ldots ,h_i$ belong to only one base
$\mu,$ while
 $h_{i+1}$ belongs to at least two bases. Applying ET1 we cut $\mu$
along the boundary $i+1$. Finally, applying ET4 we delete the
section $[1,i+1]$.

D6 ({\it Identifying closed constant sections}).

Let $\lambda$ and $\mu$ be two constant bases in $\Omega$ with
labels $a^{\varepsilon_\lambda}$ and $a^{\varepsilon_\mu}$, where
$a \in A$ and $\varepsilon_\lambda, \varepsilon_\mu \in \{1,-1\}$.
Suppose that the sections $\sigma(\lambda) = [i,i+1]$ and
$\sigma(\mu) = [j,j+1]$ are closed. Then we introduce a new
variable base $\delta$ with its dual $\Delta(\delta)$ such that
$\sigma(\delta) = [i,i+1]$, $\sigma(\Delta(\delta)) = [j,j+1]$,
$\varepsilon(\delta) = \varepsilon_\lambda,
\varepsilon(\Delta(\delta)) = \varepsilon_\mu$. After that we
transfer all bases from $\delta$ onto $\Delta(\delta)$ using ET2,
remove the bases $\delta$ and $\Delta(\delta)$, remove the item
$h_i$, and enumerate the items in a proper order. Obviously, the
coordinate group of the resulting equation is   isomorphic to the
coordinate group of the original equation.

\subsection{ Construction of the tree  $T(\Omega)$}
\label{se:5.2}

In this section we describe a branching rewrite  process for a
generalized equation $\Omega$. This process results in an
(infinite) tree $T(\Omega)$. At the end of the section we describe
infinite paths in $T(\Omega)$.

{\bf Complexity of a parametric generalized equation.}

 Denote by $\rho _A$ the number of variables $h_i$ in all active
sections of $\Omega ,$  by $n_A=n_A(\Omega )$ the number of bases
in active sections of $\Omega$,   by $\nu '$ - the number of open
boundaries in the active sections, by $\sigma '$ - the number of
closed boundaries in the active sections.

The number of closed active sections containing no bases,
precisely one base, or more than one base are denoted by $t_{A0},
t_{A1}, t_{A2}$ respectively. For a closed section $\sigma \in
\Sigma_\Omega$ denote by $n(\sigma)$, $\rho(\sigma)$  the number
of bases and, respectively, variables  in $\sigma$.

$$\rho _A= \rho_A(\Omega) = \sum_{\sigma \in
A\Sigma_\Omega}\rho(\sigma)$$
$$n_A= n_A(\Omega) = \sum_{\sigma \in
A\Sigma_\Omega} n(\sigma)$$

The {\em complexity} of a parametric  equation $\Omega $ is the
number
$$\tau = \tau (\Omega) = \sum_{\sigma \in
A\Sigma_\Omega} max\{0, n(\sigma)-2\}.$$

Notice that the entire transformation (D5) as well as  the
cleaning process (D4) do not increase complexity of equations.

Let  $\Omega $ be a parametric generalized equation. We construct
a tree $T(\Omega)$ (with associated structures), as a directed
tree oriented from a root $v_0$, starting at $v_0$ and proceeding
by induction from vertices at distance $n$ from the root to
vertices at distance $n+1$ from the root.

We start with a general description of the tree $T(\Omega)$.  For
each vertex $v$ in $T(\Omega)$ there exists a unique generalized
equation $\Omega_v$ associated with $v$. The initial equation
$\Omega$ is associated with the root $v_0$, $\Omega_{v_0} =
\Omega$. For each edge $v\rightarrow v'$ (here $v$ and $v'$ are
the origin and the terminus of the edge) there exists a unique
surjective homomorphism $\pi(v,v'):F_{R(\Omega _v )}\rightarrow
F_{R(\Omega _v' )}$ associated with $v\rightarrow v'$.

 If
 $$v\rightarrow v_1\rightarrow\ldots\rightarrow
v_s\rightarrow u$$ is a path in $T(\Omega )$, then by $\pi (v,u)$
we denote  composition of corresponding homomorphisms
$$\pi (v,u) = \pi (v,v_1)   \ldots
\pi (v_s,u).$$

The set of edges of $T(\Omega)$ is subdivided into two classes:
{\it principal} and {\it auxiliary}. Every newly constructed edge
is principle, if not said otherwise. If $v \rightarrow v^\prime$
is a principle edge then there exists a finite sequence of
elementary  or derived transformations from $\Omega_v$ to
$\Omega_{v^\prime}$ and the homomorphism $\pi(v,v')$ is
composition of the homomorphisms corresponding to these
transformations. We also assume that active [non-active] sections
in $\Omega_{v^\prime}$ are naturally inherited from $\Omega_v$, if
not said otherwise.

Suppose the tree $T(\Omega)$ is constructed by induction up to a
level $n$, and suppose  $v$ is a vertex at distance $n$ from the
root $v_0$. We describe now how to extend the tree from $v$.  The
construction of the  outgoing edges at  $v$ depends on which case
described below takes place at the vertex $v$. We always assume
that if we have Case $i,$ then all Cases $j$, with $j \leq i-1$,
do not take place at $v$. We will see from the description below
that there is an effective procedure to check whether or not a
given case takes place at a given vertex. It will be obvious for
all cases, except Case 1. We treat this case below.

{\bf Preprocessing}

Case 0.  In $\Omega_v$ we transport closed sections using D2 in
such a way that all active sections are at the left end of the
interval (the active part of the equation), then come all
non-active sections (the non-active part of the equation), then
come parametric sections (the parametric part of the equation),
and behind them  all constant sections are located (the constant
part of the equation).

 {\bf Termination conditions}

Case 1. The homomorphism $\pi (v_0,v)$ is not an isomorphism (or
equivalently, the homomorphism $\pi(v_1,v)$, where $v_1$ is the
parent of $v$, is not an isomorphism). The vertex $v$ is called a
{\it leaf} or an {\it end} vertex. There are no outgoing edges
from $v$.

\begin{lm}
There is an algorithm  to verify whether the homomorphism
$\pi(v,u)$, associated with an edge $v \rightarrow u$ in
$T(\Omega)$  is an isomorphism or not.
\end{lm} {\em Proof.}   We will see below (by a straightforward inspection of
Cases 1-15 below) that every homomorphism of the type $\pi(v,u)$
is a composition of the canonical homomorphisms corresponding to
the elementary (derived) transformations. Moreover, this
composition is effectively given. Now the result follows from
Lemma \ref{le:hom-check}.

Case 2. $\Omega _v$ does not contain active sections. The vertex
$v$ is called a {\it leaf} or an {\it end} vertex. There are no
outgoing edges from $v$.

{\bf Moving constants to the right}

Case 3. $\Omega _v$ contains a constant base $\lambda$ in an
active section such that  the section $\sigma(\lambda)$ is not
closed.

Here we  close the section $\sigma(\lambda)$ using the derived
transformation $D1$.

Case 4. $\Omega _v$ contains a constant base $\lambda$ with a
label $a \in A^{\pm 1}$ such that the section $\sigma(\lambda)$ is
closed.

Here we transport the section $\sigma(\lambda)$ to the location
right after all variable and parametric sections in $\Omega_v$
using the derived transformation D2. Then we identify  all closed
sections of the type $[i,i+1]$, which contain a constant base with
the label $a^{\pm 1}$, with the transported section
$\sigma(\lambda)$, using the derived transformation D6. In the
resulting generalized equation $\Omega_{v^\prime}$ the section
$\sigma(\lambda)$ becomes a constant section, and the
corresponding edge $(v,v^\prime)$ is auxiliary. See Fig. 10.

\begin{figure}
\centering{\mbox{\psfig{figure=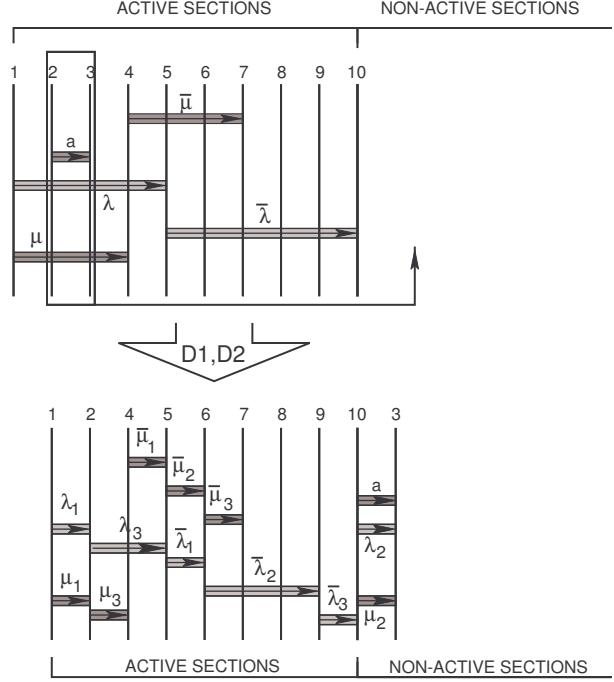}}} \caption{Case 3-4:
Moving constant bases.} \label{constants}
\end{figure}

{\bf Moving free variables to the right}

 Case 5. $\Omega _v$ contains a free variable $h_q$ in an  active section.

Here we close the section $[q,q+1]$ using D1, transport it to the
very end of the interval behind all items in $\Omega_v$ using D2.
In the resulting generalized equation $\Omega_{v^\prime}$ the
transported section  becomes a constant section, and the
corresponding edge $(v,v^\prime)$ is auxiliary.

\begin{rk} If Cases 0-5 are not possible at $v$ then the parametric
generalized equation $\Omega_v$ is in standard form.
\end{rk}

 Case 6. $\Omega _v$ contains a pair of matched bases in an
active section.

Here we perform ET3 and delete it. See Fig. 11.

\begin{figure}[here]
\centering{\mbox{\psfig{figure=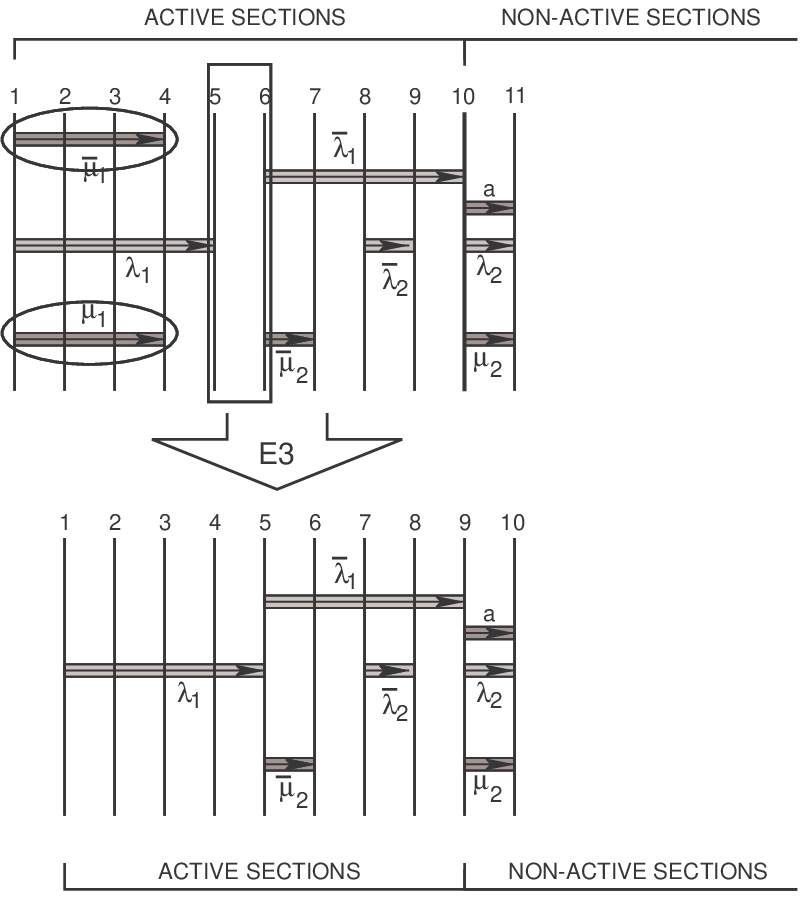}}} \caption{Case 5-6:
Trivial equations and useless variables.} \label{useless}
\end{figure}

{\bf Eliminating linear variables}

 Case 7. In $\Omega_v$ there is $h_i$ in an active section with
$\gamma _i=1$ and  such that both boundaries $i$ and $i+1$ are
closed.

Here we remove the closed section $[i,i+1]$ together with the lone
base using ET4.

Case 8. In $\Omega_v$ there is $h_i$ in an active section with
$\gamma _i=1$ and  such that  one of the boundaries $i,i+1$ is
open, say $i+1$, and the other is closed.

Here we perform ET5 and $\mu$-tie $i+1$ through the only base
$\mu$ it intersects; using ET1 we cut $\mu$ in $i+1$; and then we
delete the closed section $[i,i+1]$ by ET4.  See Fig. 12.

\begin{figure}[here]
\centering{\mbox{\psfig{figure=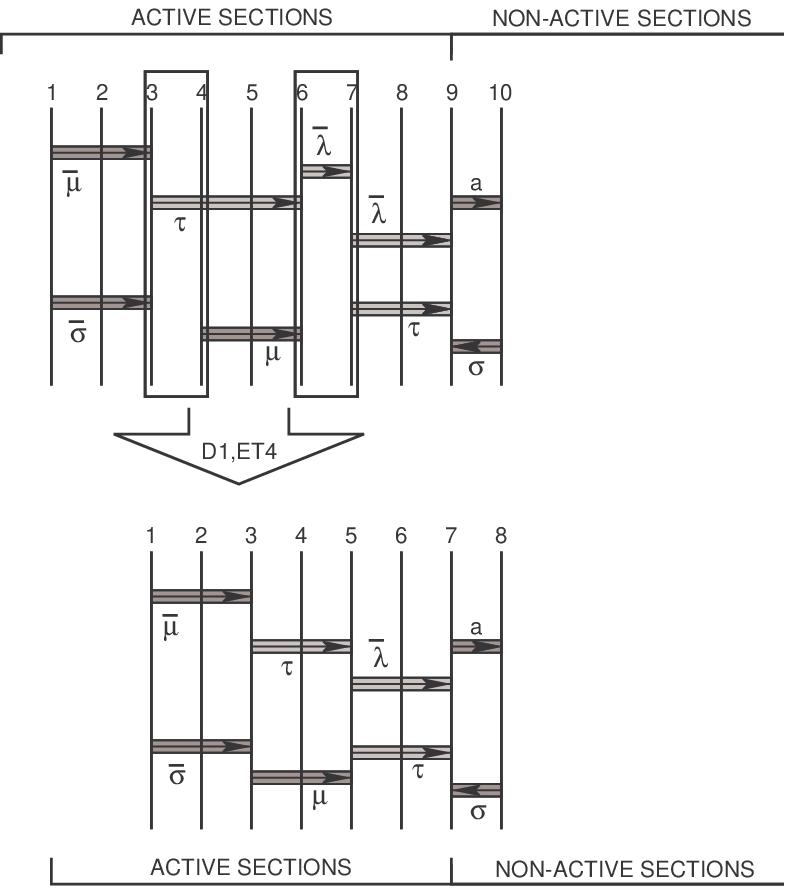}}} \caption{Case 7-10:
Linear variables.} \label{linear'}
\end{figure}

Case 9. In $\Omega_v$ there is $h_i$ in an active section with
$\gamma _i=1$ and  such that both boundaries $i$ and $i+1$ are
open. In addition, assume that there is a closed section $\sigma$
containing exactly two (not matched)  bases $\mu _1$ and $\mu _2$,
such that $\sigma = \sigma(\mu_1) = \sigma(\mu_2)$ and in the
generalized equation $\tilde {\Omega}_v$ (see the derived
transformation D3) all the bases obtained from $\mu _1,\mu _2$  by
ET1 in constructing $\tilde {\Omega}_v$ from $\Omega_v$, do not
belong to the kernel of $\tilde{\Omega}_v.$

Here, using ET5, we  $\mu _1$-tie  all the boundaries inside
$\mu_1$;
 using ET2, we transfer $\mu _2$ onto $\Delta (\mu_1)$; and
 remove $\mu _1$ together with the closed section $\sigma$  using
 ET4.

Case 10. $\Omega_v$  satisfies the first assumption of Case 9 and
does not satisfy the second one.

In this event we close the section $[i,i+1]$ using D1 and remove
it using ET4.

{\bf Tying a free boundary}

Case 11. Some boundary $i$ in  the active part of $\Omega_v$ is
free. Since we do not have Case 5 the boundary $i$ intersects at
least one base, say, $\mu$.

Here we $\mu$-tie $i$  using ET5.

{\bf Quadratic case}

Case 12. $\Omega_v$ satisfies the condition  $\gamma _i = 2$ for
each $h_i$ in the active part.

We apply the entire transformation D5.

Case 13. $\Omega_v$ satisfies the condition $\gamma _i\geq 2$ for
each $h_i$ in the active part,  and $\gamma _i>2$ or at least one
such $h_i$. In addition,  for some active base $\mu$ section
$\sigma(\mu) = [\alpha (\mu ),\beta (\mu )]$ is closed.

In this case using $ET5$, we $\mu$-tie every boundary inside
$\mu$; using $ET2$, we transfer all bases from $\mu$ to $\Delta
(\mu )$; using $ET4$, we remove the lone base $\mu$ together with
the section $\sigma(\mu)$.

Case 14.  $\Omega_v$ satisfies the condition $\gamma _i\geq 2$ for
each $h_i$ in the active part, and $\gamma _i>   2$ for at least
one such $h_i$. In addition,  some boundary $j$ in the active part
touches some base $\lambda$, intersects some base $\mu$,  and $j$
is not $\mu$-tied.

Here we $\mu$-tie $j$.

{\bf General JSJ-case}

 Case 15. $\Omega_v$ satisfies the condition $\gamma _i\geq 2$ for
each $h_i$ in the active part, and $\gamma _i>   2$ for at least
one such $h_i$.  We apply, first,  the entire transformation D5.

Here for every boundary $j$ in the active part that touches at
least one base, we $\mu$-tie $j$ by every base $\mu$ containing
$j$. This results in finitely many new vertices
$\Omega_{v^\prime}$ with principle edges $(v,v^\prime)$.

If, in addition,  $\Omega_v$ satisfies the following condition (we
called it  Case 15.1 in \cite{KMIrc}) then we construct the
principle edges as was described above, and also  construct a few
more auxiliary edges outgoing from the vertex $v$:

Case 15.1. The carrier base $\mu$ of the equation $\Omega _v$
intersects with its dual $\Delta (\mu)$.

Here we  construct an auxiliary    equation ${\hat \Omega _{ v}}$
(which does not occur in $T(\Omega)$)  as follows. Firstly, we add
a new constant section $[\rho _v+1,\rho _v+2]$ to the right of all
sections in $\Omega_v$ (in particular, $h_{\rho _v+1}$ is a new
free variable). Secondly, we introduce a new pair of bases
$(\lambda ,\Delta (\lambda))$ such that
$$\alpha(\lambda)=1, \beta (\lambda)= \beta
(\Delta (\mu)), \alpha (\Delta(\lambda))=\rho _v+1, \beta (\Delta
(\lambda))=\rho _v+2.$$
 Notice that $\Omega _v$ can be obtained from ${\hat \Omega
_{v}}$ by ET4:  deleting $\delta(\lambda)$ together with the
closed section $[\rho _v+1,\rho _v+2]$. Let
$${\hat \pi}_{v}: F_{R(\Omega _v)}\rightarrow F_{R({\hat \Omega}_{v})}$$
 be the isomorphism induced by ET4.  Case 15 still holds for ${\hat
 \Omega _{ v}}$,  but now $\lambda$ is the carrier base. Applying to
 ${\hat \Omega _{v}}$ transformations described in  Case 15, we obtain a
 list of new vertices $\Omega_{v^\prime}$ together with isomorphisms
$${\eta}_{v'}: F_{R({\hat \Omega}_{ v})} \rightarrow F_{R(\Omega _{v'})}.$$

Now for each such $v'$ we add to $T(\Omega)$ an auxiliary edge
$(v,v')$ equipped  with composition of homomorphisms $\pi(v,v') =
\eta_{v'} \circ {\hat \pi}_{v}$ and assign $\Omega_{v'}$ to the
vertex $v'$.

If none of the Cases 0-15 is possible, then we stop, and the tree
$T(\Omega)$ is constructed. In any case, the tree $T(\Omega)$ is
constructed by induction. Observe that, in general, $T(\Omega)$ is
an infinite locally finite tree.

If Case $i\ (0\leq i\leq 15)$ takes place at  a vertex $v$ then we
say that $v$ has type $i$ and write $tp(v)=i$.

\begin{lm} \label{3.1}(Lemma 3.1, \cite{Razborov3}) If $u\rightarrow
v$ is a principal edge of the tree $T(\Omega)$, then
\begin{enumerate}
\item $n_{A}(\Omega_v) \leq n_{A}(\Omega_u)$, if $tp(v_1)\not =
3,10,$ this inequality is proper if $tp(v_1)=6,7,9,13;$ \item If
$tp(v_1)=10,$ then $n_{A}(\Omega_v) \leq n_{A}(\Omega_u) + 2;$
\item $\nu^\prime(\Omega_v)
 \leq \nu^\prime(\Omega_u)$ if $tp(v_1)\leq 13$ and $tp (v_1)\not = 3,11$; \item
$\tau(\Omega_v)  \leq \tau(\Omega_u)$, if $tp (v_1)\not = 3.$
\end{enumerate}\end{lm}
{\em Proof.}   Straightforward verification.

\begin{lm} \label{3.2} Let
 $$v_1\rightarrow v_2\rightarrow\ldots \rightarrow v_r \rightarrow \ldots $$
  be an infinite path in the tree $T(\Omega ).$
Then there exists a natural number $N$ such that all the edges
$v_n \rightarrow v_{n+1}$ of this path with $n \geq N$ are
principal edges, and one of the following situations holds:
\begin{enumerate}
\item [1)] ({\bf linear case})\  $7\leq tp(v_n)\leq 10 $ for all
$n\geq N;$ \item [2)] ({\bf quadratic case}) \ $tp(v_n)=12$ for
all $n\geq N;$ \item [3)] ({\bf general JSJ case})\  $tp(v_n)=15$
for all $n\geq N$.
\end{enumerate}
\end{lm}
{\em Proof.}   Observe that starting with a generalized equation
$\Omega$ we can have Case 0 only once, afterward in all other
equations the active part is at the left, then comes the
non-active part, then - the parametric part,  and at the end - the
constant part. Obviously, Cases 1 and 2 do not occur on an
infinite path. Notice also that Cases 3 and 4 can only occur
finitely many times, namely, not more then $2t$ times where $t$ is
the number of constant bases in the original equation $\Omega$.
Therefore, there exists a natural number $N_1$ such that
$tp(v_i)\geq 5$ for all $i \geq N_1$.

Now we show that the number of vertices $v_i$ ($i \geq N$) for
which $tp (v_i)=5$ is not more than the minimal number of
generators of the group $F_{R(\Omega )}$, in particular, it cannot
be greater than $\rho +1 +|A|$, where $\rho =  \rho (\Omega).$
Indeed, if a path from the root $v_0$ to a vertex $v$ contains $k$
vertices of type 5, then $\Omega _v$ has at least $k$ free
variables in the constant part. This implies that the coordinate
group $F_{R(\Omega _v)}$ has a free group of rank $k$ as a free
factor, hence it cannot be generated by less than $k$ elements.
Since $\pi(v_0,v): F_{R(\Omega )} \rightarrow F_{R(\Omega _v)}$ is
a surjective homomorphism, the group $F_{R(\Omega )}$ cannot be
generated by less then $k$ elements. This shows that $k \leq \rho
+1 +|A|$.  It follows that there exists  a number $N_2 \geq N_1$
such that  $tp(v_i) >   5$ for every $i \geq N_2$.

Suppose $i >   N_2$. If $tp(v_i)=12,$ then it is easy to see that
$tp(v_{i+1})=6$ or $tp(v_{i+1})=12$. But if $tp(v_{i+1})=6,$ then
$tp(v_{i+2})=5$ - contradiction with $i >   N_2$. Therefore,
$tp(v_{i+1})=tp(v_{i+2})= \ldots = tp(v_{i+j}) = 12$ for every $j
>   0$ and we have  situation 2) of the lemma.

Suppose now $tp(v_{i})\neq 12$ for all $i \geq N_2$. By Lemma
\ref{3.1} $\tau(\Omega_{v_{j+1}}) \leq \tau(\Omega_{v_{j}})$ for
every principle edge $v_j \rightarrow v_{j+1}$ where $j \geq N_2$.
If $v_j \rightarrow v_{j+1}$, where $j \geq N_2$, is  an auxiliary
edge then  $tp(v_j) = 15$ and, in fact, Case 15.1 takes place at
$v_j$.  In the notation of Case 15.1  $\Omega_{v_{j+1}}$ is
obtained from ${\hat \Omega}_{v_j}$ by transformations from Case
15. In this event, both bases $\mu$ and $\Delta (\mu)$ will be
transferred from the new carrier base $\lambda $ to the constant
part, so the complexity will be decreased at least by two:
$\tau(\Omega_{v_{j+1}}) \leq \tau({\hat \Omega}_{v_j}) - 2$.
Observe also that $\tau ({\hat \Omega_{v_j}}) = \tau(\Omega_{v_j})
+ 1$. Hence $\tau(\Omega_{v_{j+1}})
 <  \tau(\Omega_{v_j}).$

It follows that there exists a number $N_3 \geq N_2$ such that
$\tau(\Omega_{v_j}) = \tau(\Omega_{v_{N_3}})$ for every $j \geq
N_3$, i.e., complexity stabilizes. Since every auxiliary edge
gives a decrease of complexity, this implies that for every $j
\geq N_3$ the edge $v_j \rightarrow v_{j+1}$ is principle.

Suppose now that $i \geq N_3$. We claim that $tp(v_i) \neq 6$.
Indeed, if $tp(v_i)=6$, then the closed section, containing the
matched bases $\mu, \Delta (\mu)$, does not contain any other
bases (otherwise the complexity of $\Omega_{v_{i+1}}$ would
decrease).  But in this event $tp(v_{i+1})=5$ which is impossible.

So  $tp (v_i)\geq 7$ for every $i \geq N_3$. Observe that ET3
(deleting match bases) is the only elementary transformation that
can produce new free boundaries. Observe also that ET3 can be
applied only in Case 6.  Since Case 6 does not occur anymore along
the path  for $i \geq N_3$, one can see that no new free
boundaries occur in equations $\Omega_{v_j}$ for $j \geq N_3$. It
follows that there exists a number $N_4 \geq N_3$ such that $tp
(v_i)\neq 11$ for every $j \geq N_4$.

Suppose now that for some $i \geq N_4$  $13\leq tp(v_i)\leq 15.$
It is easy to see from the description of these cases  that in
this event $tp(v_{i+1})\in\{6,13,14,15\}.$ Since $tp(v_{i+1})\neq
6$, this implies that $13\leq tp(v_j)\leq 15$ for every $j \geq
i$. In this case the sequence $n_A(\Omega_{v_j})$  stabilizes by
lemma \ref{3.1}. In addition, if $tp (v_j)=13,$ then $
n_A(\Omega_{v_{j+1}}) <  n_A(\Omega_{v_j}).$ Hence there exists a
number $N_5 \geq N_4$ such that $tp (v_j)\not = 13$ for all $j
\geq N_5.$

Suppose $i \geq N_5$.  There cannot be more than
$8(n_A(\Omega_{v_i})^2$ vertices of type 14 in a row starting at a
vertex $v_i$; hence there exists $j\geq i$ such that $tp(v_j)=15$.
The series of transformations ET5 in Case 15 guarantees the
inequality $tp (v_{j+1})\not = 14$; hence $tp(v_{j+1})=15,$ and we
have situation 3) of the lemma.

So we can suppose $tp (v_i)\leq 10$ for all the vertices of our
path. Then we have situation 1) of the lemma. $\Box$

\subsection {Periodized equations}

In this section we introduce a notion of a periodic structure
which allows one to describe  periodic solutions of generalized
equations. Recall that a reduced word $P$ in a free group $F$ is
called a {\it period} if it is cyclically reduced and not a proper
power. A word $w \in F$ is called $P$-{\em periodic}  if $|w| \geq
|P|$ and it is a subword of  $P^n$ for some $n$. Every
$P$-periodic word $w$ can be presented in the form
\begin{equation}
\label{2.50} w = A^rA_1
\end{equation}
where $A$  is a cyclic permutation of $P^{\pm 1}$,  $r \geq 1$, $A
= A_1 \circ A_2$, and $A_2 \neq 1$. This representation is unique
if $r \geq 2$. The number $r$ is called the {\em exponent}  of
$w$. A maximal exponent of $P$-periodic subword in a word $u$ is
called the {\em exponent of $P$-periodicity in $u$}. We denote it
$e_P(u).$

\begin{df}\label{11'}
Let $\Omega$ be a  standard generalized equation. A solution
$H:h_i \rightarrow H_i$ of $\Omega$ is called {\em periodic with
respect to a period $P$}, if for every  variable section $\sigma$
of $\Omega$ one of the following conditions hold:
\begin{enumerate}
 \item [1)] $H(\sigma)$ is $P$-periodic with exponent $r \geq 2$;
 \item [2)] $|H(\sigma)| \leq |P|$;
 \item [3)]  $H(\sigma)$ is $A$-periodic and $|A| \leq |P|$;
\end{enumerate}
Moreover, condition 1) holds at least for one such $\sigma$.
\end{df}

Let $H$ be a $P$-periodic solution of $\Omega$. Then a section
$\sigma$ satisfying 1) is called
 {\em $P$-periodic} (with
respect to $H$).

\subsubsection{Periodic  structure}
Let $\Omega$ be a parametrized generalized equation. It turns out
that every periodic  solution of $\Omega$ is a composition of a
canonical automorphism of the coordinate group $F_{R(\Omega)}$
with  either a  solution with bounded exponent of periodicity
(modulo parameters) or a solution of a "proper" equation. These
canonical automorphisms  correspond to Dehn twists of
$F_{R(\Omega)}$ which are related to the splitting of this group
(which comes from the periodic structure) over an abelian edge
group.

We fix till the end of the section a generalized equation $\Omega$
in standard form. Recall that in $\Omega$ all closed sections
$\sigma$, bases $\mu$, and variables $h_i$ belong to either the
variable part $V\Sigma$, or the parametric part $P\Sigma$, or the
constant part $C\Sigma$ of $\Omega$.

\begin{df} \label{above}
Let  $\Omega$ be a generalized equation in  standard form with no
boundary connections. A {\em periodic structure} on $\Omega$ is a
pair $\langle {\mathcal P}, R \rangle$, where
\begin{enumerate}
\item [1)]  ${\mathcal P}$ is a set consisting of some variables
$h_i$, some bases $\mu$, and some closed sections $\sigma$ from
$V\Sigma$  and such that the following conditions are satisfied:
\begin{itemize}
\item[a)] if $h_i \in {\mathcal P}$ and $h_i \in \mu$,  and
$\Delta(\mu) \in V\Sigma$, then $\mu \in {\mathcal P}$;

\item[b)] if $\mu \in {\mathcal P}$, then $\Delta(\mu) \in
{\mathcal{P}}$;

\item[c)] if $\mu \in {\mathcal P}$ and $\mu \in \sigma$, then
$\sigma \in {\mathcal P}$;

\item [d)] there exists a function ${\mathcal X}$ mapping the set
of closed sections from ${\mathcal P}$ into $\{-1, +1\}$ such that
for every $\mu, \sigma_1, \sigma_2 \in {\mathcal P}$, the
condition that $\mu \in \sigma_1$ and $\Delta(\mu) \in \sigma_2$
implies

$\varepsilon(\mu) \cdot \varepsilon(\Delta(\mu)) = {\mathcal
X}(\sigma_1) \cdot {\mathcal X}(\sigma_2)$;
\end{itemize}
\item [2)]  $R$ is an equivalence relation on a certain set
${\mathcal B}$ (defined below) such that the following conditions
are satisfied:
\begin{itemize}

 \item[e)] Notice, that for every boundary $l$ belonging to a closed section in $\mathcal P$ either
there exists a unique closed section $\sigma(l)$ in ${\mathcal P}$
containing $l$, or there exist precisely two closed section
$\sigma_{left}(l) = [i,l], \sigma_{right} =  [l,j]$ in ${\mathcal
P}$ containing $l$. The set of boundaries of the first type we
denote
 by ${\mathcal B}_1$, and of the second  type - by  ${\mathcal B}_2$. Put
$${\mathcal B} = {\mathcal B}_1  \cup \{l_{left}, l_{right}  \mid l \in {\mathcal B}_2
\}$$
 here $l_{left}, l_{right}$ are two  "formal copies" of $l$. We
 will use the following agreement: for any
base $\mu$ if $\alpha(\mu) \in {\mathcal B}_2$ then by
$\alpha(\mu)$ we mean $\alpha(\mu)_{right}$ and, similarly, if
$\beta(\mu) \in {\mathcal B}_2$ then by $\beta(\mu)$ we mean
$\beta(\mu)_{left}$.

\item[f)]  Now, we define $R$ as follows. If $\mu \in {\mathcal
P}$ then
$$\alpha(\mu) \sim_R \alpha(\Delta(\mu)), \ \ \beta(\mu)
\sim_R \beta(\Delta(\mu)) \ \ if \ \varepsilon(\mu) =
\varepsilon(\Delta(\mu))$$
$$\alpha(\mu)) \sim_R \beta(\Delta(\mu)),
 \beta(\mu)\sim_R  \alpha(\Delta(\mu)) \ \ if \ \varepsilon(\mu) = -
 \varepsilon(\Delta(\mu)).$$
\end{itemize}
\end{enumerate}
\end{df}
\begin{rk}
This definition coincides with the definition of a periodic
structure given in \cite{KMIrc} in the case of empty set of
parameters  $P\Sigma$. For a given $\Omega$ one can effectively
find all periodic structures on $\Omega$.
\end{rk}

 Let $\langle {\mathcal P}, R \rangle $ be a periodic structure of
$\Omega$. Put
$$ {\mathcal NP} = \{\mu \in B\Omega \mid \exists h_i \in {\mathcal P} \ such \
that  \ h_i \in \mu \ and  \ \Delta(\mu) \ is \ parametric \ or \
constant \}$$

Now we will show how one can associate with a $P$-periodic
solution $H$ of  $\Omega$  a periodic structure ${\mathcal P}(H,
P) = \langle {\mathcal P}, R \rangle$. We define ${\mathcal P}$ as
follows. A closed section $\sigma$ is  in ${\mathcal P}$ if and
only if $\sigma$ is $P$-periodic. A variable $h_i$ is in
${\mathcal P}$ if and only if $h_i \in \sigma$ for some $\sigma
\in {\mathcal P}$ and $d(H_i) \geq 2 d(P)$. A base $\mu$ is in
${\mathcal P}$ if and only if both $\mu$ and $\Delta(\mu)$ are in
$V\Sigma$  and one of them contains $h_i$ from ${\mathcal P}$.

Put ${\mathcal X}([i,j]) = \pm 1$ depending on whether in
(\ref{2.50}) the word $A$ is conjugate to $P$ or to $P^{-1}$.

Now let $[i,j]\in {\mathcal P}$ and $ i \leq l \leq j$. Then there
exists a subdivision $P = P_1P_2$ such that if ${\mathcal X}
([i,j]) =1$, then the word $H[i,l]$ is the end of the word
$(P^\infty)P_1$, where $P^\infty$ is the infinite word obtained by
concatenations of powers of $P$, and $H[l,j]$ is the beginning of
the word $P_2(P^\infty)$, and if ${\mathcal X}([i,j])= -1$, then
the word $H[i,l]$ is the end of the word $(P^{-1})^\infty
P_2^{-1}$ and $H[l,j]$ is the beginning of
$P_1^{-1}(P^{-1})^\infty$. Lemma 1.2.9 \cite{1} implies that the
subdivision $P=P_1P_2$ with the indicated properties is unique;
denote it by $\delta(l)$. Let us define a relation $R$ in the
following way: $R(l_1, l_2) \rightleftharpoons \delta(l_1) =
\delta(l_2)$.

\begin{lm}
\label{le:PP} Let $H$ be a periodic solution of $\Omega$. Then
${\mathcal P}(H, P)$ is a periodic structure on $\Omega$.
\end{lm}
{\em Proof.}   Let ${\mathcal P}(H, P) = \langle {\mathcal P}, R
\rangle$.
 Obviously,  ${\mathcal P}$ satisfies a) and b) from the definition \ref{above}.

Let $\mu \in {\mathcal P}$ and $\mu \in [i,j]$. There exists an
unknown $h_k \in {\mathcal P}$ such that $h_k \in \mu$ or $h_k \in
\Delta({\mu})$. If $h_k \in \mu$, then, obviously, $[i,j] \in
{\mathcal P}$. If $h_k \in \Delta(\mu)$ and $\Delta(\mu) \in [i',
j']$, then $[i', j'] \in {\mathcal P}$, and hence, the word
$H[\alpha(\Delta(\mu)), \beta(\Delta(\mu))]$ can be written in the
form $Q^{r'} Q_1$, where $Q=Q_1 Q_2$; $Q$ is a cyclic shift of the
word $P^{\pm 1}$ and $r' \geq 2$. Now let (\ref{2.50}) be a
presentation for the section $[i,j]$. Then $H[\alpha(\mu),
\beta(\mu)] = B^s B_1$, where $B$ is a cyclic shift of the word
$A^{\pm 1}$, $d(B) \leq d(P)$, $B = B_1 B_2$, and $s \geq 0$. From
the equality $H[\alpha(\mu), \beta(\mu)]^{\varepsilon(\mu)} =
H[\alpha(\Delta(\mu)),
\beta(\Delta(\mu)))]^{\varepsilon(\Delta(\mu))}$ and Lemma 1.2.9
\cite{1} it follows that $B$ is a cyclic shift of the word
$Q^{\pm1}$. Consequently, $A$ is a cyclic shift of the word
$P^{\pm 1}$, and $r \geq 2$ in (\ref{2.50}), since $d(H[i,j]) \geq
d(H[\alpha(\mu), \beta(\mu)]) \geq 2 d(P)$. Therefore, $[i,j] \in
{\mathcal P}$; i.e, part c) of the definition \ref{above} holds.

 If
$\mu \in [i_1, j_1]$, $\Delta(\mu) \in [i_2, j_2]$, and $ \mu \in
{\mathcal P}$, then the equality $\varepsilon(\mu) \cdot
\varepsilon(\Delta(\mu))$ = ${\mathcal X}([i_1, j_1]) \cdot
{\mathcal X}([i_2, j_2])$ follows from the fact that given $A^r
A_1 = B^s B_1$ and $r,s \geq 2$, the word $A$ cannot be a cyclic
shift of the word $B^{-1}$. Hence part d) also holds.

 Condition e) of the definition of a periodic
structure obviously holds.

Condition f) follows from the graphic equality $H[\alpha(\mu),
\beta(\mu)]^{\varepsilon(\mu)}=$ $H[\alpha(\Delta(\mu)),$ $
\beta(\Delta(\mu))]^{\varepsilon(\Delta(\mu))}$ and Lemma 1.2.9
\cite{1}.

This proves the lemma. \hfill $\Box$

Now let us fix a nonempty periodic structure $\langle {\mathcal
P}, R \rangle$. Item d) allows us to assume (after replacing the
variables $h_i, \ldots, h_{j-1}$ by $h_{j-1}^{-1}, \ldots,
h_i^{-1}$ on those sections $[i,j] \in {\mathcal P}$ for which
${\mathcal X}([i,j])=-1$) that $\varepsilon(\mu)=1$ for all $\mu
\in {\mathcal P}$. For a boundary $k$, we will denote by $(k)$ the
equivalence class of the relation $R$ to which it belongs.

 Let us construct
an oriented  graph $\Gamma$ whose set of vertices is the set of
$R$--equivalence classes. For each unknown $h_k$ lying on a
certain closed section from ${\mathcal P}$, we introduce an
oriented edge $e$ leading from $(k)$ to $(k+1)$ and an inverse
edge $e^{-1}$ leading from $(k+1)$ to $(k)$. This edge $e$ is
assigned the label $h(e) \rightleftharpoons h_k$ (respectively,
$h(e^{-1}) \rightleftharpoons h_k^{-1}$.) For every path
$r=e_1^{\pm 1} \ldots e_s^{\pm 1}$ in the graph $\Gamma$, we
denote by $h(r)$ its label $h(e_1^{\pm 1}) \ldots h(e_j^{\pm 1})$.
The periodic structure $\langle {\mathcal P}, R \rangle$ is called
{\em connected}, if the graph $\Gamma$ is connected. Suppose first
that $\langle {\mathcal P}, R \rangle$ is connected.  Suppose that
some boundary $k$ (between $h_{k-1}$ and $h_k$) in the variable
part of $\Omega$ is not a boundary between two bases. Since
$h_{k-1}$ and $h_k$ appear in all the basic equations together,
and there is no boundary equations, one can consider a generalized
equation $\Omega _1$ obtained from $\Omega$ by replacing the
product $h_{k-1}h_k$ in all basic equations by one variable
$h_k'$. The group $F_{R(\Omega)}$ splits as a free product of the
cyclic group generated by $h_{k-1}$ and $F_{R(\Omega _1)}.$  In
this case we can consider $\Omega _1$ instead of $\Omega$.
Therefore we  suppose now that each boundary of $\Omega $ is a
boundary between two bases.

\begin{lm} \label{2.9}
Let $H$ be a $P$-periodic solution of a generalized equation
$\Omega$,    $\langle {\mathcal P}, R \rangle = {\mathcal P}(H,
P)$; $c$ be a cycle in the graph $\Gamma$ at the vertex $(l)$;
$\delta(l)=P_1P_2$. Then there exists $n \in {\bf Z}$ such that
$H(c) = (P_2P_1)^n.$
\end{lm}

\noindent {\em Proof.}   \hspace{2mm} If $e$ is an edge in the
graph $\Gamma$ with initial vertex $V'$ and terminal vertex $V''$,
and $P = P_1'P_2', $ $P = P_1 '' P_2''$ are two subdivisions
corresponding to the boundaries from $V'$, $V''$ respectively,
then, obviously, $H(e) = P_2' P^{n_k}P_1''$ $(n_k \in {\bf Z})$.
The claim is easily  proven by multiplying together the values
$H(E)$  for all the edges $e$ taking part in the cycle $c$.

\hfill $\Box$

\begin{df}
\label{2.51}
 A generalized equation $\Omega$ is called {\em
periodized} with respect to a given  periodic structure $\langle
{\mathcal P}, R \rangle$ of $\Omega$ , if for every  two cycles
$c_1$ and $c_2$ with the same initial vertex  in the graph
$\Gamma$ , there is a relation
 $[h(c_1), h(c_2)]=1$  in
 $F_{R(\Omega)}$.
\end{df}

\subsubsection{ Case 1. Set $N\mathcal P$ is empty.}

Let $\Gamma_0$ be the subgraph of the graph $\Gamma$ having  the
same set of vertices and consisting of the edges $e$ whose labels
do not belong to ${\mathcal P}$. Choose a maximal subforest  $T_0$
in the graph $\Gamma_0$ and extend it to a maximal subforest $T$
of the graph $\Gamma$. Since $\langle {\mathcal P}, R \rangle$ is
connected by assumption, it follows that  $T$ is a tree. Let $v_0$
be an arbitrary vertex of the graph $\Gamma$ and $r(v_0, v)$ the
(unique) path from $v_0$ to $v$ all of whose vertices belong to
$T$. For every edge $e: v \rightarrow v'$ not lying in $T$, we
introduce a cycle $c_e = r(v_0, v) e (r(v_0, v'))^{-1}$. Then  the
fundamental group $\pi_1(\Gamma, v_0)$ is generated by the cycles
$c_e$ (see, for example, the proof of Proposition 3.2.1
\cite{LS}). This and decidability of the universal theory of a
free group imply that the property of a generalized equation ``to
be periodized with respect to a given periodic structure'' is
algorithmically decidable.

Furthermore, the set of elements

\begin{equation} \label{2.52}
\{h(e) \mid e \in T \} \cup \{h(c_e) \mid e \not \in T \}
\end{equation}

\noindent forms a basis of the free group with the set of
generators $\{h_k \mid h_k$ is an unknown lying on a closed
section from ${\mathcal P} \}$. If $\mu \in {\mathcal P}$, then
$(\beta(\mu)) = (\beta(\Delta(\mu)))$, $(\alpha(\mu)) =
(\alpha(\Delta(\mu)))$ by part f) from Definition \ref{above} and,
consequently, the word $h[\alpha(\mu), \beta(\mu)]
h[\alpha(\Delta(\mu)), \beta(\Delta(\mu))]^{-1}$ is the label of a
cycle $c'(\mu)$ from $\pi_1 (\Gamma, (\alpha(\mu)))$. Let $c(\mu)
= r(v_0, (\alpha(\mu)))c'(\mu) r(v_0, (\alpha(\mu)))^{-1}$. Then

\begin{equation} \label{2.53}
h(c(\mu)) = uh[\alpha(\mu), \beta(\mu)] h[\alpha(\Delta(\mu)),
\beta(\Delta(\mu))]^{-1} u^{-1},
\end{equation}

\noindent where $u$ is a certain word. Since $c(\mu) \in
\pi_1(\Gamma, v_0)$, it follows that $c(\mu) = b_\mu ( \{c_e \mid
e \not \in T \})$, where $b_\mu$ is a certain word in the
indicated generators  which can be effectively constructed (see
Proposition 3.2.1 \cite{LS}).

Let $\tilde{b}_\mu$ denote the image of the word $b_\mu$ in the
abelianization of $\pi (\Gamma ,v_0).$ Denote by $\widetilde{Z}$
the free abelian group consisting of formal linear combinations
$\sum_{e \not \in T} n_e \tilde{c}_e$ $(n_e \in {\bf Z})$, and by
$\widetilde{B}$ its subgroup generated by the elements
$\tilde{b}_\mu$ $(\mu \in {\mathcal P})$ and the elements
$\tilde{c}_e$ $(e \not \in T, \ h(e) \not \in {\mathcal P}).$ Let
$\widetilde{A} = \widetilde{Z} / \widetilde{B}$,
$T(\widetilde{A})$ the  torsion subgroups of the group
$\widetilde{A}$, and $\widetilde{Z}_1$ the preimage of
$T(\widetilde{A})$ in $\widetilde{Z}$. The group $\widetilde{Z} /
\widetilde{Z}_1$ is free; therefore, there exists a decomposition
of the form

\begin{equation} \label{2.54}
\widetilde{Z} = \widetilde{Z}_1 \oplus \widetilde{Z}_2, \
\widetilde{B} \subseteq \widetilde{Z}_1, \ (\widetilde{Z}_1 :
\widetilde{B}) < \infty .
\end{equation}

Note that it is possible to express effectively a certain basis
$\widetilde{\bar{c}}^{(1)}$,  $\widetilde{\bar{c}}^{(2)}$ of the
group $\widetilde{Z}$ in terms of the generators $\widetilde{c}_e$
so that for the subgroups $\widetilde{Z}_1$, $\widetilde{Z}_2$
generated by the sets $\widetilde{\bar{c}}^{(1)}$,
$\widetilde{\bar{c}}^{(2)}$ respectively, relation (\ref{2.54})
holds. For this it suffices, for instance, to look through the
bases one by one, using the fact that under the condition
$\widetilde{Z} = \widetilde{Z}_1 \oplus \widetilde{Z}_2$ the
relations $\widetilde{B} \subseteq \widetilde{Z}_1$,
$(\widetilde{Z}_1 : \widetilde{B}) <  \infty$ hold if and only if
the generators of the groups $\widetilde{B}$ and $\widetilde{Z}_1$
generate the same linear subspace over ${\bf Q}$, and the latter
is easily verified algorithmically. Notice, that a more economical
algorithm can be constructed by analyzing the proof of the
classification theorem for finitely generated abelian groups. By
Proposition 1.4.4  \cite{LS}, one can effectively construct a
basis $\bar{c}^{(1)}$, $\bar{c}^{(2)}$ of the free (non-abelian)
group $\pi_1(\Gamma, v_0)$ so that $\widetilde{\bar{c}}^{(1)}$,
$\widetilde{\bar{c}}^{(2)}$ are the natural images of the elements
$\bar{c}^{(1)}$, $\bar{c}^{(2)}$ in $\widetilde{Z}$.

Now assume that $\langle {\mathcal P}, R \rangle$ is an arbitrary
periodic structure of a periodized generalized equation $\Omega$,
not necessarily connected. Let $\Gamma_1, \ldots, \Gamma_r$ be the
connected components of the graph $\Gamma$. The labels of edges of
the component $\Gamma_i$ form in the equation $\Omega$ a union of
closed sections from ${\mathcal P}$; moreover, if a base $\mu \in
{\mathcal P}$ belongs to such a section, then its dual
$\Delta(\mu)$, by condition f) of Definition \ref{above}, also
possesses this property. Therefore, by taking for ${\mathcal P}_i$
the set of labels of edges from $\Gamma_i$ belonging to ${\mathcal
P}$, sections to which these labels belong, and bases $\mu \in
{\mathcal P}$ belonging to these sections, and restricting in the
corresponding way the relation $R$, we obtain a periodic connected
structure $\langle {\mathcal P}_i, R_i \rangle$ with the graph
$\Gamma_i$.

The notation $\langle {\mathcal P}', R'  \rangle$ $\subseteq$
$\langle {\mathcal P}, R \rangle$ means that  ${\mathcal P}'
\subseteq {\mathcal P},$ and the relation $R'$ is a restriction of
the relation $R$. In particular, $\langle {\mathcal P}_i, R_i
\rangle$ $\subseteq$ $\langle {\mathcal P}, R \rangle$ in the
situation described in the previous paragraph. Since $\Omega$ is
periodized, the periodic structure must be connected.

Let $e_1, \ldots, e_m$ be all the edges of the graph $\Gamma$ from
$T \setminus T_0$. Since $T_0$ is the spanning forest of the graph
$\Gamma_0$, it follows that $h(e_1), \ldots, h(e_m) \in {\mathcal
P}$. Let $F(\Omega)$ be a free group generated by the variables of
$\Omega$. Consider in the  group $F(\Omega )$ a new basis
$A\cup\bar x$ consisting of $A,$ variables not belonging to the
closed sections from $\mathcal P$ (we denote by $\bar t$ the
family of these variables), variables $\{h(e)|e\in T\}$ and words
$h(\bar c^{(1)}),h(\bar c^{(2)})$. Let $v_i$ be the initial vertex
of the edge $e_i$. We introduce new variables $\bar
u^{(i)}=\{u_{ie}|e\not\in T,\ e\not\in {\mathcal P}\},$ $\bar
z^{(i)}=\{z_{ie}| e\not\in T, e\not\in {\mathcal P}\}$ for $1\leq
i\leq m,$ as follows
\begin{equation}\label{2.59}
u_{ie}=h(r(v_0,v_i)^{-1}h(c_e)h(r(v_0,v_i)),\end{equation}
\begin{equation}
\label{2.60} h(e_i)^{-1}u_{ie}h(e_i)=z_{ie}.\end{equation}

Notice, that without loss of generality we can assume that $ v_0$
corresponds to the beginning of the period $P$.

\begin{lm} \label{2.10''} Let $\Omega $ be a consistent generalized
equation periodized with respect to a periodic structure $\langle
{\mathcal P},R\rangle  $ with empty set $N(\mathcal P )$.   Then
the following is true.
\begin{enumerate}\item [(1)]
One can choose the basis $\bar c^{(1)}$ so that  for any solution
$H$ of $\Omega$ periodic with respect to a period $P$ and
${\mathcal P}(H,P)=\langle {\mathcal P},R\rangle $ and any $c\in
\bar c^{(1)}$, $H(c)=P^n$, where $|n|< 2\rho$.

\item [(2)] In a fully residually free quotient of $F_{R(\Omega)}$
discriminated by  solutions from (1) the image of $\langle h(\bar
c^{(1)})\rangle $ is either trivial or a cyclic subgroup.

\item [(3)] Let $K$ be the subgroup of $F_{R(\Omega )}$ generated
by $\bar t$, $h(e), e\in T_0$, $h(\bar c^{(1)})$, $\bar u^{(i)}$
and $\bar z^{(i)}, i=1,\ldots ,m.$ If $|\bar c^{(2)}|=s\geq 1$,
then the group $F_{R(\Omega )}$ splits as a fundamental group of a
graph of groups with two vertices, where   one vertex group is $K$
and the other is a free abelian group generated  by $h(\bar
c^{(2)})$ and $h(\bar c^{(1)})$. The corresponding edge group is
generated by $h(\bar c^{(1)})$. The other edges are loops at the
vertex with vertex group $K$, have stable letters $h(e_i),\
i=1,\ldots ,m$, and associated subgroups $\langle \bar
u^{i}\rangle ,$ $\langle \bar z^{i}\rangle .$ If $\bar
c^{(2)}=\emptyset,$ then there is no vertex with abelian vertex
group.

\item [(4)]  Let $A\cup \bar x$ be the generators of the group
$F_{R(\Omega )}$ constructed above. If $e_i\in {\mathcal P}\cap
T$, then the mapping defined as $h(e_i)\rightarrow u_{ie}^kh(e_i)
$ ($k$ is any integer) on the generator $h(e_i)$ and fixing all
the other generators can be extended to an automorphism of
$F_{R(\Omega )}.$

\item [(5)]  If $c\in\bar c^{(2)}$ and   $c'$  is a cycle with
initial vertex $v_0$, then the mapping defined by $h(c)\rightarrow
h(c')^kh(c)$ and fixing all the other generators can be extended
to an automorphism of $F_{R(\Omega )}.$

\end{enumerate}

\end{lm}

{\it Proof}. To prove assertion (1) we have to show that each
simple cycle in the graph $\Gamma _0$ has length less than
$2\rho$. This is obvious, because the total number of edges in
$\Gamma _0$ is not more than $\rho$ and corresponding variables do
not  belong to $\mathcal P$.

(2) The image of the group $\langle h(\bar c^{(1)})\rangle $ in
$F$ is cyclic, therefore
 one of the finite number of equalities $h(c_1)^{n}=h(c_2)^m$,
 where $c_1,c_2\in c^{(1)},\ n,m<2\rho$
must hold for any solution. Therefore in a fully residually free
quotient the group generated by the image of $\langle\ h(\bar
c^{(1)})\rangle $ is a cyclic subgroup.

To prove (3) we are to study in more detail how the unknowns
$h(e_i)$ ($1 \leq i \leq m$) can participate in the equations from
$\Omega^\ast$ rewritten in the set of variables $\bar x\cup A.$

If $h_k$ does not lie on a closed section from ${\mathcal P}$, or
$h_k \not\in {\mathcal P}$, but $e \in T$ (where $h(e)=h_k$), then
$h_k$ belongs to the basis $\bar{x}\cup A$ and is distinct from
each of $h(e_1), \ldots, h(e_m)$. Now let $h(e) = h_k$, $h_k \not
\in {\mathcal P}$ and $e \not \in T$. Then $e=r_1c_er_2,$ where
$r_1,r_2$ are paths in $T$. Since $e \in \Gamma_0$, $h(c_e)$
belongs to $\langle c^{(1)}\rangle$ modulo commutation of cycles.
The vertices $(k)$ and $(k+1)$ lie in the same connected component
of the graph $\Gamma_0$, and hence they are connected by a path
$s$ in the forest $T_0$. Furthermore, $r_1$ and $sr_2^{-1}$ are
paths in the tree $T$ connecting the vertices $(k)$ and $v_0$;
consequently, $r_1 = s r_2^{-1}$. Thus, $e=sr_2^{-1}c_er_2$ and
$h_k = h(s) h(r_2)^{-1} h(c_e)h(r_2)$. The unknown $h(e_i)$ ($1
\leq i \leq m$) can occur in the right-hand side of the expression
obtained (written in the basis $\bar{x}\cup A$) only in $h(r_2)$
and at most once. Moreover, the sign of this occurrence (if it
exists) depends only on the orientation of the edge $e_i$ with
respect to the root $v_0$ of the tree $T$. If $r_2 = r_2' e_i^{\pm
1}r_2 ''$, then
 all the occurrences of the unknown $h(e_i)$ in the words
$h_k$ written in the basis $\bar{x}\cup A$, with $h_k \not\in
{\mathcal P}$, are contained in the occurrences of words of the
form $h(e_i)^{\mp 1} h((r_2')^{-1}c_er_2')h(e_i)^{\pm 1}$, i.e.,
in occurrences of the form $h(e_i)^{\mp 1} h(c) h(e_i)^{\pm 1}$,
where $c$ is a certain cycle of the graph $\Gamma$ starting at the
initial vertex of the edge $e_i^{\pm 1}$.

Therefore all the occurrences of $h(e_i),\ i=1,\ldots ,m$ in the
equations corresponding to $\mu\not\in\mathcal P$ are of the form
$h(e_i^{-1})h(c)h(e_i)$. Also, $h(e_i)$ does not occur in the
equations corresponding to $\mu\in\mathcal P$ in the basis
$A\cup\bar x.$ The system $\Omega ^*$ is equivalent to the
following system in the variables $ \bar x, \bar z^{(i)},\bar
u^{(i)},A, i=1,\ldots ,m$
 : equations (\ref{2.59}), (\ref{2.60}),

\begin{equation}\label{2.61}
[u_{ie_1},u_{ie_2}]=1,\end{equation}
\begin{equation}
[h(c_1),h(c_2)]=1, \ c_1,c_2\in c^{(1)}, c^{(2)},\end{equation}
and a system $\bar \psi (h(e),e\in T\setminus {\mathcal P}, h(\bar
c^{(1)}), \bar t,\bar z^{(i)},\bar u^{(i)},A)=1$, such that either
$h(e_i)$ or $\bar c^{(2)}$ do not occur in $\bar\psi $. Let
$K=F_{R(\bar\psi )}.$ Then to obtain $F_{R(\Omega )}$ we fist take
an HNN extension of the group $K$ with abelian associated
subgroups generated by $\bar u ^{(i)}$ and $\bar z^{(i)}$ and
stable letters $h(e_i)$, and then extend the centralizer of the
image of $\langle \bar c^{(1)}\rangle$ by the free abelian
subgroup generated by the images of $\bar c^{(2)}.$

 Statements (4) and (5) follow from (3).

$\Box$

We now introduce the notion of a {\em canonical group of
automorphisms corresponding to a connected periodic structure}.
\begin{df} \label{canon} In the case when the family of bases $N{\mathcal P}$ is empty
automorphisms described in Lemma \ref{2.10''} for $e_1,\ldots
,e_m\in T\setminus T_0$ and all $c_e$ for $e\in {\mathcal
P}\setminus T$ generate the {\em canonical group of automorphisms}
$P_0$ corresponding to a connected periodic structure.
\end{df}

\begin{lm}\label{2.10'}
Let $\Omega $ be a nondegenerate generalized equation with no
boundary connections, periodized with respect to the periodic
structure $\langle \mathcal P, R\rangle  $. Suppose that the set
$N\mathcal P$ is empty. Let  $H$ be a solution of $\Omega$
periodic with respect to a period $P$ and ${\mathcal
P}(H,P)=<{\mathcal P},R>$.
 Combining canonical automorphisms of
$F_{R(\Omega )}$
  one can get a solution $H^+$ of $\Omega$ with the property that
for any $h_k\in{\mathcal P}$ such that $H_k=P_2P^{n_k}P_1$ ($P_2$
and $P_1$ are an end and a beginning of $P$),
$H_k^+=P_2P^{n^+_k}P_1$, where $n_k, n_k^+ >  0$ and the numbers
$n_k^+$'s are bounded by a certain computable function $f_2(\Omega
,{\mathcal P},R)$. For all $h_k\not\in{\mathcal P}\ H_k=H_k^+.$
\end{lm}

{\em Proof.}   Let $\delta((k)) = P_1^{(k)} P_2^{(k)}$. Denote by
$t(\mu, h_k)$ the number of occurrences of the edge with label
$h_k$ in the cycle $c_{\mu}$, calculated taking into account the
orientation. Let

\begin{equation} \label{2.65}
H_k = P_2^{(k)} P^{n_k} P_1^{(k+1)}
\end{equation}

\noindent ($h_k$ lies on a closed section from ${\mathcal P}$),
where the equality in (\ref{2.65}) is graphic whenever $h_k \in
{\mathcal P}$. Direct calculations show that

\begin{equation} \label{2.66}
H(b_\mu) = P^{\sum_k t(\mu, h_k)(n_k+1)}.
\end{equation}
This equation implies that the vector $\{n_k\}$ is a solution to
the following system of Diophantine equations in variables
$\{z_k|h_k\in{\mathcal P}\}$:

\begin{equation}\label{2..}\sum _{h_k\in{\mathcal P}}t(\mu ,h_k)z_k+
\sum _{h_k\not\in{\mathcal P}}t(\mu ,h_k)n_k=0,
\end{equation}
$\mu\in {\mathcal P}.$ Note that the number of unknowns is
bounded, and coefficients of this system are bounded from above
($|n_k|\leq 2$ for $h_k\not\in {\mathcal P}$) by a certain
computable function of $\Omega, {\mathcal P},$ and $R$. Obviously,
$(P_2^{(k)})^{-1}H^+_kH^{-1}_kP_2^{(k)}=P^{n_k^+-n_k}$ commutes
with $H(c)$, where $c$ is a cycle such that $H(c)=P^{n_0},\ n_0<
2\rho$.

If system (\ref{2..}) has only one solution, then it is bounded.
Suppose it has infinitely many solutions. Then $(z_1,\ldots
,z_k,\ldots )$  is a composition of a bounded solution of
(\ref{2..}) and a linear combination of independent solutions of
the corresponding homogeneous system. Applying canonical
automorphisms from Lemma \ref{2.10''} we can decrease the
coefficients of this linear combination to obtain a bounded
solution $H^+$. Hence for $h_k=h(e_i),\ e_i\in{\mathcal P}$, the
value $H_k$ can be obtained by a composition of a canonical
automorphism (Lemma \ref{2.10''}) and a suitable bounded solution
$H^+$ of $\Omega$. $\Box$

\subsubsection{ Case 2. Set $N\mathcal P $ is non-empty.}

We construct an oriented graph $B\Gamma$ with the same set of
vertices as $\Gamma$. For each item $h_k\not\in\mathcal P$ such
that $h_k$ lie on a certain closed section from $\mathcal P$
introduce an edge $e$ leading from $(k)$ to $(k+1)$ and $e^{-1}$
leading from $(k+1)$ to $(k)$. For each pair of bases $\mu ,\Delta
(\mu )\in \mathcal P$ introduce an edge $e$ leading from $(\alpha
(\mu))=(\alpha (\Delta (\mu )))$ to $(\beta (\mu))=(\beta (\delta
(\mu )))$ and $e^{-1}$ leading from $(\beta (\mu))$ to $(\alpha
(\mu ))$. For each base $\mu\in N\mathcal P$ introduce an edge $e$
leading from $(\alpha (\mu)$ to $(\beta (\mu))$ and $e^{-1}$
leading from $(\beta (\mu))$ to $(\alpha (\mu ))$. denote by
$B\Gamma _0$ the subgraph with the same set of vertices and edges
corresponding to items not from ${\mathcal P}$ and bases from
$\mu\in N\mathcal P$. Choose a maximal subforest $BT_0$ in the
graph $B\Gamma _0$ and extend it to a
 maximal subforest
$BT$ of the graph $B\Gamma $. Since $\mathcal P$ is connected,
$BT$ is a tree. The proof of the following lemma is similar to the
proof of Lemma \ref{2.9}.

\begin{lm} \label{2.9b}
Let $H$ be a solution of a generalized equation $\Omega$ periodic
with respect to a period $P$, $\langle {\mathcal P}, R \rangle =
{\mathcal P}(H, P)$;  $c$ be a cycle in the graph $B\Gamma$ at the
vertex $(l)$; $\delta(l)=P_1P_2$. Then there exists $n \in {\bf
Z}$ such that $H(c) = (P_2P_1)^n.$
\end{lm}

As we did in the graph $\Gamma$, we choose a vertex $v_0$. Let
$r(v_0,v)$ be the unique path in $BT$ from $v_0$ to $v$. For every
edge $e=e(\mu ):v\rightarrow v'$  not lying in $BT$, introduce a
cycle $c_{\mu}= r(v_0,v)e(\mu)r(v_0,v')^{-1}$. For every edge
$e=e(h_k ):V\rightarrow V'$  not lying in $BT$, introduce a cycle
$c_{h_k}= r(v_0,v)e(h_k)r(v_0,v')^{-1}$.

It suffices to restrict ourselves to the case of a connected
periodic structure. If $e=e(h_k)$, we denote $h(e)=h_k$; if
$e=e(\mu )$, then $h(e)=\mu $. Let $e_1, \ldots, e_m$ be all the
edges of the graph $B\Gamma$ from $BT \setminus BT_0$. Since
$BT_0$ is the spanning forest of the graph $B{\Gamma} _0$, it
follows that $h(e_1), \ldots, h(e_m) \in {\mathcal P}$. Consider
in the free group $F(\Omega )$ a new basis $A\cup\bar x$
consisting of $A$, items  $h_k$ such that $h_k$ does not belong to
closed sections from $\mathcal P$ (denote this set by $\bar t$),
variables $\{h(e)|e\in T\}$ and words from $h( C^{(1)}),h(
C^{(2)})$, where the set $ C^{(1)}, C^{(2)}$ form a basis of the
free group $\pi (B\Gamma ,v_0)$, $C^{(1)}$ correspond to the
cycles that represent the identity in $F_{R(\Omega  )}$ (if $v$
and $v'$ are initial and terminal vertices of some closed section
in $\mathcal P$ and  $r$ and $r_1$ are different paths from $v$ to
$v'$, then $r(v_0,v)rr_1^{-1}r(v_0,v)^{-1}$ represents the
identity),  cycles $c_{\mu},\mu\in N\mathcal P$ and $c_{h_k}$,
$h_k\not\in \mathcal P$;
 and $C^{(2)}$ contains the rest of the basis of $\pi (B{\Gamma },v_0)$.

We study in more detail how the unknowns $h(e_i)$ ($1 \leq i \leq
m$) can participate in the equations from $\Omega^\ast$ rewritten
in this basis.

If $h_k$ does not lie on a closed section from ${\mathcal P}$, or
$h_k=h(e), h(\mu)=h(e) \not\in {\mathcal P}$, but $e \in T$, then
$h(\mu)$ or $h_k$  belongs to the basis $\bar{x}\cup {A}$ and is
distinct from each of $h(e_1), \ldots, h(e_m)$. Now let $h(e) =
h(\mu)$, $h(\mu) \not \in {\mathcal P}$ and $e \not \in T$. Then
$e=r_1c_er_2$, where $r_1,r_2$ are path in $BT$ from $(\alpha (\mu
))$ to $v_0$ and from $(\beta (\mu ))$ to $v_0$. Since $e \in
B\Gamma_0$, the vertices $(\alpha (\mu))$ and $(\beta (\mu))$ lie
in the same connected component of the graph $B\Gamma _0$, and
hence are connected by a path $s$ in the forest $BT_0$.
Furthermore, $r_1$ and $sr_2^{-1}$ are paths in the tree $BT$
connecting the vertices $(\alpha (\mu ))$ and $v_0$; consequently,
$r_1 = s r_2^{-1}$. Thus, $e=sr_2^{-1}c_er_2$ and $h(\mu ) = h(s)
h(r_2)^{-1} h(c_e)h(r_2)$. The unknown $h(e_i)$ ($1 \leq i \leq
m$) can occur in the right-hand side of the expression obtained
(written in the basis $\bar{x}\cup A$) only in $h(r_2)$ and at
most once. Moreover, the sign of this occurrence (if it exists)
depends only on the orientation of the edge $e_i$ with respect to
the root $v_0$ of the tree $T$. If $r_2 = r_2' e_i^{\pm 1}r_2 ''$,
then
 all the occurrences of the unknown $h(e_i)$ in the words
$h (\mu)$ written in the basis $\bar{x}\cup A$, with $h(\mu)
\not\in {\mathcal P}$, are contained in the occurrences of words
of the form $h(e_i)^{\mp 1} h((r_2')^{-1}c_er_2')h(e_i)^{\pm 1}$,
i.e., in occurrences of the form $h(e_i)^{\mp 1} h(c) h(e_i)^{\pm
1}$, where $c$ is a certain cycle of the graph $B\Gamma$ starting
at the initial vertex of the edge $e_i^{\pm 1}$. Similarly, all
the occurences of the unknown $h(e_i)$ in the words $h_k$ written
in the basis $\bar{x}, A$, with $h_k \not\in {\mathcal P}$, are
contained in  occurrences of words of the form $h(e_i)^{\mp 1}
h(c) h(e_i)^{\pm 1}$.

Therefore all the occurences of $h(e_i),\ i=1,\ldots ,m$ in the
equations corresponding to $\mu\not\in {\mathcal P}$ are of the
form $h(e_i^{-1})h(c)h(e_i)$. Also,
 cycles from $C^{(1)}$ that represent the identity and not in $B\Gamma _0$
 are basis elements themselves.
This implies
\begin{lm} \label{2n}
\begin{enumerate}
\item [(1)] Let $K$ be the subgroup of $F_{R(\Omega )}$ generated
by $\bar t$, $h(e), e\in BT_0$, $h(C^{(1)})$ and $\bar u^{(i)},
\bar z^{(i)}, i=1,\ldots ,m,$ where elements $\bar z^{(i)}$ are
defined similarly to the case of empty $\mathcal NP$.

If $|C^{(2)}|=s\geq 1$, then the group $F_{R(\Omega )}$   splits
as a fundamental group of a graph of groups with two vertices,
where one vertex group is $K$ and the other is a free abelian
group generated  by $h(C^{(2)})$ and  $ h(C^{(1)})$. The edge
group is generated by $h(C^{(1)})$. The other edges are loops at
the vertex with vertex group $K$ and have stable letters $h(e),
e\in BT\setminus BT_0$. If $C^{(2)}=\emptyset,$ then there is no
vertex with abelian vertex group.

\item [(2)] Let $H$ be a solution of $\Omega$ periodic with
respect to a period $P$ and $\langle {\mathcal P},R\rangle =
{\mathcal P}(H,P).$ Let $P_1P_2$ be a partition of $P$
corresponding to the initial vertex of $e_i$. A transformation
$H(e_i)\rightarrow P_2P_1H(e_i)$, $i\in\{1,\ldots ,m\}$,  which is
identical on all the other elements from $A, H(\bar x)$, can be
extended to another solution of $\Omega ^*$. If $c$ is a cycle
beginning at the initial vertex of $e_i$, then the transformation
$h(e_i)\rightarrow h(c)h(e_i)$ which is identical on all other
elements from $A\cup\bar x$, is an automorphism of
$F_{R(\Omega)}.$

\item [(3)]  If $c(e)\in C^{(2)}$, then the transformation
$H(c(e))\rightarrow PH(c(e))$ which is identical on all other
elements from $A, H(\bar x)$ , can be extended to another solution
of $\Omega ^*$. A transformation $h(c(e))\rightarrow h(c)h(c(e))$
which is identical on all other elements from $A\cup\bar x$, is an
automorphism of $F_{R(\Omega )}.$\end{enumerate}
\end{lm}

\begin {df} If $\Omega $ is a nondegenerate   generalized equation
periodic with respect to a connected periodic structure $\langle
\mathcal P, R\rangle  $ and the set $N\mathcal P$ is non-empty, we
consider the group ${\bar A}(\Omega )$  of transformations  of
solutions of $\Omega ^*$, where $\bar A(\Omega )$ is generated by
the transformations  defined in Lemma \ref{2n}. If these
transformations are automorphisms, the group will be denoted
$A(\Omega ).$
\end{df}

\begin{df} In the case when for a connected periodic structure $\langle \mathcal P,R\rangle  ,$ the set
$C^{(2)}$ has more than one element or $C^{(2)}$ has one element,
and $C^{(1)}$ contains a cycle formed by edges $e$ such that
variables $h_k=h(e)$ are not from $\mathcal P$,  the periodic
structure will be called {\em singular}.
\end{df} This definition coincides with the definition of singular
periodic structure given in \cite{KMIrc}) in the case of empty set
$\Lambda$.

Lemma \ref{2n} implies the following
\begin{lm}\label{2.11}
Let $\Omega $ be a nondegenerate generalized equation with no
boundary connections, periodized with respect to a singular
periodic structure $\langle \mathcal P, R\rangle  $. Let $H$ be a
solution of $\Omega$ periodic with respect to a period $P$ and
$\langle {\mathcal P},R\rangle = {\mathcal P}(H,P).$ Combining
canonical automorphisms from $A(\Omega )$ one can get a solution
$H^+$ of $\Omega ^*$ with the following properties:

1) For any $h_k\in {\mathcal P}$ such that $H_k=P_2P^{n_k}P_1$
($P_2$ and $P_1$ are an end and a beginning of $P$)
$H_k^+=P_2P^{n^+_k}P_1$, where $n_k, n_k^+\in {\mathbb Z}$;

2) For any $h_k\not\in {\mathcal P}$,  $H_k=H^+_k$;

3) For any base $\mu\not\in{\mathcal P}$, $H(\mu )=H^+(\mu )$;

4) There exists a  cycle $c$ such that $h(c)\neq 1$ in
$F_{R(\Omega )}$ but $H^+(c)=1.$


\end{lm}

Notice, that in the case described in the lemma, solution $H^+$
satisfies a proper equation. Solution $H^+$ is not necessarily a
solution of the generalized equation $\Omega $, but we will modify
$\Omega$ into a generalized equation $\Omega ({\mathcal P},BT)$.
This modification will be called the {\bf first minimal
replacement }. Equation  $\Omega ({\mathcal P},BT)$ will have the
following properties:

(1)
 $\Omega ({\mathcal P},BT)$ contains all the same parameter sections and
 closed sections which are not in $\mathcal P$, as $\Omega
 $;

(2) $H^+$ is a solution of $\Omega ({\mathcal P},BT)$;

(3) group $F_{R(\Omega ({\mathcal P},BT) )}$ is generated by the
same set of variables $h_1,\ldots ,h_{\delta}$;

(4)  $\Omega ({\mathcal P},BT)$ has the same set of bases as
$\Omega$ and possibly some new bases, but each new base is a
product of bases from $\Omega $;

(5) the mapping $h_i\rightarrow h_i$ is a proper homomorphism from
$F_{R(\Omega  )}$ onto $F_{R(\Omega ({\mathcal P},BT) )}$.

To obtain $\Omega ({\mathcal P},BT)$ we have to modify the closed
sections from $\mathcal P$.

The label of each cycle in $B\Gamma$ is a product of some bases
$\mu_1\ldots \mu _k$. Write a generalized equation
$\widetilde\Omega$ for the equations that say that $\mu _1\ldots
\mu _k=1$ for each cycle from $C^{(1)}$ representing the trivial
element  and for each cycle from $C^{(2)}$. Each $\mu _i$ is a
product $\mu_i=h_{i1}\ldots h_{it}$. Due to the first statement of
Lemma \ref{2.11}, in each product $H^+_{ij}H^+_{i,j+1}$ either
there is no cancellations between $H^+_{ij}$ and $H^+_{i,j+1}$, or
one of them is completely cancelled in the other. Therefore the
same can be said about each pair $H^+(\mu _i)H^+(\mu _{i+1})$, and
we can make a cancellation table without cutting items or bases of
$\Omega $.

Let $\widehat\Omega $ be a generalized equation obtained from
$\Omega$ by deleting bases from $\mathcal P\cup N\mathcal P$ and
items from $\mathcal P$ from the closed sections from $\mathcal
P$. Take a union of $\widetilde\Omega$ and $\widehat\Omega$ on the
disjoint set of variables, and add basic equations identifying in
$\widehat\Omega$ and $\widetilde\Omega$ the same bases that don't
belong to $\mathcal P$. This gives us $\Omega ({\mathcal P},BT).$

Suppose that
 $C^{(2)}$ for the equation $\Omega $ is either empty or contains one cycle.
 Suppose also that for each closed section from $\mathcal P$ in $\Omega $
 there exists a base $\mu$ such that the initial boundary of this section is
 $\alpha (\mu )$ and the terminal boundary is $\beta (\Delta (\mu ))$.

\begin{lm}\label{4n} Suppose that the generalized equation $\Omega $
is periodized  with respect to a non-singular periodic structure
$\mathcal P$. Then for any periodic solution $H$ of $\Omega$ we
can choose a tree BT, some set of variables $S=\{h_{j_1},\ldots ,
h_{j_s}\}$ and a solution $H^+$ of $\Omega$ equivalent to $H$ with
respect to the group of canonical transformations $\bar A(\Omega
)$ in such a way that each of the bases $\lambda _i\in BT\setminus
BT_0$ can be represented as $\lambda_i=\lambda _{i1}h_{k_i}\lambda
_{i2}$, where $h_{k_i}\in S$ and for any $h_j\in S,$ $\mid
H^+_{j}\mid < f_3\mid P\mid$, where $f_3$ is some constructible
function depending on $\Omega $. This representation gives a new
generalized equation $\Omega '$ periodic with respect to a
periodic structure $\mathcal P'$ with the same period $P$ and all
$h_{j}\in S$  considered as variables not from $\mathcal P'$. The
graph $B\Gamma '$ for the periodic structure $\mathcal P'$ has the
same set of vertices as $B\Gamma $, has empty set $C^{(2)}$ and
$BT'=BT'_0$.

Let $c$ be a cycle from $C^{(1)}$ of minimal length, then
$H(c)=P^{n_c},$ where $|n_c|\leq 2\rho$ . Using canonical
automorphisms from $A(\Omega )$ one can transform any solution $H$
of $\Omega$ into a solution $H^{+}$ such that for any $h_j\in S,$
$\mid H^{+}_j\mid\leq f_3d\mid c mid .$ Let $\mathcal P '$ be a
periodic structure, in which all $h_i\in S$ are considered as
variables not from $\mathcal P '$, then $B\Gamma '$ has empty set
$C^{(2)}$ and $BT'=BT_0'$.
\end{lm}

{\em Proof.}   Suppose first that $C^{(2)}$ is empty. We prove the
statement of the lemma by induction on the number of edges in
$BT\setminus BT_0$. It is true, when this set is empty. Consider
temporarily all   the edges in $BT\setminus BT_0$ except one edge
$e(\lambda)$ as edges corresponding to bases from $N\mathcal P$.
Then the difference between $BT_0$ and $BT$ is one edge.

Changing $H(e(\lambda))$ by a transformation from $\bar A(\Omega
)$ we can  change only $H(e')$ for  $e'\in B\Gamma$ that could be
included into $BT\setminus BT_0$ instead of $e$. For each base
$\mu\in N\mathcal P$, $H(\mu)=P_2(\mu)P^{n(\mu)}P_1(\mu)$, for
each base $\mu\in \mathcal P$,
$H(\mu)=P_2(\mu)P^{x(\mu)}P_1(\mu)$. For each cycle $c$ in
$C^{(1)}$ such that $h(c)$ represents the identity element we have
a linear equation in variables $x(\mu)$ with coefficients
depending on $n(\mu)$. We also know that this system has a
solution for arbitrary $x(\lambda)$ (where $\lambda\in BT\setminus
BT_0$) and the other $x(\nu )$ are uniquely determined by the
value of $x(\lambda)$.

If we write for each variable $h_k\in {\mathcal P},$
$H_k=P_{2k}P^{y_k}P_{1k}$, then the positive unknowns $y_k$'s
satisfy the system of equations saying that $H(\mu)=H(\Delta (\mu
))$ for bases $\mu\in\mathcal P$ and equations saying that $\mu$
is a constant for bases $\mu\in N\mathcal P$. Fixing $x(\lambda)$
we automatically fix all the $y_k$'s. Therefore at least one of
the $y_k$ belonging to $\lambda$ can be taken arbitrary.   So
there exist some elements $y_k$ which can be taken as free
variables for the second system of linear equations. Using
elementary transformations over $\mathbb Z$ we can write the
system of equations for $y_k$'s in the form:
\begin{equation}\label{linear}
\begin{array}{cccc}
n_1y_1         & 0             & \cdots                  &=m_1y_k       +C_1             \\
             & n_2y_2         &   \cdots               &=m_2y_k     +C_2   \\

             &              & \ddots                      &    \\
\vdots   & \vdots   &                  & \ddots                  \\
             & \cdots   & n_{k-1}y_{k-1}       & =m_{k-1}y_k               +C_{k-1},         \\
\end{array}
\end{equation}
where $C_1,\ldots C_k$ are constants depending on parameters, we
can suppose that they are sufficiently large positive or negative
(small constants we can treat as constants not depending on
parameters). Notice that integers $n_1,m_1,\ldots
,n_{k-1},m_{k-1}$ in this system do not depend on parameters.  We
can always suppose that all $n_1,\ldots ,n_{k-1}$ are positive.
Notice that $m_i$ and $C_i$ cannot be simultaneously negative,
because in this case it would not be a positive solution of the
system. Changing the order of the equations we can write first all
equation with $m_i,C_i$ positive, then equations with negative
$m_i$ and positive $C_i$ and, finally,  equations with negative
$C_i$ and positive $m_i$.  The system will have  the form:
\begin{equation}\label{linear1}
\begin{array}{cccc}
n_1y_1         & 0             & \cdots                  &=|m_1|y_k       +|C_1|,             \\
               &              & \ddots       &                \\

& n_ty_t         &   \cdots               &=-|m_t|y_k     +|C_t|,   \\

             &              & \ddots       &                   \\
  &    &                  & \ddots                  \\
             &    & n_{s}y_{s}       & =|m_{s}|y_k               -|C_{s}|         \\
\end{array}
\end{equation}

If the last block (with negative $C_s$) is non-empty, we can take
a minimal $y_s$ of bounded value.
 Indeed, instead of $y_s$
we can always take a remainder of the division of $y_s$ by the
product \newline$n_1\ldots n_{k-1}|m_1\ldots m_{k-1}|$, which is
less than this product (or by the product $n_1\ldots
n_{k-1}|m_1\ldots m_{k-1}|n_c$ if we wish to decrease $y_s$ by a
multiple of $n_c$). We respectively decrease $y_k$ and adjust
$y_i$'s in the blocks with positive $C_i$'s. If the third block is
not present, we decrease $y_k$ taking a remainder of the division
of $y_k$ by $n_1\ldots n_{k-1}$ (or by $n_1\ldots n_{k-1}n_c$) and
adjust $y_i$'s. Therefore for some $h_i$ belonging to a base which
can be included into $BT\setminus BT_0$,  $\mid H^+(h_i)\mid <
f_3\mid P\mid .$ Suppose this base is $\lambda$, represent
$\lambda =\lambda _1h_i\lambda _2$. Suppose $e(\lambda
):v\rightarrow v_1$ in $B\Gamma$. Let $v_2,v_3$ be the vertices in
$B\Gamma$ corresponding to the initial and terminal boundary of
$h_k$. They would be the vertices in $\Gamma$, and $\Gamma$ and
$B\Gamma$ have the same set of vertices. To obtain the graph
$B\Gamma '$ from $B\Gamma$ we have to replace $e(\lambda)$ by
three edges $e(\lambda _1):v\rightarrow v_2$,
$e(h_k):v_2\rightarrow v_3$ and $e(\lambda _2):v_3\rightarrow
v_1$. There is no path in $BT_0$ from $v_2$ to $v_3$, because if
there were such  a path $p$, then we would have the equality
$h_k=h(c_1)h(p)h(c_2),$ in $F_{R(\Omega
 )},$ where $c_1$ and $c_2$ are cycles in $B\Gamma $ beginning in
vertices $v_2$ and $v_3$ respectively. Changing $H_k$ we do not
change $H(c_1),H(c_2)$ and $H(p)$, because all the cycles are
generated by cycles in $C^{(1)}$. Therefore there are paths
$r:v\rightarrow v_2$ and $r_1:v_3\rightarrow v_1$ in $BT_0,$ and
edges $e(\lambda _1), e(\lambda _2)$ cannot be included in
$BT'\setminus BT'_0$ in $B\Gamma '$. Therefore $BT'=BT_0'$. Now we
can recall that all the edges except one in $BT\setminus BT_0$
were temporarily considered as edges in $N\mathcal P$. We managed
to decrease the number of such edges by one. Induction finishes
the proof.

If the set $C^{(2)}$ contains one cycle, we can temporarily
consider all the bases from $BT$ as parameters, and consider the
same system of linear equations for $y_i$'s. Similarly, as above,
at  least one $y_t$ can be bounded. We will bound as many  $y_i$'s
as we can. For the new periodic structure either $BT$ contains
less elements or the set $C^{(2)}$ is empty.

 The second part of the lemma follows from the remark that for
$\mu\in T$ left multiplication of $h(\mu )$ by $h(rcr^{-1})$,
where $r$ is the path in $T$ from $v_0$ to the initial vertex of
$\mu$, is an automorphism from $A(\Omega )$. $\Box$

We call a solution $H^+$ constructed in Lemma \ref{4n} a {\em
solution equivalent to $H$ with  maximal number of short
variables}.

Consider now variables from $S$ as variables not from $\mathcal
P'$, so that for the equation $\Omega $ the sets $C^{(2)}$ and
$BT'\setminus BT_0'$ are both empty. In this case we make the {\em
second minimal replacement}, which we will describe in the lemma
below.
\begin{df} \label{overl} A pair of bases $\mu ,\Delta (\mu)$ is called an
overlapping pair if $\epsilon (\mu )=1$ and $\beta (\mu )> \alpha
(\Delta (\mu ))>  \alpha (\mu )$ or $\epsilon (\mu)=-1$ and $\beta
(\mu )< \beta (\Delta (\mu ))< \alpha (\mu )$. If a closed section
begins with $\alpha (\mu)$ and ends with $\beta (\Delta (\mu ))$
for an overlapping pair of bases we call such a pair of bases a
{\em principal overlapping pair} and say that a section is in {\em
overlapping form}.\end{df}

Notice, that if $\lambda\in N\mathcal P$, then $H(\lambda )$ is
the same for any solution $H$, and we just write $\lambda$ instead
of $H(\lambda )$.

\begin{lm} \label{3n'} Suppose that for the generalized equation $\Omega '$ obtained in
Lemma \ref{4n} the sets $C^{(2)}$ and $BT'\setminus BT_0'$ are
empty, $\mathcal P '$ is a non-empty periodic structure, and each
closed section from $\mathcal P '$ has a principal overlapping
pair. Then for each base $\mu\in\mathcal P '$ there is a  fixed
presentation for $h(\mu)=\prod (parameters)$ as a  product of
elements $h(\lambda ), \lambda\in N\mathcal P$,
$h_k\not\in\mathcal P'$ corresponding to a path in $B\Gamma _0'$.
The maximal number of terms in this presentation is bounded by a
computable function of $\Omega$.
\end{lm}

 {\em Proof.}   Let $e$ be the edge in the graph $B\Gamma'$
corresponding to a base $\mu$ and suppose $e:v\rightarrow v'$.
There is a path $s$ in $BT'$ joining $v$ and $v'$ ,  and a cycle
$\bar c$ which is a product of cycles from $C^{(1)}$ such that
$h(\mu )=h(\bar c )h(s).$ For each cycle $c$ from $C^{(1)}$ either
$h(c)=1$ or $c$ can be written using only edges with labels not
from $\mathcal P'$; therefore, $\bar c$ contains only edges with
labels not from $\mathcal P '$. Therefore
\begin{equation}\label{cat}
h(\mu )=\prod (parameters)=h(\lambda _{i_1})\Pi _1\ldots h(\lambda
_{s_i})\Pi _s,\end{equation}
 where the doubles of all
$\lambda _i$ are parameters, and $\Pi _1,\ldots ,\Pi _s$ are
products of variables $h_{k_i}\not \in \mathcal P '$. $\Box$

In the equality
\begin{equation}\label{cat1}
H(\mu )=H(\lambda _{i_1})\bar \Pi _1\ldots H(\lambda
_{s_i})\bar\Pi _s,\end{equation} where $\bar\Pi _1,\ldots ,\bar\Pi
_s$ are products of $H_{k_i}$ for variables $h_{k_i}\not \in
\mathcal P '$, the cancellations between two terms in the right
side are complete because the equality corresponds to a path in
$B\Gamma _0'$. Therefore the cancellation tree for the equality
(\ref{cat1}) can be situated on a horizontal axis with intervals
corresponding to $\lambda _i$'s directed either to the right or to
the left. This tree can be drawn on a $P$-scaled axis. We call
this one-dimensional tree a $\mu$-{\em tree}. Denote by
$I(\lambda)$ the interval corresponding to $\lambda$ in the
$\mu$-tree. If $I(\mu)\subseteq\bigcup _{\lambda _i\in N{\mathcal
P}}I(\lambda _i),$ then we say that $\mu$ is covered by
parameters. In this case a generalized equation corresponding to
(\ref{cat1}) can be situated on the intervals corresponding to
bases from $N{\mathcal P}$.

We can shift the whole $\mu$-tree to the left or to the right so
that in the new situation the uncovered part becomes covered by
the bases from $N\mathcal P$. Certainly, we have to make sure that
the shift is through the interval corresponding to a cycle in $
C^{(1)}$. Equivalently, we can shift any base belonging to the
$\mu$-tree through such an interval.

If $c$ is a  cycle from $C^{(1)}$ with shortest $H(c)$, then there
is a corresponding $c$-tree. Shifting this $c$-tree to the right
or to the left through the intervals corresponding to $H(c)$
bounded number of times we can cover every $H_i$, where $h_i\in S$
by a product $H(\lambda _{j_1})\bar\Pi_1\ldots H(\lambda
_{j_t})\bar \Pi _t$, where $\bar\Pi_ 1,\ldots ,\bar\Pi _t$ are
products of values of variables not from $\mathcal P$ and $\lambda
_{j_1},\ldots \lambda _{j_t}$ are bases from $N\mathcal P$.
Combining this covering together with the covering of $H(\mu )$ by
the product (\ref{cat1}), we obtain that
 $H([\alpha (\mu ),\beta (\Delta (\mu ))])$ is almost
covered by parameters, except for the short products $\bar \Pi$.
Let $h(\mu )$ be covered by \begin{equation}\label{cat2} h(\Lambda
_1)\Pi _1\ldots ,h(\Lambda _s)\Pi _s,
\end{equation}
where $h(\Lambda _1),\ldots ,h(\Lambda _s)$ are parts completely
covered by parameters,  and $\Pi _1,\ldots ,\Pi _s$ are products
of variables not in $\mathcal P$. We also remove those bases from
$N\mathcal P$ from each $\Lambda _i$ which do not overlap with
$h(\mu )$. Denote by $f_4$ the maximal number of bases in
$N\mathcal P$ and $h_i\not \in\mathcal P$ in the covering
(\ref{cat2}).

If $\lambda _{i_1},\ldots ,\lambda _{i_s}$ are parametric bases,
then for any solution $H$ and any pair $\lambda _i,\lambda _j\in
\{\lambda _{i_1},\ldots ,\lambda _{i_s}\}$ we have either $\mid
H(\lambda _i)\mid <\mid H(\lambda _j)\mid$ or  $\mid H(\lambda
_i)\mid =\mid H(\lambda _j)\mid$ or $\mid H(\lambda _i)\mid >\mid
H(\lambda _j)\mid$. We call a {\em relationship between lengths of
parametric bases} a collection that consists of one such
inequality or equality for each pair of bases. There is only a
finite number of possible relationships between lengths of
parametric bases. Therefore we can talk about a parametric base
$\lambda$ of maximal length meaning that we consider the family of
solutions for which $H(\lambda )$ has maximal length.

\begin{lm}\label{3n}Let $\lambda _{\mu}\in{N\mathcal P}$ be a base
of max length in the covering (\ref{cat2})
 for $\mu\in\mathcal P$. If  for a solution $H$ of $\Omega ,$ and for
each closed section $[\alpha (\mu ),\beta (\Delta (\mu)]$ in
$\mathcal P$,
 min$\ \mid H[\alpha (\nu),\alpha (\Delta (\nu))]\mid\leq \mid H(\lambda _{\mu })\mid,$
 where the minimum is taken for all pairs of overlapping bases for this section,
then one can transform $\Omega$ into one of the finite number
(depending on $\Omega $) of generalized equations $\Omega
(\mathcal P)$ which do not contain  closed sections from $\mathcal
P$ but contain the same other closed sections  except for
parametric sections. The content of  closed sections from
$\mathcal P$ is transferred using bases from $N\mathcal P$ to the
parametric part. This transformation is called the {\bf second
minimal replacement}.\end{lm}

{\em Proof.}   Suppose for a closed section $[\alpha (\mu),\beta
(\Delta (\mu ))]$ that there exists a base $\lambda$ in
(\ref{cat2}) such that $\mid H(\lambda) \mid\geq min \ (H(\alpha
(\nu ),\alpha (\Delta (\nu ))),$ where the minimum is taken for
all pairs of overlapping bases for this section. We can shift the
cover  $H(\Lambda _1)\bar\Pi _1,\ldots ,H(\Lambda _s)\bar\Pi _s$
through the distance $d_1=\mid H[\alpha (\mu ),\alpha (\Delta (\mu
))]\mid .$ Consider first the case when $d_1\leq \mid H(\lambda
)\mid$ for the largest base in (\ref{cat2}). Suppose the part of
$H(\mu )$ corresponding to $\bar\Pi _i$ is not covered by
parameters. Take the first base $\lambda _j$ in (\ref{cat2}) to
the right or to the left of $\bar \Pi _i$ such that $\mid
H(\lambda _j)\mid\geq d_1$. Suppose $\lambda _j$ is situated to
the left from $\bar\Pi _i$. Shifting $\lambda _j$ to the right
through a bounded by $f_4$ multiple of $d_1$ we will cover
$\bar\Pi _i$.

Consider  now the case when $d_1>  \mid H(\lambda )\mid$, but
there exists an overlapping pair $\nu ,\Delta (\nu)$ such that
$$d_2=\mid H[\alpha (\nu ),\alpha (\Delta (\nu ))]\mid\leq \mid H(\lambda ) \mid.$$
If the part of $H(\mu )$ corresponding to $\bar\Pi _i$ is not
covered by parameters, we take the first base $\lambda _j$ in
(\ref{cat2}) to the right or to the left of $\bar \Pi _i$ such
that $\mid H(\lambda _j)\mid\geq d_2$. Without loss of generality
we can suppose that $\lambda _j$ is situated to the left of
$\bar\Pi _i$. Shifting $\lambda _j$ to the right through a bounded
by $f_4$ multiple of $d_2$  we will cover $\bar\Pi _i$.

Therefore, if the first alternative in the lemma does not take
place, we can cover the whole section $[\alpha (\mu),\beta (\Delta
(\mu ))]$ by the bases from $N{\mathcal P}$, and transform $\Omega
$ into one of the finite number of generalized equations which do
not contain the closed section $[\alpha (\mu),\beta (\Delta (\mu
))]$ and have all the other non-parametric sections the same. All
the cancellations between two neighboring terms of any equality
that we have gotten are complete, therefore the coordinate groups
of new equations are quotients of $F_{R(\Omega )}$. $\Box$

\subsection{Minimal solutions
and tree $T_0(\Omega )$} \label{se:5.3}
\subsubsection{Minimal solutions}\label{5.5.1}
Let $F = F(A\cup B)$ be a free group with basis $A\cup B$,
$\Omega$  be a generalized equation with constants from $(A\cup
B)^{\pm 1}$, and parameters $\Lambda$. Let $A(\Omega )$ be an
arbitrary  group of $(A\cup \Lambda )$-automorphisms of
$F_{R(\Omega )}$. For   solutions  $ H^{(1)}$ and $H^{(2)}$ of the
equation $\Omega$ in the group $F$ we write $H^{(1)}<  _{A(\Omega
)} H^{(2)}$ if there exists an endomorphism $\pi$ of the group $F$
which is an
\newline $(A,\Lambda )$-homomorphism, and an automorphism
$\sigma\in A(\Omega )$  such that the following conditions hold:
(1)  $\pi _{ H^{(2)}}=\sigma\pi _{ H^{(1)}}\pi $, (2) For all
active variables $d(H_k^{(1)})\leq d(H_k^{(2)})$ for all $1\leq
k\leq\rho$ and $d(H_k^{(1)})< d(H_k^{(2)})$ at least for one such
$k$.

We also define a relation $< _{cA(\Omega )}$ by the same way as $<
_{A(\Omega )}$ but with extra property: (3) for any $k,j$, if
$(H_k^{(2)})^{\epsilon}(H_j^{(2)})^{\delta}$ in non-cancellable,
then $(H_k^{(1)})^{\epsilon}(H_j^{(1)})^{\delta}$ in
non-cancellable ($\epsilon ,\delta =\pm 1$).  Obviously, both
relations are transitive.

A solution $\bar H$ of $\Omega$ is called {\em $A(\Omega
)$-minimal} if there is no any solution $ \bar H^\prime$ of the
equation $\Omega$ such that $\bar H^\prime < _{A(\Omega )} \bar
H.$
 Since the total length $\sum_{i = 1}^\rho l(H_i)$ of  a solution $\bar H$ is a
non-negative integer, every strictly decreasing chain of solutions
$\bar H >   \bar H^1 >   \ldots >   \bar H^k > _{A(\Omega )}
\ldots $ is finite. It follows that for every solution $\bar H$ of
$\Omega$ there exists a minimal solution $\bar H^0$ such that
$\bar H^0 < _{A(\Omega )}\bar H$.

\subsubsection{Automorphisms}\label{5.5.2}Assign to  some vertices $v$ of the tree $T(\Omega)$
the groups of automorphisms  of groups $F_{R(\Omega _v)}$. For
each vertex $v$ such that $tp (v)=12$ the canonical group of
automorphisms $A (\Omega  _v)$ assigned to it is the group of
automorphisms of $F_{R(\Omega _v)}$ identical on $\Lambda .$ For
each vertex $v$ such that $7\leq tp (v)\leq 10$ we assign the
group of automorphisms invariant with respect to the kernel.

For each vertex $v$ such that $tp (v)=2,$ assign the group $\bar
A_v$ generated by the groups of
 automorphisms constructed in Lemma
\ref{2n} that applied to $\Omega _v$ and all possible non-singular
periodic structures of this equation.

Let $tp (v)=15$. Apply transformation $D_3$ and consider
$\Omega=\tilde \Omega _v$. Notice that the function $\gamma _i$ is
constant when $h_i$ belongs to some closed section of $\tilde
{\Omega _v}$. Applying $D_2$, we can suppose that the section
$[1,j+1]$ is covered exactly twice. We say now that this is a
quadratic section. Assign to the vertex $v$ the group of
automorphisms of $F_{R(\Omega )}$ acting identically on the
non-quadratic part.

\subsubsection{The finite subtree $T_0(\Omega )$: cutting off long branches}\label{5.5.3}
For a generalized equation $\Omega $ with parameters we construct
a finite tree $T_0(\Omega)$. Then we show that the subtree of
$T(\Omega)$  obtained by tracing those path in $T(\Omega)$ which
actually can happen for ``short'' solutions is a subtree of
$T_0(\Omega )$.

According to Lemma \ref{3.2}, along an infinite path in
$T(\Omega)$ one can  either have $7\leq tp(v_k)\leq 10$ for all
$k$ or $tp(v_k)=12$ for all $k$, or $tp(v_k)=15$ for all $k$.
\begin{lm}\label{3.3}[Lemma 15 from \cite{KMIrc}] Let $v_1\rightarrow v_2
\rightarrow\ldots\rightarrow v_k\rightarrow\ldots $ be an infinite
path in the tree $T(\Omega )$, and $7\leq tp(v_k)\leq 10$ for all
$k$. Then among $\{\Omega _k\}$ some generalized equation occurs
infinitely many times. If $\Omega _{v_k}=\Omega _{v_l}$, then $\pi
(v_k,v_l)$ is an isomorphism invariant with respect to the kernel.
\end{lm}

\begin{lm} Let $tp(v)=12.$ If a solution $\bar H$ of $\Omega _v$
is minimal with respect to the canonical group of automorphisms,
then there is a recursive function $f_0$ such that in the sequence
\begin{equation}\label{12'}
 (\Omega _v,\bar H)\rightarrow (\Omega _{v_1},\bar H^1)\rightarrow\ldots
 \rightarrow (\Omega _{v_i},\bar H ^i),\ldots,\end{equation}
corresponding to the path in $T(\Omega _v)$ and for the solution
$\bar H$, case 12 cannot be repeated more than $f_0$ times.
\end{lm} {\em Proof.}   If $\mu $ and $\Delta\mu$ both belong to the
quadratic section, then $\mu$ is called a {\em quadratic base}.
Consider the following set of generators for $F_{R(\Omega _v )}$:
variables from $\Lambda$ and quadratic bases from the active part.
Relations in this set of generators consist of the following three
families.

1. Relations between variables in $\Lambda$.

2. If $\mu$ is an active base and $\Delta (\mu)$ is a parametric
base, and $\Delta (\mu )=h_{i}\ldots h_{i+t},$ then there is a
relation $\mu=h_i\ldots h_{i+t}$.

3. Since $\gamma _i$=2 for each $h_i$ in the active part the
product of $h_i\ldots h_j, $ where $[i,j+1]$  is a closed active
section, can be written in two different ways $w_1$ and $w_2$ as a
product of active bases. We write the relations $w_1w_2^{-1}=1.$
These relations give a quadratic system of equations with
coefficients in the subgroup generated by $\Lambda $.

When we apply the entire transformation in Case 12, the number of
variables is not increasing and the complexity of the generalized
equation is not increasing. Suppose the same generalized equation
is repeated twice in the sequence (\ref{12'}). for example,
$\Omega _j=\Omega _{j+k}$. Then $\pi (v_j,v_{j+k})$ is an
automorphism of $F_{R(\Omega _j)}$ induced by the automorphism of
the free product $\langle \Lambda\rangle  *B,$ where $B$ is a free
group generated by quadratic bases, identical on $\langle \Lambda
\rangle  $ and fixing all words $w_1w_2^{-1}$. Therefore, $\bar H
^j>  \bar H ^{j+k},$ which contradicts to the minimality of $\bar
H$. Therefore there is only a finite number (bounded by $f_0$) of
possible generalized equations that can appear in the sequence
(\ref{12'}).$\Box$

Let $\bar H$ be a solution of the equation $\Omega$ with quadratic
part
 $[1,j+1]$.If $\mu $ belongs and $\Delta\mu$ does not  belong
to the quadratic section, then $\mu$ is called a {\em
quadratic-coefficient base}. Define the following numbers:
\begin{equation}\label{2.31}
d _1(\bar H)=\sum _{i=1}^{j}d(H_i),\end{equation}
\begin{equation}\label{2.32}
d _2(\bar H)=\sum _{\mu}d(H[\alpha (\mu),\beta
(\mu)]),\end{equation} where $\mu$ is a quadratic-coefficient
base.

\begin{lm}
\label{2.8} Let $tp (v)=15$ For any solution $\bar H$ of $\Omega
_v$ there is a minimal solution $\bar H^+$, which is an
automorphic image of $\bar H$ with respect to the group of
automorphisms defined in the beginning of this section, such that
$$d_1(\bar H^+)\leq f_1(\Omega _v)\ max \ \{d_2(\bar H^+),1\},$$
where $f_1(\Omega)$ is some recursive function.  \end{lm}

{\em Proof.}   Consider instead of $\Omega _v$ equation $\Omega
=(\tilde \Omega _v)$ which does not have any boundary connections,
$F_{R(\Omega _v)}$ is isomorphic to $F_{R(\Omega )}.$ Consider a
presentation of $F_{R(\Omega _v)}$ in the set of generators
consisting of variables in the non-quadratic part and active
bases. Relations in this generating set consist of  the following
three families.

1. Relations between variables in the non-quadratic part.

2. If $\mu$ is a quadratic-coefficient base and  $\Delta (\mu
)=h_{i}\cdots h_{i+t}$ in the non-quadratic part, then there is a
relation $\mu=h_i\cdots h_{i+t}.$

3. Since $\gamma _i$=2 for each $h_i$ in the active part the
product $h_i\cdots h_j, $ where $[i,j+1]$  is a closed active
section, can be written in two different ways $w_1$ and $w_2$ as a
product of quadratic and quadratic-coefficient bases. We write the
relations $w_1w_2^{-1}=1.$

Let $\bar H$ be a solution of $\Omega _v$ minimal with respect to
the canonical group of automorphisms of the
 free product $B_1*B,$ where $B$ is a free group generated by quadratic bases,
 and $B_1$ is a subgroup of $F_{R(\Omega _v)}$ generated by variables
 in the non-quadratic part, identical on $\langle \Lambda \rangle  $ and fixing all words $w_1w_2^{-1}$.

Consider the sequence
\begin{equation}\label{122}
 (\Omega ,\bar H)\rightarrow (\Omega _{v_1},\bar H^1)\rightarrow\ldots
 \rightarrow (\Omega _{v_i},\bar H^i),\ldots .\end{equation}

Apply now the entire transformations to the quadratic section of
$\Omega.$ As in the proof of the previous lemma, each time we
apply the entire transformation, we do not increase complexity
 and
do not increase the total number of items
 in the whole interval.
Since $\bar H$ is a solution of $\Omega _v$, if the same
generalized equation  appear in this sequence $2^{4^{j^2}}+1$
times then for some $j,j+k$ we have $\bar H^j > _c\bar H^{j+k}$,
therefore the same equation can only appear a bounded number of
times. Every quadratic base (except those that become matching
bases of length 1)
 in the quadratic part
can be transferred to the non-quadratic  part with the use of some
quadratic-coefficient base as a carrier base. This means that the
length of the transferred base is equal to the length of the part
of the quadratic-coefficient carrier base, which will then be
deleted. The double of the transferred base becomes a
quadratic-coefficient base. Because there are not more than $n_A$
 bases in the active part, this would give
$$d_1(\bar H')\leq n_A d_2(\bar H'),$$ for some solution $\bar H^+$
of the equation $ \widetilde\Omega _v$. But $\bar H^+$ is obtained
from the minimal solution $\bar H$ in a bounded number of steps.
$\Box$

We call a path $v_1\rightarrow v_2 \rightarrow\ldots\rightarrow
v_k\rightarrow\ldots $ in $T(\Omega )$ for which $7\leq
tp(v_k)\leq 10$ for all $k$ or type 12  {\em prohibited} if some
generalized equation with $\rho$ variables occurs among $\{\Omega
_{v_i} \mid 1\leq i\leq \ell\}$ at least $2^{(4\rho ^2)}+1$ times.
We will define below also prohibited paths in $T(\Omega )$, for
which $tp(v_k)=15$ for all $k$. We will need some auxiliary
definitions.

Introduce a new parameter
$$\tau _v'=\tau _v+\rho -\rho _{v}',$$
where $\rho $ is the number of variables
 of the initial equation $\Omega $ and $\rho _{v}'$ the number of free
variables belonging to the non-active sections of the equation
$\Omega _v.$ We have $\rho _{v}'\leq \rho$ (see the proof of Lemma
\ref{3.2}), hence $\tau _v'\geq 0$. In addition if $v_1\rightarrow
v_2$ is an auxiliary edge, then $\tau _2'< \tau _1'.$

Define by the joint induction on $\tau _v'$ a finite subtree
$T_0(\Omega _v)$ and a natural number $s(\Omega _v)$. The tree
$T_0(\Omega _v)$ will have $v$ as a root and consist of some
vertices and edges of $T(\Omega ) $ that lie higher than $v$. Let
$\tau _v'$=0; then in $T(\Omega )$ there can not be auxiliary
edges and vertices of type 15 higher than $v.$ Hence a subtree
$T_0(\Omega _v)$ consisting of vertices $v_1$ of $T(\Omega)$ that
are higher than $v$, and for which the path from $v$ to $v_1$ does
not contain prohibited subpaths, is finite.

Let now \begin{equation}\label{so} s(\Omega
_v)=\max_w\max_{\langle {\mathcal P},R\rangle  }\{\rho
_wf_2(\Omega _w,{\mathcal P},R),\ f_4 (\Omega _w',{\mathcal
P},R)\},\end{equation} where $w$ runs through all the vertices of
$T_0(v)$ for which $tp (w)=2$,
 $\Omega _w$ contains non-trivial non-parametric sections,
 $\langle {\mathcal P},R\rangle  $ is the set of non-singular periodic structures
of the equation $\tilde {\Omega} _w$,
  $f_2$ is a function appearing
in Lemma \ref{2.10'} ($f_2$ is present only if a periodic
structure has empty set $N{\mathcal P}$) and $\Omega _w'$ is
constructed as in Lemma \ref{4n}, where $f_4$ is a function
appearing in covering \ref{cat2}.

Suppose now that $\tau _v'>  0$ and that for all $v_1$ with $\tau
_{v_1}'< \tau _v'$ the tree $T_0(\Omega _{v_1})$ and the number
$s(\Omega _{v_1})$ are already defined. We begin with the
consideration of the paths
\begin{equation}\label{3.6}
r=v_1\rightarrow v_2\rightarrow \ldots\rightarrow
v_m,\end{equation} where $tp(v_i)=15\ (1\leq i\leq m)$. We have
$\tau _{v_i}'=\tau _v'.$

Denote by $\mu _i$ the carrier base of the equation $\Omega
_{v_i}$. The path (\ref{3.6}) will be called $\mu $-reducing if
$\mu _1=\mu$ and either there are no auxiliary edges from the
vertex $v_2$ and $\mu$ occurs in the sequence $\mu _1,\ldots ,\mu
_{m-1}$ at least twice, or there are auxiliary edges
$v_2\rightarrow w_1, v_2\rightarrow w_2\ldots ,v_2\rightarrow w_k$
from $v_2$ and $\mu$ occurs in the sequence $\mu _1,\ldots ,\mu
_{m-1}$ at least $max _{1\leq i\leq k}s(\Omega _{w_i})$ times.

The path (\ref{3.6}) will be called {\em prohibited}, if it can be
represented in the form
\begin{equation}\label{3.7} r=r_1s_1\ldots r_ls_lr',\end{equation}
such that for some sequence of bases $\eta _1,\ldots ,\eta _l$ the
following three properties hold:
\begin{enumerate}
\item [1)] every base occurring at least once in the sequence $\mu
_1,\ldots ,\mu _{m-1}$ occurs at least $40n^2f_1(\Omega
_{v_2})+20n+1$ times in the sequence $\eta _1,\ldots ,\eta _l$,
where $n$ is the number of pairs of bases in equations $\Omega
_{v_i}$, \item [2)] the path $r_i$ is $\eta _i$-reducing; \item
[3)] every transfer base of some equation of path $r$ is a
transfer base of some equation of path $r'$.\end{enumerate} The
property of path (\ref{3.6}) of being prohibited is
algorithmically decidable. Every infinite path (\ref{3.6})
contains a prohibited subpath. Indeed, let $\omega$ be the set of
all bases occurring in the sequence $\mu _1,\ldots ,\mu _m,\ldots
$ infinitely many times, and $\tilde\omega$ the set of all bases,
that are transfer bases of infinitely many equations $\Omega
_{v_i}$. If one cuts out some finite part in the beginning of this
infinite path, one can suppose that all the bases in the sequence
$\mu _1,\ldots ,\mu _m,\ldots $ belong to $\omega$ and each base
that is a transfer base of at least one equation, belongs to
$\tilde\omega$. Such an infinite path for any $\mu\in\omega$
contains infinitely many non-intersecting $\mu$-reducing finite
subpaths. Hence it is possible to construct a subpath (\ref{3.7})
of this path satisfying the first two conditions in the definition
of a prohibited subpath. Making $r'$ longer, one obtains a
prohibited subpath.

 Let $T'(\Omega
_v)$ be a subtree of $T(\Omega _v)$ consisting of the vertices
$v_1$ for which the path from $v$ to $v_1$ in $T(\Omega )$
contains neither prohibited subpaths nor vertices $v_2$ with $\tau
_{v_2}'< \tau _v',$ except perhaps $v_1$. So the terminal vertices
of $T'(\Omega _v)$ are either vertices $v_1$ such that $\tau
_{v_1}'< \tau _v',$ or terminal vertices of $T(\Omega _v)$. A
subtree $T'(\Omega _v)$ can be effectively constructed.
$T_0(\Omega _v)$ is obtained by attaching of $T_0(\Omega _{v_1})$
(already constructed by the induction hypothesis) to those
terminal vertices $v_1$ of $T'(\Omega _v)$ for which $\tau
_{v_1}'< \tau _v'.$ The function $s(\Omega _v)$ is defined by
(\ref{so}). Let now $T_0(\Omega)=T_0(\Omega _{v_0}).$ This tree is
finite by construction.

\subsubsection{Paths corresponding to minimal solutions of $\Omega$ are in $T_0(\Omega )$}\label{5.5.4}
Notice, that if $tp (v)\geq 6$ and $v\rightarrow w_1,\ldots
,v\rightarrow w_m$ is the list of principal outgoing edges from
$v$, then the generalized equations $\Omega _{w_1},\ldots ,\Omega
_{w_m}$ are obtained from $\Omega _v$ by the application of
several elementary transformations. Denote by $e$ a function that
assigns a pair $(\Omega _{w_i},\bar H ^{(i)})$ to the pair
$(\Omega _v,\bar H).$ For $tp (v)=4,5$ this function is identical.

If $tp (v)=15$ and there are auxiliary edges from the vertex $v$,
then the carrier base $\mu$ of the equation $\Omega _v$ intersects
$\Delta (\mu)$. For any solution $\bar H$ of the equation $\Omega
_v$ one can construct a solution $\bar H'$ of the equation $\Omega
_{v'}$ by $H'_{\rho _v+1}=H[1,\beta (\Delta (\mu))].$ Let
$e'(\Omega _v,\bar H)=e(\Omega _{v'},\bar H').$

In the beginning of this section we assigned to  vertices $v$ of
type 12, 15, 2 and such that $7\leq tp(v)\leq 10$ of $ T(\Omega)$
the groups of   automorphisms $A(\Omega _v)$.%
 Denote by $Aut (\Omega )$ the group of  automorphisms
of $F_{R(\Omega )}$ , generated by all groups \newline $\pi
(v_0,v)A(\Omega _v)\pi (v_0,v)^{-1}$, $v\in T_0(\Omega)$. (Here
$\pi (v_0,v)$ is an isomorphism, because $tp(v)\not = 1$.) We are
to formulate the main technical result of this section. The
following proposition states that every minimal solution of a
generalized equation $\Omega$ with respect to the group $A(\Omega
)$ either factors through one of the finite family of proper
quotients of the group $F_{R(\Omega )}$ or (in the case of a
non-empty parametric part) can be transferred to the parametric
part.
\begin{prop} \label{3.4} For any solution $\bar H$ of a generalized equation
$\Omega $ there exists a  terminal vertex $w$ of the tree
$T_0(\Omega )$ having type 1 or 2,  and a solution $\bar H^{(w)}$
of a generalized equation $\Omega _w$ such that
\begin {enumerate}
\item [1)] $\pi _{\bar H}=\sigma\pi(v_0,w)\pi _{\bar H^{(w)}}\pi\
$ where $\pi$ is an endomorphism of a free group $F$, $\sigma\in
Aut (\Omega )$; \item [2)] if $tp (w)=2$ and the equation $\Omega
_{w}$ contains  nontrivial non-parametric sections, then there
exists a primitive cyclically reduced word $P$ such that $\bar
H^{(w)}$ is periodic with respect to ${\mathcal P}$ and one of the
following conditions holds:
\begin{enumerate}\item the equation $\Omega _w$ is singular with
respect  to a periodic structure
\newline ${\mathcal P}(\bar H^{(w)},P) $ and the first minimal
replacement can be applied, \item it is possible to apply the
second minimal replacement and make the family of closed sections
in $\mathcal P$ empty.
\end{enumerate}
\end{enumerate}
\end{prop}

Construct a directed tree with paths from the initial vertex
\begin{equation}\label{3.8}
(\Omega ,\bar H)=(\Omega _{v_0},\bar H^{(0)})\rightarrow (\Omega
_{v_1}, \bar H^{(1)})\rightarrow\ldots \rightarrow(\Omega
_{v_u},\bar H^{(u)})\rightarrow\ldots\end{equation} in which the
$v_i$ are the vertices of the tree $T(\Omega )$ in the following
way. Let $v_1=v_0$ and let $\bar H^{(1)}$ be some solution of
 the equation $\Omega $, minimal with respect to the group of automorphisms $A(\Omega v_0)$
with the property $\bar H\geq \bar H^{(1)}.$

Let $i\geq 1$ and suppose the term $(\Omega _{v_i},\bar H^{(i)})$
of the sequence (\ref{3.8}) has been already constructed.  If
$7\leq tp(v_i)\leq 10$ or $tp(v_i)=12$ and there exists a minimal
solution $\bar H^+$ of $\Omega _{v_i}$ such that $\bar H^+< \bar
H^{(i)}$ , then we set $v_{i+1}=v_i$, $\bar H^{(i+1)}=\bar H^+.$

If $tp(v_i)=15,\  v_i\neq v_{i-1}$ and there are auxiliary edges
from vertex $v_i$: $v_i\rightarrow w_1,\ldots ,v_i\rightarrow w_k$
(the carrier base $\mu$ intersects with its double $\Delta
(\mu)$),  then there exists a primitive word $P$ such that
\begin{equation}\label{3.9}
H^{(i)}[1,\beta (\Delta (\mu))]\equiv P^rP_1, r\geq 2,\ P\equiv
P_1P_2,\end{equation} where $\equiv$ denotes a graphical equality.
In this case the path (\ref{3.8}) can be continued along several
possible edges of $T(\Omega )$.

For each group of automorphisms assigned to vertices of type 2 in
the trees $T_0 (\Omega _{w_i})$, $i=1,\ldots ,k$ and non-singular
periodic structure including the closed section $[1,\beta (\Delta
(\mu)]$ of the equation $\Omega _{v_i}$ and corresponding to
solution $\bar H^{(i)}$ we replace $\bar H^{(i)}$ by a solution
$\bar H^{(i)+}$ with maximal number of short variables (see the
definition after Lemma \ref{4n}). This collection of short
variables can be different for different periodic structures.
Either all the variables in $\bar H^{(i)+}$ are short or there
exists a parametric base $\lambda _{max}$ of maximal length in the
covering (\ref{cat2}). Suppose there is a $\mu$-reducing path
(\ref{3.6}) beginning at $v_i$ and corresponding to $\bar
H^{(i)+}$. Let $\mu _1,\ldots ,\mu _m$ be the leading bases of
this path. Let ${\tilde H}^1=H^{(i)+},\ldots ,{\tilde H}^j$ be
solutions of the generalized equations corresponding to the
vertices of this path. If for some $\mu _i$ there is an inequality
$d({\tilde H}^j[\alpha (\mu _i),\alpha (\Delta (\mu _i))])\leq d
(\lambda _{max}),$ we set $(\Omega _{v_{i+1}},\bar
H^{(i+1)})=e'(\Omega _{v_i},\bar H^{(i)})$ and call the section
$[1,\beta (\Delta (\mu ))]$ which becomes non-active, {\em
potentially transferable}.

If there is a singular periodic structure in a  vertex of type 2
of some tree $T_0(\Omega _{w_i}), i\in\{1,\ldots ,k\},$ including
the closed section $[1,\beta (\Delta (\mu)]$ of the equation
$\Omega _{v_i}$ and corresponding to the  solution $\bar H^{(i)}$,
we also include the possibility \newline $(\Omega _{v_{i+1}},\bar
H^{(i+1)})=e'(\Omega _{v_i},\bar H^{(i)}).$

In all of the other cases we set $(\Omega _{v_{i+1}},\bar
H^{(i+1)})=e(\Omega _{v_i},\bar H^{(i)+}),$ where
 $\bar H^{(i)+}$ is a solution with maximal number of short variables and minimal solution of
$\Omega _{v_{i}}$ with respect to the canonical group of
automorphisms $P_{v_i}$ (if it exists). The path (\ref{3.8}) ends
if $tp (v_i)\leq 2.$

We will show that in the path (\ref{3.8}) $v_i\in T_0(\Omega)$. We
use induction on $\tau '$.
Suppose $v_i\not\in T_0(\Omega ),$ and let $i_0$ be the first of
such numbers. It follows from the construction of $T_0(\Omega )$
that there exists $i_1< i_0$ such that the path from $v_{i_1}$
into $v_{i_0}$ contains a subpath prohibited in the construction
of $T_2(\Omega _{v_{i_1}}).$ From the minimality of $i_0$ it
follows that this subpath goes from $v_{i_2}\ \ (i_1\leq i_2<
i_0)$ to $v_{i_0}$. It cannot be that $7\leq tp(v_i)\leq 10$ or
$tp(v_i)=12$ for all $i_2\leq i\leq i_1$, because there will be
two indices $p< q$ between $i_2$ and $i_0$ such that $\bar
H^{(p)}=\bar H^{(q)}$, and this gives a contradiction, because in
this case it must be by construction $v_{p+1}=v_p$. So
$tp(v_i)=15$ ($i_2\leq i\leq i_0$).

Suppose we have a subpath (\ref{3.6}) corresponding to the
fragment
\begin{equation}\label{3.11}
(\Omega _{v_1},\bar H^{(1)})\rightarrow (\Omega _{v_2}, \bar
H^{(2)})\rightarrow\ldots \rightarrow(\Omega _{v_m},\bar
H^{(m)})\rightarrow\ldots\end{equation} of the sequence
(\ref{3.8}). Here $v_1,v_2,\ldots,v_{m-1}$ are vertices of the
tree $T_0(\Omega)$, and for all vertices $v_i$ the edge
$v_i\rightarrow v_{i+1}$ is principal.

As before, let $\mu _i$ denote the carrier base of $\Omega
_{v_i}$, and $\omega =\{\mu _1,\ldots ,\mu _{m-1}\},$ and
$\tilde\omega $ denote the set of such bases which are transfer
bases for at least one equation in (\ref{3.11}). By $\omega _1$
denote the set of such bases $\mu $ for which either $\mu$ or
$\Delta (\mu)$ belongs to $\omega\cup\tilde\omega $; by $\omega
_2$ denote the set of all the other bases. Let
$$\alpha (\omega)=\min(\min _{\mu\in\omega _2}\alpha (\mu),j),$$
where $j$ is the boundary between active and non-active sections.
Let $X_{\mu}\circeq H[\alpha (\mu),\beta (\mu)].$ If $(\Omega
,\bar H)$ is a member of sequence (\ref{3.11}), then denote
\begin{equation}\label{3.12}
d _{\omega}(\bar H)=\sum _{i=1}^{\alpha (\omega )-1}d(H_i),
\end{equation}

\begin{equation}\label{3.13}
\psi _{\omega}(\bar H)=\sum _{\mu\in\omega
_1}d(X_{\mu})-2d_{\omega}(\bar H).\end{equation}

Every item $h_i$ of the section $[1,\alpha (\omega)]$ belongs to
at least two bases, and both bases are in $\omega _1$, hence $\psi
_{\omega}(\bar H)\geq 0.$

Consider the quadratic part of $\widetilde\Omega _{v_1}$ which is
situated to the left of $\alpha(\omega )$. The  solution $\bar
H^{(1)}$ is minimal with respect to the canonical group of
automorphisms corresponding to this vertex. By Lemma \ref{2.8} we
have
\begin{equation}
\label{3.14} d_1(\bar H^{(1)})\leq f_1(\Omega _{v_1})d_2(\bar
H^{(1)}).
\end{equation}

 Using this inequality we estimate the length of the interval
participating in the process $d_{\omega}(\bar H^{(1)})$ from above
by a product of $\psi _{\omega }$ and some function depending on
$f_1$. This will be inequality \ref{3.19}. Then we will show that
for a prohibited subpath the length of the participating interval
must be reduced by more than this figure (equalities \ref{3.27},
\ref{3.28}). This will imply that there is no prohibited subpath
in the path \ref{3.11}.

Denote by $\gamma _i(\omega)$ the number of bases $\mu\in\omega
_1$ containing $h_i$. Then
\begin{equation}\label{3.15}
\sum _{\mu\in\omega _1}d(X_{\mu}^{(1)})=\sum
_{i=1}^{\rho}d(H_i^{(1)}) \gamma _i(\omega),\end{equation} where
$\rho =\rho (\Omega _{v_1}).$ Let $I=\{i|1\leq i\leq\alpha
(\omega)-1 \&\gamma _i=2\}$ and $J=\{i|1\leq i\leq\alpha
(\omega)-1 \&\gamma _i>  2\}.$ By (\ref{3.12})
\begin{equation}\label{3.16}
d_{\omega}(\bar H^{(1)})=\sum _{i\in I}d(H_i^{(1)})+ \sum _{i\in
J}d(H_i^{(1)})= d_1(\bar H^{(1)})+\sum _{i\in
J}d(H_i^{(1)}).\end{equation} Let  $(\lambda ,\Delta(\lambda))$ be
a pair of quadratic-coefficient bases of the equation $\tilde
\Omega _{v_1}$, where $\lambda$ belongs to the nonquadratic part.
This pair can appear only from the bases $\mu\in\omega _1$. There
are two types of quadratic-coefficient bases.

{\em Type} 1. $\lambda$ is situated to the left of the boundary
$\alpha (\omega)$. Then $\lambda$ is formed by items $\{h_i|i\in
J\}$ and hence $d(X_{\lambda})\leq\sum _{i\in J}d(H_i^{(1)}).$
Thus the sum of the lengths $d(X_{\lambda})+d(X_{\Delta
(\lambda)})$ for quadratic-coefficient bases of this type is not
more than $2n\sum _{i\in J}d(H_i^{(1)}).$

{\em Type} 2. $\lambda$ is
 situated to the right of the boundary $\alpha
(\omega)$. The sum of length of the quadratic-coefficient bases of
the second type is not more than $2\sum _{i=\alpha
(\omega)}^{\rho}d(H_i^{(1)})\gamma _i(\omega).$

 We have \begin{equation}\label{3.17}
 d_2(\bar H^{(1)})\leq 2n\sum _{i\in J}d(H_i^{(1)})+2\sum _{i=\alpha (\omega)}^{\rho}d(H_i^{(1)})\gamma
_i(\omega).\end{equation} Now (\ref{3.13}) and (\ref{3.15}) imply
\begin{equation}\label{3.18}
\psi _{\omega}(\bar H^{(1)}_i)\geq \sum _{i\in J}d(H_i^{(1)})+\sum
_{i=\alpha (\omega)}^{\rho}d(H_i^{(1)})\gamma
_i(\omega).\end{equation} Inequalities (\ref{3.14}),
(\ref{3.16}),(\ref{3.17}),(\ref{3.18}) imply
\begin{equation}\label{3.19}
d_{\omega}(\bar H^{(1)})\leq max \{\psi _{\omega}(\bar
H^{(1)})(2nf_1(\Omega _{v_1})+1),\ f_1(\Omega _{v_1})\}.
\end{equation}

 From the definition of Case 15 it follows that all the words
 $H^{(i)}[1,\rho _i+1]$ are the ends of the word $H^{(1)}[1,\rho _1+1]$,
 that is
 \begin{equation}\label{3.20}
 H^{(1)}[1,\rho _1+1]\doteq U_iH^{(i)}[1,\rho _i+1].\end{equation}
 On the other hand bases $\mu\in\omega _2$  participate in these
 transformations neither as carrier bases nor as transfer bases; hence
 $H^{(1)}[\alpha (\omega ),\rho _1+1]$ is the end of the word $H^{(i)}[
 1,\rho _i+1]$, that is
 \begin{equation} \label{3.21}
 H^{(i)}[1,\rho _i+1]\doteq V_iH^{(1)}[\alpha
(\omega ),\rho _1 +1].\end{equation} So we have
\begin{equation}\label{3.22} d_{\omega}(\bar H
^{(i)})-d_{\omega}(\bar H ^{(i+1)})=
d(V_i)-d(V_{i+1})=d(U_{i+1})-d(U_{i})= d(X_{\mu
_i}^{(i)})-d(X_{\mu _{i}}^{(i+1)}).\end{equation} In particular
(\ref{3.13}),(\ref{3.22}) imply that $\psi _{\omega }(\bar
H^{(1)})=\psi _{\omega }(\bar H^{(2)})=\ldots \psi _{\omega }(\bar
H^{(m)})=\psi _{\omega }.$ Denote  the number (\ref{3.22})
 by $\delta _i$.

Let the path (\ref{3.6}) be $\mu$-reducing, that is either $\mu
_1=\mu$ and $v_2$ does not have auxiliary edges and $\mu$ occurs
in the sequence $\mu _1,\ldots ,\mu _{m-1}$ at least twice, or
$v_2$ does have auxiliary edges $v_2\rightarrow w_1,\ldots
v_2\rightarrow w_k $ and the base $\mu$ occurs in the sequence
$\mu _1,\ldots ,\mu _{m-1}$ at least $max _{1\leq i\leq k}s(\Omega
_{w_i})$ times. Estimate $d(U_m)=\sum _{i=1}^{m-1}\delta _i$ from
below. First notice that if $\mu _{i_1}=\mu _{i_2}=\mu (i_1< i_2)
$ and $\mu _i\not =\mu$ for $i_1< i< i_2$, then
\begin{equation}
\label{3.23} \sum _{i=i_1}^{i_2-1}\delta _i\geq
d(H^{i_1+1}[1,\alpha (\Delta (\mu_{i_1+1}))]).
\end{equation}
Indeed, if $i_2=i_1+1,$ then $\delta _{i_1}=d(H^{(i_1)}[1,\alpha
(\Delta (\mu))]=d(H^{(i_1+1)}[1,\alpha (\Delta (\mu))].$ If $i_2 >
i_1+1,$ then $\mu _{i_1+1}\not = \mu$ and $\mu$ is a transfer base
in the equation $\Omega _{v_{i_1+1}}.$ Hence $\delta
_{i_1+1}+d(H^{(i_1+2)}[1,\alpha (\mu)])=d(H^{(i_1+1)}[1,\alpha
(\mu _{i_1+1})]).$ Now (\ref{3.23}) follows from
$$\sum _{i=i_1+2}^{i_2-1}\delta _i\geq d(H^{(i_1+2)}[1,\alpha (\mu)]).$$
So if $v_2$ does not have outgoing auxiliary edges, that is the
bases $\mu _2$ and $\Delta(\mu _2)$ do not intersect in the
equation $\Omega _{v_2}$; then (\ref{3.23}) implies that

$$\sum _{i=1}^{m-1}\delta _i\geq d(H^{(2)}[1,\alpha (\Delta\mu _2)])\geq d(
X_{\mu _2}^{(2)})\geq d(X_{\mu}^{(2)})=d(X_{\mu}^{(1)})-\delta
_1,$$ which implies that
\begin{equation}\label{3.24}
\sum _{i=1}^{m-1}\delta
_i\geq\frac{1}{2}d(X_{\mu}^{(1)}).\end{equation}

Suppose now there are outgoing auxiliary edges from the vertex
$v_2$: $v_2\rightarrow w_1,\ldots ,v_2\rightarrow w_k$. The
equation $\Omega _{v_1}$ has some solution. Let $H^{(2)}[1,\alpha
(\Delta (\mu _2))]\doteq Q$, and $P$ a  word (in the final $h$'s)
such that $Q\doteq P^{d}$, then $X_{\mu _2}^{(2)}$ and $X_{\mu
}^{(2)}$ are beginnings of the word $H^{(2)}[1,\beta (\Delta (\mu
_2))],$ which is a beginning of $P^{\infty}$. Denote $M=max
_{1\leq j\leq k} s(\Omega _{w_j}).$

By the construction of (\ref{3.8})we either have
\begin{equation}\label{3.25}
X_{\mu}^{(2)}\doteq P^rP_1, P\doteq P_1P_2,r< M.
\end{equation}
or for each base  $\mu _i,\ i\geq 2,$ there is an inequality
$d(H^{(i)}(\alpha (\mu _i),\alpha (\Delta (\mu _i))))\geq
d(\lambda )$ and therefore
\begin{equation}\label{korova}
d(X_{\mu}^{(2)})<  M d( H^{(i)}[\alpha (\mu _i),\alpha (\Delta
(\mu _i))]).\end{equation}

Let $\mu _{i_1}=\mu _{i_2}=\mu; i_1< i_2; \mu _i\not = \mu$ for
$i_1< i< i_2.$ If
\begin{equation} \label{3.26a}
d(X_{\mu _{i_1+1}}^{(i_1+1)})\geq 2 d(P)\end{equation} and
$H^{(i_1+1)}[1,\rho _{i_1+1}+1]$ begins with a cyclic permutation
of $P^3$, then \newline $d(H^{(i_1+1)}[1,\alpha (\Delta (\mu
_{i_1+1}))])> d(X^{(2)}_{\mu})/M$. Together with (\ref{3.23}) this
gives $\sum _{i=i_1}^{i_2-1}\delta _i> d(X_{\mu }^{(2)})/M.$ The
base $\mu$ occurs in the sequence $\mu _1,\ldots ,\mu _{m-1}$ at
least $M$ times, so either (\ref{3.26a}) fails for some $i_1\leq
m-1$ or $\sum _{i=1}^{m-1}\delta _i\ (M-3)d(X^{(2)}_{\mu})/M.$

If (\ref{3.26a}) fails, then the inequality $d(X_{\mu
_i}^{(i+1)})\leq d(X_{\mu _{i+1}}^{(i+1)}) ,$ and the definition
(\ref{3.22})  imply that
$$\sum _{i=1}^{i_1}\delta _i\geq d(X_{\mu}^{(1)})-d(X_{\mu _{i_1+1}}^
{(i_1+1)})\geq (M-2)d(X_{\mu}^{(2)})/M;$$ so everything is reduced
to the second case.

Let $$\sum _{i=1}^{m-1}\delta _i\geq
 (M-3)d(X_{\mu}^{(1)})/M.$$ Notice that (\ref{3.23}) implies for $i_1=1$
$\sum _{i=1}^{m-1}\delta _i\geq d(Q)\geq d(P)$; so $\sum
_{i=1}^{m-1}\delta _i\geq max\{1,M-3\}d(X^{(2)}_{\mu})/M.$
Together with (\ref{3.25}) this implies $\sum _{i=1}^{m-1}\delta
_i\geq \frac{1}{5}d(X_{\mu}^{(2)})=\frac{1}{5}
(d(X_{\mu}^{(1)})-\delta _1).$ Finally,
\begin{equation}\label{3.26}
\sum _{i=1}^{m-1}\delta _i\geq\frac{1}{10} d(X_{\mu}^{(1)}).
\end{equation}
Comparing (\ref{3.24}) and (\ref{3.26}) we can see that for the
$\mu$-reducing path (\ref{3.6}) inequality (\ref{3.26}) always
holds.

Suppose now that the path (\ref{3.6}) is prohibited; hence it can
be represented in the form (\ref{3.7}). From  definition
(\ref{3.13}) we have $\sum _{\mu \in\omega _1}d(X_{\mu}^{(m)})\geq
\psi _{\omega}$; so at least for one base $\mu\in\omega _1$ the
inequality $d(X_{\mu}^{(m)})\geq\frac{1}{2n}\psi _{\omega}$ holds.
Because $X_{\mu}^{(m)}\doteq (X_{\Delta (\mu)}^{(m)})^{\pm 1},$ we
can suppose that $\mu\in\omega\cup\tilde{\omega}.$ Let $m_1$ be
the length of the path $r_1s_1\ldots r_ls_l$ in (\ref{3.7}). If
$\mu\in\tilde{\omega}$ then  by the third part of the definition
of a prohibited path there exists $m_1\leq i\leq m$ such that
$\mu$ is a transfer base of $\Omega _{v_i}$. Hence, $d(X_{\mu
_i}^{(m_1)})\geq d(X_{\mu _i}^{(i)})\geq d(X_{\mu}^{(i)})\geq
d(X_{\mu}^{(m)})\geq\frac{1}{2n}\psi _{\omega}.$ If
$\mu\in\omega$, then take $\mu$ instead of $\mu _i$. We proved the
existence of a base $\mu\in\omega$ such that
\begin{equation}\label{3.27}
d(X_{\mu }^{(m_1)})\geq\frac{1}{2n}\psi _{\omega}.\end{equation}
 By the definition of a prohibited path, the inequality
$d(X_{\mu}^{(i)})\geq d(X_{\mu}^{(m_1)}) (1\leq i\leq m_1),$
(\ref{3.26}), and (\ref{3.27}) we obtain
\begin{equation}\label{3.28}
\sum _{i=1}^{m_1-1}\delta _i\geq max\{\frac{1}{20n}\psi
_{\omega},1\}(40n^2f_1+20n+1).
\end{equation}

By (\ref{3.22}) the sum in the left part of the inequality
(\ref{3.28}) equals $d_{\omega}(\bar H^{(1)})-d_{\omega}(\bar
H^{(m_1)});$ hence
$$d_{\omega}(\bar H^{(1)})\geq max \{\frac{1}{20n}\psi _{\omega},1\}
(40n^2f_1 +20n +1) ,$$ which contradicts (\ref{3.19}).

This contradiction was obtained from the supposition that there
are prohibited paths (\ref{3.11}) in the path (\ref{3.8}). Hence
(\ref{3.8}) does not contain prohibited paths. This implies that
$v_i\in T_0(\Omega )$ for all $v_i$ in (\ref{3.8}).
 For all $i$ $v_i\rightarrow v_{i+1}$ is an edge of a finite tree. Hence
the path (\ref{3.8}) is finite. Let $(\Omega _{w},\bar H^{w})$ be
the final term of this sequence. We show that $(\Omega _{w},\bar
H^{w})$ satisfies all the  properties formulated in the lemma.

The first property is obvious.

Let $tp (w)=2$ and let $\Omega _w$ have non-trivial non-parametric
part. It follows from the construction of (\ref{3.8}) that if
$[j,k]$ is a non-active section for $\Omega _{v_i}$ then
$H^{(i)}[j,k]\doteq H^{(i+1)}[j,k]\doteq\ldots H^{(w)}[j,k]$.
Hence (\ref{3.9}) and the definition of $s(\Omega _v)$ imply that
the word $h_1\ldots h_{\rho _w}$ can be subdivided into subwords
$h[i_1,i_2],\ldots ,h[i_{k-1},i_k],$ such that for any $a$ either
$H^{(w)}$ has length 1, or $h[i_a,i_{a+1}]$ does not participate
in basic and coefficient equations, or $H^{(w)}[i_a,i_{a+1}]$ can
be written as
\begin{equation}
\label{3.29} H^{(w)}[i_a,i_{a+1}]\doteq P_a^rP_a';\ P_a\doteq
P_a'P_a''; r\geq max_{\langle \mathcal P,R\rangle  }\ max \{\rho
_{w}f_2(\Omega _w,P,R), f_4(\Omega _w')\},\end{equation} where
$P_a$ is a primitive word, and $\langle \mathcal P,R\rangle  $
runs
 through all the periodic
structures of ${\tilde \Omega} _w$  such that either one of them
is singular or for a  solution with maximal number of short
variables with respect to the group of extended automorphisms all
the closed sections are potentially transferable. The proof of
Proposition \ref{3.4} will be completed after we prove the
following statement.
\begin{lm}\label{33} If $tp (w)=2$ and every closed section belonging to a
periodic structure $\mathcal P$ is potentially transferable (the
definition is given in the construction of $T_0$ in case 15), one
can apply the second minimal replacement and get a finite number
(depending on periodic structures containing this section in the
vertices of type 2 in the trees $T_0(w_i), i=1,\ldots ,m$ ) of
possible generalized equations containing the same closed sections
not from $\mathcal P$ and not containing closed sections from
$\mathcal P$.\end{lm}

{\em Proof.}   From the definition of  potentially transferable
section it follows that after finite number of transformations
depending on $f_4(\Omega '_u,\mathcal P)$, where $u$ runs through
the vertices of type 2 in the trees $T_0(w_i), i=1,\ldots ,m$, we
obtain a cycle that is shorter than or equal to $d(\lambda
_{max}).$ This cycle is exactly $h[\alpha (\mu _i),\alpha (\Delta
(\mu _i)]$ for the base $\mu _i$ in the $\mu$-reducing subpath.
The rest of the proof of Lemma \ref{33} is a repetition of the
proof of Lemma \ref{3n}. $\Box$

\subsubsection{The decomposition tree $T_{dec}(\Omega )$}\label{5.5.5}
We can define now a decomposition tree $T_{dec}(\Omega ).$ To
obtain $T_{dec}(\Omega )$ we add some edges to the terminal
vertices of type 2 of $T_0(\Omega )$. Let $v$ be a vertex of type
2 in $T_0(\Omega )$. If there is no periodic structures in $\Omega
_v$ then this is a terminal vertex of $T_{dec}(\Omega )$. Suppose
there exists a finite number of combinations of different periodic
structures ${\mathcal P}_1,\ldots ,{\mathcal P}_s$ in $\Omega _v$.
If some ${\mathcal P}_i$ is singular, we consider a generalized
equation $\Omega _{u({\mathcal P}_1,\ldots ,{\mathcal P}_s)}$
obtained from $\Omega _v({\mathcal P}_1,\ldots ,{\mathcal P}_s)$
by the first minimal replacement corresponding to ${\mathcal
P}_i$. We also draw the edge $v\rightarrow u=u({\mathcal
P}_1,\ldots ,{\mathcal P}_s)$. This vertex $u$ is a terminal
vertex of $T_{dec}(\Omega )$. If all ${\mathcal P}_1,\ldots
,{\mathcal P}_s$ in $\Omega _v$ are not singular, we can suppose
that for each periodic structure ${\mathcal P}_i$ with period
$P_i$ some values of variables in ${\mathcal P}_i$ are shorter
than $2|P_i|$ and values of some other variables are shorter than
$f_3(\Omega _v)|P_i|,$ where $f_3$ is the function from Lemma
\ref{4n}. Then we apply the second minimal replacement. The
resulting generalized equations $\Omega _{u_1},\ldots ,\Omega
_{u_t}$ will have empty non-parametric part. We draw the edges
$v\rightarrow u_1,\ldots , v\rightarrow u_t$ in $T_{dec}(\Omega
)$. Vertices $u_1,\ldots ,u_t$ are terminal vertices of
$T_{dec}(\Omega )$.

\subsection{The solution tree $T_{sol}(\Omega,\Lambda)$}\label{se:5.5}

Let $\Omega = \Omega(H)$ be a generalized  equation in variables
$H$ with the set of bases $B_\Omega = B \cup \Lambda$.  Let
$T_{dec}(\Omega)$ be  the tree constructed in Subsection
\ref{5.5.5} for a generalized equation $\Omega$ with parameters
$\Lambda$.

  Recall that in a leaf-vertex $v$  of $T_{dec}(\Omega)$  we have the
coordinate group $F_{R(\Omega_v)}$ which is a proper homomorphic
image of $F_{R(\Omega)}$.  We define  a new  transformation $R_v$
(we call it {\it leaf-extension}) of the tree $T_{dec}(\Omega)$ at
the leaf vertex $v$. We take the union of two trees
$T_{dec}(\Omega)$ and $T_{dec}(\Omega_v)$ and identify the
vertices $v$ in both trees (i.e., we extend the tree
$T_{dec}(\Omega)$ by gluing the tree $T_{dec}(\Omega_v)$ to the
vertex $v$).  Observe that if the equation $\Omega_v$ has
non-parametric  non-constant sections  (in this event we call $v$
{\em a terminal vertex}), then $T_{dec}(\Omega_v)$ consists of a
single vertex, namely $v$.

  Now we construct a solution tree $T_{sol}(\Omega)$ by induction
starting at $T_{dec}(\Omega)$.  Let $v$ be  a leaf  non-terminal
vertex of $T^{(0)} = T_{dec}(\Omega)$. Then we apply the
transformation $R_v$ and obtain a new tree $T^{(1)} =
R_v(T_{dec}(\Omega))$. If there exists a leaf  non-terminal vertex
$v_1$ of $T^{(1)}$, then we apply the transformation $R_{v_1}$,
and so on. By induction we construct a strictly increasing
sequence of trees

\begin{equation}
\label{eq:5.3.1} T^{(0)} \subset  T^{(1)} \subset \ldots \subset
T^{(i)}  \subset \ldots .\end{equation} This sequence is finite.
Indeed, suppose to the contrary that the sequence is infinite and
hence  the union $T^{(\infty)}$ of this sequence is an infinite
tree in which every vertex has a finite degree. By  Konig's lemma
there is  an infinite branch $B$  in $T^{(\infty)}$. Observe that
along any infinite  branch in $T^{(\infty)}$  one has to encounter
infinitely many proper epimorphisms. This contradicts the fact
that $F$ is equationally Noetherian.

Denote the union of the sequence of the trees (\ref{eq:5.3.1}) by
$T_{sol}(\Omega,\Lambda)$. We call $T_{sol}(\Omega,\Lambda)$ the
{\em solution tree of $\Omega$ with parameters $\Lambda$}.
 Recall that with every edge $e$ in $T_{dec}(\Omega)$ (as well as in
$T_{sol}(\Omega,\Lambda)$)
  with the initial vertex $v$ and
 the terminal  vertex $w$ we associate an epimorphism
 $$\pi_e: F_{R(\Omega_v)} \rightarrow  F_{R(\Omega_v)}.$$
 It follows that every connected (directed) path $p$ in the graph gives
rise to a
 composition of homomorphisms which we denote by $\pi_p$. Since
$T_{sol}(\Omega,\Lambda)$
 is a tree the path $p$ is completely defined by its initial and
terminal vertices $u, v$;
 in this case we  sometimes write $\pi_{u,v}$ instead of $\pi_p$.
  Let $\pi_v$ be the homomorphism
 corresponding to the path from the initial vertex $v_0$ to a given
vertex $v$,
  we call it the {\it canonical epimorphism} from $F_{R(\Omega)}$ onto
 $F_{R(\Omega_v)}$.

 Also,
 with some vertices $v$ in the tree
$T_{dec}(\Omega)$, as well as in the tree
$T_{sol}(\Omega,\Lambda)$, we associate groups of canonical
automorphisms $A(\Omega _v)$ or extended automorphisms  $\bar
A(\Omega _v)$ of the coordinate group $F_{R(\Omega_v)}$ which, in
particular, fix  all variables  in the non-active part of
$\Omega_v$.  We can suppose that the group  $\bar A(\Omega _v)$ is
associated to every vertex, but for some vertices it is trivial.
Observe also, that canonical epimorphisms map parametric parts
into parametric parts (i.e., subgroups generated by variables in
parametric parts).

Recall that  writing $(\Omega,U)$ means that $U$ is a solution of
$\Omega$. If $(\Omega,U)$ and $\mu \in B_{\Omega}$, then by
$\mu_U$ we denote the element
\begin{equation}
\label{eq:5.3.2} \mu_U = [u_{\alpha(\mu)} \ldots u_{\beta(\mu)-1}
]^{\varepsilon(\mu)}. \end{equation}
 Let $B_U = \{\mu_U \mid \mu \in B\}$ and
$\Lambda_U = \{\mu_U \mid \mu \in \Lambda\}$. We refer to  these
sets as the set of values of bases from $B$ and the set of values
of parameters from $\Lambda$ with respect to the solution $U$.
Notice, that the value  $\mu_U$ is given in (\ref{eq:5.3.2}) as a
value of one fixed word mapping  $$P_\mu(H) = [h_{\alpha(\mu)}
\ldots h_{\beta(\mu)-1} ]^{\varepsilon(\mu)}. $$ In vector
notation we can write that
  $$B_U =  P_B(U), \ \ \  \Lambda_U = P_{\Lambda}(U),$$
  where $P_B(H)$ and $P_{\Lambda}(H)$ are corresponding word mappings.

  The following result explains the name of the
tree $T_{sol}(\Omega,\Lambda)$.

\begin{theorem}
\label{th:5.3.1} Let $\Omega = \Omega(H,\Lambda)$ be a generalized
equation in variables $H$ with parameters $\Lambda .$  Let
$T_{sol}(\Omega,\Lambda)$ be the solution tree for $\Omega$ with
parameters.
  Then the following conditions hold.

1.  For any solution $U$ of the generalized equation $\Omega$
there exists a path \newline $v_0, v_1, \ldots, v_n = v$  in
$T_{sol}(\Omega,\Lambda)$ from the root vertex $v_0$ to a terminal
vertex $v$, a sequence of canonical automorphisms $\sigma =
(\sigma_0, \ldots, \sigma_n ), \sigma_i \in  A(\Omega _{v_i}),$
and a solution $U_v$ of the generalized equation $\Omega _v$ such
that the solution $U$ (viewed as a homomorphism $F_{R(\Omega
)}\rightarrow F$) is equal to the following composition of
homomorphisms
\begin{equation}
 \label{eq:5.3.3}
 U=\Phi_{\sigma,U_v} =
\sigma_0 \pi_{v_0,v_1}\sigma_{1} \ldots  \pi_{v_{n-1},v_n}\sigma
_n U_v.\end{equation}

2. For any path $v_0, v_1, \ldots, v_n = v$  in
$T_{sol}(\Omega,\Lambda)$ from the root vertex $v_0$ to a terminal
vertex $v$, a sequence of canonical automorphisms $\sigma =
(\sigma_0, \ldots, \sigma_n ), \sigma_i \in  A(\Omega _{v_i}),$
and a solution $U_v$ of the generalized equation $\Omega_v$, $\Phi
_ {\sigma,U_v}$ gives a solution of the group equation $\Omega^* =
1$; moreover, every solution of $\Omega^* =1$ can be obtained this
way.

 3.  For each  terminal vertex $v$ in $T_{sol}(\Omega,\Lambda)$
there exists a word mapping $Q_v(H_v)$ such that for any solution
$U_v$ of $\Omega_v$ and any solution $U = \Phi_{\sigma,U_v}$ from
(\ref{eq:5.3.3}) the values of the parameters $\Lambda$ with
respect to $U$ can be written as $\Lambda_U = Q_v(U_v)$ (i.e.,
these values do not depend on $\sigma$)  and the word $Q_v(U_v)$
is reduced as written.
\end{theorem}
{\em Proof.}   Statements 1 and 2 follow from the construction of
the tree $T_{sol}(\Omega,\Lambda)$. To verify  3 we need to invoke
the argument above this theorem which  claims that the canonical
automorphisms associated with generalized equations in
$T_{sol}(\Omega,\Lambda)$ fix  all variables  in the parametric
part and, also, that  the  canonical epimorphisms map variables
from the parametric part into themselves.

The set of homomorphisms having form (\ref{eq:5.3.3}) is called a
{\em fundamental sequence}.

\begin{theorem} \label{Ase} For any finite system  $S(X)=1$ over a free group $F$, one  can find effectively a
finite family of nondegenerate triangular quasi-quadratic systems
$U_1,\ldots, U_k$ and word mappings $p_i: V_F(U_i) \rightarrow
V_F(S)$ $(i = 1, \ldots,k)$ such that for every $b \in V_F(S)$
there exists $i$ and $c \in V_F(U_i)$ for which $b = p_i(c)$, i.e.
$$
V_F(S) = p_1(V_F(U_1)) \cup \ldots \cup p_k(V_F(U_k))
$$
and all sets $p_i(V_F(U_i))$ are irreducible; moreover, every
irreducible component of $V_F(S)$ can be obtained as a closure of
some $p_i(V_F(U_i))$ in the Zariski topology.
\end{theorem}

{\em Proof.}    Each solution of the system $S(X)=1$ can be
obtained as $X=p_i(Y_i),$ where $Y_i$ are variables of $\Omega
=\Omega _i$ for a finite number of generalized equations. We have
to show that all solutions of $\Omega ^*$ are solutions of some
NTQ system. We can use Theorem \ref{th:5.3.1} without parameters.
In this case $\Omega _v$ is an empty equation with non-empty set
of variables. In other words $F_{R(\Omega _v)}=F*F(h_1,\ldots
,h_{\rho}).$ To each of the branches of $T_{sol}$ we assign an NTQ
system from the formulation of the theorem. Let $\Omega _w$ be a
leaf vertex in $T_{dec}.$ Then $F_{R(\Omega _w)}$ is a proper
quotient of $F_{R(\Omega )}$. Consider the path $v_0, v_1, \ldots,
v_n = w$ in $T_{dec}(\Omega)$ from the root vertex $v_0$ to a
terminal vertex $w$. All the groups $F_{R(\Omega _{v_i})}$ are
isomorphic. There are the following four possibilities.

1. $tp (v_{n-1})= 2$. In this case there is a singular periodic
structure on $\Omega _{v_{n-1}}$.  By Lemma \ref{2.10''},
$F_{R(\Omega _{v_{n-1}})}$ is a fundamental group of a graph of
groups with  one vertex group $K$,  some free abelian vertex
groups, and some edges defining HNN extensions of $K$. Recall that
making the first minimal replacement we first replaced
$F_{R(\Omega _{v_{n-1}})}$  by a finite number of proper quotients
in which the edge groups corresponding to abelian vertex groups
and HNN extensions are maximal cyclic in $K$. Extend the
centralizers of the edge groups of $\Omega _{v_{n-1}}$
corresponding to HNN extensions by stable letters $t_1,\ldots
,t_k$. This new group that we denote by $N$ is the coordinate
group of a quadratic equation  over $F_{R(\Omega _w)}$ which has a
solution in $F_{R(\Omega _w)}$.

In all the other cases  $tp (v_{n-1})\neq 2$.

2. There were no auxiliary edges from vertices $v_0, v_1, \ldots,
v_n = w$, and if one of the Cases 7--10 appeared at one of these
vertices, then it only appeared a bounded (the boundary depends on
$\Omega _{v_0}$) number of times in the sequence. In this case we
replace $F_{R(\Omega )}$ by $F_{R(\Omega _w)}$

3. $F_{R(\Omega _w)}$ is a term in a free decomposition of
$F_{R(\Omega _{v_{n-1}})}$ ( $\Omega _w$ is a kernel of a
generalized equation $\Omega _{v_{n-1}}$). In this case we also
consider $F_{R(\Omega _w)}$ instead of $F_{R(\Omega)}$.

4. For some $i$ $tp(v_i)=12$ and the path $v_i,\ldots ,v_n=w$ does
not contain vertices of type $ 7-10, 12$ or $15$. In this case
$F_{R(\Omega)}$ is the coordinate group of a quadratic equation.

5. The path $v_0, v_1, \ldots, v_n = w$ contains  vertices of type
15.  Suppose $v_{i_j},\ldots ,v_{i_j+k_j},$ $j=1,\ldots ,l$ are
all blocks of  consecutive vertices of type 15 This means that
$tp(v_{i_j+k_j+1})\neq 15$ and $i_j+k_j+1<i_{j+1}$. Suppose also
that none of the previous cases takes place. To each $v_{i_j}$ we
assigned a quadratic equation and a group of canonical
automorphisms corresponding to this equation. Going alon the path
$v_{i_j},\ldots ,v_{i_j+k_j},$ we take minimal solutions
corresponding to some non-singular periodic structures. Each such
structure corresponds to a representation of $F_{R(\Omega
_{v_{i_j}})}$ as an HNN extension. As in the case of a singular
periodic structure, we can suppose that the edge groups
corresponding to  HNN extensions are maximal cyclic and not
conjugated in $K$. Extend the centralizers of the edge groups
corresponding to HNN extensions by stable letters $t_1,\ldots
,t_k$.  Let $N$ be the new group. Then $N$ is the coordinate group
of a quadratic system of equations over $F_{R(\Omega
_{v_{i_j+k_j+1}})}$. Repeating this construction for each
$j=1,\ldots ,l$, we construct NTQ system over $F_{R(\Omega _w)}$.

Since $F_{R(\Omega _w)}$ is a proper quotient of $F_{R(\Omega )}$,
the theorem can  now be proved by induction.

$\Box$

\begin{theorem}\label{ge}For any finitely generated group G and a free group $F$ the set
$Hom(G,F)$ [$Hom _F(G,F)$]  can be effectively described by
   a finite rooted tree oriented from the
  root, all vertices except for the root vertex are labelled by
  coordinate groups of generalized equations. Edges from the root vertex correspond
  to a finite number of homomorphisms from $G$ into coordinate groups
of generalized equations. Leaf vertices are labelled by free
groups. To each vertex
   group we assign the group of canonical
  automorphisms.
Each edge (except for the edges from the root) in this tree  is
labelled by a quotient map, and all quotients are proper. Every
homomorphism from G to F can be written as a composition of the
homomorphisms corresponding to edges, canonical automorphisms of
the groups assigned to vertices, and  some homomorphism [retract]
from a free group in a leaf vertex into $F$.\end{theorem}

\subsection{Cut Equations}\label{se:cut}
In the proof of the implicit function theorems it will be
convenient to use a modification of the notion of a generalized
equation.  The following definition provides a framework for such
a modification.

\begin{df}\label{df:cut}
A cut equation $\Pi = ({\mathcal E}, M, X, f_M, f_X)$ consists of
a set of intervals $\mathcal E$, a set of variables $M$, a set of
parameters $X$, and two  labeling functions $$f_X: {\mathcal E}
\rightarrow F[X], \ \ \ f_M:  {\mathcal E} \rightarrow F[M] .$$
For an interval $\sigma \in {\mathcal E}$ the image  $f_M(\sigma)
= f_M(\sigma)(M)$ is a reduced word in variables $M^{\pm 1}$ and
constants from $F$, we call it a {\it partition}  of
$f_X(\sigma)$.
\end{df}

Sometimes we write $\Pi = ({\mathcal E}, f_M, f_X)$ omitting $M$
and $X$.

\begin{df}
A solution of a cut equation $\Pi = ({\mathcal E}, f_M, f_X)$ with
respect to an $F$-homomorphism $\beta : F[X] \rightarrow F$  is an
$F$-homomorphism $\alpha : F[M] \rightarrow F$ such that: 1) for
every $\mu \in M$ $\alpha(\mu)$ is a reduced non-empty word; 2)
for every reduced word $f_M(\sigma)(M)\ (\sigma \in {\mathcal E})$
the replacement $m \rightarrow \alpha(m) \ (m \in M) $ results in
a word $f_M(\sigma)(\alpha(M))$ which is  a  reduced word as
written and  such that $f_M(\sigma)(\alpha(M))$ is graphically
equal to the reduced form of $\beta(f_X(\sigma))$; in particular,
the following diagram is commutative.
\begin{center}
\begin{picture}(100,100)(0,0)
\put(50,100){$\mathcal E$} \put(0,50){$F(X)$} \put(100,50){$F(M)$}
\put(50,0){$F$} \put(47,97){\vector(-1,-1){30}}
\put(53,97){\vector(1,-1){30}} \put(3,47){\vector(1,-1){30}}
\put(97,47){\vector(-1,-1){30}} \put(16,73){$f_X$}
\put(77,73){$f_M$} \put(16,23){$\beta$} \put(77,23){$\alpha$}
\end{picture}
\end{center}
\end{df}

If $\alpha: F[M] \rightarrow F$ is a solution of a cut equation
$\Pi = ({\mathcal E}, f_M, f_X)$ with respect to an
$F$-homomorphism $\beta : F[X] \rightarrow F$, then we write
$(\Pi, \beta,\alpha)$ and refer to $\alpha$ as a \emph{solution
of} $\Pi$ {\it modulo} $\beta$. In this event,   for a given
$\sigma \in {\mathcal E}$ we say that $f_M(\sigma)(\alpha(M))$ is
a {\it partition} of $\beta(f_X(\sigma))$. Sometimes we also
consider homomorphisms $\alpha:F[M] \rightarrow F$, for which the
diagram above is still commutative, but cancellation may occur in
the words $f_M(\sigma)(\alpha(M))$. In this event we refer to
$\alpha$ as a {\em group} solution of $\Pi$ with respect to
$\beta$.

\begin{lm}
\label{le:cut}
 For a generalized equation $\Omega(H)$ one can effectively construct a cut
equation
 $\Pi_{\Omega} = ({\mathcal E}, f_X, f_M)$ such that  the following
conditions hold:
\begin{enumerate}
 \item [(1)] $ X$ is a partition of the whole interval $[1,\rho_{\Omega}]$ into
disjoint
  closed subintervals;

 \item  [(2)] $M$ contains the set of variables $H$;

 \item [(3)] for any solution $U = (u_1, \ldots, u_\rho)$ of $\Omega$  the cut
equation
 $\Pi_{\Omega}$ has a solution
 $\alpha$ modulo the canonical homomorphism $\beta_U: F(X) \rightarrow
F$
 ($\beta_U(x) = u_i u_{i+1} \ldots u_j$ where $i,j$ are,
correspondingly,
  the left and the right  end-points of the interval $x$);

 \item [(4)] for any solution $(\beta,\alpha)$ of the cut equation $\Pi_{\Omega}$
the
 restriction of $\alpha$ on $H$ gives a solution of the generalized
equation
 $\Omega$.
 \end{enumerate}
 \end{lm}

{\it Proof.}
 We begin with defining the sets $X$ and $M$. Recall that  a
closed interval of $\Omega$ is a union of closed sections of
$\Omega$.  Let $X$ be an arbitrary partition of the whole interval
$[1,\rho_{\Omega}]$ into closed subintervals (i.e., any two
intervals in $X$ are disjoint and the union of $X$ is the whole
interval $[1,\rho_{\Omega}]$).

 Let $B$ be a set of representatives of dual bases of $\Omega$, i.e.,
for every base
  $\mu$  of $\Omega$ either $\mu$ or $\Delta(\mu)$ belongs to $B$, but
not both.
 Put $M = H \cup B$.

Now let $\sigma \in X$.  We denote  by $B_{\sigma}$  the set of
all bases over $\sigma$ and by  $H_\sigma$ the set of all items in
$\sigma$. Put $S_\sigma = B_{\sigma} \cup H_\sigma.$  For $e \in
S_\sigma$ let $s(e)$ be the interval $[i,j]$, where $i < j$ are
the endpoints of $e$.  A sequence  $P = (e_1, \ldots,e_k)$ of
elements from $S_\sigma$  is called a {\it partition} of $\sigma$
if $s(e_1) \cup \cdots \cup s(e_k) = \sigma$ and $s(e_i) \cap
s(e_j) = \emptyset$ for $i \neq j$.  Let ${\rm Part}_\sigma$ be
the set of all partitions of $\sigma$. Now put
 $${\mathcal E} = \{P \mid P \in {\rm Part}_\sigma, \sigma \in X\}.$$
Then for every $P \in {\mathcal E}$
 there exists one and only one $\sigma \in X$ such that $P \in
{\rm Part}_\sigma$.
 Denote this $\sigma$ by $f_X(P)$. The map $f_X: P \rightarrow f_X(P)$
is a
 well-defined function from ${\mathcal E}$ into $F(X)$.

 Each partition $P = (e_1, \ldots,e_k) \in {\rm Part}_\sigma$ gives rise to a
word
  $w_P(M) = w_1 \ldots w_k$  as follows.  If $e_i \in H_\sigma$ then
$w_i = e_i$.
  If $e_i = \mu \in B_\sigma$ then $w_i = \mu^{\varepsilon(\mu)}$.
  If $e_i = \mu$ and $\Delta(\mu) \in B_\sigma$ then
  $w_i = \Delta(\mu)^{\varepsilon(\mu)}$. The map $f_M(P) = w_P(M)$ is a
  well-defined function from ${\mathcal E}$ into $F(M)$.

  Now set $\Pi_{\Omega} = ({\mathcal E}, f_X, f_M)$.  It is not hard to see
from the
  construction that the cut equation $\Pi_\Omega$ satisfies all the
required properties.
  Indeed, (1) and (2) follow directly from the construction.

  To verify (3), let's consider a
  solution $U = (u_1, \ldots, u_{\rho_\Omega})$ of $\Omega$. To define
  corresponding functions $\beta_U$ and $\alpha$, observe that the
function $s(e)$
  (see above) is defined for every $e \in X \cup M$.
   Now  for $\sigma \in X$ put
 $\beta_U(\sigma) = u_i\ldots u_j$, where $s(\sigma) = [i,j]$, and for $m
\in M$
 put $\alpha(m) =  u_i\ldots u_j$, where $s(m) = [i,j]$.  Clearly,
$\alpha$ is a
 solution of $\Pi_\Omega$ modulo $\beta$.

 To verify (4) observe that if $\alpha$ is a solution of $\Pi_\Omega$ modulo
$\beta$,
  then the restriction of $\alpha$ onto the subset $H \subset M$ gives a
solution
 of the generalized equation $\Omega$. This follows from the
construction
 of the words $w_p$ and the fact that the words $w_p(\alpha(M))$ are
reduced as
 written (see definition of a solution of a cut equation). Indeed, if a
base
 $\mu$ occurs in a partition $P \in {\mathcal E}$, then there is a partition
 $P^\prime \in {\mathcal E}$  which is obtained from $P$ by replacing $\mu$
by the
 sequence $h_i\ldots h_j$. Since there is no cancellation in words
$w_P(\alpha(M))$ and
 $w_{P^\prime}(\alpha(M))$, this implies that
 $\alpha(\mu)^{\varepsilon(\mu)} = \alpha(h_i\ldots h_j)$. This shows
that
 $\alpha_H$ is a solution of $\Omega$.
 \hfill $\Box$

\begin{theorem}\label{th:cut}
Let $S(X,Y,A)) = 1$ be a system of  equations over $F = F(A)$.
Then one can effectively construct a finite set of cut equations
 $${\mathcal CE}(S) = \{\Pi_i \mid  \Pi_i =({\mathcal E}_i, f_{X_i},  f_{M_i}),
i = 1 \ldots,k \}$$
  and a finite set of tuples of words $\{Q_i(M_i) \mid i = 1,
\dots,k\}$  such that:
\begin{enumerate}
\item for every equation $\Pi_i =({\mathcal E}_i, f_{X_i},
f_{M_i}) \in {\mathcal CE}(S)$,
 one has $X_i = X$ and   $f_{X_i}({\mathcal E}_i) \subset  X^{\pm 1}$;

\item  for  any  solution $(U,V)$ of $S(X,Y,A) = 1$ in $F(A)$,
there exists a number $i$
 and a tuple of words $P_{i,V}$ such that the  cut
equation $\Pi_i \in {\mathcal CE}(S)$  has a solution $\alpha: M_i
\rightarrow F$ with respect to the $F$-homomorphism  $\beta_U:
F[X] \rightarrow F$ which is induced by the map $ X \rightarrow
U$. Moreover, $U = Q_i(\alpha(M_i))$, the word $Q_i(\alpha(M_i))$
is reduced as written,  and  $V = P_{i,V}(\alpha(M_i))$;

\item  for any $\Pi_i \in {\mathcal CE}(S)$ there exists a tuple
of words $P_{i,V}$ such that for  any solution (group solution)
$(\beta, \alpha)$ of  $\Pi_i$ the pair $(U,V),$ where $U =
Q_i(\alpha(M_i))$ and $V = P_{i,V}(\alpha(M_i)),$ is a solution of
$S(X,Y) = 1$ in $F$.
\end{enumerate}
\end{theorem}

{\it Proof.}
 Let $S(X,Y) = 1$ be a system of  equations over a free group $F$.
In  Subsection \ref{se:parametric}    we have constructed a set of
initial parameterized generalized equations  ${\mathcal
GE}_{par}(S) = \{\Omega_1, \ldots, \Omega_r\}$ for $S(X,Y) = 1$
with respect to the set of parameters $X$.
 For each $\Omega \in {\mathcal GE}_{par}(S)$ in Section \ref{se:5.5} we
constructed  the finite tree $T_{sol}(\Omega)$ with respect to
parameters $X$. Observe that parametric part
$[j_{v_0},\rho_{v_0}]$ in the root equation $\Omega =
\Omega_{v_0}$ of the tree $T_{sol}(\Omega)$ is  partitioned into a
disjoint union of closed sections  corresponding to $X$-bases and
constant bases (this follows from the construction of the initial
equations in the set  ${\mathcal GE}_{par}(S)$). We label every
closed section $\sigma$ corresponding to a variable  $x \in X^{\pm
1}$ by $x$, and every constant section corresponding to a constant
$a$ by $a$.  Due to our construction of the tree $T_{sol}(\Omega)$
moving along a branch $B$ from the initial vertex $v_0$ to a
terminal vertex $v$, we transfer all the bases from the active and
non-active parts into parametric parts until, eventually,  in
$\Omega_v$ the whole interval consists of the parametric part.
 Observe also that, moving along $B$ in the parametric part,  we neither
introduce
  new closed sections nor  delete any. All we do is we split
(sometimes) an item in a closed  parametric section into two new
ones. In any event we keep the same label of the section.

Now for a terminal vertex $v$ in $T_{sol}(\Omega)$  we construct a
cut equation $\Pi^\prime_v = ({\mathcal E}_v, f_{X_v}, f_{M_v})$
as in Lemma \ref{le:cut} taking the set of all  closed sections of
$\Omega_v$ as the partition $X_v$.  The set of cut equations
 $$ {\mathcal CE}^\prime(S) = \{\Pi^\prime_v \mid \Omega \in {\mathcal
GE}_{par}(S), v
 \in VTerm(T_{sol}(\Omega))\}$$
satisfies all the requirements of the theorem except $X_v$ might
not be equal to $X$. To satisfy this condition we adjust slightly
the equations $\Pi_v^\prime$.

To do this, we denote by  $l:X_v \rightarrow X^{\pm 1} \cup A^{\pm
1}$ the labelling function on the set of closed sections of
$\Omega_v$. Put $\Pi_v = ({\mathcal E}_v, f_{X}, f_{M_v})$ where
$f_X$ is the composition of $f_{X_v}$ and $l$. The set of  cut
equations
 $$ {\mathcal CE}(S) = \{\Pi_v \mid  \Omega \in {\mathcal GE}_{par}(S), v
 \in VTerm(T_{sol}(\Omega))\}$$
satisfies all the conditions of the theorem. This follows from
Theorem \ref{th:5.3.1} and from Lemma \ref{le:cut}. Indeed, to
satisfy 3) one can take the words $P_{i,V}$ that correspond to a
minimal solution of $\Pi_i$, i.e., the words $P_{i,V}$ can be
obtained from a given particular way to transfer all bases from
$Y$-part onto $X$-part.

\hfill $\Box$

The next result shows that for every cut equation $\Pi$ one can
effectively and canonically  associate a generalized equation
$\Omega_{\Pi}$.

   For  every  cut equation $\Pi= ({\mathcal E}, X, M, f_X, f_M)$  one
  can  canonically associate a generalized equation $\Omega_{\Pi}(M,X)$
as follows.
  Consider the following word

$$ V = f_X(\sigma_1)f_M(\sigma_1)  \ldots f_X(\sigma_k)f_M(\sigma_k). $$
 Now we are going to mimic the construction of the generalized equation
in Lemma \ref{le:4.2}.
  The set of boundaries $BD$ of $\Omega_{\Pi}$ consists of
positive integers
 $1, \ldots, |V| +1$. The set of bases $BS$ is union of the following
sets:

a)  every letter $\mu$ in  the word $V$. Letters  $X^{\pm 1} \cup
M^{\pm 1}$ are variable bases, for every two different occurrences
$\mu^{\varepsilon_1}, \mu^{\varepsilon_2}$ of  a letter  $\mu \in
X^{\pm 1} \cup M^{\pm 1}$ in $V$ we say that these bases  are dual
and they have the same orientation if $\varepsilon_1\varepsilon_2
= 1$, and different orientation otherwise. Each occurrence of a
letter  $a \in A^{\pm 1}$ provides a constant base with the label
$a$.  Endpoints of these bases correspond to their positions in
the word $V$ (see  Lemma
 \ref{le:R1}).

 b)   every pair of subwords $f_X(\sigma_i), f_M(\sigma_i)$ provides a
pair of dual
 bases $\lambda_i, \Delta(\lambda_i)$, the base $\lambda_i$ is located
above the
 subword $f_X(\sigma_i)$, and  $\Delta(\lambda_i)$ is located above
 $f_M(\sigma_i)$ (this defines the endpoints of the bases).

 Informally, one can visualize the generalized equation $\Omega_{\Pi}$
as
 follows.   Let ${\mathcal E} = \{\sigma_1, \ldots, \sigma_k\}$ and let
  ${\mathcal E}^\prime = \{ \sigma^\prime \mid \sigma \in {\mathcal E}\}$
   be another disjoint copy of the set ${\mathcal E}$. Locate intervals from
${\mathcal
   E} \cup {\mathcal E}^\prime$ on a segment $I$ of a straight line from
left to the
   right in the following order $\sigma_1, \sigma_1^\prime,
   \ldots, \sigma_k, \sigma_k^\prime$; then  put bases over $I$
according to the word
   $V$.
 The next result summarizes the discussion above.
 \begin{lm}
  \label{le:5.4.3}
 For  every  cut equation $\Pi= ({\mathcal E}, X, M, f_x, f_M)$,  one
  can  canonically associate a generalized equation $\Omega_{\Pi}(M,X)$
such that
   if $\alpha_\beta : F[M] \rightarrow F$ is a solution of the cut
 equation $\Pi$, then the maps $\alpha: F[M] \rightarrow F$ and $\beta:
F[X]
 \rightarrow F$ give rise to a solution of the group equation (not
generalized!)
 $\Omega_{\Pi}^* = 1$ in such a way that  for every $\sigma \in {\mathcal E}$
 $f_M(\sigma)(\alpha(M))$ is a reduced word which is graphically equal
to
 $\beta(f_X(\sigma)(X))$, and vice versa.
  \end{lm}

\section{Definitions and elementary properties of liftings}

In this section we give necessary definitions for the further
discussion of liftings of equations and inequaliies into
coordinate groups.

Let $G$ be a group and let $S(X) = 1$ be a system of equations
over $G$. Recall that by $G_{S}$ we denote the quotient group
$G[X]/ncl(S)$, where $ncl(S)$ is the normal closure of $S$ in
$G[X]$. In particular, $G_{R(S)} = G[X]/R(S)$ is the coordinate
group defined by $S(X) = 1$. The radical $R(S)$ can be described
as follows. Consider a set of $G$-homomorphisms
 $$\Phi_{G,S} = \{\phi \in Hom_G(G[S],G)  \mid
\phi(S) = 1\}.$$ Then  \[R(S) = \left\{\begin{array}{ll}\cap_{\phi
\in \Phi_{G,S}} \ker \phi & \mbox{if $\Phi_{G,S} \neq
\emptyset$}\\ G[X] & \mbox{otherwise} \end{array} \right. \]

  Now we put these definitions in a more general framework.
  Let  $H$ and $K$ be $G$-groups   and $M \subset H$. Put
  $$\Phi_{K,M} = \{\phi \in Hom_G(H,K)  \mid
\phi(M) = 1\}.$$ Then the following subgroup
  is termed the {\it  $G$-radical of $M$
  with respect to $K$:}
 \[Rad_K(M) = \left\{\begin{array}{ll}\cap_{\phi \in
\Phi_{K,M}} \ker \phi , & \mbox{if $\Phi_{K,M} \neq \emptyset$,}\\
G[X] & \mbox{otherwise.} \end{array} \right. \] Sometimes, to
emphasize that $M$ is a subset of $H$, we write $Rad_K(M,H)$.
Clearly, if $K = G$, then $R(S) = Rad_G(S,G[X])$.

Let $$H_K^* = H/Rad_K(1).$$ Then $H_K^*$ is either a $G$-group or
trivial.  If $H_K^* \neq 1$, then it is  $G$-separated by $K$. In
the case $K = G$ we omit $K$ in the notation above and simply
write $H^*$. Notice that $$(H/ncl(M))^*_K \simeq H/Rad_K(M),$$ in
particular, $(G_S)^* = G_{R(S)}$.

\begin{lm}
Let $\alpha: H_1 \rightarrow H_2$ be a $G$-homomorphism and
suppose $\Phi = \{\phi:H_2 \rightarrow K \}$ be a separating
family of $G$- homomorphisms. Then
$$ \ker \alpha = \bigcap \{\ker(\alpha \circ \phi) \ \mid \ \phi \in \Phi \} $$
\end{lm}
{\it Proof}. Suppose $h \in H_1$ and $h \not \in \ker(\alpha).$
Then $\alpha(h) \neq 1$ in $H_2$. Hence there exists $\phi \in
\Phi$ such that $\phi(\alpha(h)) \neq 1$. This shows that $ \ker
\alpha \supset \bigcap \{\ker(\alpha \circ \phi) \ \mid \ \phi \in
\Phi \} $. The other inclusion is obvious.  $\Box$

\begin{lm}
\label{le:6.2} Let $H_1$, $H_2$, and $K$  be $G$-groups.
\begin{enumerate}
\item  Let  $\alpha: H_1 \rightarrow H_2$ be a $G$-homomorphism
and let $H_2$  be
 $G$-separated by $K$. If $M \subset \ker \alpha$,  then
 $Rad_K(M) \subseteq \ker \alpha$.
\item Every $G$-homomorphism $\phi: H_1 \rightarrow H_2$ gives
rise to a unique homomorphism $$\phi^*:(H_1)_K^* \rightarrow
(H_2)_K^*$$ such that $\eta_2 \circ \phi = \phi^* \circ \eta_1 $,
where $\eta_i: H_i \rightarrow H_i^*$ is the canonical
epimorphism.
\end{enumerate}
\end{lm}
{\it Proof}. 1. We have $$Rad_K(M,H_1) = \bigcap \{\ker \phi \mid
\phi: H_1 \rightarrow_G K  \ \wedge \ \phi(M) = 1\}  \subseteq
\bigcap \{\ker(\alpha \circ \beta) \mid \beta : H_2 \rightarrow_G
K\} = \ker \alpha.$$ 2. Let $\alpha:H_1 \rightarrow (H_2)_K^*$ be
the composition of the following homomorphisms $$ H_1
\stackrel{\phi}{\rightarrow} H_2 \stackrel{\eta_2}{\rightarrow}
(H_2)_K^*. $$ Then by assertion 1 $Rad_K(1,H_1) \subseteq \ker
\alpha$, therefore $\alpha$ induces the canonical $G$-homomorphism
$\phi^*:(H_1)_K^* \rightarrow (H_2)_K^*$. $\Box$

\begin{lm}
\label{le:6.3}
\begin{enumerate}
\item The canonical map $\lambda:G \rightarrow G_S$ is an
embedding $\Longleftrightarrow$ $S(X) = 1$ has a solution in some
$G$-group $H$. \item The canonical map $\mu:G \rightarrow
G_{R(S)}$ is an embedding $\Longleftrightarrow$ $S(X) = 1$ has a
solution in some $G$-group $H$ which is $G$-separated by $G$.
\end{enumerate}
\end{lm}
{\it Proof}.  1. If $S(x_1, \ldots,x_m) = 1$ has a solution
$(h_1,\ldots,h_m)$ in some $G$-group $H$, then the
$G$-homomorphism $x_i \rightarrow h_i, \ (i =1,\ldots,m)$ from
$G[x_1,\ldots,x_m]$ into $H$ induces a homomorphism  $\phi: G_S
\rightarrow H$. Since $H$ is a $G$-group all non-trivial elements
from $G$ are also non-trivial in the factor-group $G_S$, therefore
$\lambda:G \rightarrow G_S$ is an embedding. The converse is
obvious.

2. Let $S(x_1, \ldots,x_m) = 1$ have a solution $(h_1,\ldots,h_m)$
in some $G$-group $H$ which is $G$-separated by $G$. Then there
exists the canonical $G$-homomorphism $\alpha: G_S \rightarrow H$
defined as in the proof of the first assertion. Hence $R(S)
\subseteq \ker \alpha$ by Lemma \ref{le:6.2},  and $\alpha$
induces a homomorphism from $G_{R(S)}$ into $H$, which is monic on
$G$. Therefore $G$ embeds into $G_{R(S)}$. The converse is
obvious. $\Box$

Now we apply Lemma \ref{le:6.2} to coordinate groups of nonempty
algebraic sets.
\begin{lm}
\label{le:6.4} Let subsets $S$ and $T$ from $ G[X]$ define
non-empty algebraic sets in a group $G$. Then every
$G$-homomorphism $\phi:G_S \rightarrow G_T$ gives rise to a
$G$-homomorphism $\phi^*: G_{R(S)} \rightarrow G_{R(T)}$.
\end{lm}
{\em Proof.}   The result follows from Lemma \ref{le:6.2}  and
Lemma \ref{le:6.3}.

Now we are in a position to give the following
\begin{df}
Let $S(X) = 1$ be a system of  equations over a group $G$ which
has a solution in $G$.  We say that a system of equations  $T(X,Y)
= 1$ is compatible with $S(X) = 1$ over $G$ if for every solution
$U$ of $S(X) = 1$ in $G$ the equation $T(U,Y) = 1$ also has a
solution in $G$, i.e., the algebraic set $V_G(S)$ is a projection
of the algebraic set $V_G(S\cup T).$
\end{df}
The next proposition describes compatibility of two equations in
terms of their coordinate groups.
\begin{prop}
\label{prop:6.1} Let $S(X) = 1$ be a system of  equations over a
group $G$ which has a solution in $G$. Then $T(X,Y) = 1$ is
compatible with $S(X) = 1$ over $G$ if and only if $G_{R(S)}$ is
canonically embedded into $G_{R(S\cup T)},$ and every
$G$-homomorphism $\alpha: G_{R(S)} \rightarrow G$ extends to a
$G$-homomorphisms $\alpha^\prime: G_{R(S\cup T)} \rightarrow G$.
\end{prop}
{\it Proof}. Suppose first that $T(X,Y) = 1$ is compatible with
$S(X) = 1$ over $G$ and suppose that $V_G(S) \neq \emptyset.$ The
identity map $X \rightarrow X$ gives rise to a $G$-homomorphism $$
\lambda: G_S \longrightarrow G_{S\cup T} $$ (notice that both
$G_S$ and $G_{S\cup T}$ are $G$-groups by Lemma \ref{le:6.3}),
which by Lemma \ref{le:6.4} induces a $G$-homomorphism $$
\lambda^*: G_{R(S)} \longrightarrow G_{R(S\cup T)}. $$ We claim
that $\lambda^*$ is an embedding. To show this  we need to prove
first the  statement about the extensions of homomorphisms. Let
$\alpha:G_{R(S)} \rightarrow G$ be  an arbitrary $G$-homomorphism.
It follows that $ \alpha(X)$ is a solution of $S(X) = 1$ in $G$.
Since $T(X,Y) = 1$ is compatible with $S(X) = 1$ over $G$, there
exists a solution, say $\beta(Y)$,  of  $T(\alpha(X),Y) = 1$ in
$G$. The  map $$ X \rightarrow \alpha(X), Y \rightarrow \beta(Y)
$$ gives rise to  a $G$-homomorphism $ G[X,Y] \rightarrow G$,
which induces a $G$-homomorphism $\phi: G_{S\cup T} \rightarrow
G$. By Lemma \ref{le:6.4} $\phi$ induces a $G$-homomorphism $$
\phi^* : G_{R(S\cup T)} \longrightarrow G.$$ Clearly, $\phi^*$
makes the following diagram to commute.
\bigskip
\begin{center}
\begin{picture}(100,100)(0,0)
\put(0,100){$G_{R(S)}$} \put(100,100){$G_{R(S\cup T)}$}
\put(0,0){$G$} \put(15,103){\vector(1,0){80}}
\put(5,93){\vector(0,-1){78}} \put(95,95){\vector(-1,-1){80}}
\put(-10,50){$\alpha$} \put(55,45){$\phi^*$}
\put(50,108){$\lambda^*$}
\end{picture}
\end{center}
Now to prove  that $\lambda^*$ is an embedding, observe that
$G_{R(S)}$ is $G$-separated by $G$. Therefore for every
non-trivial $h \in G_{R(S)}$ there exists a $G$-homomorphism
$\alpha: G_{R(S)} \rightarrow G$ such that $\alpha(h) \neq 1$. But
then $\phi^*(\lambda^*(h)) \neq 1$ and consequently $h \not \in
\ker \lambda^*$. The converse statement is obvious. $\Box$

\smallskip
Let $S(X) = 1$ be a system of  equations over $G$ and suppose
$V_G(S) \neq \emptyset$. The canonical embedding $X \rightarrow
G[X]$ induces the canonical map $$\mu: X \rightarrow G_{R(S)}.$$
We are ready to formulate the main definition.
\begin{df} Let $S(X) = 1$ be a system of  equations over $G$ with
$V_G(S) \neq \emptyset$ and let $\mu :X \rightarrow G_{R(S)}$ be
the canonical map.  Let a system  $T(X,Y) = 1$ be compatible with
$S(X) = 1 $ over $G$. We say that $T(X,Y) = 1$ admits a lift to a
generic point of $S = 1$ over $G$ (or, shortly, $S$-lift over $G$)
if $T(X^\mu,Y) = 1$ has a solution in $G_{R(S)}$ (here $Y$ are
variables and $X^\mu$ are constants from $G_{R(S)}$).
\end{df}
\begin{lm}
\label{le:oger} Let $T(X,Y) = 1$ be compatible with $S(X) = 1 $
over $G$. If $T(X,Y) = 1$ admits an $S$-lift, then the identity
map $Y \rightarrow Y$ gives rise to  a canonical
$G_{R(S)}$-epimorphism from $G_{R(S\cup T)}$ onto the coordinate
group of $T(X^\mu,Y) = 1$ over $G_{R(S)}$:
 $$  \psi^* : G_{R(S\cup T)} \rightarrow
G_{R(S)}[Y]/Rad_{G_{R(S)}}(T(X^\mu,Y)).$$ Moreover, every solution
$U$
 of  $T(X^\mu,Y) = 1$ in $G_{R(S)}$ gives rise to a $G_{R(S)}$-homomorphism $\phi_U:  G_{R(S\cup T)} \rightarrow G_{R(S)}$, where
$\phi_U(Y) = U$.
\end{lm}
{\em Proof.}   Observe that the following chain of isomorphisms
hold: $$ G_{R(S\cup T)} \simeq_G G[X][Y]/Rad_G(S\cup T) \simeq_G
G[X][Y]/Rad_G(Rad_G(S,G[X])\cup T)$$ $$ \simeq_G
(G[X][Y]/ncl(Rad_G(S,G[X]) \cup T))^*  \simeq_G
(G_{R(S)}[Y]/ncl(T(X^\mu,Y)))^*. $$
 Denote by $\overline{G_{R(S)}}$ the canonical image of $G_{R(S)}$ in
 $(G_{R(S)}[Y]/ncl(T(X^\mu,Y)))^* .$

 Since $Rad_{G_{R(S)}}(T(X^\mu,Y))$ is a normal subgroup in
$G_{R(S)}[Y]$ containing
 $T(X^\mu,Y)$  there exists  a canonical $G$-epimorphism
 $$ \psi:  G_{R(S)}[Y]/ncl(T(X^\mu,Y)) \rightarrow
G_{R(S)}[Y]/Rad_{G_{R(S)}}(T(X^\mu,Y)).$$
 By Lemma \ref{le:6.2} the homomorphism $\psi$  gives rise to a
canonical
  $G$-homomorphism
   $$ \psi^*:  (G_{R(S)}[Y]/ncl(T(X^\mu,Y)))^* \rightarrow
(G_{R(S)}[Y]/Rad_{G_{R(S)}}(T(X^\mu,Y)))^*.$$
  Notice that the group $G_{R(S)}[Y]/Rad_{G_{R(S)}}(T(X^\mu,Y))$ is the
  coordinate group of the system $T(X^\mu,Y)= 1$ over  $G_{R(S)}$  and
  this system has a solution in  $G_{R(S)}$. Therefore this group is a
  $G_{R(S)}$-group and  it is $G_{R(S)}$-separated by
  $G_{R(S)}$. Now since $G_{R(S)}$ is the coordinate group of $S(X) = 1$
over $G$ and
  this system has a solution in $G$, we see that $G_{R(S)}$ is $G$-separated by $G$.
  It follows that the group $G_{R(S)}[Y]/Rad_{G_{R(S)}}(T(X^\mu,Y))$  is
$G$-separated
   by $G$. Therefore
  $$G_{R(S)}[Y]/Rad_{G_{R(S)}}(T(X^\mu,Y)) =
(G_{R(S)}[Y]/Rad_{G_{R(S)}}(T(X^\mu,Y)))^*.$$
  Now we can see that
  $$ \psi^* : G_{R(S\cup T)} \rightarrow
G_{R(S)}[Y]/Rad_{G_{R(S)}}(T(X^\mu,Y))$$
  is a $G$-homomorphism which maps the subgroup $\overline{G_{R(S)}}$
from
  $G_{R(S \cup T)}$   onto the subgroup $G_{R(S)}$ in
  $G_{R(S)}[Y]/Rad_{G_{R(S)}}(T(X^\mu,Y))$.

  This shows  that $\overline{G_{R(S)}} \simeq_G G_{R(S)}$ and
$\psi^*$ is a  $G_{R(S)}$-homomorphism. If $U$ is a solution of
$T(X^\mu,Y) = 1$ in  $G_{R(S)}$, then there exists a
$G_{R(S)}$-homomorphism
 $$\phi_U:
G_{R(S)}[Y]/Rad_{G_{R(S)}}(T(X^\mu,Y)) \rightarrow G_{R(S)}.$$
  such that $\phi_U(Y) = U$.  Obviously, composition of $\phi_U$ and
$\psi^*$
  gives a  $G_{R(S)}$-homomorphism from $G_{R(S\cup T)}$ into
$G_{R(S)}$,
 as desired. $\Box$.

The next result characterizes lifts in terms of the coordinate
groups of the corresponding equations.
\begin{prop}
\label{pr:Slift}
 Let $S(X) = 1$ be an equation over $G$ which has a solution in
$G$. Then for an arbitrary  equation  $T(X,Y) = 1$ over $G$  the
following conditions are equivalent:
\begin{enumerate}
\item [1)] $T(X,Y) = 1$ is compatible with $S(X) = 1$ and $T(X,Y)
= 1$ admits $S$-lift over $G$; \item [2)] $G_{R(S)}$ is a retract
of $G_{R(S,T)}$, i.e., $G_{R(S)}$ is a subgroup of $G_{R(S,T)}$
and there exists a $G_{R(S)}$-homomorphism $G_{R(S,T)} \rightarrow
G_{R(S)}.$
\end{enumerate}
\end{prop}
{\em Proof.}   1) $\Longrightarrow$ 2). By Proposition
\ref{prop:6.1} $G_{R(S)}$ is a subgroup of $G_{R(S,T)}$. Moreover,
$T(X^\mu,Y) = 1$ has a solution in $G_{R(S)}$, so by Lemma
\ref{le:oger} there exists a $G_{R(S)}$- homomorphism $G_{R(S,T)}
\rightarrow G_{R(S)},$ i.e., $G_{R(S)}$ is a retract of
$G_{R(S,T)}$.

2) $\Longrightarrow$ 1). If $\phi: G_{R(S,T)} \rightarrow
G_{R(S)}$ is a retract then every $G$-homomorphism $\alpha:
G_{R(S)}  \rightarrow G$ extends to a $G$-homomorphism $\alpha
\circ \phi :  G_{R(S,T)} \rightarrow G$. It follows from
 Proposition \ref{prop:6.1} that $T(X,Y) = 1$ is compatible with $S(X) =
1$ and
 $\phi$ gives a solution of $T(X^\mu,Y) = 1$ in $G_{R(S)}$, as desired.
\hfill$\Box$

One can ask whether it is possible to lift a system of equations
and inequalities into a generic point of some equation $S = 1$?
This is the question that we are going to address below. We start
with very general definitions.
\begin{df}
Let $S(X) = 1$ be an equation over a group $G$ which has a
solution in $G$.  We say that a formula $\Phi(X,Y)$ in the
language $L_A$ is compatible with $S(X) = 1$ over $G$, if for
every solution $\bar{a}$ of $S(X) = 1$ in $G$ there exists a tuple
$\bar{b}$ over $G$ such that the formula $\Phi(\bar a, \bar b )$
is true in $G$, i.e., the algebraic set $V_G(S)$ is a projection
of the truth  set of the formula $\Phi(X,Y) \ \wedge \ (S(X) = 1).
$
\end{df}
\begin{df}
Let a formula $\Phi(X,Y)$ be compatible with $S(X) = 1 $ over $G$.
We say that $\Phi(X,Y)$ admits a lift to a generic point of $S =
1$ over $G$ (or shortly $S$-lift over $G$), if $\exists Y
\Phi(X^\mu,Y)$ is true in $G_{R(S)}$ (here $Y$ are variables and
$X^\mu$ are constants from $G_{R(S)}$).
\end{df}

\begin{df} Let $S(X) = 1$ be an equation over $G$ which
has a solution in  $G$, and let $T(X,Y) = 1$ be compatible with
$S(X)=1$. We say that an equation $T(X,Y) = 1$ admits a complete
$S$-lift if every formula $T(X,Y) = 1 \ \& \ W(X,Y) \neq 1$, which
is compatible with $S(X) = 1$
 over $G$,  admits an $S$-lift.
\end{df}

\section{Implicit function theorem: lifting  solutions  into
 generic points}
Now we are ready to formulate and prove the main results of this
paper, Theorems \ref{1,2,3,4}, \ref{reg}, and \ref{tqe}. Let
$F(A)$ be a free non-abelian group.
\begin{theorem}\label{1,2,3,4}
Let $S(X,A)=1$ be a regular standard  quadratic equation over
$F(A)$. Every equation $T(X,Y,A) = 1$ compatible with $S(X,A) = 1$
admits a complete $S$-lift.
\end{theorem}
 We divide the proof of this theorem into two parts: for orientable  $S(X,A) = 1$,
 and for  a  non-orientable one.

\subsection{Basic Automorphisms of Orientable Quadratic  Equations}
\label{se:7.2.5}

In this section, for a finitely generated fully residually free
group $G$ we introduce some particular $G$-automorphisms of a free
$G$-group $G[X]$ which fix a given standard orientable quadratic
word with coefficients in $G$. Then we describe some cancellation
properties of  these automorphisms.

 Let $G$ be a group and let $S(X) = 1$ be a regular standard orientable
quadratic  equation over $G:$
\begin{equation}
\label{5}
 \prod_{i=1}^{m}z_i^{-1}c_iz_i\prod_{i=1}^{n}[x_i,y_i] d^{-1} = 1,
 \end{equation}
where $c_i, d$ are non-trivial constants from $G$,  and
  $$X = \{x_i, y_i, z_j \mid i = 1, \dots, n, j = 1, \dots, m\}$$
is the set of  variables. Observe that if $n=0$, then $m\geq 3$ by
definition of a regular quadratic equation (Definition
\ref{regular}). Sometimes we omit $X$ and write simply $S = 1$.
Denote by
$$C_S =  \{c_1, \ldots, c_m, d\}$$
 the set of constants which occur in the equation $S = 1$.

 Below we  define a {\em basic sequence}
$$\Gamma = (\gamma_1, \gamma_2, \ldots, \gamma_{K(m,n)})$$
 of $G$-automorphisms of the free $G$-group $G[X]$,  each of which
 fixes  the element
$$S_0 = \prod_{i=1}^{m}z_i^{-1}c_iz_i\prod_{i=1}^{n}[x_i,y_i] \in G[X].$$  We
assume that each $\gamma \in \Gamma$    acts identically on all
the generators from $X$ that are not mentioned in the description
of $\gamma$.

\medskip \noindent
Let $m \geq 3, n = 0$. In this case $K(m,0) = m-1.$ Put

\smallskip
$\gamma _{i} \ \ \ :\ z_i\rightarrow
z_i(c_i^{z_i}c_{i+1}^{z_{i+1}}), \ \ \ z_{i+1}\rightarrow
z_{i+1}(c_i^{z_i}c_{i+1}^{z_{i+1}})$, \ \ \ for $i=1,\dots ,m-1$.

\medskip \noindent
Let $m = 0, n\geq 1$. In this case $K(0,n) = 4n-1.$ Put

\medskip
$\gamma _{4i-3}:\ y_i\rightarrow x_iy_i, $ \ \ \ for $i=1,\dots
,n$;

\smallskip
$\gamma _{4i-2}:\ x_i\rightarrow y_ix_i, $ \ \ \ for $i=1,\dots
,n$;

\smallskip
$\gamma _{4i-1}:\ y_i\rightarrow x_iy_i, $ \ \ \  for $i=1,\dots
,n$;

\smallskip
$\gamma _{4i} \ \ \ :\ x_i\rightarrow (y_ix_{i+1}^{-1})^{-1}x_i,\
\ \ y_i\rightarrow y_i^{y_ix_{i+1}^{-1}},\ \ \  x_{i+1}\rightarrow
x_{i+1}^{y_ix_{i+1}^{-1}}, \ \ \ y_{i+1}\rightarrow
(y_ix_{i+1}^{-1})^{-1}y_{i+1}$,
\smallskip
  for  $i=1,\dots ,n-1$.

 \medskip \noindent
 Let $m \geq 1, n\geq 1$. In this case $K(m,n) = m + 4n-1.$ Put

\smallskip
$\gamma _{i} \ \ \ :\ z_i\rightarrow
z_i(c_i^{z_i}c_{i+1}^{z_{i+1}}), \ \ \ z_{i+1}\rightarrow
z_{i+1}(c_i^{z_i}c_{i+1}^{z_{i+1}})$, \ \ \ for $i=1,\dots ,m-1$;

\smallskip
 $\gamma _{m}\ \ \ :\ \ \ z_m\rightarrow z_m(c_m^{z_m}x_1^{-1}),\ \ \
x_1\rightarrow x_1^{c_m^{z_m}x_1^{-1}},\ \ \  y_1\rightarrow
(c_m^{z_m}x_1^{-1})^{-1}y_1;$

\medskip
$\gamma _{m + 4i-3}:\ y_i\rightarrow x_iy_i, $ \ \ \ for
$i=1,\dots ,n$;

\smallskip
$\gamma _{m + 4i-2}:\ x_i\rightarrow y_ix_i, $ \ \ \ for
$i=1,\dots ,n$;

\smallskip
$\gamma _{m + 4i-1}:\ y_i\rightarrow x_iy_i, $ \ \ \  for
$i=1,\dots ,n$;

\smallskip
$\gamma _{m + 4i} \ \ \ :\ x_i\rightarrow
(y_ix_{i+1}^{-1})^{-1}x_i,\ \ \ y_i\rightarrow
y_i^{y_ix_{i+1}^{-1}},\ \ \  x_{i+1}\rightarrow
x_{i+1}^{y_ix_{i+1}^{-1}}, \ \ \ y_{i+1}\rightarrow
(y_ix_{i+1}^{-1})^{-1}y_{i+1}$,
\smallskip
  for  $i=1,\dots ,n-1$.

\medskip

It is easy to check that each $\gamma \in \Gamma$ fixes the word
$S_0$ as well as
 the word $S$. This shows that  $\gamma$ induces a $G$-automorphism on the group
$G_S = G[X]/{\rm ncl}(S)$. We denote  the induced automorphism
again by $\gamma$, so $\Gamma \subset Aut_G(G_S)$.  Since $S = 1$
is regular, $G_S = G_{R(S)}$. It follows
 that composition of any product of automorphisms from
$\Gamma$ and  a particular  solution $\beta$ of $S = 1$  is again
a solution of $S = 1$.

Observe, that in the case $m \neq 0, n\neq 0$  the basic sequence
of automorphisms $\Gamma$
 contains the basic automorphisms from the other two cases. This allows
us, without loss of
 generality, to formulate some of the  results below only for the case $K(m,n)= m + 4n - 1$.
 Obvious adjustments provide the proper argument in the other cases.
 From now on we  order elements of the set $X$
 in the following  way
 $$z_1<\ldots <z_m<x_1<y_1<\ldots <x_n<y_n.$$
For a word $w \in F(X)$ we denote by $v(w)$ the {\em leading}
variable (the highest variable with respect to the order
introduced above) that occurs in $w$. For $v = v(w)$ denote by
$j(v)$ the following number

 \[ j(v) = \left\{\begin{array}{ll}
                  m+4i, & \mbox{if $v = x_i$ or $v = y_i$ and $i < n$,}\\
                  m+4i-1,  & \mbox{if $v = x_i$ or $v = y_i$ and $i =                  n$,}\\
                  i, &   \mbox{if $v = z_i$ and $n\neq 0$,} \\
                  m-1, &  \mbox{if $v = z_m$, n= 0.}
                 \end{array}
                   \right. \]

The following lemma describes the action of powers of basic
automorphisms from
 $\Gamma$ on $X$. The proof is obvious, and we omit it.

\begin{lm}
\label{le:7.1.29} Let $\Gamma = (\gamma_1, \ldots,
\gamma_{m+4n-1})$ be the basic sequence of automorphisms and $p$
be a  positive integer. Then the following hold: \bea \gamma
_{i}^p \ \ \ &:& \ z_i\rightarrow
z_i(c_i^{z_i}c_{i+1}^{z_{i+1}})^p, \ \ \ z_{i+1}\rightarrow
z_{i+1}(c_i^{z_i}c_{i+1}^{z_{i+1}})^p, \\
&& \qquad\qquad  \mbox{for}\;\; i=1,\dots
,m-1;\\
\gamma _{m}^p \ \  &:&\ z_m\rightarrow z_m(c_m^{z_m}x_1^{-1})^p,\
\ \ x_1\rightarrow x_1^{(c_m^{z_m}x_1^{-1})^p},\ \ \
y_1\rightarrow
(c_m^{z_m}x_1^{-1})^{-p}y_1;\\
\gamma _{m+4i-3}^p &:&\ y_i\rightarrow x_i^py_i,  \ \ \
\mbox{for}\;\;
i=1,\dots ,n;\\
\gamma _{m+4i-2}^p &:&\ x_i\rightarrow y_i^px_i,  \ \ \
\mbox{for}\;\;
i=1,\dots ,n;\\
\gamma _{m+4i-1}^p &:&\ y_i\rightarrow x_i^py_i,  \ \ \
\mbox{for}\;\;
i=1,\dots ,n;\\
\gamma _{m+4i}^p \ \ \ &:&\ x_i\rightarrow
(y_ix_{i+1}^{-1})^{-p}x_i,\ y_i\rightarrow
y_i^{(y_ix_{i+1}^{-1})^p},\\
 &&\ x_{i+1}\rightarrow x_{i+1}^{(y_ix_{i+1}^{-1})^p},\ \ \
y_{i+1}\rightarrow
 (y_ix_{i+1}^{-1})^{-p}y_{i+1},\\
&& \qquad\qquad \mbox{for}\;\; i=1,\dots ,n-1. \eea
 \end{lm}

The $p$-powers of elements that occur in Lemma  \ref{le:7.1.29}
play an important part in what follows, so we describe them in  a
separate definition.
\begin{df}
\label{de:A(gamma)} Let $\Gamma = (\gamma_1, \ldots,
\gamma_{m+4n-1})$ be the basic sequence of automorphism for $S=
1$. For every $\gamma \in \Gamma$ we define the leading term
$A(\gamma)$ as follows:

\smallskip $A(\gamma _{i}) = c_i^{z_i}c_{i+1}^{z_{i+1}}$, \ \ \ for
$i=1,\dots ,m-1$;

\smallskip
 $A(\gamma _{m}) = c_m^{z_m}x_1^{-1};$

\smallskip
 $A(\gamma _{m+4i-3}) = x_i, $ \ \ \ for $i=1,\dots ,n$;

\smallskip $A(\gamma _{m+4i-2}) =  y_i, $ \ \ \ for $i=1,\dots ,n$;

\smallskip $A( \gamma _{m+4i-1}) =  x_i, $ \ \ \ for $i=1,\dots ,n$;

\smallskip $A(\gamma _{m+4i}) = y_ix_{i+1}^{-1} $,  for $i=1,\dots
,n-1$.

\end{df}


Now we introduce  vector notations for automorphisms of particular
type.

 Let ${\mathbb{N}}$ be the set of all positive integers and  ${\mathbb{N}}^k$
 the set of all $k$-tuples of elements from ${\mathbb{N}}$.
 For $s \in {\mathbb{N}}$ and $p \in {\mathbb{N}}^k$ we say that the tuple $p$ is
{\it $s$-large} if every coordinate of $p$  is greater then  $s$.
Similarly, a subset $P \subset  {\mathbb{N}}^k$ is {\em $s$-large}
if every tuple in $P$ is $s$-large.   We say that the set $P$ is
{\em unbounded} if for any $s \in  {\mathbb{N}}$ there exists an
$s$-large tuple in $P$.

Let $\delta = (\delta_1, \ldots, \delta_k)$ be a sequence of
$G$-automorphisms of the group $G[X]$, and $p = (p_1, \ldots,p_k)
\in {\mathbb{N}}^k$.  Then by $\delta^p$ we denote the following
automorphism of $G[X]$: $$\delta^p = \delta_1^{p_1} \cdots
\delta_k^{p_k}.$$

\begin{notation}
  Let $\Gamma = (\gamma_1, \dots, \gamma_K)$ be the
 basic sequence of automorphisms for $S = 1$. Denote by
$\Gamma_{\infty}$ the
 infinite periodic sequence with period $\Gamma$, i.e.,
 $\Gamma_{\infty} = \{\,\gamma_i\,\}_{i \geq 1}$ with $\gamma_{i+K}= \gamma_i$.
 For $j \in {\mathbb{N}}$ denote by $\Gamma_j$ the initial segment of
 $\Gamma_{\infty}$ of length $j$. Then for a given $j$ and  $p \in
{\mathbb{N}}^j$ put
 $$ \phi_{j,p} =\stackrel{\leftarrow}{\Gamma}_j^{\stackrel{\leftarrow}{p}} =\gamma_j^{p_j}\gamma_{j-1}^{p_{j-1}} \cdots \gamma_1^{p_1}.$$ Sometimes
we omit $p$ from $ \phi_{j,p}$ and write simply $\phi_j$.
\end{notation}

{\bf Agreement.} {\it From now on we fix an arbitrary positive
multiple $L$ of the number $K = K(m,n)$,  a $2$-large tuple $p \in
{\mathbb{N}}^L$, and the automorphism $\phi = \phi_{L,p}$ {\rm
(}as well as all the automorphism $\phi_j$, $j \leq L${\rm )}. }

\begin{df}
The leading term $A_j=A(\phi_j)$ of the automorphism $\phi_j$ is
defined to be the
 cyclically  reduced form of the word

 \[ \left\{\begin{array}{ll}
                  A(\gamma_j)^{\phi_{j- 1}}, & \mbox{if $j\leq K$, $j \neq  m+4i-1$ \ for \ any \
                  $i = 1, \ldots,n;$},\\
                  y_i^{-\phi _{j-2}} A(\gamma _{j})^{\phi _{j-1}}y_i^
                  {\phi _{j-2}},
                   & \mbox{if $j = m+4i-1$ \ for \ some \
                   $i = 1, \ldots,n;$}\\
  A_r^{\phi _{sK}} &\text{if $j=r+sK,\ r\leq K,  s \in {\mathbb{N}}$.}
\end{array}
                   \right. \]
\end{df}

\begin{lm}
\label{le:Apower} For every $j \leq L$ the element $A(\phi_j)$ is
not a proper power in $G[X]$. \end{lm}

{\it Proof.}
 It is easy to check that $A(\gamma_s)$ from  Definition  \ref{de:A(gamma)}
 is not a proper power for $s = 1, \dots, K.$ Since  $A(\phi_j)$ is
the  image of some $A(\gamma_s)$ under an automorphism of $G[X]$
it is not a proper power in $G[X]$. \hfill $\Box$

For words $w, u, v \in G[X]$,  the notation $\ig{w}{u}{v}$ \ means
that $w = u \circ w^\prime \circ v$ for some $w^\prime \in G[X]$,
where the length of elements and reduced form defined as in the
free product $G*\left<X\right>$. Similarly, notations $\igl{w}{u}$
\  and   $ \igr{w}{v}$ \ mean that $w = u \circ w^\prime$ and $w =
w^\prime \circ v$. Sometimes we write $\ig{w}{u}{\ast}$ \ or
$\ig{w}{\ast}{v}$ \ when the corresponding words are irrelevant.


 If $n$ is a positive integer and $w \in G[X]$, then by $Sub_n(w)$ we denote the
set of all $n$-subwords of $w$, i.e.,
$$Sub_n(w) = \{u \ \mid \ |u| = n \ and \  w = w_1 \circ u \circ w_2 \ for \ some \
w_1, w_2 \in G[X] \}.$$ Similarly, by $SubC_n(w)$ we denote all
$n$-subwords of the {\it cyclic} word $w$. More generally, if $W
\subseteq G[X]$, then
$$Sub_n(W) = \bigcup_{w \in W} Sub_n(w), \ \ \ SubC_n(W) = \bigcup_{w \in W}
SubC_n(w).$$ Obviously, the set $Sub_i(w)$ ($SubC_i(w)$) can be
effectively reconstructed from $Sub_n(w)$  ($SubC_n(w)$) for $i
\leq n$.

 In the following  series of lemmas we write  down explicit
expressions for images of elements of $X$ under the automorphism
 $$\phi_K = \gamma_K^{p_K} \cdots \gamma_1^{p_1}, \ \ \ K = K(m,n).$$
 These lemmas are very easy and straightforward, though tiresome in terms of
 notations. They provide  basic data needed to  prove the implicit function theorem.
 All elements that occur in the lemmas below can be viewed as
 elements (words) from the free group $F(X \cup C_S)$. In particular, the notations $\circ$,
 $\raisebox{1ex}{\ig{w}{u}{v}}$ \ , and $Sub_n(W)$ correspond to the  standard length function on $F(X \cup C_S)$.
Furthermore, until the end of this section we assume that the
elements $c_1, \dots, c_m$ are {\em pairwise different}.

\begin{lm}
\label{le:7.1.zforms}  Let $m \neq 0$,  $K = K(m,n)$, $p = (p_1,
\ldots, p_K)$ be a 3-large tuple, and
 $$\phi_K = \gamma_K^{p_K} \cdots \gamma_1^{p_1}.$$
 The following statements hold.
\begin{enumerate}
\item  [(1)] All automorphisms from $\Gamma$, except for
$\gamma_{i-1}, \gamma_{i}$ {\rm (}if  defined\/{\rm )},  fix
$z_i$, $i=1,\dots, m$. It follows that
 $$z_i^{\phi_K} \ldots = z_i^{\phi_{i}} \ (i = 1, \ldots,
 m-1).$$

\item [(2)]
 Let ${\tilde z}_i = z_i^{\phi_{i-1}}$ \  ($i = 2, \ldots, m$), $ \tilde z_1=z_1.$ Then

\medskip
${\tilde z}_i  = \ig{ z_i \circ ( c_{i-1}^{{\tilde z}_{i-1}} \circ
c_i^{{z_i}})^{p_{i-1}} }{z_iz_{i-1}^{-1}}{c_iz_i} $ \ \ \ {\rm
(}$i = 2, \ldots, m${\rm )}.

\item  [(3)] The reduced forms of the leading terms of the
corresponding automorphisms are listed below:

\medskip
 $A_1 = \ig{c_1^{z_1} \circ c_2^{z_2}}{z_1^{-1}c_1}{c_2z_2}$ \  \
 \  {\rm (}$m\geq 2${\rm )},

\medskip
$A_2=(c_2^{z_2}x_1^{-1})^{\phi
_1}=A_1^{-p_1}c_2^{z_2}A_1^{p_1}x_1^{-1}\ (n\neq 0, m=2),$

\medskip
$A_2=A_1^{-p_1}c_2^{z_2}A_1^{p_1}c_3^{z_3}\ (n\neq 0, m>2),$

\medskip
$SubC_3(A_1) = \{z_1^{-1}c_1z_1, \ c_1z_1z_2^{-1}, \
z_1z_2^{-1}c_2, \ z_2^{-1}c_2z_2, \ c_2z_2z_1^{-1}, \
z_2z_1^{-1}c_1\};$

\medskip
 $A_{i} = \ig{A_{i-1}^{-p_{i-1}}}{z_i^{-1}c_i^{-1}}{c_{i-1}z_{i-1}}\
 c_i^{z_i} \ \ig{A_{i-1}^{p_{i-1}}}{z_{i-1}^{-1}c_{i-1}^{-1}}{c_iz_i}
   \ig{c_{i+1}^{z_{i+1}}}{z_{i+1}^{-1}}{c_{i+1}z_{i+1}}$\ ,
 \newline( $i=3,\dots ,m-1$),

\medskip
$SubC_3(A_{i}) = SubC_3(A_{i-1})^{\pm 1} \cup
\{c_{i-1}z_{i-1}z_i^{-1}, \ z_{i-1}z_i^{-1}c_i, \ z_i^{-1}c_iz_i,
\ c_iz_iz_{i-1}^{-1}$, $$ \ z_iz_{i-1}^{-1}c_{i-1}^{-1}, \
c_iz_iz_{i+1}^{-1}, \ z_iz_{i+1}^{-1}c_{i+1}, \
z_{i+1}^{-1}c_{i+1}z_{i+1}, \ c_{i+1}z_{i+1}z_i^{-1}, \
z_{i+1}z_i^{-1}c_i^{-1} \};$$

\medskip
$A_{m} =  (c_m^{z_m}x_1^{-1})^{\phi _{m-1}}
=\ig{A_{m-1}^{-p_{m-1}}}{z_m^{-1}c_m^{-1}}{c_{m-1}z_{m-1}} \
c_m^{z_m} \
\ig{A_{m-1}^{p_{m-1}}}{z_{m-1}^{-1}c_{m-1}^{-1}}{c_mz_m} \
x_1^{-1}$ \ \ \ ($n \neq 0, m\geq 3$),

\medskip
$SubC_3(A_{m}) = SubC_3(A_{m-1})^{\pm 1} \cup
\{c_{m-1}z_{m-1}z_m^{-1}, \ z_{m-1}z_m^{-1}c_m, $ $$\
z_{m-1}^{-1}c_mz_m, \ c_mz_mz_{m-1}^{-1}, \ c_mz_mx_1^{-1},
z_mx_1^{-1}z_m^{-1}, \ x_1^{-1}z_m^{-1}c_m^{-1} \}.$$

\item   [(4)] The reduced forms of $z_i^{\phi_{i-1}},
z_i^{\phi_{i}}$ are listed below:

\medskip
 $z_1^{\phi_K} = z_1^{\phi _{1}} =c_1  \ig{z_1 c_2^{z_2}}{z_1z_2^{-1}}{c_2z_2}
\ig{A_1^{p_1-1}}{z_1^{-1}c_1}{c_2z_2}\ \ (m\geq 2)$,

\medskip
$SubC_3(z_1^{\phi_K}) = \{ \ c_1z_1z_2^{-1}, \ z_1z_2^{-1}c_2, \
z_2^{-1}c_2z_2, \ c_2z_2z_1^{-1}, \ z_2z_1^{-1}c_1, \
z_1^{-1}c_1z_1\};$

\medskip
$z_i^{\phi_{i-1}} = z_i
\ig{A_{i-1}^{p_{i-1}}}{z_{i-1}^{-1}c_{i-1}^{-1}}{c_iz_i} \ =\tilde
z_i$,

\medskip
$ z_{i}^{\phi _{K}} =   z_{i}^{\phi _{i}}   = c_i z_i \
\ig{A_{i-1}^{p_{i-1}}}{z_{i-1}^{-1}c_{i-1}^{-1}}{c_iz_i} \
c_{i+1}^{z_{i+1}}
\ig{A_{i}^{p_{i}-1}}{z_i^{-1}c_i^{-1}}{c_{i+1}z_{i+1}} $ \ \
\newline ( $i=2,\dots ,m-1$),

\medskip
$Sub_3(z_i^{\phi_K}) = SubC_3(A_{i-1}) \cup SubC_3(A_{i}) \cup \{
c_iz_iz_{i-1}^{-1}, \ z_iz_{i-1}^{-1}c_{i-1}^{-1}, $
 $$\ c_iz_iz_{i+1}^{-1}, \ z_iz_{i+1}^{-1}c_{i+1}, \ z_{i+1}^{-1}c_{i+1}z_{i+1}, \
c_{i+1}z_{i+1}z_i^{-1}, \ z_{i+1}z_i^{-1}c_i^{-1} \} \ ;$$

\medskip
$z_m^{\phi_K} = z_m^{\phi_{m-1}} =  z_m
\ig{A_{m-1}^{p_{m-1}}}{z_{m-1}^{-1}c_{m-1}^{-1}}{c_mz_m}\ $, \ \ \
($n = 0$),

\medskip
$Sub_3(z_m^{\phi_K})_{(when \ n = 0)} = SubC_3(A_{m-1}) \cup \{
z_mz_{m-1}^{-1}c_{m-1}^{-1} \} \ ;$

\medskip
$z_m^{\phi_K} = z_m^{\phi_{m}} =  c_m z_{m}
\ig{A_{m-1}^{p_{m-1}}}{z_{m-1}^{-1}c_{m-1}^{-1}}{c_mz_m}\ x_1^{-1}
\ \ig{A_{m}^{p_{m}-1}}{z_m^{-1}c_m^{-1}}{z_mx_1^{-1}}$\ \ \ ($n
\neq 0, m\geq 2$),

\medskip
$Sub_3(z_m^{\phi_K}) = Sub_3(z_m^{\phi_K})_{(when \ n = 0)} \cup
SubC_3(A_m)\cup \{c_mz_mx_1^{-1}, \ z_mx_1^{-1}z_m^{-1}, \
x_1^{-1}z_m^{-1}c_m^{-1} \}.$

\item [(5)] The elements $z_i^{\phi _K}$ have the following
properties:

\medskip
$z_i^{\phi _K}=c_iz_i\hat z_i \ \ \ (i=1,\dots ,m-1)$, where $\hat
z_i$ is a word in the alphabet $\{ c_1^{z_1}, \ldots,
c_{i+1}^{z_{i+1}} \}$ which begins with $c_{i-1}^{-z_{i-1}}$, if
$i\neq 1$, and with $c_2^{z_2}$, if $i=1$;

\medskip
$z_m^{\phi _K}=z_m\hat z_m  \ \ \ (n=0)$, where $\hat z_m$ is a
word in the alphabet $\left\{ c_1^{z_1}, \ldots,
c_{m}^{z_{m}}\right \};$

\medskip
$z_m^{\phi _K}=c_mz_m\hat z_m  \ \ \ (n \neq 0)$, where $\hat z_m$
is a word in the alphabet $\{ c_1^{z_1}, \ldots, c_{m}^{z_{m}},
x_1 \};$

 \medskip
 Moreover, if $m\geq 3$ the word $(c_m^{z_m})^{\pm 1}$ occurs in  $z_i^{\phi _K} \ \ \ (i = m-1,m)$
 only  as a part of the subword $(\prod _{i=1}^m c_{i}^{z_{i}})^{\pm 1}.$
\end{enumerate}
\end{lm}

{\it Proof.}
 (1) is obvious. We prove (2) by induction. For $i\geq
2$

\medskip
${\tilde z}_i = z_i^{\phi_{i-1}} =
z_i^{\gamma_{i-1}^{p_{i-1}}\phi_{i-2 }}$

\medskip\noindent
Therefore, by induction,

\medskip
${\tilde z}_i = z_i ( c_{i-1}^{{\tilde z}_{i-1}}
c_i^{{z_i}})^{p_{i-1}} = z_i \circ ( c_{i-1}^{{\tilde z}_{i-1}}
\circ c_i^{{z_i}})^{p_{i-1}}. $

 Now we prove (3) and (4) simultaneously. Let $m\geq 2.$ By the straightforward verification one has:

\medskip
$A_1 = \ig{c_1^{z_1} \circ c_2^{z_2}}{z_1^{-1}}{z_2}$;

\medskip
$z_1^{\phi_1} = z_1^{\gamma _1^{p_1}} =
z_1(c_1^{z_1}c_2^{z_2})^{p_1} = \ig{c_1 \circ z_1 \circ c_2^{z_2}
\circ A_1^{p_1-1}}{c_1}{z_2}$.

Denote by cycred $(w)$ the cyclically reduced form of $w$.

\medskip
$A_i =$ cycred $( \left (c_i^{z_i}c_{i+1}^{z_{i+1}} \right
)^{\phi_{i-1}})= \ig {c_i^{{\tilde z}_i} \circ
c_{i+1}^{z_{i+1}}}{z_i^{-1}}{z_{i+1}}\ (i\leq m-1).$

\medskip \noindent
Observe that in the notation above

\medskip
${\tilde z}_i =  z_iA_{i-1}^{p_{i-1}} \ (i\geq 2).$

\medskip \noindent
This shows that we can rewrite  $A(\phi_{i})$ as follows:

\medskip
$A_i =  A_{i-1}^{-p_{i-1}}\circ c_i^{z_i}\circ
  A_{i-1}^{p_{i-1}} \circ
  c_{i+1}^{z_{i+1}}$,

\medskip \noindent
 beginning with $z_i^{-1}$ and ending with $z_{i+1}$ \ \ ( $i=2,\dots ,m-1$);

\medskip
$A_{m} =$ cycred $(c_m^{\tilde z_m}x_1^{-1}) = c_m^{\tilde
z_m}x_1^{-1} = A_{m-1}^{-p_{m-1}} \circ c_m^{z_m}\circ
A_{m-1}^{p_{m-1}} \circ x_1^{-1}\ (m\geq 2).$

\medskip \noindent
beginning with $z_m^{-1}$ and ending with $x_1^{-1}$  \ \ ( $n
\neq 0$).

\medskip
$z_i^{\phi_{i-1}} = \left (z_i
(c_{i-1}^{z_{i-1}}c_i^{z_i})^{p_{i-1}} \right )^{\phi_{i-2}} = z_i
(c_{i-1}^{{\tilde z}_{i-1}}c_i^{z_i})^{p_{i-1}} = z_i \circ
A_{i-1}^{p_{i-1}} $,

\medskip \noindent
beginning with $z_i$ and ending with $z_i$;

\medskip
$z_i^{\phi_{i}} = \left (z_i(c_i^{z_i}c_{i+1}^{z_{i+1}})^{p_{i}}
\right )^{\phi_{i-1}} = {\tilde z}_i(c_i^{{\tilde
z}_i}c_{i+1}^{z_{i+1}})^{p_{i}} = c_i \circ {\tilde z}_i \circ
c_{i+1}^{z_{i+1}} \circ (c_i^{{\tilde
z}_i}c_{i+1}^{z_{i+1}})^{p_{i}-1} = $

\medskip
$c_i \circ z_i \circ A_{i-1}^{p_{i-1}} \circ c_{i+1}^{z_{i+1}}
\circ A_{i}^{p_{i}-1}$,

\medskip \noindent
beginning with $c_i$ and ending with $z_{i+1}$ \ \  ($i = 2,
\ldots, m-1$);

\medskip
$z_m^{\phi_{m}} = \left (z_m(c_m^{z_m}x_1^{-1})^{p_{m}} \right
)^{\phi_{m-1}} =  {\tilde z}_m(c_m^{{\tilde z}_m}x_1^{-1})^{p_{m}}
= c_m {\tilde z}_m x_1^{-1} (c_m^{{\tilde z}_m}x_1^{-1})^{p_{m}-1}
= $

\medskip
$c_m \circ z_{m} \circ A_{m-1}^{p_{m-1}}\circ x_1^{-1} \circ
A_{m}^{p_{m}-1}$\ \  \ ($n \neq 0$),

\medskip \noindent
beginning with $c_m$ and ending with $x_1^{-1}.$  This proves the
lemma.

(5) Direct verification using formulas in (3) and (4). \hfill
$\Box$

In the following two lemmas  we describe the  reduced expressions
of the elements $x_1^{\phi_K}$ and $y_1^{\phi_K}$.
\begin{lm}
\label{le:7.1.x1formsm0}
 Let $m = 0$,  $K = 4n-1$, $p = (p_1, \ldots, p_K)$ be a $3$-large tuple, and
 $$\phi_K = \gamma_K^{p_K} \ldots \gamma_1^{p_1}.$$

\begin{enumerate}
\item  [(1)] All automorphisms from $\Gamma $, except for  $\gamma
_{2}, \gamma _{4} $, fix $x_1$,  and all automorphisms from
$\Gamma $, except  for $\gamma _{1}, \gamma_3, \gamma _{4} $,  fix
$y_1$. It follows that

 \medskip
$ x_1^{\phi_K} = x_1^{\phi_4}, \ y_1^{\phi_K} = y_1^{\phi_4} \ \ \
(n \geq 2).$

\item [(2)] Below we list the reduced forms of the leading terms
of the corresponding automorphisms  (the words on the right are
reduced as written)

\medskip
$A_1 = x_1$;

\medskip
$A_2 =  x_1^{p_1}y_1 = A_1^{p_1}\circ y_1 $ ;

\medskip
$A_3 = \ig{A_2^{p_2-1}}{x_1^2}{x_1y_1} \ x_1^{p_1+1}  y_1$,
\medskip
$SubC_3(A_3) = SubC_3(A_2) = \{x_1^3, \  x_1^2y_1, \ x_1y_1x_1,
y_1x_1^2 \} \ ;$

\medskip
$A_4 =  \ig{\left(\ig{A_2^{p_2}}{x_1^2}{x_1y_1} \
x_1\right)^{p_3}}{x_1^2}{y_1x_1} \ig{A_2}{x_1^2}{x_1y_1}\
x_2^{-1}\ (n\geq 2)$,

\medskip
$SubC_3(A_4) = SubC_3(A_2) \cup \{x_1y_1x_2^{-1}, \
y_1x_2^{-1}x_1, \ x_2^{-1}x_1^2\}\ (n\geq 2).$

 \item [(3)] Below we list reduced forms of $x_1^{\phi _j}, y_1^{\phi _j}$ for $j = 1,
 \ldots, 4$:

\medskip
$x_1^{\phi_1} = x_1$;

\medskip
$y_1^{\phi_1} = x_1^{p_1} y_1; $

\medskip
$x_1^{\phi_2} =  \ig{A_2^{p_2}}{x_1^2}{x_1y_1}\  x_1;$

\medskip
$y_1^{\phi_2} = x_1^{p_1}y_1;$

\medskip
$x_1^{\phi_3} =  x_1^{\phi_2} = \ig{A_2^{p_2}}{x_1^2}{x_1y_1}\
x_1=_{({\rm when}\ n=1)}x_1^{\phi _K}$;

\medskip
$Sub_3(x_1^{\phi_K})_{(when \ n = 1)} = SubC_3(A_2)$;

\medskip
$y_1^{\phi_3} =   \ig{(\ig{A_2^{p_2}}{x_1^2}{x_1y_1} \
x_1)^{p_3}}{x_1^2}{x_1y_1} \ x_1^{p_1}y_1=_{({\rm when}\
n=1)}y_1^{\phi_K} ;$

\medskip
$Sub_3(y_1^{\phi_K})_{(when \ n = 1)} = SubC_3(A_2);$

 \medskip
 $x_1^{\phi_4} = x_1^{\phi_K} = \ig{A_4^{-(p_4-1)}}{x_2y_1^{-1}}{x_1^{-2}} \  x_2
 \ \ig{A_2^{-1}}{y_1^{-1} x_1^{-1}}{x_1^{-2}}  \ig{(x_1^{-1} \
 \ig{A_2^{-p_2}}{y_1^{-1}x_1^{-1}}{x_1^{-2}})^{p_3-1}}{x_1^{-1}y_1^{-1}x_1^{-1}}{x_1^{-2}}$
 \  $(n \geq 2)$,

 \medskip
 $Sub_3(x_1^{\phi_K}) = SubC_3(A_4)^{-1} \cup SubC_3(A_2)^{-1} \cup
 \{x_1^{-2}x_2, \ x_1^{-1}x_2y_1^{-1}, \ x_2y_1^{-1}x_1^{-1}, \
 x_1^{-3}, $
 $$\ x_1^{-2}y_1^{-1}, x_1^{-1}y_1^{-1}x_1^{-1} \}\ \ \ (n\geq 2);$$

\medskip
$y_1^{\phi_4} =  \ig{A_4^{-(p_4-1)}}{x_2y_1^{-1}}{x_1^{-1}} \ x_2
\ \ig{A_4^{p_4}}{x_1^2}{y_1x_2^{-1}}\ \ \ (n\geq 2), $

\medskip
$Sub_3(y_1^{\phi_K}) = SubC_3(A_4)^{\pm 1} \cup \{x_1^{-2}x_2, \
x_1^{-1}x_2x_1, \ x_2x_1^2 \} \ \ \ (n\geq 2).$
\end{enumerate}
 \end{lm}
{\it Proof.} (1) follows directly from definitions.

To show (2) observe that

\medskip
$A_1 = A(\gamma_1) = x_1$;

\medskip
$x_1^{\phi_1} = x_1$;

\medskip
$y_1^{\phi_1} = x_1^{p_1}y_1 = A_1^{p_1} \circ y_1.$

\medskip\noindent
Then \bea A_2 &=& {\rm cycred} (A(\gamma_2)^{\phi_1}) ={\rm
cycred} (y_1^{\phi_1}) = x_1^{p_1} \circ
y_1 = A_1^{p_1} \circ y_1;\\
x_1^{\phi_2} &=& (x_1^{\gamma_2^{p_2}})^{\gamma_1^{p_1}}=
(y_1^{p_2}x_1)^{\gamma_1^{p_1}}
= (x_1^{p_1}y_1)^{p_2} x_1 =A_2^{p_2} \circ x_1;\\
y_1^{\phi_2} &=& (y_1^{\gamma_2^{p_2}})^{\gamma_1^{p_1}}=
y_1^{\gamma_1^{p_1}} = x_1^{p_1}y_1 = A_2.\eea Now \bea A_3 &=&
{\rm cycred}(y_1^{-\phi_1} A(\gamma_3)^{\phi_2}y_1^{\phi_1})= {\rm
cycred}( (x_1^{p_1}y_1)^{-1} x_1^{\phi_2} (x_1^{p_1}y_1))={\rm
cycred} (
(x_1^{p_1}y_1)^{-1}(x_1^{p_1}y_1)^{p_2}x_1(x_1^{p_1}y_1))\\
& =&   (x_1^{p_1}y_1)^{p_2-1} x_1^{p_1+1}y_1\\
& =&  A_2^{p_2-1}\circ A_1^{p_1+1} \circ y_1. \eea

\medskip\noindent
It follows that

\medskip
$x_1^{\phi_3} = (x_1^{\gamma_3^{p_3}})^{\phi_2} = x_1^{\phi_2}$;

\medskip
$y_1^{\phi_3} = (y_1^{\gamma_3^{p_3}})^{\phi_2}
(x_1^{p_3}y_1)^{\phi_2} = (x_1^{\phi_2})^{p_3}y_1^{\phi_2}
(A_2^{p_2} \circ x_1)^{p_3}\circ
 A_2.$

 \medskip\noindent
 Hence
\[
A_4 ={\rm cycred} (A(\gamma_4)^{\phi_3}) ={\rm cycred}
((y_1x_2^{-1})^{\phi_3})={\rm cycred} (y_1^{\phi_3}x_2^{-\phi_3})
= (A_2^{p_2} \circ x_1)^{p_3}\circ
 A_2\circ x_2^{-1}.
 \]
 Finally:
\bea x_1^{\phi_4} &=& (x_1^{\gamma_4^{p_4}})^{\phi_3}  =\left
((y_1x_2^{-1})^{-p_4}x_1 \right )^{\phi_3} =
 \left ((y_1x_2^{-1})^{\phi_3})\right )^{-p_4}x_1^{\phi_3}= A_4^{-p_4}A_2^{p_2} \circ x_1\\
& = &  A_4^{-(p_4-1)} \circ x_2
 \circ A_2^{-1} \circ (x_1^{-1} \circ
A_2^{-p_2})^{p_3-1}\\
y_1^{\phi_4} &=& (y_1^{\gamma_4^{p_4}})^{\phi_3} =(y_1^{(y_1x_2^{-1})^{p_4}})^{\phi_3}\\
& = &  \left ((y_1x_2^{-1})^{\phi_3}\right ) ^{-p_4} y_1^{\phi_3}
\left
((y_1x_2^{-1})^{\phi_3}\right )^{p_4}\\
& = & A_4^{-p_4} y_1^{\phi_3}  A_4^{p_4}= A_4^{-(p_4-1)} A_4^{-1}
y_1^{\phi_3} A_4^{p_4} =\\ && \qquad\qquad A_4^{-(p_4-1)} \circ
x_2 \circ  A_4^{p_4} . \eea This proves the lemma. \hfill $\Box$

\begin{lm}
\label{le:7.1.x1formsmneq0}
 Let $m \neq 0$, $n \neq 0$,  $K = m + 4n-1$, $p = (p_1, \ldots, p_K)$ be a $3$-large tuple, and
 $$\phi_K = \gamma_K^{p_K} \cdots \gamma_1^{p_1}.$$

\begin{enumerate}
\item [(1)] All automorphisms from $\Gamma $ except  for $\gamma_{m},
\gamma _{m+2}, \gamma _{m+4} $ fix $x_1$; and all automorphisms
from $\Gamma $ except for $\gamma_{m}, \gamma _{m+1},
\gamma_{m+3}, \gamma _{m+4} $ fix $y_1$.  It follows that

\medskip
$x_1^{\phi_K} = x_1^{\phi_{m+4}}, \  y_1^{\phi_K}
=y_1^{\phi_{m+4}} \ \ \ (n \geq 2).$

\item [(2)] Below we list the reduced forms of the leading terms  of the
corresponding automorphisms  (the words on the right are reduced
as written)

$ A_{m+1} =  x_1,$

$A_{m+2} = y_1^{\phi_{m+1}}
=\ig{A_{m}^{-p_{m}}}{x_1z_m^{-1}}{c_mz_m} \ x_1^{p_{m+1}} y_1,$

$SubC_3(A_{m+2}) = SubC_3(A_{m})^{-1} \cup \{c_mz_mx_1, \
z_mx_1^2, \ x_1^3, \ x_1^2y_1, \ x_1y_1x_1, \ y_1x_1z_m^{-1} \};$

$ A_{m+3} = \ig{A_{m+2}^{p_{m+2}-1}}{x_1z_m^{-1}}{x_1y_1}
 \ig{A_{m}^{-p_{m}}}{x_1z_m^{-1}}{c_mz_m}\  x_1^{p_{m+1}+1}
  y_1, $

$SubC_3(A_{m+3}) = SubC_3(A_{m+2});$

$ A_{m+4} = \ig{A_{m}^{-p_{m}}}{x_1z_m^{-1}}{c_mz_m}  \circ \left
(x_1^{p_{m+1}} y_1 \
\ig{A_{m+2}^{p_{m+2}-1}}{x_1z_m^{-1}}{x_1y_1}\ig{A_{m}^{-p_{m}}}{x_1z_m^{-1}}{c_mz_m}\
x_1\right )^{p_{m+3}}
 x_1^{p_{m+1}} y_1  x_2^{-1}$

 $ (n \geq 2),$

$ SubC_3(A_{m+4}) = SubC_3(A_{m+2}) \cup \{x_1y_1x_2^{-1}, \
 y_1x_2^{-1}x_1, \ x_2^{-1}x_1z_m^{-1} \}\
 (n
\geq 2). $
 \item [(3)]Below we list reduced forms of $x_1^{\phi _j}, y_1^{\phi _j}$ for
$j=m,\dots ,m+4$ and their expressions via the leading terms:

$x_1^{\phi _{ m}}=  A_{m}^{-p_{m}}\circ x_1 \circ A _{m}^{p_{m}}$,
\medskip

$y_1^{\phi _{m}}=  A_{m}^{-p_{m}} \circ y_1$,
\medskip
$x_1^{\phi _{m+1}}=x_1^{\phi _{m}}$,

\medskip
$y_1^{\phi _{ m+1}}= A _{m}^{-p_{m}} \circ x_1^{p_{m+1}} \circ
y_1$,

\medskip
$x_1^{\phi _{ m+2}}=_{(when \ n = 1)}{x_1^{\phi_K}}
=\ig{A_{m+2}^{p_{m+2}}}{x_1z_m^{-1}}{x_1y_1}\ig{A_{m}^{-p_{m}}}{x_1z_m^{-1}}{c_mz_m}
\ x_1 \ \ig{A_{m}^{p_{m}}}{z_m^{-1}c_m^{-1}}{z_mx_1^{-1}}\ $,

\medskip
$Sub_3(x_1^{\phi_K})_{(when \ n = 1)} = SubC_3(A_{m+2}) \cup
SubC_3(A_{m}) \cup \{z_mx_1z_m^{-1}, \ x_1z_m^{-1}c_m^{-1} \};$

\medskip
$y_1^{\phi _{m+2}}=y_1^{\phi _{m+1}}$,
\medskip
$x_1^{\phi _{m+3}}=x_1^{\phi _{m+2}}$,
\medskip
\begin{multline*} y_1^{\phi _{ m+3}}=_{(when \ n
= 1)} {y_1^{\phi_K}} =\ig{A_{m}^{-p_{m}}}{x_1z_m^{-1}}{c_mz_m}  \\
\left (x_1^{p_{m+1}}y_1
\ig{A_{m+2}^{p_{m+2}-1}}{x_1z_m^{-1}}{x_1y_1}
\ig{A_{m})^{-p_{m}}}{x_1z_m^{-1}}{c_mz_m} \ x_1\right )^{p_{m+3}}
x_1^{p_{m+1}} y_1.\end{multline*}

\medskip
\vspace{.2cm} $Sub_3(y_1^{\phi_K})=_{(when \ n =
1)}Sub_3(y_1^{\phi_{m+3}})= SubC_3(A_{m+2});$

\medskip
\begin{multline*} x_1^{\phi_{m+4}} = {x_1^{\phi_K}}_{(when \ n \geq
2)}= \ig{A_{m+4}^{-p_{m+4}+1}}{x_2y_1^{-1}}{z_mx_1^{-1}} \
x_2y_1^{-1}x_1^{-p_{m+1}} \circ \\ \left( x_1^{-1} \
\ig{A_{m}^{p_{m}}}{z_m^{-1}c_m^{-1}}{z_mx_1^{-1}}
\ig{A_{m+2}^{-p_{m+2}}}{y_1^{-1}x_1^{-1}}{z_mx_1^{-1}}\ y_1^{-1}
x_1^{-p_{m+1}}\right )^{p_{m+3}-1}
\ig{A_{m}^{p_{m}}}{z_m^{-1}c_m^{-1}}{z_mx_1^{-1}} \ \ \ (n \geq
2),\end{multline*}

\medskip\vspace{.2cm}
$Sub_3(x_1^{\phi_K}) = SubC_3(A_{m+2})^{-1} \cup \{z_mx_1^{-1}x_2,
\ x_1^{-1}x_2y_1^{-1}, \ x_2y_1^{-1}x_1^{-1} \} \ (n\geq 2);$

\vspace{.2cm}
 $y_1^{\phi _{m+4}} = {y_1^{\phi_K}}_{(when \ n \geq
2)} =\ig{A_{m+4}^{-(p_{m+4}-1)}}{x_2y_1^{-1}}{z_mx_1^{-1}} \ x_2 \
\ig{A_{m+4}^{p_{m+4}}}{x_1z_m^{-1}}{y_1x_2^{-1}} \ \ \ (n\geq 2)$,

\medskip
$Sub_3(y_1^{\phi_K}) = SubC_3(A_{m+4})^{\pm 1} \cup
\{z_mx_1^{-1}x_2, \ x_1^{-1}x_2x_1, x_2x_1z_m^{-1}\}$ \ \  $(n
\geq 2).$
\end{enumerate}
 \end{lm}

{\it Proof.}
 Statement (1) follows immediately from definitions of
automorphisms of $\Gamma$.

\medskip \noindent
We prove formulas in the second and third statements
simultaneously using Lemma \ref{le:7.1.zforms}:

\medskip
$x_1^{\phi _{m}}= \left (x_1^{(c_m^{{z}_m}x_1^{-1})^{p_{m}}}
\right )^{\phi_{m-1}}= x_1^{A(\phi_{m})^{p_{m}}} =  A
_{m}^{-p_{m}}\circ x_1 \circ A _{m})^{p_{m}}$,

\medskip\noindent
beginning with $x_1$ and ending with $x_1^{-1}$.

\medskip
$y_1^{\phi _{m}}= \left ((c_m^{{ z}_m}x_1^{-1})^{-p_{m}}y_1\right
)^{\phi_{m-1}}= A(\phi_{m})^{-p_{m}} \circ y_1$,
\medskip\noindent
beginning with $x_1$ and ending with $y_1$. Now
 $A_{m+1}=$  cycred $(
  A(\gamma_{m+1})^{\phi_{m}} = x_1^{\phi_{m}} = A
_{m}^{-p_{m}}\circ x_1 \circ A _{m}^{p_{m}}.$

\medskip\noindent
$A_{m+1}=x_1.$

\medskip
$x_1^{\phi _{m+1}}=x_1^{\phi _{m}}$,

\medskip
$y_1^{\phi _{ m+1}}= \left (y_1^{\gamma_{m+1}^{p_{m+1}}}\right
)^{\phi_{m}} = (x_1^{p_{m+1}}y_1)^{\phi_{m}}=
(x_1^{\phi_{m}})^{p_{m+1}}y_1^{\phi_{m}} = A_{m}^{-p_{m}} \circ
x_1^{p_{m+1}} \circ y_1$,

\medskip\noindent
beginning with $x_1$ and ending with $y_1$, moreover, the element
that cancels in reducing

$A_{m+1}^{p_{m+1}} A_{m}^{-p_{m}} y_1$ is equal to
$A_{m}^{p_{m}}.$

\medskip
$A_{m+2} =$ cycred $(A(\gamma_{m+2})^{\phi_{m+1}}) =$ cycred
$(y_1^{\phi_{m+1}})= A_{m}^{-p_{m}} \circ x_1^{p_{m+1}} \circ
y_1$,

\medskip\noindent
 beginning with $x_1$ and ending with $y_1$.

\medskip
$x_1^{\phi _{ m+2}}= \left (x_1^{\gamma_{m+2}^{p_{m+2}}}\right
)^{\phi_{m+1}} =  (y_1^{\phi _{ m+1}})^{p_{m+2}} x_1^{\phi _{
m+1}}$

\medskip
$= A_{m+2}^{p_{m+2}} \circ  A_{m}^{-p_{m}}\circ x_1 \circ A
_{m}^{p_{m}} =$

\medskip
$A_{m}^{-p_{m}} \circ \left (x_1^{p_{m+1}} \circ y_1 \circ A
_{m+2}^{p_{m+2}-1}  \circ  A_{m}^{-p_{m}}\circ x_1\right ) \circ
A_{m}^{p_{m}}$,

\medskip \noindent
beginning with $x_1$ and ending with $x_1^{-1}$;

\medskip
$y_1^{\phi _{m+2}}=y_1^{\phi _{m+1}}.$

 \medskip
 $A_{m+3} =$ cycred $ (y_1^{-\phi_{m+1}} x_1^{\phi_{m+2}} y_1^{\phi_{m+1}})  =A_{m+2}^{p_{m+2}-1}\circ A_{m}^{-p_{m}}\circ
 x_1^{p_{m+1}+1}
  \circ y_1$,

\medskip\noindent
beginning with $x_1$ and ending with $y_1$;

\medskip
$x_1^{\phi _{m+3}}=x_1^{\phi _{m+2}}$,

\medskip
$y_1^{\phi _{ m+3}}= (x_1^{\phi _{ m+2}})^{p_{m+3}} y_1^{\phi _{
m+1}} = $

\medskip
$A_{m}^{-p_{m}} \circ \left (x_1^{p_{m+1}} \circ y_1 \circ A
_{m+2}^{p_{m+2}-1}  \circ  A_{m}^{-p_{m}}\circ x_1\right
)^{p_{m+3}} \circ  x_1^{p_{m+1}} \circ y_1$,

\medskip\noindent
 beginning with $x_1$ and ending with $y_1$. Finally, for $n\geq
 2$,

\medskip
$A_{m+4} =$cycred $ (A(\gamma_{m+4})^{\phi_{m+3}})=$ cycred
$((y_1x_2^{-1})^{\phi_{m+3}}) = y_1^{\phi_{m+3}}x_2^{-1}=
y_1^{\phi_{m+3}}\circ x_2^{-1}$,

\medskip\noindent
beginning with $x_1$ and ending with $x_2^{-1}$;

\medskip
$x_1^{\phi_{m+4}} = \left ((y_1x_2^{-1})^{-p_{m+4}}x_1\right
)^{\phi_{m+3}} = \left ((x_2 y_1^{-\phi_{m+3}}\right
)^{p_{m+4}}x_1^{\phi_{m+3}} = $

\medskip
\begin{multline*} \left ((x_2
y_1^{-\phi_{m+1}}(x_1^{\phi_{m+2}})^{-p_{m+3}}\right
)^{p_{m+4}}x_1^{\phi_{m+2}} = (x_2y_1^{-\phi
_{m+3}})^{p_{m+4}-1}\circ x_2\circ y_1^{-1} \circ x_1^{-p_{m+1}}
\circ \\ \left( x_1^{-1} \circ A_{m}^{p_{m}} \circ
A_{m+2}^{-p_{m+2}} \circ y_1^{-1} \circ x_1^{-p_{m+1}}\right
)^{p_{m+3}-1} \circ A_{m}^{p_{m}},\end{multline*}

\medskip\noindent
beginning with $x_2$ and ending with $x_1^{-1}$, moreover,  the
element that is cancelled out is $x_1^{\phi _{m+2}}$. Similarly,

\medskip
$y_1^{\phi _{m+4}}=(x_2y_1^{-\phi _{m+3}})^{p_{m+4}}y_1^{\phi
_{m+3}}(y_1^{\phi _{m+3}}x_2^{-1})^{p_{m+4}} =$

\medskip
$(x_2y_1^{-\phi _{m+3}})^{p_{m+4}-1} \circ x_2 \circ (y_1^{\phi
_{m+3}}x_2^{-1})^{p_{m+4}}=A_{m+4}^{-(p_{m+4}-1)} \circ x_2 \circ
A_{m+4}^{p_{m+4}}$,
\medskip
\noindent  beginning with $x_2$ and ending with $x_2^{-1}$,
moreover, the element that is cancelled out is $y_1^{\phi
_{m+3}}$.

This proves the lemma. \hfill $\Box$

In the following lemmas  we describe the  reduced expressions of
the elements $x_i^{\phi_j}$ and $y_i^{\phi_j}$ for $i\geq 2$.
\begin{lm}
\label{le:7.1.xiforms}
 Let $n \geq 2$, $K = K(m,n)$, $p = (p_1, \ldots,p_K)$ be a 3-large tuple,  and
 $$\phi_K = \gamma_K^{p_K} \ldots \gamma_1^{p_1}.$$
 Then for any $i , n \geq i \geq 2$,  the following holds:

\begin{enumerate}
\item  [(1)] All automorphisms from $\Gamma $, except for
$\gamma_{m+4(i-1)}, \gamma_{m+4i-2},
 \gamma_{m+4i}$ fix $x_i$, and all automorphisms from $\Gamma $, except for $\gamma_{m+4(i-1)},
 \gamma_{m+4i-3},
\gamma_{m+4i-1}, \gamma_{m+4i}$ fix $y_i$. It follows that

\medskip
$x_i^{\phi_K} = x_i^{\phi_{K-1}} = \ldots = x_i^{\phi_{m+4i}}, $

\medskip
$y_i^{\phi_K} = y_i^{\phi_{K-1}} = \ldots = y_i^{\phi_{m+4i}}. $

\item  [(2)] Let $\tilde  y_{i}=y_{i}^{\phi _{ m+4i-1}}$. Then
${\tilde y}_{i} = \ig{{\tilde y}_{i}}{x_iy_{i-1}^{-1}}{x_iy_i}\ $
where (for $i = 1$) we assume that  $y_0 = x_1^{-1}$ for $m = 0$,
and $y_0 = z_m$ for $m \neq 0$;

 \item [(3)] Below we list the reduced forms of the leading terms
of the corresponding automorphisms.   Put $q_j = p_{m+4(i-1)+j}$
for $j = 0, \ldots, 4$. In the formulas below we assume that $y_0
= x_1^{-1}$ for $m = 0$, and $y_0 = z_m$ for $m \neq 0$.

\medskip
$A_{m+4i-4}  = \ig{{\tilde
y}_{i-1}}{x_{i-1}y_{i-2}^{-1}}{x_{i-1}y_{i-1}
}\circ
x_i^{-1} \ $,

\medskip $SubC_3(A_{m+4i-4}) = Sub_3({\tilde y}_{i-1}) \cup
\{x_{i-1}y_{i-1}x_i^{-1}, \ y_{i-1}x_i^{-1}x_{i-1}, \
x_i^{-1}x_{i-1}y_{i-2}^{-1} \} ;$

\medskip
$A_{m+4i-3} =  x_i; $

\medskip
$A_{m+4i-2}=
\ig{A_{m+4i-4}^{-q_0}}{x_iy_{i-1}^{-1}}{y_{i-2}x_{i-1}^{-1}}
 \  x_i^{q_1} y_i$,
\medskip
$$SubC_3(A_{m+4i-2}) = SubC_3(A_{m+4i-4}) \cup
\{y_{i-2}x_{i-1}^{-1}x_i, \ x_{i-1}^{-1}x_i^2, \ x_i^2y_i, \
x_iy_ix_i, \ y_ix_iy_{i-1}^{-1}, x_i^3 \} ;$$

\medskip
$A_{m+4i-1}= \ig{A_{m+4i-2}^{q_2-1}}{x_iy_{i-1}^{-1}}{x_iy_i}
\ig{A_{m+4i-4}^{-q_0}}{x_iy_{i-1}^{-1}}{y_{i-2}x_{i-1}^{-1}} \
 x_i^{q_1+1}y_i, $

\medskip
$SubC_3(A_{m+4i-1}) = SubC_3(A_{m+4i-2});$

\item [(4)] Below we list the reduced forms of elements
$x_i^{\phi_{m+4(i-1)+j}}, y_i^{\phi_{m+4(i-1)+j}}$
 for $j = 0, \ldots, 4.$ Again, in the formulas below we assume that $y_0
= x_1^{-1}$ for $m = 0$, and $y_0 = z_m$ for $m \neq 0$.

\medskip
$x_i^{\phi _{ m+4i-4}} = A_{m+4i-4}^{-q_0}\circ x_i\circ A
_{m+4i-4}^{q_0}$,

\medskip
 $y_i^{\phi _{ m+4i-4}} = A_{m+4i-4}^{-q_0}\circ y_i$,

\medskip
$x_i^{\phi _{m+4i-3}}=x_i^{\phi_{m+4i-4}}$,

\medskip
$y_i^{\phi _{ m+4i-3}}= A_{m+4i-4}^{-q_0}\circ x_i^{q_1} \circ
y_i$,

\medskip
$x_i^{\phi _{
m+4i-2}}=\ig{A_{m+4i-2}^{q_2}}{x_iy_{i-1}^{-1}}{x_iy_i}
\ig{A_{m+4i-4}^{-q_0}}{x_iy_{i-1}^{-1}}{y_{i-2}x_{i-1}^{-1}}\ x_i\
\ig{A_{m+4i-4}^{q_0}}{x_{i-1}y_{i-2}^{-1}}{y_{i-1}x_i^{-1}}\ $,

\medskip
$y_i^{\phi_{m+4i-2}} = y_i^{\phi_{m+4i-3}}, $

\medskip
$x_i^{\phi_{m+4i-1}} = x_i^{\phi_{m+4i-2}} =_{(when \ i = n)}
x_i^{\phi_K}$,

\medskip
$Sub_3(x_i^{\phi_K}) =_{(when \ i = n)} \ SubC_3(A_{m+4i-2})\ \cup
\ SubC_3(A_{m+4i -4})^{\pm 1}\  \cup\ $ $$
\{y_{i-2}x_{i-1}^{-1}x_i, \ x_{i-1}^{-1}x_ix_{i-1}, \
x_ix_{i-1}y_{i-2}^{-1} \} ;$$

\medskip
$y_i^{\phi_{m+4i-1}}= {\tilde y}_i =_{(when \ i = n)} y_i^{\phi_K}
=$\newline $
\ig{A_{m+4i-4}^{-q_0}}{x_iy_{i-1}^{-1}}{y_{i-2}x_{i-1}^{-1}} \left
(x_i^{q_1} y_i \ \ig{A_{m+4i-2}^{q_2-1}}{x_iy_{i-1}^{-1}}{x_iy_i}
\ig{A_{m+4i-4}^{-q_0}}{x_iy_{i-1}^{-1}}{y_{i-2}x_{i-1}^{-1}} \ x_i
\right )^{q_3} \ x_i^{q_1} y_i $,

\medskip\begin{multline*}
Sub_3({\tilde y}_i) = SubC_3(A_{m+4i-2}) \cup
SubC_3(A_{m+4i-4})^{- 1} \cup \{y_{i-2}x_{i-1}^{-1}x_i, \
x_{i-1}^{-1}x_i^2,\\ \ x_i^3, \ x_iy_ix_i, \ y_ix_iy_{i-1}^{-1},
\ x_i^2y_i \}\end{multline*}

\medskip \begin{multline*}
x_i^{\phi_{m+4i}} =_{(when \ i \neq n)} x_i^{\phi_K}=
\ig{A_{m+4i}^{-q_4+1}}{x_{i+1}y_i^{-1}}{y_{i-1}x_i^{-1}}  \
x_{i+1} \circ y_i^{-1}x_i^{-q_1} \circ \\
\circ  \left ( x_i^{-1} \
\ig{A_{m+4i-4}^{q_0}}{x_{i-1}y_{i-2}^{-1}}{y_{i-1}x_i^{-1}}
\ig{A_{m+4i-2}^{-q_2+1}}{y_{i}^{-1}x_i^{-1}}{y_{i-1}x_i^{-1}} \
y_i^{-1}x_i^{-q_1} \right )^{q_3-1}\
\ig{A_{m+4i-4}^{q_0}}{x_{i-1}y_{i-2}^{-1}}{y_{i-1}x_i^{-1}} \
,\end{multline*}

\medskip\begin{multline*}
Sub_3(x_i^{\phi_K}) = SubC_3(A_{m+4i})^{-1} \cup
SubC_3(A_{m+4i-2})^{-1} \cup SubC_3(A_{m+4i-4}) \\
\cup \{y_{i-1}x_i^{-1}x_{i+1}, \ x_i^{-1}x_{i+1}y_i^{-1}, \
x_{i+1}y_i^{-1}x_i^{-1}, \ y_i^{-1}x_i^{-2}, \ x_i^{-3}, \
x_i^{-2}x_{i-1},\\ \ x_i^{-1}x_{i-1}y_{i-2}^{-1}, \
y_{i-1}x_i^{-1}y_i^{-1}, \ x_i^{-1}y_i^{-1}x_i^{-1}
\};\end{multline*}

\medskip
\begin{multline*} y_i^{\phi_{m+4i}} = _{(when \ i\neq n)}y_i^{\phi_K} =\ig{A_{m+4i}^{-q_4+1}}{x_{i+1}y_i^{-1}}{y_{i-1}x_i^{-1}}  \
x_{i+1} \ \ig{{\tilde y}_i}{x_iy_{i-1}^{-1}}{x_iy_i} \
x_{i+1}^{-1} \
\ig{A_{m+4i}^{q_4-1}}{x_iy_{i-1}^{-1}}{y_ix_{i+1}^{-1}}\
,\end{multline*}

\medskip
\begin{multline*}Sub_3(y_i^{\phi_K}) = SubC_3(A_{m+4i})^{\pm 1} \cup Sub_3({\tilde
y}_i) \cup \{y_{i-1}x_i^{-1}x_{i+1},\\ \ x_i^{-1}x_{i+1}x_i, \
x_{i+1}x_iy_{i-1}^{-1}, \ x_iy_ix_{i+1}^{-1}, \
y_ix_{i+1}^{-1}x_i, \
x_{i+1}^{-1}x_iy_{i-1}^{-1}\}.\end{multline*}

\item [(5)] $A_j$=
\[ \left\{\begin{array}{ll}
                  A(\gamma_j)^{\phi_{j- 1}}, & \mbox{if $j \neq  m+4i-1,\ m+4i-3$ \ for \ any \
                  $i = 1, \ldots,n;$},\\
                  A_{m+4i-4}^{p_{m+4i-4}}A(\gamma _j)^{\phi
                  _{j-1}}
                  A_{m+4i-4}^{-p_{m+4i-4}},& \mbox{if $j = m+4i-3$ \ for \ some \
                   $i = 1, \ldots,n\ (m+4i-4\neq 0);$}\\
                  y_i^{-\phi _{j-1}} A(\gamma _{j})^{\phi _{j-1}}y_i^
                  {\phi _{j-1}},
                   & \mbox{if $j = m+4i-1$ \ for \ some
                   $i = 1, \ldots,n.$ }
                 \end{array}
                   \right. \]
\end{enumerate}
\end{lm}
{\it Proof.} Statement (1) is obvious. We prove  statement (2) by
induction on $i \geq 2$. Notice that by Lemmas
\ref{le:7.1.x1formsm0} and \ref{le:7.1.x1formsmneq0}
  ${\tilde y}_1 =  y_1^{\phi_{m+3}}$ begins with $x_1z_m^{-1}$  and ends  with $x_1y_1$.
   Now let $i \geq  2$. Then denoting exponents by $q_i$ as in (3), we have

\medskip
${\tilde y}_i  = y_i^{\phi_{m+4i-1}}=
(x_i^{q_3}y_i)^{\phi_{m+4i-2}} = \left ((y_i^{q_2}x_i)^{q_3}y_i
\right )^{\phi_{m+4i-3}} = $

\medskip
$ \left ( \left ((x_i^{q_1}y_i)^{q_2}x_i\right
)^{q_3}x_i^{q_1}y_i\right )^{\phi_{m+4i-4}}.$

\medskip \noindent
 Before we  continue, and to avoid huge formulas,  we compute separately
$x_i^{\phi_{m+4i-4}}$ and $y_i^{\phi_{m+4i-4}}$:

\medskip
$x_i^{\phi_{m+4i-4}} = \left (x_i^{(y_{i-1}x_i^{-1})^{q_0}}\right
)^{\phi_{m+4(i-1)-1}} = x_i^{({\tilde y}_{i-1}x_i^{-1})^{q_0}} = $
$\ig{(x_i {\tilde y}_{i-1}^{-1})^{q_0} \circ x_i \circ ({\tilde
y}_{i-1}x_i^{-1})^{q_0}}{x_iy_{i-1}^{-1}}{y_{i-1}x_i^{-1}}\ $,

\medskip \noindent
  by induction (by Lemmas \ref{le:7.1.x1formsm0} and
\ref{le:7.1.x1formsmneq0}  in the case $i = 2$).

\medskip
$y_i^{\phi_{m+4i-4}} =  \left ((y_{i-1}x_i^{-1})^{-q_0} y_i \right
)^{\phi_{m+4(i-1)-1}} =  ({\tilde y}_{i-1}x_i^{-1})^{-q_0}y_i
=(x_i \circ {\tilde y}_{i-1}^{-1})^{q_0} \circ y_i $,
\medskip \noindent
beginning with $x_iy_{i-1}^{-1}$ and ending with
$x_{i-1}^{-1}y_i$. It follows that

\medskip
$ (x_i^{q_1}y_i )^{\phi_{m+4i-4}} =  (x_i {\tilde
y}_{i-1}^{-1})^{q_0}  x_i^{q_1}   ({\tilde
y}_{i-1}x_i^{-1})^{q_0}(x_i  {\tilde y}_{i-1}^{-1})^{q_0}
 y_i = $
$ (x_i {\tilde y}_{i-1}^{-1})^{q_0}  \circ x_i^{q_1} \circ y_i $,

\medskip \noindent
beginning with $x_iy_{i-1}^{-1}$ and ending with $x_iy_i$.  Now
looking at the formula

\medskip
${\tilde y}_i  = \left ( \left ((x_i^{q_1}y_i)^{q_2}x_i\right
)^{q_3}x_i^{q_1}y_i\right )^{\phi_{m+4i-4}}$

\medskip \noindent
it is obvious  that ${\tilde y}_i $ begins with $x_iy_{i-1}^{-1}$
and ends with $x_iy_i$, as required.

Now we prove statements (3) and (4) simultaneously.

\medskip
$A_{m+4i-4} =$cycred $ ((y_{i-1}x_i^{-1})^{\phi_{m+4(i-1)-1}})
={\tilde y}_{i-1}\circ x_i^{-1}$,
\medskip \noindent
beginning with $x_{i-1}$ and ending with  $x_i^{-1}$. As we have
observed in proving (2)

\medskip
$x_i^{\phi_{m+4i-4}} = (x_i {\tilde y}_{i-1}^{-1})^{q_0} \circ x_i
\circ ({\tilde y}_{i-1}x_i^{-1})^{q_0} = $

\medskip
$A_{m+4i-4}^{-q_0}\circ x_i \circ A_{m+4i-4}^{q_0}, $

\medskip \noindent
beginning with $x_i$  and ending with $x_i^{-1}$.

\medskip
$y_i^{\phi_{m+4i-4}} =   (x_i \circ {\tilde y}_{i-1}^{-1})^{q_0}
\circ y_i  = A_{m+4i-4}^{-q_0} \circ y_i$,
\medskip \noindent
beginning with $x_i$  and ending with $y_i$. Now

\medskip
$A_{m+4i-3} = $ cycred $( x_i^{\phi_{m+4i-4}}) = x_i$,
\medskip \noindent
beginning with $x_i$  and ending with $x_i$.

\medskip
$x_i^{\phi_{m+4i-3}} = x_i^{\phi_{m+4i-4}}$,

\medskip
$y_i^{\phi_{m+4i-3}} =  (x_i^{q_1}y_i)^{\phi_{m+4i-4}} = $
 $A_{m+4i-4}^{-q_0} x_i^{q_1}
A_{m+4i-4}^{q_0} A(\phi_{m+4i-4})^{-q_0}  y_i = $

\medskip
$ A_{m+4i-4}^{-q_0}\circ x_i^{q_1}  \circ y_i$,
\medskip \noindent
 beginning with $x_i$ and ending with $y_i$.

Now

\medskip
$A_{m+4i-2} = y_i^{\phi_{m+4i-3}}$,

\medskip
$x_i^{\phi_{m+4i-2}} =  (y_i^{q_2}x_i)^{\phi_{m+4i-3}} = $
\medskip
$A_{m+4i-2}^{q_2}\circ A_{m+4i-4}^{-q_0}\circ x_i \circ
A_{m+4i-4}^{q_0}, $

\medskip\noindent
beginning with $x_i$ and ending with $x_i^{-1}$. It is also
convenient to rewrite $x_i^{\phi_{m+4i-2}}$ (by rewriting the
subword $A_{m+4i-2}$) to show its cyclically reduced form:

\medskip
$x_i^{\phi_{m+4i-2}} = A_{m+4i-4}^{-q_0}\circ \left ( x_i^{q_1}
\circ y_i  \circ  A_{m+4i-2}^{q_2-1}\circ A_{m+4i-4}^{-q_0}\circ
x_i \right ) \circ A_{m+4i-4}^{q_0}. $

\medskip
$y_i^{\phi_{m+4i-2}} = y_i^{\phi_{m+4i-3}}. $

\medskip \noindent
Now we can write down the next set of formulas:

\medskip
$A_{m+4i-1}=$cycred $(
y_i^{-\phi_{m+4i-3}}x_i^{\phi_{m+4i-2}}y_i^{\phi_{m+4i-3}}) =
$

\medskip
cycred $(A_{m+4i-2}^{-1} A_{m+4i-2}^{q_2} A_{m+4i-4}^{-q_0} x_i
A_{m+4i-4}^{q_0} A_{m+4i-2})=$

\medskip
$ A_{m+4i-2}^{q_2-1} \circ A_{m+4i-4}^{-q_0}  \circ x_i^{q_1+1}
\circ y_i$,
\medskip \noindent
beginning with $x_i$ and ending with $y_i$,

\medskip
$x_i^{\phi_{m+4i-1}} = x_i^{\phi_{m+4i-2}}$,
\medskip
$ y_i^{\phi_{m+4i-1}} = {\tilde y}_i = (x_i^{q_3}y_i)^{
\phi_{m+4i-2}}  = (x_i^{ \phi_{m+4i-2}})^{q_3}y_i^{ \phi_{m+4i-2}}
= $

\medskip \noindent
substituting the cyclic decomposition of $x_i^{ \phi_{m+4i-2}}$
from above one has

\medskip
$= A_{m+4i-4}^{-q_0} \circ \left (x_i^{q_1} \circ y_i \circ
A_{m+4i-2}^{q_2-1} \circ  A_{m+4i-4}^{-q_0} \circ x_i \right
)^{q_3} \circ  x_i^{q_1} \circ y_i .$

\medskip \noindent
beginning with $x_i$ and ending with $y_i$.

\medskip \noindent
Finally

\medskip
$A_{m+4i} =$cycred $ ((y_ix_{i+1}^{-1})^{\phi_{m+4i-1}}) = {\tilde
y}_i \circ x_{i+1}^{-1}$,
\medskip \noindent
beginning with $x_i$ and ending with  $x_{i+1}^{-1}$.

\medskip
$x_i^{\phi_{m+4i}} = \left ((y_ix_{i+1}^{-1})^{-q_4}x_i\right
)^{\phi_{m+4i-1}} = ({\tilde
y}_ix_{i+1}^{-1})^{-q_4}x_i^{\phi_{m+4i-1}} = $

\medskip
$A_{m+4i}^{-q_4+1}x_{i+1}{\tilde y}_i^{-1} x_i^{\phi_{m+4i-1}} = $

\medskip
$A_{m+4i}^{-q_4+1} \circ x_{i+1} \circ  \left (
(x_i^{\phi_{m+4i-2}})^{q_3-1} y_i^{\phi_{m+4i-2}} \right )^{-1} $

\medskip \noindent
observe that computations similar to that for $
y_i^{\phi_{m+4i-1}}$ show that
\begin{multline*}\left ( (x_i^{\phi_{m+4i-2}})^{q_3-1} y_i^{\phi_{m+4i-2}} \right
)^{-1} = \\ \left (A_{m+4i-4}^{-q_0} \circ \left (x_i^{q_1} \circ
y_i \circ A_{m+4i-2}^{q_2-1} \circ A_{m+4i-4}^{-q_0} \circ x_i
\right )^{q_3-1} \circ x_i^{q_1} \circ y_i\right )^{-1}.
\end{multline*}

\noindent Therefore
\begin{multline*}x_i^{\phi_{m+4i}} = A_{m+4i}^{-q_4+1} \circ x_{i+1} \circ
\\
\circ  \left (A_{m+4i-4}^{-q_0} \circ \left (x_i^{q_1} \circ y_i
\circ A_{m+4i-2}^{q_2-1} \circ A_{m+4i-4}^{-q_0} \circ x_i \right
)^{q_3-1} \circ x_i^{q_1} \circ y_i\right )^{-1} ,\end{multline*}

\medskip \noindent
beginning with $x_{i+1}$ and ending with  $x_{i}^{-1}$.

\medskip
$y_i^{\phi_{m+4i}} = \left ( y_i^{(y_ix_{i+1}^{-1})^{q_4}}\right
)^{\phi_{m+4i-1}} = (x_{i+1}{\tilde y}_i^{-1})^{q_4} {\tilde y}_i
({\tilde y}_ix_{i+1}^{-1})^{q_4} =$

\medskip
$A_{m+4i}^{-q_4+1} \circ x_{i+1} \circ {\tilde y}_i \circ
x_{i+1}^{-1} \circ  A_{m+4i}^{q_4-1}$,
\medskip \noindent
beginning with $x_{i+1}$ and ending with  $x_{i+1}^{-1}$.

(5) If $j\neq m+4i-1, m+4i-3$, $A_j$ is either
$(c_{j}^{z_{j}^{\phi j}}c_{j+1}^{z_{j+1}^{\phi j}})$ or $y_i^{\phi
_j}$ ($j=m+4i-3$) or $(y_{i}x_{i+1}^{-1})^{\phi _j}$ ($j=m+4i-1$).
In all these cases $A_j=A(\gamma _j)^{\phi _{j-1}}.$ The formulas
for the other two cases can be found in the proof of Statement
(2). This finishes the proof of the lemma. \hfill $\Box$

\begin{lm}
\label{le:7.1.2words}  Let $m \geq 1$, $K = K(m,n)$, $p = (p_1,
\ldots,p_K)$ be a 3-large tuple,  $\phi_K = \gamma_K^{p_K} \cdots
\gamma_1^{p_1}$,  and $X^{\pm \phi_K}  = \{x^{\phi_K} \mid x \in
X^{\pm 1} \}.$
 Then the following holds:

\medskip
\begin{enumerate}
 \item [(1)] $
 Sub_2(X^{\pm \phi_K})  = \left\{
 \begin{array}{ll}    c_jz_j,\ z_j^{-1}c_j & (1 \leq j \leq m), \\
z_jz_{j+1}^{-1}    & (1 \leq j \leq m-1), \\
z_mx_1^{-1},\ z_mx_1 &  (if \ m \neq 0, n \neq 0),  \\
x_i^2, \ x_iy_i,\ y_ix_i &  (1 \leq i \leq n), \\
x_{i+1}y_{i}^{-1},  \ x_i^{-1}x_{i+1},\ x_{i+1}x_i & (1 \leq i
\leq n-1) \ \end{array} \right\}^{\pm1}$

\medskip \noindent moreover, the word $z_j^{-1}c_j$, as well as $c_jz_j$,
 occurs only as a part of the subword $(z_j^{-1}c_jz_j)^{\pm 1}$
 in $x^{\phi_K}$ $( x \in X^{\pm 1})$;

\item [(2)] $ Sub_3(X^{\pm \phi_K})  =$ \begin{multline*} \left\{
\begin{array}{ll}
z_{j}^{-1}c_{j}z_{j}, & (1 \leq j \leq m), \\
 c_jz_jz_{j+1}^{-1},  \
z_jz_{j+1}^{-1}c_{j+1},  &  (1 \leq j \leq m-1),  \\
z_jz_{j+1}^{-1}c_{j+1}^{-1}, &  (2 \leq j \leq m-1),\\
y_1x_1^2, \ & (m = 0, n = 1), \\
 x_2^{-1}x_1^2, \ x_2x_1^2,  & (m = 0, n \geq 2) \\
 c_m^{-1}z_mx_1, & (m= 1, n\neq 0)\\
   c_mz_mx_1^{-1},  \ z_mx_1^{-1}z_m^{-1}, \
z_mx_1^2, \  z_mx_1^{-1}y_1^{-1},  & (\ m
\neq 0, n \neq 0), \\
\ c_mz_mx_1, & (\ m\geq 2
,n \neq 0), \\
   z_mx_1^{-1}x_2, \ z_mx_1^{-1}x_2^{-1},  &  (m \neq 0, n \geq
2),  \\
  c_1^{-1}z_1z_2^{-1}, & (m \geq 2),  \\
     x_i^3, \ x_i^2y_i,\ \ x_iy_ix_i,  &  (1 \leq i \leq n),  \\
  x_i^{-1}x_{i+1}x_i,  \ y_ix_{i+1}^{-1}x_i, \
x_iy_ix_{i+1}^{-1},  &  (1 \leq i \leq n-1),  \\
 x_{i-1}^{-1}x_i^2,  \ y_ix_{i}y_{i-1}^{-1},  &  (2 \leq i \leq n), \\
 y_{i-2}x_{i-1}^{-1}x_i^{-1}, \ y_{i-2}x_{i-1}^{-1}x_i &  (3 \leq
i \leq n).\   \end{array} \right\}^{\pm 1}.\end{multline*}

\item [(3)] for any 2-letter word $uv \in  Sub_2(X^{\pm \phi_K}) $
  one has
  $$Sub_2(u^{\phi_K}v^{\phi_K}) \subseteq Sub_2(X^{\pm \phi_K})\cup\{c_i^2\},
   \ \ \ Sub_3(u^{\phi_K}v^{\phi_K}) \subseteq Sub_3(X^{\pm
   \phi_K})\cup\{ c_i^2z_i\}.$$\end{enumerate}

\end{lm}
 {\it Proof.} (1) and (2) follow by  straightforward inspection of the reduced
forms of elements $x^{\phi_K}$ in Lemmas \ref{le:7.1.zforms},
\ref{le:7.1.x1formsm0},  \ref{le:7.1.x1formsmneq0}, and
\ref{le:7.1.xiforms}.

To prove (3) it suffices for every word $uv \in Sub_2(X^{\pm
\phi_K}) $ to write down  the product $u^{\phi_K}v^{\phi_K}$
(using formulas from the lemmas mentioned above),  then make all
possible cancellations and check whether   3-subwords  of the
resulting word all lie in $Sub_3(X^{\pm \phi_K})$. Now we do the
checking one by one for all possible 2-words from $Sub_2(X^{\pm
\phi_K}) $.

\begin{enumerate}
\item [1)] For $uv  \in \{  c_jz_j, \   z_j^{-1}c_j\}$ the
checking is obvious and we omit it.
 \item [2)] Let $uv =z_jz_{j+1}^{-1}$. Then there are three cases to consider:
\begin{enumerate}

\item [2.a)]  Let $j \leq m-2$, then

\medskip
$(z_jz_{j+1}^{-1})^{\phi_K}=
\ig{z_j^{\phi_K}}{*}{c_{j+1}z_{j+1}}\ig{z_{j+1}^{-\phi_K}}{z_{j+2}^{-1}c_{j+2}^{-1}}{*}
\ $,
\medskip
in this case there is no cancellation in $u^{\phi_K}v^{\phi_K}$.
All 3-subwords of $u^{\phi_K}$ and $v^{\phi_K}$ are obviously in
$Sub_3(X^{\pm \phi_K})$. So one needs only to check the new
3-subwords which arise "in between"  $u^{\phi_K}$ and $v^{\phi_K}$
(below we will check only subwords of this type). These subwords
are $c_{j+1}z_{j+1}z_{j+2}^{-1}$ and
$z_{j+1}z_{j+2}^{-1}c_{j+2}^{-1}$  which both lie in $Sub_3(X^{\pm
\phi_K})$.

\item [2.b)] Let $j = m-1$ and $n \neq 0$. Then

\medskip
$(z_{m-1}z_{m}^{-1})^{\phi_K}=
\ig{z_{m-1}^{\phi_K}}{*}{c_mz_m}\ig{z_{m}^{-\phi_K}}{x_1z_m^{-1}}
{*} \ $,
\medskip
 again, there is no cancellation in this case and the words "in
between" are $c_mz_mx_1$ and $z_mx_1z_m^{-1}$, which are in
$Sub_3(X^{\pm \phi_K})$.

\item [2.c)] Let $j = m-1$ and $n = 0$. Then ( below we put
$\cdot$ at the place where the corresponding   initial segment of
$u^{\phi_K}$ and the corresponding terminal segment of
$v^{\phi_K}$ meet)

\medskip
$(z_{m-1}z_{m}^{-1})^{\phi_K} = z_{m-1}^{\phi_K} \cdot
z_{m}^{-\phi_K}=
c_{m-1}z_{m-1}A_{m-4}^{p_{m-4}}c_m^{z_m}A_{m-1}^{p_{m-1}-1} \cdot
A_{m-1}^{-p_{m-1}}z_m^{-1} =$

\noindent (cancelling $A_{m-1}^{p_{m-1}-1}$ and substituting for
$A_{m-1}^{-1}$ its expression via the leading terms)

\medskip
$= c_{m-1}z_{m-1}A_{m-4}^{p_{m-4}}c_m^{z_m} \cdot (c_m^{-z_m}
A_{m-4}^{-p_{m-4}} c_{m-1}^{-z_{m-1}}A_{m-4}^{p_{m-4}})z_m^{-1} =$

\medskip
$ z_{m-1}\ \ig{A_{m-4}^{p_{m-4}}}{z_{m-2}^{-1}}{*} \ z_m^{-1}.$

\medskip
Here $z_{m-1}^{\phi _K}$ is completely cancelled.

\item [3.a)] Let $n=1.$ Then $(z_mx_1^{-1})^{\phi _K}=$

\medskip
$c_m z_mA_{m-1}^{p_{m-1}}x_1^{-1}A_{m}^{p_{m}-1}\cdot
A_{m}^{-p_{m}}x_1^{-1}A_{m}^{p_{m}}A_{m+2}^{p_{m+2}}=$

\medskip
$c_mz_mA_{m-1}^{p_{m-1}}x_1^{-1}\cdot
x_1A_{m-1}^{-p_{m-1}}c_m^{-z_m}A_{m-1}^{p_{m-1}}x_1^{-1}A_{m}^{p_{m}}A_{m+2}^{p_{m+2}}=$

\medskip
$\ig
{z_mA_{m-1}^{p_{m-1}}x_1^{-1}A_{m}^{p_{m}}A_{m+2}^{p_{m+2}}}{z_mz_{m-1}^{-1}}{}\
$, and $z_m^{\phi _K}$ is completely cancelled.

\item [3.b)] Let $n>1$. Then \begin{multline*}(z_mx_1^{-1})^{\phi
_K}=c_mz_mA_{m-1}^{p_{m-1}}x_1^{-1}A_{m}^{p_{m}-1}\\
A_{m}^{-p_{m}}(x_1^{-1}A_{m}^{p_{m}}A_{m+2}^{-p_{m+2}+1}y_1^{-1}x_1^{-p_{m+1}})^{-p_{m+3}+1}x_1^{p_{m+1}}y_1x_2^{-1}
A_{m+4}^{p_{m+4}-1}=\\ c_mz_mA_{m-1}^{p_{m-1}}x_1^{-1}
A_{m}^{-1}(x_1^{-1}A_{m}^{p_{m}}A_{m+2}^{-p_{m+2}+1}y_1^{-1}x_1^{-p_{m+1}})^{-p_{m+3}+1}
x_1^{p_{m+1}}y_1x_2^{-1}A_{m+4}^{p_{m+4}-1}=\\
c_mz_mA_{m-1}^{p_{m-1}}x_1^{-1}\cdot
x_1A_{m-1}^{-p_{m-1}}c_m^{-z_m}A_{m-1}^{p_{m-1}}
(x_1^{-1}A_{m}^{p_{m}}A_{m+2}^{-p_{m+2}+1}y_1^{-1}x_1^{-p_{m+1}})^{-p_{m+3}+1}\\
x_1^{p_{m+1}}y_1x_2^{-1}A_{m+4}^{p_{m+4}-1}=\ig{z_mA_{m-1}^{p_{m-1}}}{z_mz_{m-1}^{-1}c_{m-1}^{-1}}{}
\ ,\end{multline*} and $z_m^{\phi _K}$ is completely cancelled.

\item [4.a)]  Let $n=1$. Then $(z_mx_1)^{\phi
_K}=z_mA_{m-1}^{p_{m-1}}x_1^{-1}A_{m}^{p_{m}-1}\cdot
A_{m+2}^{p_{m+2}}A_{m}^{-p_{m}}x_1A_{m}^{p_{m}}= $

\medskip
$ \ig{z_mA_{m-1}^{p_{m-1}}**}{z_mz_{m-1}^{-1}c_{m-1}^{-1}}{*}\ $,
and $z_m^{\phi _K}$ is completely cancelled.

\item [4.b)] Let $n>1.$ Then $(z_mx_1)^{\phi _K}=\ig{z_m^{\phi
_K}}{*}{z_mx_1^{-1}}\ig{x_1^{\phi _K}}{x_2y_1^{-1}}{*}\ .$

\item [5.a)] Let $n=1.$ Then $x_1^{2\phi
_K}=A_{m+2}^{p_{m+2}}A_{m}^{-p_{m}}x_1A_{m}^{p_{m}}\cdot
A_{m+2}^{p_{m+2}}A_{m}^{-p_{m}}x_1A_{m}^{p_{m}}=$

\medskip
$A_{m+2}^{p_{m+2}}A_{m}^{-p_{m}}x_1A_{m}^{p_{m}}\cdot
(A_{m}^{-p_{m}}x_1^{p_{m+1}}y_1)A_{m+2}^{p_{m+2}-1}A_{m}^{-p_{m}}x_1A_{m}^{p_{m}}=$

\medskip
$A_{m+2}^{p_{m+2}}\ig{A_{m}^{-p_{m}}x_1}{*}{z_mx_1}\cdot
x_1^{p_{m+1}}y_1**\ .$

\item [5.b)] Let $n>1$. Then $x_1^{2\phi _K}=\ig{x_1^{\phi
_K}}{}{z_mx_1^{-1}}\ig{x_1^{\phi _K}}{x_2y_1^{-1}}{}\ .$

\item [6.a)] Let $1<i<n.$

Then $x_i^{2\phi
_K}=A_{m+4i}^{-q_4+1}x_{i+1}y_i^{-1}x_i^{-q_1}(x^{-1}A_{m+4i-4}^{q_0}A_{m+4i-2}^{-q_2+1}y_i^{-1}x_i^{-q_1})^{q_3-1}$

\medskip
$\ig{A_{m+4i-4}^{q_0}}{} {y_{i-1}x_i^{-1}} \cdot
\ig{A_{m+4i}^{-q_4+1}}{x_{i+1}y_i^{-1}}{}**.$

\item [6.b)] $x_{n}^{2\phi
_K}=A_{m+4n-2}^{q_2}A_{m+4n-4}^{-q_0}x_nA_{m+4n-4}^{q_0}\cdot
A_{m+4n-2}^{q_2}A_{m+4n-4}^{-q_0}x_nA_{m+4n-4}^{q_0}=$

\medskip
\begin{multline*}A_{m+4n-2}^{q_2}A_{m+4n-4}^{-q_0}x_nA_{m+4n-4}^{q_0}\cdot
A_{m+4n-4}^{-q_0}x_n^{q_1}y_n
A_{m+4n-2}^{q_2-1}A_{m+4n-4}^{-q_0}x_nA_{m+4n-4}^{q_0}= \\
A_{m+4n-2}^{q_2}\ig{A_{m+4n-4}^{-q_0}x_n}{}{x_{n-1}x_n}\cdot
x_n^{q_1}**.\end{multline*}

\item [7.a)] If $n=1$. Then $(x_1y_1)^{\phi
_K}=A_{m+2}^{p_{m+2}}A_{m}^{-p_{m}}x_1\cdot x_1^{p_{m+1}}**.$

\medskip
\item  [7.b)] If $n>1$. Then $(x_1y_1)^{\phi _K}=\ig{x_1^{\phi
_K}}{}{z_{m}x_1^{-1}}\ig{y_1^{\phi _K}}{x_2y_1^{-1}}{}\ .$

\medskip
\item [7.c)] If $1<i<n.$ Then $(x_iy_i)^{\phi _K}=\ig{x_i^{\phi
_K}}{}{y_{i-1}x_i^{-1}}\ig{y_i^{\phi _K}}{x_{i+1}y_i^{-1}}{}\ .$

\medskip
\item [7.d)] $(x_ny_n)^{\phi _K}=\ig{x_n^{\phi
_K}}{}{x_{n-1}^{-1}x_n}\ig{y_n^{\phi _K}}{x_{n}^2}{}\ .$

\item [8a)] If $n=1$. Then $(y_1x_1)^{\phi _K}=\ig{y_1^{\phi
_K}}{}{x_1y_1}\ig{x_1^{\phi _K}}{x_1z_m^{-1}}{}\ .$

\item [8.b)] If $n>1$. Then $(y_1x_1)^{\phi
_K}=A_{m+4}^{-p_{m+4}+1}x_2A_{m+4}^{p_{m+4}}\cdot
A_{m+4}^{-p_{m+4}+1}x_2y_1^{-1}x_1^{-p_{m+1}}\circ **=$

\medskip
$A_{m+4}^{-p_{m+4}+1}x_2A_{m+4}\cdot
x_2y_1^{-1}x_1^{-p_{m+1}}\circ **=$

\medskip
$A_{m+4}^{-p_{m+4}+1}x_2A_{m}^{-p_{m}}
(x_1^{p_{m+1}}y_1A_{m+2}^{p_{m+2}-1}A_{m}^{-p_{m}}x_1)^{p_{m+3}}
x_1^{p_{m+1}}y_1x_2^{-1}\cdot$

\medskip
$x_2y_1^{-1}x_1^{-p_{m+1}}( )^{p_{m+3}-1}A_{m}^{p_{m}}=$

\medskip
$A_{m+4}^{-p_{m+4}+1}x_2A_{m}^{-p_{m}}
(x_1^{p_{m+1}}y_1A_{m+2}^{p_{m+2}-1}\ig{A_{m}^{-p_{m}}x_1}{}
{z_mx_1}\ig{A_{m}}{z_m^{-1}c_m^{-1}}{}\ .$

\item[8.c)] $(y_nx_n)^{\phi _K}=\ig{y_n^{\phi
_K}}{}{x_ny_n}\ig{x_n^{\phi _K}}{x_ny_{n-1}}{}\ .$

\item [9.a)] If $n=2$, then
 $(x_2y_1^{-1})^{\phi _K}=A_{m+6}^{q_2}A_{m+4}^{-1}.$
\item [9.b)] If $n>2$,
 $1<i<n$. Then
$(x_iy_{i-1}^{-1})^{\phi
_K}=\ig{A_{m+4i}^{-q_4+1}}{x_{i+1}y_i^{-1}}{y_{i-1}x_i^{-1}}\
x_{i+1} \circ y_i^{-1}x_i^{-q_1} \circ$

\medskip
$\circ  \left ( x_i^{-1} \
\ig{A_{m+4i-4}^{q_0}}{x_{i-1}y_{i-2}^{-1}}{y_{i-1}x_i^{-1}}
\ig{A_{m+4i-2}^{-q_2+1}}{x_{i-1}y_{i-2}^{-1}}{y_{i-1}x_i^{-1}}\
y_i^{-1}x_i^{-q_1} \right )^{q_3-1}$

\medskip
$\ig{A_{m+4i-4}^{q_0}}{x_{i-1}y_{i-2}^{-1}}{y_{i-1}x_i^{-1}} \cdot
A_{m+4i-4}^{-q_0+1}\circ x_i\circ \tilde y_{i-1}\circ x_i^{-1}
\ig{A_{m+4i-4}^{q_0-1}}{x_{i-1}y_{i-2}^{-1}}{y_{i-1}x_i^{-1}}=$

\medskip
$\ig{A_{m+4i}^{-q_4+1}}{x_{i+1}y_i^{-1}}{y_{i-1}x_i^{-1}} \
x_{i+1} \circ y_i^{-1}x_i^{-q_1} \circ$

\medskip
$  \left ( x_i^{-1} \
\ig{A_{m+4i-4}^{q_0}}{x_{i-1}y_{i-2}^{-1}}{y_{i-1}x_i^{-1}}
\ig{A_{m+4i-2}^{-q_2+1}}{x_{i-1}y_{i-2}^{-1}}{y_{i-1}x_i^{-1}} \
y_i^{-1}x_i^{-q_1} \right )^{q_3-1}\cdot$

\medskip$
x_i^{-1}\ig{A_{m+4i-4}^{q_0-1}}{x_{i-1}y_{i-2}^{-1}}{y_{i-1}x_i^{-1}}\
.$ \item [9.c)] $(x_ny_{n-1}^{-1})^{\phi
_K}=\ig{A_{m+4n-2}^{q_2}}{}{x_ny_n}\ig{A_{m+4n-4}}{x_ny_{n-1}^{-1}}{}\
\ .$

\item [10.a)] Let $n=2$, then $ (x_1^{-1}x_2)^{\phi _K}=$

\noindent
 $A_{m}^{-p_{m}}(x_1^{p_{m+1}}y_1A_{m+2}^{p_{m+2}-1}
A_{m}^{-p_{m}}x_1)^{p_{m+3}-1}x_1^{p_{m+1}}y_1x_2^{-1}A_{m+4}^{p_{m+4}-1}
A_{m+6}^{p_{m+6}}A_{m+4}^{-p_{m+4}}x_2A_{m+4}^{p_{m+4}}=$

\noindent $A_{m}^{-p_{m}}(x_1^{p_{m+1}}y_1A_{m+2}^{p_{m+2}-1}
A_{m}^{-p_{m}}x_1)^{p_{m+3}-1}x_1^{p_{m+1}}y_1x_2^{-1}A_{m+4}^{p_{m+4}-1}
$

\noindent
$(A_{m+4}^{-p_{m+4}}x_2^{p_{m+5}}y_2)^{p_{m+6}}A_{m+4}^{-p_{m+4}}x_2A_{m+4}^{p_{m+4}}=$

\noindent $A_{m}^{-p_{m}}(x_1^{p_{m+1}}y_1A_{m+2}^{p_{m+2}-1}
A_{m}^{-p_{m}}x_1)^{p_{m+3}-1}x_1^{p_{m+1}}y_1x_2^{-1}\cdot $

\noindent $ A_{m+4}^{-1}x_2^{p_m+5}y_2(A_{m+4}^{-p_{m+4}}
x_2^{p_{m+5}}y_2)^{p_{m+6}-1}A_{m+4}^{-p_{m+4}}x_2A_{m+4}^{p_{m+4}}=\ig{A_{m}^{-p_{m}}}{}{c_mz_m}\
$

\noindent $
\ig{x_1^{-1}A_{m}^{p_{m}}}{x_1^{-1}z_m^{-1}}{}A_{m+2}^{-p_{m+2}+1}y_1^{-1}x_1^{-p_{m+1}}A_{m}^{p_{m}}x_2^{p_{m+5}}
y_2(A_{m+4}^{-p_{m+4}}x_2^{p_{m+5}}y_2)^{p_{m+6}-1}$

\noindent
 $A_{m+4}^{-p_{m+4}}x_2A_{m+4}^{p_{m+4}}.$

\item [10.b)] If $1<i<n-1$, then $(x_i^{-1}x_{i+1})^{\phi
_K}=\ig{x_i^{-\phi _K}}{}{y_{i}x_{i+1}^{-1}}\ig{x_{i+1}^{\phi
_K}}{x_{i+2}y_{i+1}^{-1}}{}.$

\item [10.c)] Similarly to 10.a) we get $(x_{n-1}^{-1}x_n)^{\phi
_K}=$

\noindent $\ig{A_{2n+4n-8}^{-p_{m+4n-8}}}{}{y_{n-3}x_{n-2}^{-1}}\
\cdot\
\ig{x_{n-1}^{-1}A_{m+4n-8}^{p_{m+4n-8}}}{x_{n-1}^{-1}x_{n-2}}{}A_{m+4n-6}^{p_{m+4n-6}+1}**.$

\item [11.a)] If $1<i<n-1$, then $(x_{i+1}x_i)^{\phi _K}=$

\noindent
$A_{m+4i+4}^{-q_8+1}x_{i+2}y_{i+1}^{-1}x_{i+1}^{-q_5}\left(x_{i+1}^{-1}A_{m+4i}^{q_4}A_{m_4i+2}^{-q_6+1}y_{i+1}^{-1}x_{i+1}^{-q_5}\right
)^{q_7-1}A_{m+4i}^{q_4} $

\noindent
$A_{m+4i}^{-q_4+1}x_{i+1}y_{i}^{-1}x_{i}^{-q_1}\left(x_{i}^{-1}A_{m+4i-4}^{q_0}A_{m_4i-2}^{-q_2+1}y_{i}^{-1}x_{i}^{-q_1}\right
)^{q_3-1}A_{m+4i-4}^{q_0}=$

\noindent
$A_{m+4i+4}^{-q_8+1}x_{i+2}y_{i+1}^{-1}x_{i+1}^{-q_5}\left(x_{i+1}^{-1}A_{m+4i}^{q_4}A_{m_4i+2}^{-q_6+1}y_{i+1}^{-1}x_{i+1}^{-q_5}\right
)^{q_7-1}A_{m+4i}$

\noindent
$x_{i+1}y_{i}^{-1}x_{i}^{-q_1}\left(x_{i}^{-1}A_{m+4i-4}^{q_0}A_{m_4i-2}^{-q_2+1}y_{i}^{-1}x_{i}^{-q_1}\right
)^{q_3-1}A_{m+4i-4}^{q_0}=$

\noindent
$A_{m+4i+4}^{-q_8+1}x_{i+2}y_{i+1}^{-1}x_{i+1}^{-q_5}\left(x_{i+1}^{-1}A_{m+4i}^{q_4}A_{m_4i+2}^{-q_6+1}y_{i+1}^{-1}x_{i+1}^{-q_5}\right
)^{q_7-1}$

\noindent
$A_{m+4i-4}^{-q_0}x_i^{q_1}y_iA_{m+4i-2}^{q_2-1}\ig{A_{m+4i-4}^{-q_0}x_i}{}{x_{i-1}^{-1}x_i}\ig
{A_{m+4i-4}^{q_0}}{x_{i-1}y_{i-2}^{-1}}{}.$

\item [11.b)] If $n>2$, then $(x_2x_1)^{\phi _K}= \ **
\ig{A_{m}^{-q_0}x_1}{}{z_mx_1}\ig{A_{m}^{q_0}}{z_m^{-1}c_m^{-1}}{}.$

\item [11.c)] $(x_nx_{n-1})^{\phi
_K}=A_{m+4n-2}^{q_6}A_{m+4n-4}^{-q_4}x_nA_{m+4n-4}^{q_4}\cdot
A_{m+4n-4}^{-q_4+1}x_ny_{n-1}^{-1}x_{n-1}^{-q_1}$

\noindent
 $(x_{n-1}^{-1}A_{m+4n-8}^{q_0}A_{m+4n-6}^{-q_2+1}y_{n-1}^{-1}x_{n-1}^{-q_1}
 )^{q_3-1}A_{m+4n-8}^{q_0}=$

 \noindent
 $ \ **\ig{A_{m+4n-8}^{-q_0}x_{n-1}}{}{x_{n-2}^{-1}x_{n-1}}\ \cdot
 \ig{A_{m+4n-8}^{q_0}}{x_{n-2}y_{n-3}^{-1}}{}.$

 \item [11.d)] Similarly, if $n=2$, then $(x_2x_1)^{\phi _K}=\
 **\ig{A_{m}^{-p_{m}}x_1}{}{z_mx_1}\ig
 {A_{m}^{p_{m}}}{z_m^{-1}c_m^{-1}}{}.$

\end{enumerate}

\end{enumerate}

This proves the lemma. \hfill $\Box$

{\bf Notation.}
 \ {\it  Denote by $Y$ the following set of words
\begin{itemize}
 \item  [1)] if $n \neq 0$ and for $n=1, m\neq 1$, then
 $$Y = \{x_i,y_i,c_j^{z_j} \mid  i=1,\dots ,n,\ j=1,\dots ,m\}.$$
 \item   [2)] if $n=0$,  denote the element
$c_1^{z_1}\ldots
 c_m^{z_m}\in F(X\cup C_S)$ by a new letter $d$, then
 $$Y = \{c_1^{z_1}, \ldots,   c_{m-1}^{z_{m-1}},\ d\}.$$
 a reduced word in this alphabet is a word that does not contain
 subwords $(c_1^{-z_1}d)^{\pm 1}$ and
$(dc_m^{-z_m})^{\pm 1}$; \item [3)] if $n=1,m=1$, then

$$Y=\{A_1, x_1,y_1\},$$
 a reduced word in this alphabet is a word that does not contain
subwords $(A_1x_1)^{\pm 1}$.

\end{itemize} }

\begin{lm} \label{le:7.1.2words.new}  Let $m \geq 3,\  n=0, $ $K =K(m,0)$. Let  $p = (p_1, \ldots,p_K)$ be a 3-large tuple,  $\phi_K
= \gamma_K^{p_K} \ldots \gamma_1^{p_1}$,  and $X^{\pm \phi_K}=
\{x^{\phi_K} \mid x \in X^{\pm 1} \}.$ Then the following holds:

\bi \item[(1)]  Every element from  $X^{\phi _K}$ can be uniquely
presented as a reduced product of elements and their inverses from
the set
$$X\cup \{c_1,\dots ,c_{m-1}, d\}$$
  Moreover:
 \begin{itemize}
 \item  all elements $z_i^{\phi _K},i\neq m$ have the form $z_i^{\phi _K}=c_iz_i\hat z_i$, where
 $\hat z_i$ is a reduced word in the alphabet
 $Y$,  \item $z_m^{\phi _K}=z_{m}\hat z_m$, where $\hat z_m$ is a
reduced word in the alphabet $Y$.
\end{itemize}

 When
viewing elements from  $X^{\phi _K}$ as elements in $$F(X\cup
\{c_1,\dots ,c_{m-1},d\},$$ the following holds:

\item[(2)] \noindent
 $
 Sub_2(X^{\pm \phi_K})  = \left\{
 \begin{array}{ll}    c_jz_j\ & (1 \leq j \leq m), \\
  z_j^{-1}c_j, \ z_jz_{j+1}^{-1}    & (1 \leq j \leq m-1), \\
z_2d,\  dz_{m-1}^{-1} &
 \end{array} \right\}^{\pm1}$

\medskip \noindent Moreover:
\begin{itemize}
 \item the word $z_mz_{m-1}^{-1}$ occurs  only in the beginning
of $z_m^{\phi _K}$ as a part of the subword \newline$
z_mz_{m-1}^{-1}c_{m-1}^{-1}z_{m-1}$
 \item  the words $z_2d,\
dz_{m-1}^{-1}$ occur  only as  parts of  subwords
$$(c_1^{z_1}c_2^{z_2})^2dz_{m-1}^{-1}c_{m-1}^{-1}z_{m-1}c_{m-1} $$
and $(c_1^{z_1}c_2^{z_2})^2d$. \end{itemize}

\item[(3)] \noindent $ Sub_3(X^{\pm \phi_K})  = \left\{
\begin{array}{ll}

 z_{j}^{-1}c_{j}z_{j},\ c_jz_jz_{j+1}^{-1}, \ z_jz_{j+1}^{-1}c_{j+1}^{-1}, \
  &  (1 \leq j \leq m-1),  \\
z_jz_{j+1}^{-1}c_{j+1}  &  (1 \leq j \leq m-2),\\
c_2z_2d,\ z_2dz_{m-1}^{-1},\  dz_{m-1}^{-1}c_{m-1}^{-1},\ &
\end{array} \right\}^{\pm 1}.$
\ei
\end{lm}

{\it Proof.} The lemma follows from Lemmas \ref {le:7.1.zforms}
and \ref{le:7.1.2words}
 by replacing all the products $c_1^{z_1}\ldots
c_m^{z_m}$
 in subwords of  $X^{\pm \phi_K}$
 by the letter $d$.
\hfill $\Box$

 {\bf Notation.} Let $m \neq 0, K = K(m,n), p=(p_1,\dots
,p_K)$ be a 3-large tuple, and  $\phi_K = \gamma _K^{p_K} \ldots
\gamma _1^{p_1}.$ Let ${\mathcal W}$ be the set of words in $F(X
\cup C_S)$ with the following properties: \begin{enumerate} \item
If $v\in W$ then $ Sub_3(v) \subseteq Sub_3(X^{\pm \phi _K}),
Sub_2(v) \subseteq Sub_2(X^{\pm \phi _K});$ \item Every subword
$x_i^{\pm 2}$ of $v\in W$ is contained in a subword $x_i^{\pm 3};$
\item Every subword $c_1^{\pm z_1}$ of $v\in W$ is contained in
$(c_1^{z_1}c_2^{z_2})^{\pm 3}$ when $m\geq 2$ or in
$(c_1^{z_1}x_1^{-1})^{\pm  3}$ when $m=1$; \item Every subword
$c_m^{\pm z_m}\ (m\geq 3)$ is contained in $(\prod _{i=1}^{m}
c_i^{z_i})^{\pm 1}.$ \item  Every subword  $c_2^{\pm z_2}$ of
$v\in W$ is contained either in $(c_1^{z_1}c_2^{z_2})^{\pm 3}$ or
as the central occurrence of $c_2^{\pm z_2}$ in
$(c_2^{-z_2}c_1^{-z_1})^{ 3}c_2^{\pm z_2}
(c_1^{z_1}c_2^{z_2})^{3}$ or in
$(c_1z_1c_2^{z_2}(c_1^{z_1}c_2^{z_2})^{3})^{\pm 1}$.
\end{enumerate}

\begin{df}
 The following words are called {\em elementary periods}:
$$x_i, \  \ \ c_1^{z_1}c_2^{z_2} \ (if \ m \geq 2), \ \ c_1^{z_1}x_1^{-1} \  (if \
m=1).$$
  We call the squares (cubes) of elementary periods or
 their inverses elementary squares (cubes).
\end{df}

{\bf Notation.}
 \ {\it
 \begin{itemize}
 \item  [1)] Denote by  ${\mathcal W}_{\Gamma}$  the set of
all subwords of words in $\mathcal W$.
 \item   [2)] Denote by  $\bar{\mathcal W}_{\Gamma}$  the set of all  words $v \in {\mathcal
W}_{\Gamma}$ that are freely reduced forms of products of elements
from $Y^{\pm 1}$. In this case we say that these elements $v$ are
(group) words in the alphabet $Y$.
\end{itemize}
}

If $U$ is a set of words in alphabet $Y$  we denote by
$Sub_{n,Y}(U)$ the set of subwords of length $n$ of words from $U$
in alphabet $Y$.

\begin{lm}
\label{le:xyu} Let $v\in {\mathcal W_{\Gamma}}$. Then  the
following holds:
\begin{enumerate}
\item  [(1)] If $v$ begins and ends with an elementary square and
contains no elementary cube, then $v$ belongs to the following
set:

$$
 \left\{
\begin{array}{ll}
x_{i-2}^2y_{i-2}x_{i-1}^{-1}x_ix_{i-1}y_{i-2}^{-1}x_{i-2}^{-2}, \
x_i^2y_ix_iy_{i-1}^{-1}x_{i-1}^{-2}, &   m \geq 2, n\neq 0\\
x_{i-2}^2y_{i-2}x_{i-1}^{-1}x_i^2,\
x_{i-2}^2y_{i-2}x_{i-1}^{-1}x_iy_{i-1}^{-1}x_{i-1}^{-2}, \ & \\
 (c_2^{-z_2} c_1^{-z_1})^2 c_2^{z_2} (c_1^{z_2} c_2^{z_2})^2, & \\
 (c_1^{z_2} c_2^{z_2})^2  c_3^{z_3} \dots c_i^{z_i}
c_{i-1}^{-z_{i-1}} \dots c_3^{-z_3}
        (c_2^{-z_2} c_1^{-z_1})^2 \quad (i\ge3),& \\
 (c_1^{z_2} c_2^{z_2})^2  c_3^{z_3} \dots c_m^{z_m} x_1
c_m^{-z_m} \dots c_3^{-z_3}
        (c_2^{-z_2} c_1^{-z_1})^2, & \\
(c_1^{z_2} c_2^{z_2})^2  c_3^{z_3} \dots c_m^{z_m} x_1^2, & \\
 (c_1^{z_2} c_2^{z_2})^2  c_3^{z_3} \dots c_m^{z_m} x_1^{-1} y_1^{-1} x_1^{-2}, & \\
 (c_1^{z_2} c_2^{z_2})^2  c_3^{z_3} \dots c_m^{z_m} x_1^{-1}
x_2^{-1} x_1
        c_m^{-z_m} \dots c_3^{-z_3} (c_2^{-z_2} c_1^{-z_1})^2, & \\
 (c_1^{z_2} c_2^{z_2})^2  c_3^{z_3} \dots c_m^{z_m} x_1^{-1} x_2^2, & \\
 (c_1^{z_2} c_2^{z_2})^2  c_3^{z_3} \dots c_m^{z_m} x_1^{-1} x_2 y_1^{-1} x_1^{-2}, & \\
     &   \\
(c_1^{z_1}c_2^{z_2})^2dc_{m-1}^{-z_{m-1}}\ldots
(c_2^{-z_2}c_1^{-z_1})^2,\ z_mc_{m-1}^{-z_{m-1}}\ldots
(c_2^{-z_2}c_1^{-z_1})^2,& m\geq 3, n=0,\\
(c_1^{z_2} c_2^{z_2})^2  c_3^{z_3} \dots c_i^{z_i}
c_{i-1}^{-z_{i-1}} \dots c_3^{-z_3}
        (c_2^{-z_2} c_1^{-z_1})^2 \quad (i\ge3),& \\
(c_2^{-z_2} c_1^{-z_1})^2 c_2^{z_2} (c_1^{z_2} c_2^{z_2})^2, & \\
     &   \\
x_1^2y_1(x_1c_1^{-z_1})^2, \
(c_1^{z_1}x_1^{-1})^{2}x_2(x_1c_1^{-z_1})^2,\ & m = 1, n\geq 2, \\
(x_1c_1^{-z_1})^{2}x_1^2,\ x_1^2y_1x_2^{-1}(x_1c_1^{-z_1})^2,
x_2^{-2}(x_1c_1^{-z_1})^2,
& \\

 & \\
x_{i-2}^2y_{i-2}x_{i-1}^{-1}x_ix_{i-1}y_{i-2}^{-1}x_{i-2}^{-2}, \
x_i^2y_ix_iy_{i-1}^{-1}x_{i-1}^{-2},  & m = 0, \ n > 1, \\
x_{i-2}^2y_{i-2}x_{i-1}^{-1}x_iy_{i-1}^{-1}x_{i-1}^{-2}, \
x_1^2y_1x_2^{-1}x_1^2,\ x_1^{-2}x_2^{-1}x_1^2, & \\A_1^2,\
A_1^{-2}x_1^2, \ A_1^{-2}x_1A_1^2,\ x_1^2y_1A_1^{-2}, & m=1,\ n=1.
  \end{array} \right\}^{\pm 1},$$

\item  [(2)] If $v$ does not contain two elementary squares and
begins (ends) with an elementary square, or contains  no
elementary squares, then $v$ is  a subword of either one of the
words above or  of one of the words in $\{x_1^2y_1x_1, \
x_2^2y_2x_2\}$ for $m=0$.
\end{enumerate}
\end{lm}

{\it Proof.}
 Straightforward verification using the description of the set
$Sub_3(X^{\pm \phi_K})$ from Lemma \ref{le:7.1.2words}. \hfill
$\Box$

\begin{df}
Let $Y$ be an alphabet and $E$  a set of words of length at least
2 in $Y$. We say that an occurrence of a word $w \in Y \cup E$ in
a word $v$ is {\em maximal} relative to $E$ if it is not contained
 in any other (distinct from $w$)  occurrence of a word from $E$ in $v$.

We say that a set of words $W$ in the alphabet $Y$ admits {\em
Unique Factorization Property (UF)} with respect to $E$ if  every
word $w \in W$ can be uniquely presented as a product
  $$w = u_1 \ldots u_k$$
  where $u_i$ are maximal occurrences of words from $Y \cup E$. In
  this event the decomposition above is called {\em irreducible}.
\end{df}

\begin{lm}
\label{le:UF} Let  $E$  be a set of words of length $\geq 2$ in an
alphabet $Y$. Suppose that $W$ is a set of words in the alphabet
$Y$ such that
 if $w_1w_2w_3$ is a subword of a word from $W$ and
  $w_1w_2, w_2w_3 \in E$ then $w_1w_2w_3 \in E.$
   Then $W$ admits (UF) with respect to $E$.

\end{lm}
 {\em Proof.} Obvious.

\begin{df}
\label{de:nielsen} Let $Y$ be an alphabet, $E$  a set of words of
length at least 2 in $Y$ and $W$ a set of words in $Y$ which
admits (UF) relative to $E$.  An automorphism $\phi \in Aut F(Y)$
satisfies the Nielsen property with respect to $W$ with exceptions
$E$ if for any word $z \in Y \cup E$ there exists a decomposition
 \begin{equation}
\label{eq:LMR} z^\phi  = L_z \circ M_z \circ R_z,
 \end{equation}
  for some words $L_z, M_z, R_z   \in F(Y)$
  such that for any
  $u_1,u_2  \in Y \cup E$ with $u_1u_2 \in Sub(W) \smallsetminus
  E$ the words $L_{u_1} \circ M_{u_1}$ and $M_{u_2} \circ R_{u_2}$
occur as written  in the reduced form of $u_1^\phi u_2^\phi$.

If an automorphism $\phi$ satisfies the Nielsen property with
respect to $W$ and $E$, then for each word  $z \in Y \cup E$ there
exists a unique decomposition (\ref{eq:LMR}) with maximal length
of $M_z$. In this event we call $M_z=M_{\phi ,z}$  the {\em
middle} of $z^\phi$ with respect to $\phi$.
\end{df}

\begin{lm}
 \label{le: Nielsen-FU}
 Let $W$ be a set of words in the alphabet $Y$ which admits (UF) with respect to a set of words $E$.
 If an automorphism $\phi
\in Aut F(Y)$ satisfies the Nielsen property with respect to $W$
with exceptions $E$ then for every $w \in W$ if $w = u_1 \ldots
u_k$ is the irreducible decomposition of $w$ then the words
$M_{u_i}$ occur as written (uncancelled) in the reduced form of
$w^{\phi}$.
\end{lm}
  {\it Proof} follows directly from definitions.

Set

\medskip
 $ T(m,1) = \{c_s^{z_s} (s=1,\dots ,m),
\prod _{i=1}^{m}c_i^{z_i}x_1\prod _{i=m}^1c_i^{-z_i},\}^{\pm 1}$,
$m\geq 2.$

\medskip
$T(m,2) = T(m,1) \cup \{ \
\prod_{i=1}^{m}c_i^{z_i}x_1^{-1}x_2x_1\prod _{i=m}^{1}c_i^{-z_i},\
y_1x_2^{-1}x_1\prod _{i=m}^{1}c_i^{-z_i}, \ \prod
_{i=1}^{m}c_i^{z_i}x_1^{-1}y_1^{-1}\}^{\pm 1}$,

\medskip\noindent
if $n \geq 3$ then put

\medskip

$T(m,n) = T(m,1) \cup \{ \ \prod
_{i=1}^{m}c_i^{z_i}x_1^{-1}x_2^{-1},\ \prod
_{i=1}^{m}c_i^{z_i}x_1^{-1}y_1^{-1} \}^{\pm 1} \cup T_1(m,n), $

\medskip
where

\medskip
$T_1(m,n) = $

$\{ y_{n-2}x_{n-1}^{-1}x_nx_{n-1}y_{n-2}^{-1},\
y_{r-2}x_{r-1}^{-1}x_r^{-1}, \ y_{r-1}x_{r}^{-1}y_{r}^{-1},\
 y_{n-1}x_n^{-1}x_{n-1}y_{n-2}^{-1} \ \  (n > r \geq 2) \}^{\pm 1}.$

\medskip
Now, let
 $$E(m,n) = \bigcup_{i \geq 2} Sub_i(T(m,n))\cap\bar {\mathcal W}_{\Gamma},\ \ E(m,0)=\emptyset, \ E(1,1)=\emptyset.$$

\begin{lm}
\label{main} Let $m \neq 0, n \neq 0, K = K(m,n), p=(p_1,\dots
,p_K)$ be a 3-large tuple.  Then the following holds:
\begin{enumerate}
 \item [(1)] Let  $w\in X\cup E(m,n)$, $v = v(w)$ be the leading
variable of $w$, and $j = j(v)$ (see notations at the beginning of
Section \ref{se:7.2.5}). Then the period $A_j^{p_j-1}$ occurs in
$w^{\phi_K}$ and  each occurrence of $A_j^2$ in $w^{\phi _K}$ is
contained in some occurrence of $A_j^{p_j-1}.$ Moreover, no square
$A_k^2$ occurs in $w$ for $k > j$.
 \item [(2)]  The automorphism $\phi_K$ satisfies
the Nielsen property with respect to $\bar{\mathcal W}_{\Gamma}$
with exceptions $E(m,n)$. Moreover, the following conditions hold:
  \begin{enumerate}
\item $M_{x_j} = A_{m+4r-8}^{-p_{m+4r-8}+1}x_{r-1}$, for $j \neq
n$.
 \item $M_{x_n} = x_n^{q_1}  \circ y_n  \circ A_{m+4n-2}^{q_2-1}\circ
A_{m+4n-4}^{-q_0}\circ x_n $
 \item $M_{y_j} = y_j^{\phi_K}$, for $j < n$.
 \item $M_{y_n} = \left (x_n^{q_1} y_n \
\ig{A_{m+4n-2}^{q_2-1}}{x_ny_{n-1}^{-1}}{x_ny_n}
\ig{A_{m+4n-4}^{-q_0}}{x_ny_{n-1}^{-1}}{y_{n-2}x_{n-1}^{-1}} \ x_n
\right )^{q_3} \ x_n^{q_1} y_n$.
  \item   $M_w=w^{\phi_K}$ for any
$w \in E(m,n)$ except for the following words:
 \begin{itemize}
 \item  $w_1=y_{r-2}x_{r-1}^{-1}x_r^{-1}, 3\leq r\leq n-1$,
  $w_2=y_{r-1}x_{r}^{-1}y_{r}^{-1}, 2\leq r\leq n-1$,
 \item $w_3=y_{n-2}x_{n-1}^{-1}x_n$,
 $w_4=y_{n-2}x_{n-1}^{-1}x_ny_{n-1}^{-1},$ $w_5=y_{n-2}x_{n-1}^{-1}x_nx_{n-1}^{-1}y_{n-2}^{-1},$
 $w_6=y_{n-2}x_{n-1}^{-1}x_nx_{n-1},$
 $w_7=y_{n-2}x_{n-1}^{-1}x_n^{-1}$, $w_8=y_{n-1}x_n^{-1},$ $w_9=x_{n-1}^{-1}x_n,$
$w_{10}=x_{n-1}^{-1}x_ny_{n-1}^{-1},$
 $w_{11}=x_{n-1}^{-1}x_nx_{n-1}y_{n-2}^{-1}.$
\end{itemize}
  \item  The only letter that may occur in a word from ${\mathcal W}_{\Gamma}$ to the left of a subword
 $w\in\{w_1,\ldots ,w_8\}$ ending with $y_i$ ($i=r-1, r-2,n-1,n-2,\ i\geq 1$) is $x_{i}.$
The maximal number $j$ such that $L_{w}$
 contains $A_j^{p_j-1}$ is $j=m+4i-2$, and $R_{w_1}=R_{w_2}=1$,

 \end{enumerate}

\end{enumerate}

\end{lm}

{\em Proof.}   We first exhibit the formulas for $u^{\phi _K}$,
where $u\in \bigcup_{i \geq 2} Sub_i(T_1(m,n)).$

(1.a) Let $i<n$. Then
$$(x_iy_{i-1}^{-1})^{\phi
_{m+4i}}=(x_iy_{i-1}^{-1})^{\phi
_K}=\ig{A_{m+4i}^{-q_4+1}}{x_{i+1}y_i^{-1}}{y_{i-1}x_i^{-1}} \
x_{i+1} \circ y_i^{-1}x_i^{-q_1} \circ$$

\medskip
$\circ  \left ( x_i^{-1} \
\ig{A_{m+4i-4}^{q_0}}{x_{i-1}y_{i-2}^{-1}}{y_{i-1}x_i^{-1}}
\ig{A_{m+4i-2}^{-q_2+1}}{x_{i-1}y_{i-2}^{-1}}{y_{i-1}x_i^{-1}} \
y_i^{-1}x_i^{-q_1} \right )^{q_3-1}\
\ig{A_{m+4i-4}^{q_0}}{x_{i-1}y_{i-2}^{-1}}{y_{i-1}x_i^{-1}} $

\medskip
$\cdot A_{m+4i-4}^{-q_0+1}\circ x_i\circ \tilde y_{i-1}\circ
x_i^{-1}
\ig{A_{m+4i-4}^{q_0-1}}{x_{i-1}y_{i-2}^{-1}}{y_{i-1}x_i^{-1}}$

\medskip
$=\ig{A_{m+4i}^{-q_4+1}}{x_{i+1}y_i^{-1}}{y_{i-1}x_i^{-1}} \
x_{i+1} \circ y_i^{-1}x_i^{-q_1} \circ$

\medskip
$\circ  \left ( x_i^{-1} \
\ig{A_{m+4i-4}^{q_0}}{x_{i-1}y_{i-2}^{-1}}{y_{i-1}x_i^{-1}}
\ig{A_{m+4i-2}^{-q_2+1}}{x_{i-1}y_{i-2}^{-1}}{y_{i-1}x_i^{-1}} \
y_i^{-1}x_i^{-q_1} \right )^{q_3-1}\cdot
x_i^{-1}\ig{A_{m+4i-4}^{q_0-1}}{x_{i-1}y_{i-2}^{-1}}{y_{i-1}x_i^{-1}}\
.$

\vspace{.2cm} (1.b) Let $i=n$.  Then
 $$(x_ny_{n-1}^{-1})^{\phi
_{m+4n-1}}= (x_ny_{n-1}^{-1})^{\phi
_K}=\ig{A_{m+4n-2}^{q_2}}{x_ny_{n-1}^{-1}}{x_ny_n}
\ig{A_{m+4n-4}^{-1}}{x_ny_{n-1}^{-1}}{y_{n-2}x_{n-1}^{-1}}\ .$$

\noindent Here $y_{n-1}^{-\phi _K}$ is completely cancelled.

\vspace{.2cm} (2.a) Let $i<n-1.$ Then

$(x_{i+1}x_iy_{i-1}^{-1})^{\phi _K}=(x_{i+1}x_iy_{i-1}^{-1})^{\phi
_{m+4i+4}}=A_{m+4i+4}^{-q_8+1}\circ x_{i+2}\circ y_{i+1}^{-1}\circ
x_{i+1}^{-q_5}\circ $

\medskip
$\left (x_{i+1}^{-1}\circ A_{m+4i}^{q_4}\circ
A_{m+4i+2}^{-q_6+1}\circ y_{i+1}^{-1}x_{i+1}^{-q_5}\right
)^{q_7-1}A_{m+4i-4}^{-q_0}\circ x_i^{q_1}y_i\circ
A_{m+4i-2}^{q_2-1}\circ A_{m+4i-4}^{-1}.$

\noindent Here $(x_iy_{i-1}^{-1})^{\phi _{m+4i+4}}$ was completely
cancelled.

\vspace{.2cm} (2.b) Similarly, $(x_iy_{i-1}^{-1})^{\phi
_{m+4i+3}}$ is completely cancelled in
$(x_{i+1}x_iy_{i-1}^{-1})^{\phi _{m+4i+3}}$ and

$(x_{i+1}x_iy_{i-1}^{-1})^{\phi _{m+4i+3}}= A_{m+4i+2}^{q_6}\circ
A_{m+4i}^{-q_4}\circ x_{i+1}\circ
A_{m+4i-4}^{-q_0}A_{m+4i-2}^{q_2-1}\circ A_{m+4i-4}^{-1}.$

\vspace{.2cm} (2.c) $(x_n^{-1}x_{n-1}y_{n-2}^{-1})^{\phi
_{m+4n-1}}=A_{m+4n-4}^{-q_4}\circ x_n^{-1}\circ
A_{m+4n-4}^{q_4}\circ  A_{m+4n-2}^{-q_6+1}\circ y_n^{-1}\circ
x_n^{-q_5}\circ A_{m+4n-8}^{-q_0}\circ x_{n-1}^{q_1}\circ
y_{n-1}\circ A_{m+4n-6}^{q_2-1}\circ A_{m+4n-8}^{-1} $,

\noindent and $(x_{n-1}y_{n-2}^{-1})^{\phi _{m+4n-1}}$ is completely
cancelled.

\vspace{.2cm}
\medskip
(3.a) $(y_ix_iy_{i-1}^{-1})^{\phi _{m+4i}}= A_{m+4i}^{-q_4+1}\circ
x_{i+1}\circ A_{m+4i-4}^{-q_0}\circ x_i^{q_1}\circ y_i\circ
A_{m+4i-2}^{q_2-1}\circ A_{m+4i-4}^{-1}$,

\noindent and $(x_iy_{i-1}^{-1})^{\phi _{m+4i}}$ is completely
cancelled.

\vspace{.2cm}
\medskip
(3.b)$(y_nx_ny_{n-1}^{-1})^{\phi _{K}}=y_n^{\phi _{K}}\circ
(x_ny_{n-1}^{-1})^{\phi _{K}}.$

\vspace{.2cm}
\medskip
(3.c) $(y_{n-1}x_{n}^{-1}x_{n-1}y_{n-2}^{-1})^{\phi
_K}=A_{m+4n-4}\circ A_{m+4n-2}^{-q_6+1}\circ y_n^{-1}\circ
x_n^{-q_5}\circ A_{m+4n-8}^{-q_0}\circ x_{n-1}^{q_1}\circ
y_{n-1}\circ A_{m+4n-6}^{q_2-1}\circ A_{m+4n-8}^{-1} $,

\noindent and $y_{n-1}^{\phi _K}$ and $(x_{n-1}y_{n-2}^{-1})^{\phi
_K}$ are completely cancelled.

(4.a) Let $n\geq 2$.  $(x_1c_m^{-z_m})^{\phi
_{m+4i}}=(x_1c_m^{-z_m})^{\phi _K}=\left
(\ig{A_{m+4}^{-q_4+1}}{x_{2}y_1^{-1}}{c_m^{z_m}x_1^{-1}} \ x_{2}
\circ y_1^{-1}x_1^{-q_1} \circ\right.$

\medskip
$ \circ  \left.\left ( x_1^{-1}\circ A_{m}^{q_0}\circ
A_{m+2}^{-q_2+1}\circ y_1^{-1}\circ x_1^{-q_1}\right
)^{q_3-1}\circ A_{m}^{q_0} \right )\cdot \left (A_{m}^{-q_0}\circ
x_1^{-1}\circ A_{m}^{q_0-1}\right )=$

$A_{m+4}^{-q_4+1}\circ x_2\circ y_1^{-1}\circ x_1^{-q_1}\circ
\left ( x_1^{-1}\circ A_{m}^{q_0}\circ A_{m+2}^{-q_2+1}\circ
y_1^{-1}\circ x_1^{-q_1}\right )^{q_3-1}\circ x_1^{-1}\circ
A_{m}^{q_0-1}.$

Let $n=1$.

$(x_1z_m^{-c_m})^{\phi _K}=A_{m}^{-p_{m}}\circ x_1^{p_{m+1}}\circ
y_1\circ A_{m+2}^{p_{m+2}-1}\circ A_{m}^{-1}$,
$(y_1x_1z_m^{-c_m})^{\phi _K}=y_1^{\phi _K}\circ
(x_1z_m^{-c_m})^{\phi _K}.$

(4.b) $(x_1c_m^{-z_m})^{\phi _K}$ is completely cancelled in
$x_2^{\phi _K}$ and for $n>2$:

 $(x_2x_1c_m^{-z_m})^{\phi _K}=A_{m+8}^{-q_8+1}\circ x_3\circ y_2^{-1}\circ x_3^{-q_5}\circ$

$(x_3^{-1}\circ A_{m+4}^{q_4}\circ A_{m+6}^{-q_6+1}\circ
y_2^{-1}\circ x_3^{-q_5})^{q_7-1}\circ A_{m}^{-q_0}\circ
x_1^{q_1}\circ y_1\circ A_{m+2}^{q_2-1}\circ A_{m}^{-1} $

and for $n=2$:

$(x_2x_1c_m^{-z_m})^{\phi _K}= A_{m+6}^{q_6}\circ
A_{m+4}^{-q_4}\circ x_i\circ A_{m}^{-q_0}\circ x_1^{q_1}\circ
y_1\circ A_{m+2}^{q_2-1}\circ A_{m}^{-1}.$

 (4.c)  The cancellation between
$(x_2x_1c_m^{-z_m})^{\phi _K}$ and $c_{m-1}^{-z_{m-1}}$ is the
same as the cancellation between $A_{m}^{-1}$ and
$c_{m-1}^{-z_{m-1}^{\phi _K}}$, namely,

$A_{m}^{-1}c_{m-1}^{-z_{m-1}^{\phi _K}}= \left ( x_1\circ
A_{m-1}^{-p_{m-1}}\circ c_m^{-z_m}\circ A_{m-1}^{p_{m-1}}\right )$

$\left (A_{m-1}^{-p_{m-1}+1}\circ c_m^{-z_m}\circ
A_{m-4}^{-p_{m-4}}\circ c_{m-1}^{-z_{m-1}}\circ
A_{m-4}^{p_{m-4}}\circ  c_m^{z_m}\circ A_{m-1}^{p_{m-1}-1}\right
)=x_1A_{m-1}^{-1}$,

and $c_{m-1}^{-z_{m-1}^{\phi _K}}$ is completely cancelled.

(4.d)  The cancellations between $(x_2x_1c_m^{-z_m})^{\phi _K}$
(or between $(y_1x_1c_m^{-z_m})^{\phi _K}$) and $\prod
_{i=m-1}^1c_{i}^{-z_{i}^{\phi _K}}$ are the same as the
cancellations between $A_{m}^{-1}$ and $\prod
_{i=m-1}^1c_{i}^{-z_{i}^{\phi _K}}$ namely, the product $\prod
_{i=m-1}^1c_{i}^{-z_{i}^{\phi _K}}$ is completely cancelled and
\[
A_{m}^{-1}\prod _{i=m-1}^1c_{i}^{-z_{i}^{\phi _K}}= x_1\prod
_{i=m}^1c_i^{-z_i}.
\]

 Similarly one can write expressions for $u^{\phi _K}$ for all
$u\in E(m,n).$ The first statement of the lemma  now follows from
these formulas.

 Let us verify the second statement.
Suppose $w\in E(m,n)$  is a maximal subword from $E(m,n)$ of  a
 word $u$ from ${\mathcal W}_{\Gamma}$.
 If $w$ is a subword
of a word in $T(m,n)$, then either $u$ begins with $w$ or $w$ is
the leftmost subword of  a word in $T(m,n).$ All the words in
$T_1(m,n)$ begin with some $y_j$, therefore  the only possible
letters in $u$ in front of $w$ are $x_{j}^2$.

We have $x_{j}^{\phi _K}x_{j}^{\phi _K}w^{\phi _K }=x_{j}^{\phi _K
}\circ x_{j}^{\phi _K}\circ w^{\phi _K}$ if $w$ is a two-letter
word, and $x_{j}^{\phi _K}x_{j}^{\phi _K }w^{\phi _K }=x_{j}^{\phi
_K }\circ x_{j}^{\phi _K} w^{\phi _K}$ if $w$ is more than a
two-letter word. In this last case there are some cancellations
between $x_{j}^{\phi _K}$ and  $w^{\phi _K}$, and the middle of
$x_j$ is the non-cancelled part of $x_j$ because  $x_j$ as a
letter not belonging to $E(m,n)$  appears only in $x_j^n$.

We still have to consider all letters that can appear to the right
of $w$, if $w$ is the end of some word in $T_1(m,n)$ or
$w=y_{n-1}x_n^{-1}x_{n-1}$, $w=y_{n-1}x_n^{-1}$. There are the
following possibilities: \bi \item[(i)] $w$ is an end of $
y_{n-2}x_{n-1}^{-1}x_nx_{n-1}y_{n-2}^{-1};$

\item[(ii)] $w$ is an end of $y_{r-2}x_{r-1}^{-1}x_r^{-1}, r<i$;

\item[(iii)] $w$ is an end of   $y_{n-2}x_{n-1}^{-1}y_{n-1}^{-1}$.
\ei Situation (i) is equivalent to  the situation when $w^{-1}$ is
the beginning of the word $
y_{n-2}x_{n-1}^{-1}x_nx_{n-1}y_{n-2}^{-1}$, we have considered
this case already. In the situation (ii)  the only possible word
to the right of $w$ will be left end of
$x_{r-1}y_{r-2}^{-1}x_{r-2}^{-2}$ and $w^{\phi _K }x_{r-1}^{\phi
_K}y_{r-2}^{-\phi _K}x_{r-2}^{-2\phi _K }=w^{\phi _K }\circ
x_{r-1}^{\phi _K}y_{r-2}^{-\phi _K}\circ x_{r-2}^{-2\phi _K},$ and
$w^{\phi _K }x_{r-1}^{\phi _K}=w^{\phi _K }\circ x_{r-1}^{\phi
_K}.$ In the situation (iii) the first two letters to the right of
$w$ are $x_{n-1}x_{n-1}$, and $w^{\phi _K}x_{n-1}^{\phi
_K}=w^{\phi _K}\circ x_{n-1}^{\phi _K}.$

There is no cancellation in the words $ (c_j^{z_j})^{\phi_K}\circ
(c_{j+1}^{\pm z_{j+1}})^{\phi_K},
 (c_m^{z_m})^{\phi_K}\circ  x_1^{\pm\phi_K}, \  x_1^{\phi_K}\circ x_1^{\phi_K}. $
 For  all the
other occurrences of $x_i$ in the words from ${\mathcal
W}_{\Gamma}$, namely for occurrences in $x_i^n,\ x_i^2y_i$, we
have $(x_i^2y_i)^{\phi _K}=x_i^{\phi _K}\circ x_i^{\phi _K}\circ
y_i^{\phi _k}$ for $i<n$.

In the case $n=i$, the bold subword of the word

\noindent $x_n^{\phi_{_K}} = A_{m+4n-4}^{-q_0}\circ {\bf \left (
x_n^{q_1}  \circ y_n  \circ A_{m+4n-2}^{q_2-1}\circ
A_{m+4n-4}^{-q_0}\circ x_n \right ) }\circ A_{m+4n-4}^{q_0} $

\noindent is $M_{x_n}$ for $\phi _K$, and the bold subword in the
word

\noindent $y_n^{\phi_K}= \ig{{
A_{m+4n-4}^{-q_0}}}{x_ny_{n-1}^{-1}}{y_{n-2}x_{n-1}^{-1}}
{\bf\left (x_n^{q_1} y_n \ \ig{{\bf A_{m+4n-2}^{q_2-1}}}{{\bf
x_ny_{n-1}^{-1}}}{{\bf x_ny_n}} \ig{{\bf A_{m+4n-4}^{-q_0}}}{{\bf
x_ny_{n-1}^{-1}}}{{\bf y_{n-2}x_{n-1}^{-1}}} \ x_n \right )^{q_3}
\ x_n^{q_1} y_n }$,

\noindent is $M_{y_n}$ for $\phi _K.$ \hfill $\Box$

\begin{lm}
\label{le:prozaik}
 The following statements hold:
  \begin{enumerate}
 \item [(1)] Let  $u \in E(m,n)$. If $B^2$ occurs as a subword in
 $u^{\phi_K}$ for some cyclically reduced word $B$ ($B\neq c_i$) then $B$ is a
 power of a cyclic permutation of a period $A_j, j = 1, \ldots K$.

 \item[(2)] Let  $u \in \bar{\mathcal{W}}_{\Gamma}$. If $B^2$ occurs as a subword in
 $u^{\phi_K}$ for some cyclically reduced word $B$ ($B\neq c_i$) then $B$ is a
 power of a cyclic permutation of a period $A_j, j = 1, \ldots K$.
  \end{enumerate}
\end{lm}
{\em Proof.} 1) follows from the formulas (1.a)-(4.d) from Lemma
\ref{main}.

2) We may assume that $w$ does not contain an elementary square. In
this case $w$ is a subword of a word from Lemma \ref{le:xyu}. Now
the result follows from the formulas (1.a)-(4.d) from Lemma
\ref{main}.

\hfill $\Box$

\medskip
{\bf Notation.} \ {\it 1) Denote by ${\mathcal W}_{\Gamma ,L}$ the
least set of words in the alphabet $Y$  that contains ${\bar
{\mathcal W}}_{\Gamma },$ is closed under taking subwords, and is
$\phi _K$-invariant.

 2) Let $\bar {\mathcal W}_{\Gamma ,L}$ be union of ${\mathcal
W}_{\Gamma ,L}$ and the set of all initial subwords of $z_i^{\phi
_{Kj}}$ which are of the form $c_i^j\circ z_i \circ w,$ where
$w\in{\mathcal W}_{\Gamma ,L}.$}

\medskip
{\bf Notation.} \ {\it  Denote by $Exc$ the following set of words
in the alphabet $Y$. \begin{enumerate} \item If $m > 2, n\geq 2$,
then
 $$Exc = \{ c_1^{-z_1}c_i^{-z_i}c_{i-1}^{-z_{i-1}}, \
 c_1^{-z_1}x_1c_m^{-z_m}, \ c_1^{-z_1}x_jy_{j-1}^{-1}\}.$$
 \item If $m > 2, n=1$, then
 $$Exc = \{ c_1^{-z_1}c_i^{-z_i}c_{i-1}^{-z_{i-1}}, \
 c_1^{-z_1}x_1c_m^{-z_m}\}.$$
\item If $m = 2, n\geq 2$, then
 $$Exc = \{
 c_1^{-z_1}x_1c_m^{-z_m}, \ c_1^{-z_1}x_jy_{j-1}^{-1}\}.$$
\item If $m = 2, n=1$, then
 $$Exc = \{
 c_1^{-z_1}x_1c_m^{-z_m}\}.$$
 \item If $m =1, n\geq 2$, then
 $$Exc = \{
  c_1^{-z_1}x_jy_{j-1}^{-1}\}.$$
\item If $m = 0, n\geq 2$, then
 $$Exc = \{y_1x_1x_i, x_1x_iy_{i-1}^{-1}, 2\leq i\leq n\}.$$
 \end{enumerate}}

\begin{lm}
\label{lm:W-Gamma-L} The following holds:

\begin{enumerate}
 \item [(1)] $Sub_{3,Y}({\mathcal W}_{\Gamma ,L})=Sub_{3,Y}(X^{\pm\phi _K})\cup Exc$.
\item [(2)] Let $v \in {\mathcal W}_{\Gamma ,L}$ be a word that
begins and ends with an elementary square and does not contain any
elementary cubes. Then either $v \in {\bar {\mathcal W}}_{\Gamma }$
or $v = v_1v_2$  for some words  $v_1,v_2\in {\bar {\mathcal
W}}_{\Gamma }$ described below:
\begin{enumerate}
\item for $m>2,\ n\geq 2$,
 $$v_1 \in \{v_{11}=(c_1^{z_1}c_2^{z_2})^2 \prod_{i = 3}^{m} c_i^{z_i}x_1x_2x_1 \prod_{i = m}^{1}c_i^{-z_i},
 \ v_{12}= x_1^2y_1x_1\prod_{i = m}^{1} c_i^{-z_i} \}$$
 and
  $$v_2 \in \{  v_{2i}=c_i^{-z_i}\ldots c_3^{-z_3}(c_2^{-z_2}c_1^{-z_1})^2, \ u_{2,1}=x_1c_m^{-z_m}\ldots c_3^{-z_3}
  (c_2^{-z_1}c_1^{-z_1})^2, \
  u_{2,j}=x_jy_{j-1}^{-1}x_{j-1}^2 \};$$

  \item for $m=2,\ n\geq 2$,
 $$v_1 \in \{v_{11}=(c_1^{z_1}c_2^{z_2})^2x_1x_2x_1 \prod_{i = m}^{1}c_i^{-z_i},
 \ v_{12}= x_1^2y_1x_1\prod_{i = m}^{1} c_i^{-z_i} \}$$
 and
  $$v_2 \in \{
   u_{2,1}=x_1
  (c_2^{-z_1}c_1^{-z_1})^2, \
  u_{2,j}=x_jy_{j-1}^{-1}x_{j-1}^2 \};$$

  \item for $m>2,\ n=1$,
 $$v_1 \in \{
 v_{12}= x_1^2y_1x_1\prod_{i = m}^{1} c_i^{-z_i} \}$$
 and
  $$v_2 \in \{  v_{2i}=c_i^{-z_i}\ldots c_3^{-z_3}(c_2^{-z_2}c_1^{-z_1})^2, \ u_{2,1}=x_1c_m^{-z_m}\ldots c_3^{-z_3}
  (c_2^{-z_1}c_1^{-z_1})^2 \};$$

 \item for $m=2,\ n=1$,
 $$v_1 \in \{ v_{12}= x_1^2y_1x_1\prod_{i = m}^{1} c_i^{-z_i} \}$$
 and
  $$v_2 \in \{
   u_{2,1}=x_1
  (c_2^{-z_1}c_1^{-z_1})^2 \};$$
  \item for $m=1,\ n\geq 2$,
 $$v_1 \in \{v_{11}=(c_1^{z_1}x_1^{-1})^2x_2x_1 c_1^{-z_1},
 \ v_{12}= x_1^2y_1x_1c_1^{-z_1} \}$$
 and
  $$v_2 \in \{
     u_{2,j}=x_jy_{j-1}^{-1}x_{j-1}^2 \}.$$

\end{enumerate}
 \item  [(3)] If $v \in {\mathcal W}_{\Gamma ,L}$ and either
 $v$ does not contain two elementary squares and
begins (ends) with an elementary square, or $v$ contains  no
elementary squares, then either $v$ is  a subword of one of the
words from [2] or (for $m=0$) $v$ is  a subword of one of the words
$ x_1^2y_1x_1,  x_2^2y_2x_2$.

\item [(4)] Automorphism $\phi _K$ satisfies Nielsen property with
respect ${\mathcal W}_{\Gamma ,L}$ with exceptions $E(m,n).$
\end{enumerate}
\end{lm}

{\em Proof.} Let $T=Kl.$ We will consider only the case $m\geq 2,\
n\geq 2$. We will prove all the statements of the lemma by
simultaneous induction on $l$. If $l=1$, then $T=K$ and the lemma
is true. Suppose now that
$$Sub _{3,Y}(\bar {\mathcal W}_{\Gamma}^{\phi _{T-K}})=Sub
_{3,Y}(\bar {\mathcal W}_{\Gamma})\cup Exc.$$ Formulas in the
beginning of the proof of Lemma \ref{main} show that  $$Sub
_{3,Y}(E(m,n)^{\pm\phi _K})\subseteq Sub _{3,Y}(\bar {\mathcal
W}_{\Gamma}).$$ By the third statement for $Sub (\bar{\mathcal
W}_{\Gamma }^{\phi _{T-K}})$ the automorphism $\phi _K$ satisfies
the Nielsen property with exceptions $E(m,n)$. Let us verify that
new 3-letter subwords do not occur "between" $u^{\phi _K}$ for
$u\in T_1(m,n)$ and the power of the corresponding $x_i$ to the
left and right of it. All the cases are similar to the following:

$$(x_nx_{n-1}y_{n-2}^{-1})^{\phi _K}\cdot x_{n-2}^{\phi _K}\ldots
\ig {A_{m+4n-10}^{-q+1}}{*}{y_{n-3}x_{n-2}^{-1}}\cdot
x_{n-1}^{-1}\ig{A_{m+4n-8}^{q_0-1}}{x_{n-2}}{*}\ .$$

Words $$ (v_1v_2)^{\phi _K}$$ produce the subwords from $Exc$.
Indeed, $[(x_2x_1\prod _{i=m}^1c_i^{-z_i})]^{\phi _{Kj}}$ ends
with $v_{12}$ and $v_{12}^{\phi _K}$ ends with $v_{12}.$
Similarly, $v_{2,j}^{\phi _K}$ begins with $v_{2,j+1}$ for $j<m$
and with $u_{2,1}$ for $j=m$. And $u_{2,j}^{\phi _K}$ begins with
$u_{2,j+1}$ for $j<n$ and with $u_{2,j}$ for $j=n$.

This and the second part of Lemma \ref{le:7.1.2words} finish the
proof. \hfill $\Box$

According to the definition of $\bar{\mathcal W}_{\Gamma ,L}$,
this set contains words which are written in the alphabet
$Y^{\pm1}$ as well as extra words $u$ of the form $(c_i^j z_i
w)^{\pm1}$ or $(z_i w)^{\pm1}$ whose $Y^{\pm1}$-representation is
spoiled at the start or at the end of $u$. For those $u \in
\bar{\mathcal W}_{\Gamma ,L}$ which are written in the alphabet
$Y^{\pm1}$, Lemma \ref{le:UF} gives a unique representation as the
product $u_1 \dots u_k$ where $u_i \in Y^{\pm1} \cup E(m,n)$ and
the occurrences of $u_i$ are maximal. We call this representation
a {\em canonical decomposition of $u$}. For $u\in \bar{\mathcal
W}_{\Gamma , L}$ of the form $(c_i^j z_i w)^{\pm1}$ or $(z_i
w)^{\pm1}$ we define the canonical decomposition of $u$ as follows
: $u=c_i\dots c_iz_iu_1 \dots u_k$ where $u_i \in Y^{\pm1} \cup
E(m,n)$. Clearly, we can consider the Nielsen property of
automorphisms with exceptions $E(m,n)$ relative to this extended
notion of canonical decomposition. Below the Nielsen property is
always assumed in this sense.

\begin{lm}\label{nielL} The automorphism $\phi _{K}$ satisfies Nielsen
property with respect to $\bar{\mathcal W}_{\Gamma ,L}$ with
exceptions $E(m,n).$ The set $\bar{\mathcal W}_{\Gamma ,L}$ is
$\phi _K$-invariant.\end{lm} {\em Proof.} The first statement
follows from Lemmas \ref{main} and  \ref{lm:W-Gamma-L}. For the
second statement notice that if $c_i^jz_iw\in\bar{\mathcal
W}_{\Gamma ,L}$, then $c_i^{c_i^jz_iw}=w^{-1}\circ c_i^{z_i}\circ
w\in {\mathcal W}_{\Gamma ,L}$ and $c_i^{(c_i^jz_iw)^{\phi
_K}}=w^{-\phi _K}\circ c_i^{z_i^{\phi _K}}\circ w^{\phi _K}\in
{\mathcal W}_{\Gamma ,L}$, therefore $c_i^{j}{z_i^{\phi _K}}\circ
w^{\phi _K}\in \bar{\mathcal W}_{\Gamma ,L}$. \hfill$\Box$

 Let  $W \in G[X]$.
 We say that a word  $U \in G[X]$
 {\it occurs } in $W$ if $W = W_1 \circ U \circ W_2$ for some $W_1, W_2 \in G[X]$.
An occurrence of $U^q$ in $W$ is called {\it maximal} with respect
to a property $P$ of words if $U^q$ is not a part of any
occurrence of $U^r$ with $q < r$ and which satisfies $P$. We say
that an occurrence of $U^q$ in $W$ is $t$-{\it stable} if $q \geq
1$ and $W = W_1 \circ U^t U^q U^t \circ W_2,\ t\geq 1,$ (it
follows that $U$ is cyclically reduced). If $t=1$ it is {\em
stable}. Maximal stable occurrences $U^q$ will play an important
part in what follows. If $(U^{-1})^q$
 is a stable occurrence of $U^{-1}$ in $W$ then, sometimes, we say
 that  $U^{-q}$ is a stable occurrence of $U$ in $W$. Two given occurrences
$U^q$ and $U^p$ in a word $W$ are {\it disjoint} if they do not
have a common letter as subwords of  $W$. Observe that if integers
$p$ and $q$ have different signs then any two occurrences of $A^q$
and $A^p$ are disjoint. Also, any two different maximal stable
occurrences of powers of $U$ are disjoint. To explain the main
property of stable occurrences of powers of $U$, we need the
following definition.
 We say that a given occurrence of $U^q$ {\em occurs correctly} in a given occurrence
  of $U^p$ if  $|q |\leq |p|$ and for these
occurrences $U^q$ and $U^p$ one has $U^p = U^{p_1} \circ U^q \circ
U^{p_1}$. We say, that two
 given non-disjoint occurrences of $U^q, U^p$ {\it overlap correctly} in $W$ if
 their common subword occurs correctly in each of them.

A cyclically reduced word $A$ from $G[X]$ which is not a proper
power and does not belong to $G$ is called {\it a period}.

\begin{lm}
\label{le:stable-overlap} Let $A$ be a period in $G[X]$ and $W \in
G[X]$. Then any two stable occurrences of powers of $A$ in $W$ are
either disjoint or they overlap correctly. \end{lm}
  {\it Proof.} Let $A^q$, $A^p$ ($q
\leq p$) be two non-disjoint stable occurrences of powers of $A$
in $W$. If they overlap  incorrectly then $A^2 = u \circ A \circ
v$ for some elements $u,v \in G[X]$. This implies that $A = u
\circ v = v \circ u$ and hence $u$ and $v$ are (non-trivial)
powers of some element in $G[X]$. Since $A$ is  not a proper power
it follows that $u = 1$ or $v = 1$ - contradiction. This shows
that $A^q$ and $A^p$ overlap correctly.  \hfill $\Box$

Let $W \in G[X]$ and ${\mathcal O} = {\mathcal O}(W,A) =\{A^{q_1},
\ldots, A^{q_k}\}$ be a set of pair-wise disjoint stable
occurrences of powers of a period $A$ in $W$ (listed according to
their appearance  in $W$ from the left to the right). Then
${\mathcal O}$ induces an ${\mathcal O}$-decomposition of $W$ of
the following form:
\begin{equation}
\label{eq:new-O-decomp} W = B_1 \circ A^{q_1} \circ \cdots \circ
B_k \circ A^{q_k} \circ B_{k+1}\end{equation}

For example, let $P$ be a property of words (or just a property of
occurrences in  $W$) such that if two powers of $A$ (two
occurrences of powers of $A$ in $W$) satisfy $P$ and overlap
correctly then their union also satisfies $P$. We refer to such
$P$ as  {\it preserving correct overlappings}. In this event, by
${\mathcal O}_P = {\mathcal O}_P(W,A)$ we denote the uniquely
defined set of all maximal stable occurrences of powers of $A$ in
$W$ which satisfy the property $P$. Notice, that occurrences in
${\mathcal O}_P$ are pair-wise disjoint by Lemma
\ref{le:stable-overlap}. Thus, if $P$ holds on every power of $A$
then ${\mathcal O}_P(W,A) = {\mathcal O}(W,A)$ contains all
maximal stable occurrences of powers of $A$ in $W$. In this case,
the decomposition (\ref{eq:new-O-decomp}) is unique and it is
called the {\it canonical (stable)} $A$-decomposition of $W$.

The following example provides another property $P$ that will be
in use later. Let $N$ be a positive integer and let $P_N$ be the
property of $A^q$ that $|q| \geq N$. Obviously, $P_N$ preserves
correct overlappings. In this case the set ${\mathcal O}_{P_N}$
provides the so-called {\it canonical $N$-large }
$A$-decompositions of $W$ which are also uniquely defined.

\begin{df}
 Let
$$W = B_1 \circ A^{q_1} \circ \cdots \circ B_k
\circ A^{q_k} \circ B_{k+1}$$ be the  decomposition
(\ref{eq:new-O-decomp}) of $W$ above. Then   the numbers
$$ \max_A(W) = \max \{q_i \mid i = 1, \dots, k\}, \ \ \
 \min_A(W) = \min \{q_i \mid i = 1, \dots, k\}$$ are called,
correspondingly,  the {\em upper} and the {\em lower} $A$-bounds
of $W$.
\end{df}

\begin{df}
Let $A$ be a period in $G[X]$ and $W \in G[X]$.   For a positive
integer $N$ we say that the $N$-large $A$-decomposition of $W$
 $$W = B_1 \circ
A^{q_1} \circ \cdots \circ B_k \circ A^{q_k} \circ B_{k+1}$$
 has $A$-size $(l,r)$ if $\min_A(W) \geq l$ and
 $\max_A(B_i) \leq r$ for every $i = 1, \dots, k$.
 \end{df}

 Let ${\mathcal A} = \{A_1, A_2,
\ldots, \}$ be a sequence of periods from $G[X]$.   We say that a
word $W \in G[X]$ has ${\mathcal A}$-rank $j$ (${\rm
rank}_{\mathcal A}(W) = j$) if $W$ has a stable occurrence of
$(A_j^{\pm 1})^q$ ($q \geq 1$) and $j$ is maximal with this
property. In this event, $A_j$ is called the {\em ${\mathcal
A}$-leading term} (or just the {\em leading term}) of $W$
(notation $LT_{{\mathcal A}}(W) = A_j$ or $LT(W) = A_j$).

We now fix an arbitrary  sequence ${\mathcal A}$ of periods  in
the group $G[X]$. For a period $A = A_j$ one can consider
canonical $A_j$-decompositions of a word $W$ and define the
corresponding $A_j$-bounds and $A_j$-size. In this case
 we, sometimes, omit $A$ in the writings and simply write $max_j(W)$ or
$min_j(W)$ instead of $max_{A_j}(W)$,  $min_{A_j}(W)$.

In the case when   ${\rm rank}_{{\mathcal A}}(W) = j$  the
canonical $A_j$-decomposition of $W$ is called the {\em canonical}
${\mathcal A}$-decomposition of $W$.

Now we turn to an analog of ${\mathcal O}$-decompositions of $W$
with respect to "periods" which are not necessarily cyclically
reduced words. Let $U = D^{-1} \circ A \circ D$, where  $A$ is a
period. For  a set ${\mathcal O} = {\mathcal O}(W,A) = \{A^{q_1},
\ldots, A^{q_k}\}$ as above  consider the ${\mathcal
O}$-decomposition of a word $W$
 \begin{equation} \label{eq:A-dec}
 W = B_1 \circ A^{q_1} \circ \cdots \circ B_k \circ A^{q_k} \circ B_{k+1}
 \end{equation}
  Now it can be rewritten in the form:
 $$ W = (B_1D) (D^{-1} \circ A^{q_1} \circ D) \ldots (D^{-1}B_kD) (D^{-1}\circ A^{q_k} \circ D)
 (D^{-1} B_{k+1}).$$
Let $\varepsilon _i, \delta _i=sgn(q_i).$ Since every occurrence
of $A^{q_i}$ above is stable, $B_1 = {\bar B}_1 \circ
A^{\varepsilon _1}$, $B_i = (A^{\delta _{i-1}} \circ {\bar B}_i
\circ A^{\varepsilon _i})$, $B_{k+1} = A^{\delta _k} \circ {\bar
B}_{k+1}$ for suitable words ${\bar B}_i$.  This shows that the
decomposition above can be written as
$$W = ({\bar B}_1 A^{\varepsilon _1} D) (D^{-1} A^{q_1} D) \ldots (D^{-1}
A^{\delta _{i-1}} {\bar B}_i  A^{\varepsilon _i} D) \ldots (D^{-1}
A^{q_k}  D)
 (D^{-1} A^{\delta _k}{\bar B}_{k+1}) = $$
$$({\bar B}_1 D) (D^{-1}  A^{\varepsilon _1} D) (D^{-1} A^{q_1} D) \ldots (D^{-1}
A^{\delta _{i-1}}D)( D^{-1}{\bar B}_i D) (D^{-1} A^{^{\varepsilon
_i}} D) \ldots $$ $$(D^{-1} A^{q_k} D)
 (D^{-1} A^{\delta _k} D)(D^{-1} {\bar B}_{k+1})$$
  $$= ({\bar B}_1D) (U^{\varepsilon _1}) (U^{q_1}) \ldots (U^{\delta _{k-1}})( D^{-1}{\bar B}_k D)
 (U^{\varepsilon _k}) (U^{q_k}) (U^{\delta _k})(D^{-1} {\bar B}_{k+1}).$$
 Observe, that the cancellation between parentheses  in the
decomposition above does not exceed the length  $d = |D|$ of $D$.
Using notation $w = u \circ_d v$ to indicate that the cancellation
between $u$ and $v$ does not exceed the number $d$, we can rewrite
the decomposition above in the following form:
$$W = ({\bar B}_1D) \circ_{d}  U^{\varepsilon _1} \circ_{d} U^{q_1} \circ_{d}  U^{\delta _1} \circ_d \cdots
\circ_d U ^{\varepsilon _k} \circ_{d} U^{q_k} \circ_{d} U^{\delta
_k} \circ_{d} (D^{-1} {\bar B}_{k+1}),$$
 hence
 \begin{equation}
 \label{eq:U-dec}
 W = D_1 \circ_d U^{q_1} \circ_d \cdots
\circ_d D_k \circ_d U^{q_k} \circ_d D_{k+1}, \end{equation}
 where
$D_1 = {\bar B}_1D, \ D_{k+1} = D^{-1} {\bar B}_{k+1},  \ D_i=
D^{-1}{\bar B}_i D \ (2 \leq i \leq k)$, and the occurrences
$U^{q_i}$ are $(1,d)$-stable. (We similarly define $(t,d)$-stable
occurrences.) We will refer to this decomposition of $W$ as  {\em
$U$-decomposition with respect to ${\mathcal O}$} (to get a
rigorous definition of $U$-decompositions one has to replace in
the definition of the ${\mathcal O}$-decomposition of $W$ the
period $A$ by $U$ and $\circ$ by $\circ_{|D|}$). In the case when
an $A$-decomposition of $W$ (with respect to ${\mathcal O}$) is
unique then the corresponding $U$-decomposition of $W$  is also
unique, and in this event  one can easily rewrite
$A$-decompositions of $W$ into $U$-decomposition and vice versa.

We  summarize the discussion above in the following lemma.

\begin{lm} \label{le:A-U}
Let $A \in G[X]$ be a period and $U = D^{-1} \circ A \circ D \in
G[X]$. Then for a word $W \in G[X]$ if

 $$W = B_1 \circ A^{q_1} \circ \cdots \circ B_k \circ A^{q_k} \circ
 B_{k+1}$$
is a stable $A$-decomposition of $W$ then
 $$W = D_1 \circ_d U^{q_1} \circ_d \cdots
\circ_d D_k \circ_d U^{q_k} \circ_d D_{k+1}$$
 is a stable $U$-decomposition of $W$, where $D_i$ are defined as in (\ref{eq:U-dec}). And vice versa.

\end{lm}


 From now on we fix the following  set of leading terms
$${\mathcal A}_{L,p}
= \{ A_j \mid j \leq L,  \phi = \phi_{L,p} \}$$ for a given
multiple $L$ of $K = K(m,n)$ and a given tuple $p$.

\begin{df}
Let $W \in G[X]$
 and $N$ be a positive integer.   A word of the type $A_s$ is termed
 the  $N$-large leading term $LT_N(W)$  of the word $W$ if
$A_s^q$ has a stable occurrence in $W$ for some  $q \geq N$, and
$s$ is maximal with this property. The number $s$ is called the
$N$-rank of $W$ ($s = {\rm rank}_N(W), s\geq 1$).
\end{df}

In Lemmas  \ref{le:7.1.zforms}, \ref{le:7.1.x1formsm0},
\ref{le:7.1.x1formsmneq0}, and   \ref{le:7.1.xiforms}  we described
precisely the leading terms $A_j$ for  $j = 1, \ldots,K$.    It is
not easy to describe precisely $A_j$ for an arbitrary $j > K$.  So
we are not going to do it  here, instead, we chose a  compromise by
introducing a modified version $A_j^*$  of $A_j$ which is not
cyclically reduced, in general, but which  is ``more cyclically
reduced" then the initial word $A_j$. Namely, let $L$ be a multiple
of $K$ and  $1\leq j\leq K.$ Define
 $$A^*_{L+j}= A^*(\phi _{L+j}) = A_j^{\phi _L}.$$

\begin{lm}
\label{le:cyclicAA*}  Let $L$ be a multiple of $K$ and  $1\leq j\leq
K.$  Let $p=(p_1,\dots ,p_n)$ be $N+3$-large tuple. Then
 $$A^*_{L+j}  = R^{-1} \circ A_{L+j} \circ R$$
for some word $R \in F(X \cup C_S)$ such that   $rank(R)\leq L-K+j
+2$ and $|R| < |A_{L+j}|.$
\end{lm}

{\it Proof.} First, let $L=K.$ Consider elementary periods
$x_i=A_{m+4i-3}$ and $A_1=c_1^{z_1}c_2^{z_2}$. For $i\neq n$,
$x_i^{2\phi _K}=x_i^{\phi _k}\circ x_i^{\phi _K}$.  For $i=n$,
$$A^*(\phi _{K+m+4n-3})=R^{-1}\circ A _{K+m+4n-3}\circ R,$$
where $R=A_{m+4i-4}^{p_{m+4n-4}}$, therefore ${\rm
rank}_N(R)=m+4n-4.$ For the other elementary period,
$(c_1^{z_1}c_2^{z_2})^{2\phi _K}=(c_1^{z_1}c_2^{z_2})^{\phi
_K}\circ (c_1^{z_1}c_2^{z_2})^{\phi _K}.$

Any other $A_j$ can be written in the  form $A_j=u_1\circ v_1\circ
u_2\circ v_2\circ u_3$, where $v_1, v_2$ are the first and the
last elementary squares in $A_j$, which are parts of big powers of
elementary periods.  The Nielsen property of $\phi _K$ implies
that the word $R$ for  $A^*(\phi _{K+j})$  is  the word that
cancels between $(v_2u_3)^{\phi _K}$ and $(u_1v_1)^{\phi _K}.$ It
definitely has $N$-large rank  $\leq K$, because the element
$(v_2u_3u_1v_1)^{\phi _K}$ has $N$-large rank $\leq K$. To give an
exact bound for the rank of $R$ we consider all possibilities for
$A _j$:
\begin{enumerate}
\item $A_{i}$ begins with $z_i^{-1}$ and ends with $z_{i+1}$,
$i=1,\dots ,m-1$, \item $A_{m}$ begins with $z_m^{-1}$ and ends
with $x_1^{-1}$, \item $A_{m+4i-4}$ begins with
$x_{i-1}y_{i-2}^{-1}x_{i-2}^{-2}$, if $i=3,\dots n$, and ends with
$x_{i-1}^2y_{i-1}x_{i}^{-1}$ if $i=2,\dots ,n$, If $i=2$ it begins
with $x_1\Pi _{j=m}^1c_j^{-z_j}(c_2^{-z_2}c_1^{-z_1})^2$.

\item $A_{m+4i-2}$  and $A_{m+4i-1}$ begins with
$x_iy_{i-1}^{-1}x_{i-1}^{-2}$ and ends with $x_i^2y_{i}$ if
$i=1,\dots ,n.$

\end{enumerate}

Therefore, $A_{i}^{\phi _K}$ begins with $z_{i+1}^{-1}$ and ends
with $z_{i+2}$, $i=1,\dots ,m-2$, and is cyclically reduced.

 $A_{m-1}^{\phi _K}$ begins with $z_m^{-1}$ and ends with $x_1$, and is cyclically reduced,
 $A_{m}^{\phi _K}$ begins with $z_m^{-1}$ and ends with $x_1^{-1}$ and is cyclically reduced.

 We have already considered $A_{m+4i-3}^{\phi _K}$.

Elements $A_{m+4i-4}^{\phi _K}, A_{m+4i-2}^{\phi _K},
A_{m+4i-1}^{\phi _K}$ are not cyclically reduced. By Lemma
\ref{main}, for $A^*_{K+m+4i-4}$, one has
$R=(x_{i-1}y_{i-2}^{-1})^{\phi _K}$ ($rank(R)=m+4i-4$); for
$A^*_{K+m+4i-2}$, and $A^*_{K+m+4i-2}$, $R=(x_{i}y_{i-1}^{-1})^{\phi
_K}$ ($rank(R)=m+4i$).

This proves the statement of the Lemma for $L=K$.

We can suppose by induction that $A^*_{L-K+j}=R^{-1}\circ
A_{L-K+j}\circ R$, and $rank(R)\leq L-2K+j+2$. The cancellations
between $A_{L-K+j}^{\phi _K}$ and $R^{\phi _K}$ and between
$A_{L-K+j}^{\phi _K}$ and $A_{L-K+j}^{\phi _K}$ correspond to
cancellations in words $u^{\phi _K}$, where $u$ is a word in
${\mathcal W}_{\Gamma}$ between two elementary squares. These
cancellations  are in rank $\leq K$, and the statement of the lemma
follows. \hfill $\Box$

 \begin{lm}
 \label{le:Phiwords1}
Let  $W \in F(X\cup C_S)$ and  $A = A_j= LT_N(W)$, and $A^* =
R^{-1} \circ A \circ R$. Then $W$ can be presented in the form
\begin{equation}
\label{eq:Phiform} W = B_1 \circ_d A^{* q_1} \circ_d B_2 \circ_d
\cdots \circ_d B_k \circ _d A^{* q_k} \circ_d B_{k+1}
 \end{equation}
where  $A^{* q_i}$ are maximal stable $N$-large occurrences of $A^*$
 in $W$ and $d \leq |R|$.  This presentation  is unique and it is called
 the
canonical $N$-large $A^*$-decomposition of $W$.
\end{lm}
{\it Proof.} The result follows from existence and  uniqueness of
the canonical $A$-decompositions. Indeed, if
 $$W = B_1  \circ A^{ q_1} \circ  B_2 \circ \cdots  \circ B_k \circ A^{q_k}
\circ B_{k+1} $$ is the canonical  $A$-decomposition  of $W$, then
 $$ (B_1R)(R^{-1}AR)^{ q_1}(R^{-1} B_2R) \cdots  (R^{-1}B_kR) (R^{-1}AR)^{q_k}
 (R^{-1}B_{k+1}) $$
is the canonical $A^*$-decomposition of $W$. Indeed, since every
$A^{q_i}$ is a stable occurrence, then every $B_i$ starts with $A$
(if $i \neq 1$) and ends with $A$ (if $i = k+1$). Hence $R^{-1}
B_iR= R^{-1} \circ B_i \circ R$.

Conversely if
 $$ W = B_1 A^{* q_1} B_2 \cdots B_k A^{* q_k} B_{k+1}$$
  is an $A^*$-representation of $W$ then
 $$W = (B_1R^{-1})  \circ A^{ q_1} \circ  (RB_2R^{-1}) \circ \cdots  \circ (RB_kR^{-1})
  \circ A^{q_k}
\circ (RB_{k+1}) $$ is the canonical  $A$-decomposition for $W$.
\hfill $\Box$

Let $\phi$ be an automorphism of $F(X\cup C)$ which satisfies the
Nielsen property with respect to a set $W$ with exceptions $E$. In
Definition \ref{de:nielsen}, we have introduced the notation
$M_{\phi ,w}$ for the middle of $w$ with respect to $\phi$ for $w\in
Y \cup E$. We now introduce a similar notation for any $w \in
Sub(W)$ denoting by $\bar{M}_{\phi ,w}$ the maximal noncancelled
part of $w^\phi$ in the words $(u w v)^\phi$ for all $uwv \in W$
with $w \neq u^{-1}, v^{-1}$. Observe that, in general,
$\bar{M}_{\phi ,w}$ may be empty while this cannot hold for $M_{\phi
,w}$. If $\bar{M}_{\phi ,w}$ is nonempty then we represent $w^\phi$
as
$$ w^\phi = \bar{L}_{\phi ,w} \circ \bar{M}_{\phi ,w} \circ \bar{R}_{\phi ,w }.$$

\begin{lm}\label{R}
Let $L = lK$,  $l >0$, $p$  a $3$-large tuple.
\begin{itemize}
\item[(1)] If $E$ is closed under taking subwords then
$\bar{M}_{\phi ,w}$ is nonempty whenever the irreducible
 decomposition of $w$ has length at least 3.
  \item[(2)] $\bar{M}_{\phi_L, (A^2)}$ is nonempty for an
elementary period $A$.
  \item [(3)]The automorphism $\phi_L=\phi _{L,p}$ has the Nielsen property
with respect to $\bar{\mathcal W}_{\Gamma,L}$ with exceptions
$E(m,n)$. For $w \in X\cup E(m,n)$ and $l >1$, the middle $M_{\phi_L
,w}$ can be described in the following way. Let
$$
    M_{\phi_K ,w} = f \circ A^r \circ g \circ B^s \circ h
$$
where $A^r$ and $B^s$ are the first and the last maximal occurrences
of elementary powers in $M_{\phi_K ,w}$. Then $M_{\phi_L , w}$
contains $\bar{M}_{\phi_{L-K}, (A^r g B^s)}$ as a subword.

 \item  [(4)] If $i < j \leq L$ then $A_j^2$ does not occur in
 $A_i$.
\end{itemize}
\end{lm}

{\it Proof.} To prove (1) observe  that  if  $w = u_1 u_2 u_3$, $u_i
\in Y\cup E$, is the irreducible decomposition of $w$ then
$\bar{M}_{\phi ,w}$ should contain $M_{\phi ,u_2}$.

The middles $M_{\phi _K,x}$ of elements from $X$ and from $E(m,n)$
contain  big powers of some $A_j$, where $j=1,\dots ,K$ and,
therefore, big powers of elementary periods. Therefore, statements
(2) and (3) can be proved by the simultaneous  induction on $l$.
Notice that for $l=1$ both statements follow from Lemma \ref{nielL}.

The statement (4) follows from Lemmas
\ref{le:7.1.zforms}--\ref{le:7.1.x1formsmneq0}.

 \hfill $\Box$

\begin{lm}\label{le:properties-gamma-words}
Let $L= lK>0$, $1\le ir \le K$, $t\ge2$, $p$  a $3$-large tuple,

\begin{enumerate}\item [(1)]
and
$$
    w = u \circ A_r^s \circ v
$$
be a $t$-stable occurrence of $A_r^s$ in a word $w \in
\bar{\mathcal W}_{\Gamma,L}$. Let $A^*_{r+L} = R^{-1} \circ
A_{r+L} \circ R$ and $d = |R|$. Then
$$
    w^{\phi_L} = u^{\phi_L} \circ_d (A_{r+L}^*)^s \circ_d v^{\phi_L}
$$
where the occurrence of $(A_{r+L}^*)^s$ is $(t-2,d)$-stable.

\item [(2)] Let  $W \in \bar{\mathcal W}_{\Gamma,L}$, and
$A^*_{L+r} = R^{-1} \circ A_{L+r} \circ R$ and $d = |R|$.

 If $t\ge 2$ and
$$
W = D_1 \circ A_r^{q_1} \circ D_2 \dots \circ A_r^{q_k} \circ
D_{k+1}
$$
is a $t$-stable $A_r$-decomposition of $W$ then
$$
W^{\phi_L} = D_1^{\phi_L} \circ_d (A_{L+r}^*)^{q_1} \circ_d
D_2^{\phi_L} \circ_d \dots \circ_d (A_{L+r}^*)^{q_k} \circ_d
D_{k+1}^{\phi_L}
$$
is a $(t-2,d)$-stable $A_{L+r}^*$-decomposition of $W^{\phi_L}$.
\end{enumerate}
\end{lm}

{\em Proof.} (1) Clearly, we can assume $t=2$ without loss of
generality. Suppose first that $A_r$ is not an elementary period.
Then the canonical decomposition of $A_r$ is of length $\ge3$ and
thus $\bar{M}_{\phi_L} (A_r)$ is nonempty by Lemma \ref{R}(a).
This implies that $u^{\phi_L}$ ends with $\bar{M}_{\phi_L} (A_r)
\bar{R}_{\phi_L} (A_r)$, and thus the cancellation between
$u^{\phi_L}$ and $(A_{r+L}^*)^r$ is the same as in the product
$A_{r+L}^* \cdot A_{r+L}^*$. Similarly, the same is the
cancellation between $(A_{r+L}^*)^r$ and $v^{\phi_L}$ and the
statement of lemma follows.

If $A_r$ is an elementary period, a slightly more careful analysis
is needed. We first consider the image of $w$ under $\phi_K$. If
$r=1$ one of the images $A_1^{\pm \phi_K}$ of the periods
$A_1^{\pm1}$ in the occurrence of $A^{s+4\text{sgn}(r)}$ in $w$
(i.e.\ the first or the last one) may be completely cancelled in
$w^{\phi_K}$, but all the others have nonempty noncancelled
contributions in $w^{\phi_K}$. Then an easy application of Lemma
\ref{R} (with $L$ replaced with $L-K$) gives the result, and this is
the case when only $(t-2,d)$-stability can be stated. If $A_r$ is an
elementary period of the form $x_j$, a similar argument applies but
with no possibility of completely cancelled period $A_r^{\pm1}$
under $\phi_K$.

(2) follows from (1). \hfill $\Box$

\begin{lm}\label{elper}
Let   $A_{j_1}\ldots ,A_{j_k}$, $k \geq 0$,  be elementary periods,
$1 \leq j_1, \ldots, j_k \leq K$.  If $w \in \bar{\mathcal
W}_{\Gamma,L}$ and

\begin{equation}\label{eldec}
 w^{\phi_K}=\tilde w_0\circ_{d_{j_1}}A^{*q_1}_{j_1+K}\circ
_{d_{j_1}}\tilde w_1\ldots\circ _{d_{j_k}}A^{*q_k}_{j_k+K}\tilde
w_k,\end{equation}
 where $q_i\geq 5$, ${\tilde w}_i$ does not contain an elementary square, and $d_{j_i} = |R_{j_i}|$,
 where $A_{j_i+K}^*=R_{j_i}^{-1}\circ A_{j_i+K}\circ R_{j_i}$ (see Lemma \ref{le:cyclicAA*}), $i = 1,
 \ldots,k$, then
 $$w=w_0\circ A_{j_1}^{q_1}\circ w_1\circ\ldots \circ A_{j_k}^{q_k}\circ w_k,$$
where $w_i^{\phi _K}=\tilde w_i,\ i=0,\ldots ,k.$
\end{lm}
{\em Proof.}  {\em Case 1)}. Suppose that $w$  does not contain an
elementary square.

In this case either $w \in \bar{\mathcal{W}}_{\Gamma}$ or $w =
v_1v_2$ for some words $v_1, v_2 \in \bar{\mathcal{W}}_{\Gamma}$
which are described in Lemma \ref{lm:W-Gamma-L}.

{\em Claim 1}. If $w^{\phi _K}$ contains $B^s$ for some cyclically
reduced  word $B\neq c_i, i = 1, \ldots,m$, and $s \geq 2$, then $B$
is a power of a cyclic permutation of some uniquely defined period
$A_i, i=1,\ldots ,K$.

It suffices to consider the case $s = 2$. Notice that for $w \in
\bar{\mathcal{W}}_{\Gamma}$ the claim follows from Lemma
\ref{le:prozaik}. Now observe that if $w = v_1v_2$ for $v_1,v_2 \in
\bar{\mathcal{W}}_{\Gamma}$ then $w^{\phi_K} = v_1^{\phi_K} \circ
v_2^{\phi_K}$ and "illegal" squares do not occur on the boundary
between $v_1^{\phi_K}$ and $v_2^{\phi_K}$ (direct inspection).

{\em Claim 2}.   $w^{\phi _K}$ does not contain  $(E^{\phi _K})^2$,
where $E$ is an elementary period.
  By Claim 1 $w^{\phi _K}$ contains  $(E^{\phi _K})^2$,
where $E$ is an elementary period, if and only if $E^{\phi _K}$ is a
power of a cyclic permutation of some period $A_i, i=1,\ldots ,K$.
So it suffices to show that $E^{\phi _K}$ is not a power of a cyclic
permutation of some period $A_i, i=1,\ldots ,K$.  To this end we
list below all the words $E^{\phi_K}$:

 by Lemma \ref{le:7.1.zforms}
$$A_1^{\phi _K}=A_1^{-p_1+1}c_2^{-z_2}A_1^{p_1}A_2^{-p_2+1}
c_3^{-z_3}A_1^{-p_1+1}c_2^{z_2}A_1^{p_1-1}c_3^{z_3}A_2^{p_2-1}
(m\geq 2);$$
 by Lemma \ref{le:7.1.xiforms}
$$x_i^{\phi _{K }}=\ig{A_{m+4i-2}^{q_2}}{x_iy_{i-1}^{-1}}{x_iy_i}
\ig{A_{m+4i-4}^{-q_0}}{x_iy_{i-1}^{-1}}{y_{i-2}x_{i-1}^{-1}}\ x_i\
\ig{A_{m+4i-4}^{q_0}}{x_{i-1}y_{i-2}^{-1}}{y_{i-1}x_i^{-1}}\ (i\neq
n);$$
 and (direct computation from Lemmas  \ref{le:7.1.zforms}
and \ref{le:7.1.xiforms}).

$$(c_1^{z_1}x_1^{-1})^{\phi _K}=z_1x_1^{p_2}y_1A_3^{p_3}A_1^{-p_1}x_1(A_1^{-p_1}x_1^{p_2}y_1A_3^{p_3}A_1^{-p_1}x_1)^{p_4-2}x_1^{p_2}y_1x_2^{-1}(y_1^{\phi
_4}x_2^{-1})^{p_5-1}\ \ (n>1),$$

$$(c_1^{z_1}x_1^{-1})^{\phi
_K}=z_1x_1^{-1}A_1^{p_1}A_3^{-p_3+1}y_1^{-1}x_1^{-p_2}A_1^{p_1}\ \
(n=1).$$

The claim follows by comparing the formulas above with the
corresponding formulas for $A_j$ (Lemmas  \ref{le:7.1.zforms} -
\ref{le:7.1.xiforms}).

Now the claim 2 implies the lemma since in this case the
decomposition (\ref{eldec}) for the $w^{\phi_K}$ is of  the form
 $w^{\phi_K}= \tilde{w}_0$ and $w = w_0$, as required.

{\em Case 2}.  $w^{\phi _K}$ contains  $(E^{\phi _K})^2$, where $E$
is an elementary period. By the case 1) $w$ has a non-trivial
decomposition of the form
 $$w=w_0\circ A_{j_1}^{q_1}\circ w_1\circ\ldots \circ A_{j_k}^{q_k}\circ w_k,$$
 where $q_i \geq 2$, and $w_i$ does not have squares of elementary
 periods. Consider the $A_r$-decomposition of
 $w$ where $r = \max \{j_1, \ldots, j_k\}$:
  $$w = D_1 \circ A_r^{q_1} \circ \ldots A_r^{q_s} \circ
  D_{s+1},$$
  where $D_i$ does not contain a square of an elementary period.
  It follows from the case 1 that this decomposition is at least
  $3$-large canonical stable $A_r$-decomposition of $w$. Indeed,
  if $E_1$ and $E_2$ are two distinct elementary periods  then $E_1^{s\phi
_K}$ does not contain a cyclically reduced part of $E_2^{2\phi _K}$
as a subword (see the formulas above). So in the canonical
$3$-stable $A^*_{r+K}$-decomposition of $w^{\phi_K}$ the powers
$A_{r+K}^{*q_i}$ come from the corresponding powers of $A_r$.
 By  Lemma
\ref{le:properties-gamma-words}
  $$w^{\phi_K} = D_1^{\phi_K} \circ_d A_{K+r}^{*q_1} \circ_d \ldots A_{K+r}^{*q_s} \circ_d
  D_{s+1}^{\phi_K},$$
is the canonical stable $A_{j+K}^*$-decomposition of $w^{\phi _K}$
that contains all the occurrences of powers of $A_{j+K}^*$ in the
decomposition (\ref{eldec}). Now by induction on the maximal rank
of elementary periods which squares appear in the words $D_i$  we
can finish the proof.
  \hfill$\Box$
\begin{lm}\label{le:A-A^*}
Let $L= lK>0$, $1\le r \le K$,   $A^*_{r+L} = R^{-1} \circ A_{r+L}
\circ R$ and $d = |R|$. Then the following holds for every $w \in
\bar{\mathcal W}_{\Gamma,L}$:\begin{enumerate} \item [(1)] Suppose
there is a decomposition
$$ w^{\phi_K} = \tilde{u} \circ_f (A_{r+K}^*)^s \circ_f
    \tilde{v},$$
where $s\geq 5$ and  and the cancellation between $\tilde u$ and
$A_{r+K}^*$ (resp., between
 $A_{r+K}^*$ and $\tilde v$) is not more than $f$ which is the maximum of the corresponding $d$ and length of
 the part of  $A_{r+K}^*$ before the first stable occurrence of an
 elementary power (resp., after the last stable occurrence of an
 elementary power). Then
$$
    w = u \circ A_r^s \circ v, \quad u^{\phi_K} = \tilde{u}, \quad v^{\phi_K} = \tilde{v}.
$$
\item [(2)] Let
$$
W^{\phi_L} = \tilde{D}_1 \circ_d (A_{L+r}^*)^{q_1} \circ_d
\tilde{D}_2 \circ_d \dots \circ_d (A_{L+r}^*)^{q_k} \circ_d
\tilde{D}_{k+1}
$$
be a $(1,d)$-stable $3$-large $A_{L+r}^*$-decomposition of
$W^{\phi_L}$.
 Then $W$ has a stable  $A_r$-decomposition
$$
W = D_1 \circ A_r^{q_1} \circ D_2 \dots \circ A_k^{q_k} \circ
D_{k+1}
$$
where $D_i^{\phi_L} = \tilde{D}_i$. \end{enumerate}

\end{lm}

{\em Proof.}  (1) If $A_r$ is an elementary period, the statement
follows from Lemms \ref{elper}. Otherwise represent $A_r$ as $A_r=
A_{j_1}^{q_1}\circ w_1\circ A_{j_2}^{q_2}\circ w_2,$ where
$A_{j_1}^{q_1}$ and $A_{j_2}^{q_2}$ are the first and the last
maximal elementary powers (each $A_i$ begins with an elementary
power).

Then $$w^{\phi _K}={\tilde u}\circ _d(A_{j_1+K}^{q_1}\circ
_dw_1^{\phi _K}\circ _d A_{j_2+K}^{q_2}\circ _d w_2^{\phi
_K})^{s-1}\circ _dA_{j_1+K}^{q_1}\circ _dw_1^{\phi _K}\circ _d
A_{j_2+K}^{q_2}\circ _d w_2^{\phi _K}\circ _f\tilde v.$$ Since
$\phi _K$ is a monomorphism, by Lemma \ref{elper} we obtain
$$w=u\circ (A_{j_1}^{q_1}\circ w_1\circ A_{j_2}^{q_2}\circ
w_2)^{s-1}\circ  A_{j_1}^{q_1}\circ w_1\circ A_{j_2}^{q_2}\circ
w_2v,$$ where $u^{\phi _K}=\tilde u, \ v^{\phi _K}=\tilde v .$ We
will show that $w_2v=w_2\circ v$. Indeed, $w_2$ is either
$c_i^{z_i}, i\geq 3$, or $y_{i-1}x_i^{-1}$, or $y_i$.  If there is
a cancellation between $w$ and $v$, then $v$ must respectively
begin either with $c_i^{-z_i}$, or $x_i$ or $y_i^{-1}$ and the
image of this letter when $\phi _K$ is applied to $v$ must be
almost completely cancelled. It follows from Lemma \ref{main} that
this does not happen. Therefore  $w = u \circ A_r^s \circ v,$ and
(1) is proved.

(2) For  $L=K$ statement (1) implies statement (2). We now use
induction on $l$ to prove (2).

Suppose \begin{equation}\label{LLL}
w^{\phi _L}=\tilde u\circ
_dA_{r+L}^*\circ _d\tilde v.\end{equation}

Represent $A^*_{r+K}$ as
$$A^*_{r+K}=w_0\circ A_{i_1}^{q_1}\circ w_1\circ
A_{i_2}^{q_2}\circ w_2,$$ where $A_{i_1}^{s_1}$ and
$A_{i_2}^{s_2}$ are the first and the last maximal occurrences of
elementary powers.

Then
$$w^{\phi _L}=\tilde uw_0^{\phi _{L-K}}\circ _d(A_{i_1}^{s_1\phi _{L-K}}\circ
_dw_1^{\phi _{L-K}}\circ _dA_{i_2}^{s_2\phi _{L-K}}\circ
_d(w_2w_0)^{\phi _{L-K}})^{s-1}\circ _d$$ $$A_{i_1}^{s_1\phi
_{L-K}}\circ _dw_1^{\phi _{L-K}}\circ _dA_{i_2}^{s_2\phi
_{L-K}}\circ _d(w_2)^{\phi _{L-K}}\tilde v.$$ By the assumption of
induction
$$w^{\phi _K}=\hat uw_0\circ (A_{i_1}^{s_1}\circ w_1\circ
A_{i_2}^{s_2}\circ (w_2w_0))^{s-1}\circ A_{i_1}^{s_1}\circ
w_1\circ A_{i_2}^{s_2}\circ (w_2\hat v),$$ where $\hat u^{\phi
_{L-K}}=\tilde u, \hat v^{\phi _{L-K}}=\tilde v.$ Therefore
$$w^{\phi _K}=\hat u\circ _f A_{r+K}^{*s}\circ _f\hat v.$$ By
statement (1), $w=u\circ A_r^s\circ v,$ where  $u^{\phi _K}=\hat
u,\ v^{\phi _K}=\hat v.$ Therefore (\ref{LLL}) implies that
$w=u\circ A_r^s\circ v,$ where $u^{\phi _L}=\tilde u, v^{\phi
_L}=\tilde v.$ This implies (2) for $L$.
 \hfill$\Box$

\begin{cy} \label{cy:middles}\begin{enumerate}\item [(1)]
  Let $m \neq 0, n \neq 0, K = K(m,n), p=(p_1,\dots
,p_K)$ be a 3-large tuple, $L = Kl$.  Then  for any $u \in Y \cup
E(m,n)$ the element $M_{\phi _L, u}$  contains $A_j^q$ for some
$j> L-K$ and $q > p_j - 3$.

\item [(2)] For any  $x \in X$ if $rank(x^{\phi_L})$ = j then
every
 occurrence of $A_j^2$ in $x^{\phi_L}$ occurs inside some occurrence
 of $A_j^{N-3}$.
\end{enumerate}\end{cy} {\em Proof.} (1) follows from the formulas for $M_u$
with respect to ${\phi_K}$ in Lemma \ref{main} and Lemma
\ref{le:properties-gamma-words}.\hfill$\Box$

\begin{cy}\label{co:12-new}
Let $u, v \in  {\mathcal W}_ {\Gamma , L}$. If  the canceled
subword
 in the product  $u^{\phi_{K}}v^{\phi_{K}}$ does not contain
$A_j^l$ for some $j \leq K$ and $l \in \mathbb{Z}$ then the canceled
subword in the product $u^{\phi_{K+L}}v^{\phi_{K+L}}$ does not
contain the subword $A_{L+j}^l$.
\end{cy}

\begin{lm}\label{le:7.1.gammawords}

  Suppose $p$ is an $(N+3)$-large tuple, $\phi _j=\phi _{jp}$. Let $L$ be a
multiple of $K$. Then:

\begin{enumerate}\item [(1)]
\begin{enumerate}
   \item $x_i^{\phi _j}$ has a canonical $N$-large
$A^*_j$-decomposition of size $(N,2)$ if  either $j\equiv
m+4(i-1)(mod \ K)$, or $j \equiv m+4i-2 (mod \ K)$, or $j \equiv
m+4i (mod \ K)$. In all other cases $rank(x_i^{\phi _j}) < j$.
    \item $ y_i^{\phi _j}$ has a canonical $N$-large $A^*_j$-decomposition
 of size $(N,2)$ if either $j\equiv m+4(i-1)(mod \ K)$, or $j \equiv m+4i-3(mod \
K)$, or $j \equiv m+4i-1(mod \ K)$, or $j \equiv  m+4i\ (mod \
K).$ In all other cases $rank(y_i^{\phi _j}) < j$.
   \item $z_i^{\phi _{j}}$ has a canonical $N$-large
   $A^*_j$-decomposition of size $(N,2)$
 if $j\equiv i\ (mod \ K)$ and either $1 \leq i \leq m-1$ or $i = m$
 and $n \neq 0$. In all other
cases $rank(z_i^{\phi _j}) < j$.
    \item if $n = 0$ then $z_m^{\phi _{j}}$ has a canonical $N$-large
    $A^*_j$-decomposition of size $(N,2)$
 if $j\equiv m-1\ (mod \ K)$.  In all other
cases $rank(z_m^{\phi _j}) < j$.
 \end{enumerate}

\item [(2)] If $j=r+L$, $0<r\leq K,\ (w_1\ldots w_k)\in
Sub_k(X^{\pm\gamma_K \ldots\gamma_{r+1}})$ then either $(w_1\ldots
w_k)^{\phi _j}=(w_1\ldots w_k)^{\phi _{j-1}},$ or $ (w_1\ldots
w_k)^{\phi _j}$ has a canonical $N$-large $A^*_j$-decomposition.
In any case, $ (w_1\ldots w_k)^{\phi _j}$   has a canonical
$N$-large $A^*_s$-decomposition in some rank $s,$ $j-K+1\leq s\leq
j.$
\end{enumerate}\end{lm}

{\em Proof.} (1) Consider  $y_i^{\phi _{L+m+4i}}:$

$$y_i^{\phi _{L+m+4i}}=(   x_{i+1}^{\phi _{L}} y_i^{-\phi _{L+
m+4i-1}})^{q_4-1}x_{i+1}^{\phi _{L}} (y_i^{\phi _{ L+m+4i-1}}
x_{i+1}^{-\phi _{L}})^{q_4},$$

In this case $A^*(\phi _{L+m+4i})=x_{i+1}^{\phi
_{L+m+4i-1}}y_i^{-\phi _{L+m+4i-1}}$.

To write a formula for  $x_i^{\phi _{L+m+4i}}$, denote $\tilde
y_{i-1}=y_{i-1}^{\phi _{L+m+4i-5}},\ \bar x_i=x_i^{\phi _{L}}, \
\bar y_i=y_i^{\phi _{L}}$. Then

\medskip
\begin{multline*} x_i^{\phi _ {L+m+4i}}=( \bar x_{i+1} y_i^{-\phi _
{L+m+4i-1}})^{q_4-1} \bar x_{i+1}\\ ( ((
 \bar x_i\tilde  y_{i-1}^{-1})^{q_0}    \bar x_i^{q_1}    \bar y_i)^{q_2-1} (
   \bar x_i\tilde
y_{i-1}^{-1})^{q_0}    \bar x_i^{q_1+1}    \bar y_i)^{-q_3+1} \bar
y_i^{-1} \bar x_i^{-q_1}(\tilde y_{i-1} \bar
x_i^{-1})^{q_0}.\end{multline*}

 Similarly we consider
 $ z_i^{\phi _{L+i}}$.

(2)  If in a word $(w_1\ldots w_k)^{\phi _j}$ all the powers of
$A_j^{p_j}$ are cancelled (they can only cancel completely and the
process of cancellations does not depend on $p$) then if we
consider an $A_j^*$-decomposition of $(w_1\ldots w_k)^{\phi _j}$,
all the powers of $A_j^*$ are also completely cancelled. By
construction of the automorphisms $\gamma _j$, this implies that
$(w_1\ldots w_k)^{\gamma _j\phi _{j-1}}= (w_1\ldots w_k)^{\phi
_{j-1}}.$

$\Box$

\subsection{Generic  Solutions of Orientable Quadratic Equations}
\label{se:7.2}

Let $G$ be a finitely generated fully residually free group and $S
= 1$ a standard quadratic orientable equation over $G$ which has a
solution in $G$. In this section we effectively construct
discriminating  sets of solutions of $S = 1$ in  $G$. The main
tool in this construction is an embedding
 $$\lambda: G_{R(S)} \rightarrow G(U,T)$$
of the coordinate group $G_{R(S)}$ into a group $G(U,T)$ which is
obtained from  $G$ by finitely many extensions of centralizers.
There is a nice  set $\Xi_P$ (see Section \ref{se:2-5}) of
discriminating $G$-homomorphisms from $G(U,T)$ onto $G$. The
restrictions of homomorphisms from $\Xi_P$
 onto  the image $G_{R(S)}^\lambda$ of $G_{R(S)}$ in $G(U,T)$
 give a discriminating  set of
 $G$-homomorphisms from $G_{R(S)}^\lambda$ into $G$, i.e.,
solutions of $S = 1$ in $G$. This idea was introduced in
\cite{KMNull}  to describe the radicals of quadratic equations.

It has been shown in \cite{KMNull} that the coordinate groups of
non-regular standard quadratic equations $S = 1$ over $G$ are
already extensions of centralizers of $G$, so in this case we can
immediately put $G(U,T) = G_{R(S)}$ and the result follows. Hence
we can assume from the beginning that $S= 1$ is regular.

Notice, that all regular quadratic  equations have solutions in
general position, except for the equation $[x_1,y_1][x_2,y_2] = 1$
(see Section \ref{se:2-7}).

 For the equation $[x_1,y_1][x_2,y_2]
= 1$ we do the following trick. In this case we view the
coordinate group $G_{R(S)}$  as the coordinate group of the
equation $[x_1,y_1]= [y_2,x_2]$ over the group of constants $G
\ast F(x_2,y_2)$. So the commutator $[y_2,x_2] = d$ is a
non-trivial constant  and the new equation is of the form $[x,y] =
d$, where all solutions are in general position. Therefore, we can
assume that  $S = 1$ is  one of the following types (below $d,c_i$
are nontrivial elements from $G$):
\begin{equation}\label{eq:st7}
\noindent \prod_{i=1}^{n}[x_i,y_i] = 1, \ \ \ n \geq 3;
\end{equation}
\begin{equation}\label{eq:st8}
\prod_{i=1}^{n}[x_i,y_i] \prod_{i=1}^{m}z_i^{-1}c_iz_i d = 1,\ \ \
n \geq 1, m \geq 0;
\end{equation}
\begin{equation}\label{eq:st9}
 \prod_{i=1}^{m}z_i^{-1}c_iz_i d = 1,\ \ \  m \geq 2,
\end{equation}
and it has a solution in $G$ in general position.

Observe, that since $S= 1$ is regular then Nullstellenzats  holds
for $S = 1$, so $R(S) = ncl(S)$ and $G_{R(S)} = G[X]/ncl(S) =
G_S$.

For a group $H$ and an element $u \in H$ by  $H(u,t)$ we denote
the  extension of the centralizer $C_H(u)$ of $u$:
$$H(u,t) = \langle H,t \ | \ t^{-1}xt = x  \ \ (x \in C_H(u)) \rangle.$$
If
 $$ G = G_1  \leq G_1(u_1,t_1) = G_2 \leq  \ldots \leq
G_n(u_n,t_n) = G_{n+1}
$$
   is a chain of extensions of centralizers of elements $u_i \in G_i$,
then we denote the resulting group $G_{n+1}$ by $G(U,T)$, where $U
= \{u_1, \ldots, u_n\}$ and $T = \{t_1, \ldots, t_n\}$.

  Let   $\beta: G_{R(S)} \rightarrow G$ be  a solution of the equation
  $S(X) = 1$ in the group $G$ such that
 $$x_i^\beta = a_i, y_i^\beta = b_i, z_i^\beta = e_i. $$
Then
 $$d = \prod_{i=1}^{m}e_i^{-1}c_ie_i \prod_{i=1}^{n}[a_i,b_i].$$
Hence we can rewrite the equation $S = 1$  in the following form
(for appropriate $m$ and $n$):
\begin{equation}\label{6-general}
\prod_{i=1}^{m}z_i^{-1}c_iz_i\prod_{i=1}^{n}[x_i,y_i]  =
\prod_{i=1}^{m}e_i^{-1}c_ie_i \prod_{i=1}^{n}[a_i,b_i].
\end{equation}

\begin{prop}
\label{prop:hom-lambda} Let $S = 1$ be a regular quadratic
equation (\ref{6-general})  and $\beta:G_{R(S)} \rightarrow G$ a
solution of $S= 1$ in $G$ in a general position. Then one can
effectively construct a sequence of extensions of centralizers
$$ G = G_1  \leq G_1(u_1,t_1) = G_2 \leq  \ldots \leq
G_n(u_n,t_n) = G(U,T)$$
 and a $G$-homomorphism $\lambda_\beta : G_{R(S)} \rightarrow G(U,T)$.
 \end{prop}
{\it Proof.}  By induction we define a sequence of  extensions of
centralizers and a sequence of group homomorphisms in the
following way.

{\it Case: $m \neq 0, n = 0$.} In this event for each $i = 1,
\ldots, m-1$ we define by induction  a pair $(\theta_i, H_i)$,
consisting of a group $H_i$ and a $G$-homomorphism $\theta_i:G[X]
\rightarrow H_i$.

Before we will go into formalities let us explain the idea that
lies behind this. If $z_1 \rightarrow e_1, \ldots, z_m \rightarrow
e_m$ is a solution of an equation
\begin{equation} \label{eq:zz}
z_1^{-1}c_1z_1 \ldots z_m^{-1}c_mz_m = d,
\end{equation}
then  transformations
\begin{equation}\label{eq:ztr2}
e_i \rightarrow e_i(c_i^{e_i}c_{i+1}^{e_{i+1}})^q, \ \ e_{i +1}
\rightarrow e_{i+1}(c_i^{e_i}c_{i+1}^{e_{i+1}})^q,  \ \ e_j
\rightarrow e_j  \ \ \ (j \neq i, i+1),
\end{equation}
produce a new solution of the equation (\ref{eq:zz}) for an
arbitrary integer $q$. This solution is composition of the
automorphism $\gamma _{i}^q$ and the solution $e$. To avoid
collapses under cancellation of the periods
$(c_i^{e_i}c_{i+1}^{e_{i+1}})^q$ (which is an important part of
the construction of the discriminating set of homomorphisms
$\Xi_P$ in Section \ref{se:2-5}) one might want to have number $q$
as big as possible, the best way would be to have $q = \infty$.
Since there are no infinite powers in $G$, to realize this idea
one should go outside the group $G$ into a bigger group,
 for example, into an ultrapower $G^\prime$ of $G$, in which
 a non-standard power, say $t$,  of the element $c_i^{e_i}c_{i+1}^{e_{i+1}}$
 exists. It is not hard to see that the subgroup $\langle G,t\rangle \leq G^\prime$
 is an extension of the centralizer $C_G(c_i^{e_i}c_{i+1}^{e_{i+1}})$ of
 the element $c_i^{e_i}c_{i+1}^{e_{i+1}}$ in  $G$. Moreover, in the group
 $\langle G,t\rangle $
 the transformation  (\ref{eq:ztr2}) can be described as
\begin{equation}\label{eq:ztr3}
e_i \rightarrow e_it, \ \ e_{i +1} \rightarrow e_{i+1}t, \ \ e_j
\rightarrow e_j  \ \ \ (j \neq i, i+1),
\end{equation}
Now, we are going to construct formally the subgroup  $\langle
G,t\rangle $ and the corresponding homomorphism using
(\ref{eq:ztr3}).

 Let $H$ be an arbitrary group and $\beta: G_{S} \rightarrow H$
a homomorphism. Composition of the canonical projection $G[X]
\rightarrow G_S$ and $\beta$ gives a homomorphism $\beta_0: G[X]
\rightarrow H$.   For $i =0$ put
$$H_0 = H, \ \ \ \theta _0 = \beta_0$$
 Suppose now, that a group $H_{i}$ and
a homomorphism $\theta _{i}: G[X] \rightarrow H_i$ are already
defined. In this event we define $H_{i+1}$ and $\theta_{i+1}$ as
follows
$$
H_{i+1}=<H_{i},r_{i+1}|[C_{H_{i}}(c_{i+1}^{z_{i+1}^{\theta
_{i}}}c_{i+2}^{z_{i+2}^{\theta_{i}}}),r_{i+1}]=1>, $$
$$z_{i+1}^{\theta_{i+1}}=z_{i+1}^{\theta_{i}}r_{i+1}, \ \
z_{i+2}^{\theta_{i+1}}=z_{i+2}^{\theta_{i}}r_{i+1}, \ \
z_j^{\theta_{i+1}}= z_j^{\theta_{i}}, \ \ \ (j \neq i+1, i+2).
$$
By induction we constructed a series of extensions of centralizers
$$
G = H_0 \leq H_1 \leq \ldots \leq H_{m-1} = H_{m-1}(G)
$$
and a homomorphism
$$\theta_{m-1,\beta} = \theta _{m-1}: G[X] \rightarrow H_{m-1}(G).$$
 Observe, that,
  $$c_{i+1}^{z_{i+1}^{\theta
_{i}}}c_{i+2}^{z_{i+2}^{\theta_{i}}} =
c_{i+1}^{e_{i+1}r_{i}}c_{i+2}^{e_i}$$
 so  the element $ r_{i+1}$ extends the
centralizer of the element $c_{i+1}^{e_{i+1}r_{i}}c_{i+2}^{e_i}$.
In particular, the following equality holds in the group
$H_{m-1}(G)$ for each $i = 0, \ldots, m-1$:
\begin{equation}\label{eq:r}
[r_{i+1},c_{i+1}^{e_{i+1}r_{i}}c_{i+2}^{e_i}] = 1.
\end{equation}
(where $r_0 = 1$). Observe also, that
\begin{equation}
\label{eq:z} z_1^{\theta_{m-1}} =  e_1r_1, \ \ z_i^{\theta_{m-1}}
=  e_ir_{i-1}r_{i}, \ \ z_m^{\theta_{m-1}} =  e_mr_{m-1} \ \ \
(0<i<m).
\end{equation}

From (\ref{eq:r}) and (\ref{eq:z}) it readily follows that
\begin{equation}
\label{eq:homz} (\prod_{i=1}^{m}z_i^{-1}c_iz_i ) ^{\theta_{m- 1}}
= \prod_{i=1}^{m}e_i^{-1}c_ie_i,
\end{equation}
so   $\theta_{m-1}$ gives rise to a homomorphism (which we again
denote by $\theta_{m-1}$ or $\theta_\beta$)
 $$\theta_{m-1} : G_S \longrightarrow H_{m-1}(G). $$
  Now we iterate the construction one more time replacing $H$ by $H_{m-1}(G)$
  and $\beta$ by $\theta_{m-1}$ and put:
  $$H_\beta(G) = H_{m-1}(H_{m-1}(G)), \ \ \
  \lambda_\beta = \theta_{\theta_{m-1}}: G_S \longrightarrow H_\beta(G).$$
The group $H_\beta(G)$ is union of a chain of extensions of
centralizers which starts at the group $H$.

If $H = G$ then  all the homomorphisms
 above are $G$-homomorphisms. Now  we can write
  $$H_\beta(G) = G(U,T)$$
  where $U = \{u_1, \ldots, u_{m-1},
\bar u_1, \ldots, \bar u_{m-1}\}$, $T = \{r_1, \ldots, r_{m-1},
\bar r_1, \ldots, \bar r_{m-1}\}$ and   $ \bar u_i, \bar r_i$ are
the corresponding elements when
 we iterate the construction:
 $$
 u_{i+1} =
c_{i+1}^{e_{i+1}r_{i}}c_{i+2}^{e_{i+2}}, \ \  \bar u_{i+1} =
c_{i+1}^{e_{i+1}r_{i}r_{i+1}\bar
r_i}c_{i+2}^{e_{i+2}r_{i+1}r_{i+2}}.$$

{\it Case:  $m = 0, n > 0$. } In this case $S = [x_1,y_1] \ldots
[x_n,y_n]d^{-1}$. Similar to the case above we start with the
principal automorphisms. They consist of two Dehn's twists:
\begin{equation}\label{t1}
x \rightarrow y^px, \ \ y \rightarrow y;
\end{equation}
\begin{equation}\label{t2}
x  \rightarrow x, \ \ y \rightarrow x^py;
\end{equation}
which fix the commutator $[x,y]$, and the third transformation
which ties two consequent commutators
$[x_i,y_i][x_{i+1},y_{i+1}]$:
\begin{equation}\label{t3}
x_i \rightarrow (y_ix_{i+1}^{-1})^{-q}x_i, \ \ y_i \rightarrow
(y_ix_{i+1}^{- 1})^{-q}y_i(y_ix_{i+1}^{-1})^q,
\end{equation}
$$x_{i+1} \rightarrow (y_ix_{i+1}^{-1})^{-q}x_{i+1}(y_ix_{i+1}^{- 1})^q, \ \
y_{i+1} \rightarrow (y_ix_{i+1}^{-1})^{-q}y_{i+1}.$$

 Now we define
by induction on $i$,  for $i = 0, \ldots, 4n-1$,   pairs
$(G_i,\alpha_i)$ of groups $G_i$ and $G$-homomorphisms $\alpha_i:
G[X] \rightarrow G_i$. Put
 $$G_0 = G, \ \ \alpha_0  = \beta.$$
For each commutator $[x_i,y_i] $ in $S = 1$ we  perform
consequently three Dehn's twists (\ref{t2}), (\ref{t1}),
(\ref{t2}) (more precisely, their analogs for an extension of a
centralizer) and an analog of the connecting transformation
(\ref{t3}) provided the next commutator exists. Namely, suppose
$G_{4i}$ and $\alpha _{4i}$ have been already defined. Then
$$G_{4i+1}=<G_{4i},t_{4i+1}| [C_{G_{4i}}(x_{i+1}^{\alpha
_{4i}}),t_{4i+1}]=1>;$$ $$y_{i+1}^{\alpha
_{4i+1}}=t_{4i+1}y_{i+1}^{\alpha _{4i}}, \ \ s^{\alpha
_{4i+1}}=s^{\alpha _{4i}} \ \ (s \neq y_{i+1}).$$
$$G_{4i+2}=<G_{4i+1},t_{4i+2}| [C_{G_{4i+1}}(y_{i+1}^{\alpha
_{4i+1}}),t_{4i+2}]=1>;$$ $$x_{i+1}^{\alpha
_{4i+2}}=t_{4i+2}x_{i+1}^{\alpha _{4i+1}}, \ \ s^{\alpha
_{4i+2}}=s^{\alpha _{4i+1}} \ \ \ (s \neq x_{i+1};$$
 $$G_{4i+3}=<G_{4i+2},t_{4i+3}| [C_{G_{4i+2}}(x_{i+1}^{\alpha
_{4i+2}}),t_{4i+3}]=1>;$$ $$y_{i+1}^{\alpha
_{4i+3}}=t_{4i+3}y_{i+1}^{\alpha _{4i+2}}, \ \ s^{\alpha
_{4i+3}}=s^{\alpha _{4i+2}} \ \ \ (s \neq y_{i+1});$$
$$G_{4i+4}=<G_{4i+3},t_{4i+4}| [C_{G_{4i+3}}(y_{i+1}^{\alpha
_{4i+3}}x_{i+2}^{-\alpha _{4i+3}}),t_{4i+4}]=1>;$$
$$x_{i+1}^{\alpha _{4i+4}}=t_{4i+4}^{-1}x_{i+1}^{\alpha _{4i+3}},
y_{i+1}^{\alpha _{4i+4}}=y_{i+1}^{\alpha _{4i+3}t_{4i+4}},
x_{i+2}^{\alpha _{4i+4}}=x_{i+2}^{\alpha _{4i+3}t_{4i+4}},$$ $$
y_{i+2}^{\alpha _{4i+4}}=t_{4i+4}^{-1}y_{i+2}^{\alpha _{4i+3}}, \
\ s^{\alpha _{4i+4}}=s^{\alpha _{4i+3}} \ \  (s \neq
x_{i+1},y_{i+1},x_{i+2},y_{i+2}).$$
 Thus we have defined groups
$G_i$ and  mappings $\alpha_i$ for all $i = 0, \ldots 4n-1$.    As
above, the straightforward verification shows that the mapping
$\alpha_{4n-1}$ gives rise to a $G$-homomorphism $\alpha_{4n-1}:
G_S \longrightarrow G_{4n-1}.$ We repeat now the above
construction once more  time with $G_{4n-1}$ in the place of
$G_0$, $\alpha _{4n-1}$ in the place of $\beta$, and  $\bar t_j$
in the place of $ t_j$. We denote the corresponding groups and
homomorphisms by ${\bar G_i}$ and ${\bar \alpha_i}: G_S
\rightarrow {\bar G_i}$.

 Put
 $$G(U,T) = \bar G_{4n-1}, \ \ \ \lambda_\beta = \bar \alpha_{4n-1},$$
By induction we have constructed  a $G$-homomorphism
$$\lambda_\beta : G_S \longrightarrow G(U,T).$$

{\it Case:}  $m > 0, n > 0$.  In this case we combine the two
previous cases together. To this end we take the group $H_{m-1}$
 and the homomorphism $\theta_{m-1}: G[X] \rightarrow H_{m-1}$
 constructed in the first case and put them as the input for
 the construction in the second case. Namely, put
$$G_0=<H_{m-1},r_{m}| [C_{H_{m-1}}(c_m^{z_m^{\theta_{m-1}}}
x_1^{-\theta_{m- 1}}),r_{m}]=1>,$$
 and define the homomorphism $\alpha_0$ as follows
 $$z_m^{\alpha _0}=z_m^{\theta
_{m-1}}r_{m}, \ \ x_1^{\alpha _0}=a_1^{r_{m}}, \ \ y_1^{\alpha
_0}=r_{m}^{-1}b_1, \ \  s^{\alpha _0}=s^{\theta _{m- 1}} \ \ (s\in
X, s \neq z_m, x_1,y_1).$$
  Now we  apply the construction from the second case.
 Thus we have defined groups
$G_i$ and  mappings $\alpha_i: G[X] \rightarrow G_i$ for all $i =
0, \dots, 4n-1$. As above, the straightforward verification shows
that the mapping $\alpha_{4n-1}$ gives rise to a $G$-homomorphism
$\alpha_{4n-1}: G_S \longrightarrow G_{4n-1}.$

We repeat now the above construction once more time with
$G_{4n-1}$ in place of $G_0$ and  $\alpha _{4n-1}$ in place of
$\beta$. This results in a group $\bar G_{4n-1}$ and a
homomorphism $\bar \alpha_{4n-1}: G_S \rightarrow \bar G_{4n-1}$.

Put
 $$G(U,T) = \bar G_{4n-1}, \ \ \ \lambda_\beta = \bar \alpha_{4n-1}.$$
We have constructed  a $G$-homomorphism
$$\lambda_\beta : G_S \longrightarrow G(U,T).$$

 We proved the proposition for all three types of equations (\ref{eq:st7}),
 (\ref{eq:st8}), (\ref{eq:st9}), as required.

 \hfill $\Box$

\begin{prop}
\label{prop:hom-lambda-monic} Let $S = 1$ be a regular quadratic
equation (\ref{eq:st2})
 and $\beta:G_{R(S)} \rightarrow G$ a
solution of $S= 1$ in $G$ in a general position. Then  the
homomorphism $\lambda_\beta : G_{R(S)} \rightarrow G(U,T)$ is a
monomorphism.
\end{prop}
{\it Proof.} In the proof of this proposition we use induction on
the atomic rank of the equation in the same way as in the proof of
Theorem 1 in \cite{KMNull}.

Since all the intermediate groups are also fully residually free
 by induction it suffices to prove the following:

1. $n=1, m=0$; prove that $\psi =\alpha _3$ is an embedding of
$G_S$ into $G_3$.

2. $n=2, m=0$;  prove that $\psi =\alpha _4$ is a monomorphism on
$H=<G,x_1,y_1>.$

3. $n=1, m=1$; prove that $\psi=\alpha _{3}\bar \alpha _0$ is a
monomorphism on $H=<G,z_1>.$

4. $n=0, m\geq 3$; prove that $\theta _{2}\bar \theta _2$ is an
embedding of $G_S$ into $\bar H_2$.

Now we consider all these cases one by one.

 Case 1. Choose an
arbitrary nontrivial element $h\in G_S$. It can be written in the
form $$ h = g_1\ v_1(x_1,y_1)\ g_2\ v_2(x_1,y_1) \ g_3 \ldots
v_n(x_1,y_1)\ g_{n+1}, $$ where $1 \neq v_i(x_1,y_1) \in
F(x_1,y_1)$ are words in $x_1, y_1$, not belonging to the subgroup
$\langle [x_1,y_1]\rangle ,$ and $1 \neq g_i \in G, g_i\not\in
<[a,b]>$ (with the exception of $g_1$ and $g_{n+1}$, they could be
trivial). Then \beq \label{66} h^\psi = g_1\ v_1(t_3t_1a,t_2b)\
g_2\ v_2(t_3t_1a,t_2b) \ g_3 \ldots v_n(t_3t_1a,t_2b)\ g_{n+1}.
\eeq
 The group $G(U,T)$ is obtained from $G$ by three
HNN-extensions (extensions of centralizers), so every element in
$G(U,T)$ can be rewritten to its reduced form by making finitely
many pinches. It is easy to see that the leftmost occurrence of
either $t_3$ or $t_1$ in the product (\ref{66}) occurs in the
reduced form of $h^\psi$ uncancelled.

Case 2. $x_1\rightarrow t_4^{-1}t_2a_1,\ y_1\rightarrow
t_4^{-1}t_3t_1b_1t_4,\ x_2\rightarrow t_4^{-1}a_2t_4, \
y_2\rightarrow t_4^{-1}b_2.$ Choose an arbitrary nontrivial
element $h\in H=G*F(x_1,y_1)$. It can be written in the form $$ h
= g_1\ v_1(x_1,y_1)\ g_2\ v_2(x_1,y_1) \ g_3 \ldots v_n(x_1,y_1)\
g_{n+1}, $$ where $1 \neq v_i(x_1,y_1) \in F(x_1,y_1)$ are words
in $x_1, y_1$, and $1 \neq g_i \in G.$ (with the exception of
$g_1$ and $g_{n+1}$, they could be trivial). Then \beq \label{7}
h^\psi = g_1\ v_1(t_4^{-1}t_2a,(t_3t_1b)^{t_4})\ g_2\
v_2(t_4^{-1}t_2a,(t_3t_1b)^{t_4}) \ g_3 \ldots
v_n(t_4^{-1}t_2a,(t_3t_1b)^{t_4})\ g_{n+1}. \eeq
 The group $G(U,T)$ is obtained from $G$ by four
HNN-extensions (extensions of centralizers), so every element in
$G(U,T)$ can be rewritten to its reduced form by making finitely
many pinches. It is easy to see that the leftmost occurrence of
either $t_4$ or $t_1$ in the product (\ref{7}) occurs in the
reduced form of $h^\psi$ uncancelled.

Case 3. We have an equation $c^z[x,y]=c[a,b]$, $z\rightarrow
zr_1\bar r_1,\ x\rightarrow (t_2a^{r_1})^{\bar r _1},\
y\rightarrow \bar r_1^{-1}t_3t_1r_1^{-1}b,$ and $[r_1,ca^{-1}]=1,\
[\bar r_1, (c^{r_1}a^{-r_1}t_2^{-1})]=1.$ Here we can always
suppose, that $[c,a]\not =1$, by changing a solution, hence
$[r_1,\bar r_1]\not =1.$ The proof for this case is a repetition
of the proof of Proposition 11 in \cite{KMNull}. \hfill$\Box$

Case 4. We will consider the case when $m=3$; the general case can
be considered similarly. We have an equation
$c_1^{z_1}c_2^{z_2}c_3^{z_3}=c_1c_2c_3$, and can suppose
$[c_i,c_{i+1}]\neq 1.$

We will prove that $\psi=\theta _2\bar\theta _1$ is an embedding.
The images of $z_1,z_2,z_3$ under $\theta _2\bar\theta _1$ are the
following:
$$z_1\rightarrow c_1r_1\bar r_1,\ z_2\rightarrow c_2r_1r_2\bar
r_1,\ z_3\rightarrow c_3r_2,$$ where
$$
[r_1,c_1c_2]=1,\ [r_2, c_2^{r_1}c_3]=1,\ [\bar r_1,
c_1^{r_1}c_2^{r_1r_2}]=1.$$

Let $w$ be a reduced word in  $G*F(z_i, i=1,2,3),$ which does not
have subwords $c_1^{z_1}$. We will prove that if $w^{\psi}=1$ in
$\bar H_1$, then $w\in N,$ where  $N$ is  the normal closure of
the element $c_1^{z_1}c_2^{z_2}c_3^{z_3}c_3^{-1}c_2^{-1}c_1^{-1}.$
We use induction on the number of occurrences of $z_1^{\pm 1}$ in
$w$. The induction basis is obvious, because homomorphism $\psi$
is injective on the subgroup $<F,z_2,z_3>.$

Notice, that the homomorphism $\psi$ is also injective on the
subgroup $K=<z_1z_2^{-1}, z_3, F>.$

Consider $\bar H_1$ as an HNN-extension by letter $\bar r_1$.
Suppose $w^{\psi}=1$ in $\bar H_1$. Letter $\bar r_1$ can
disappear in two cases: 1) $w\in KN,$  2) there is a pinch between
$\bar r_1^{-1}$ and $\bar r_1$ (or between $\bar r_1$ and $\bar
r_1^{-1}$) in $w^{\psi}.$ This pinch corresponds to some element
$z_{1,2}^{-1}uz^{\prime}_{1,2}$ (or $z_{1,2}u(z^{\prime}
_{1,2})^{-1}$), where $z_{1,2}, z^{\prime}_{1,2}\in\{z_1,z_2\}.$

In the first case $w^{\psi}\neq 1$, because $w\in K$ and $w\not\in
N$.

In the second case, if the pinch happens in
$(z_{1,2}u(z^{\prime}_{1,2})^{-1})^{\psi}$, then
$z_{1,2}u(z^{\prime}_{1,2})^{-1}\in KN,$ therefore it has to be at
least one pinch that corresponds to
$(z_{1,2}^{-1}uz^{\prime}_{1,2})^{\psi}$. We can suppose, up to a
cyclic shift of $w$, that $z_{1,2}^{-1}$ is the first letter, $w$
does not end with some $z_{1,2}^{\prime\prime}$, and $w$ cannot be
represented as
$z_{1,2}^{-1}uz^{\prime}_{1,2}v_1z^{\prime\prime}_{1,2}v_2,$ such
that $z_{1,2}^{\prime}v_1\in KN.$ A pinch can only happen if
$z_{1,2}^{-1}uz^{\prime}_{1,2}\in <c_1^{z_1}c_2^{z_2}>$.
Therefore, either $z_{1,2}=z_1,$ or $z_{1,2}^{\prime}=z_1$, and
one can replace $c_1^{z_1}$ by  $c_1c_2c_3c_3^{-z_3}c_2^{-z_2}$,
therefore replace $w$ by $w_1$ such that $w=uw_1$, where $u$ is in
the normal closure of the element
$c_1^{z_1}c_2^{z_2}c_3^{z_3}c_3^{-1}c_2^{-1}c_1^{-1},$ and apply
induction.
  \hfill $\Box$

The embedding $\lambda_\beta : G_S \rightarrow G(U,T)$ allows one
to construct effectively discriminating sets of solutions in $G$
of the equation $S = 1$. Indeed, by the construction above  the
group $G(U,T)$ is union of the following  chain of length $2K =
2K(m,n)$ of extension of centralizers:
 $$
  G = H_0 \leq H_1 \ldots \leq H_{m-1} \leq G_0 \leq G_1 \leq \ldots
  \leq G_{4n-1} = $$
   $$ = \bar H_0 \leq \bar H_1 \leq \ldots \leq
 \bar H_{m-1 }=  \bar G_0 \leq \ldots \leq \bar G_{4n-1}
 = G(U,T).$$
 Now, every $2K$-tuple $p \in {\mathbb N}^{2K}$ determines  a
$G$-homomorphism
 $$\xi_p: G(U,T) \rightarrow G.$$
  Namely, if $Z_i$ is the $i$-th term of the chain above then $Z_i$ is an extension of
  the centralizer of some element $g_i \in Z_{i-1}$ by a stable letter $t_i$. The $G$-homomorphism
$\xi_p$  is defined as composition
 $$\xi_p= \psi_1 \circ \ldots \circ \psi_K$$
 of homomorphisms
$\psi_i:Z_i \rightarrow Z_{i-1}$ which are identical on $Z_{i-1}$
 and such that $t_i^{\psi_i} = g_i^{p_i}$, where $p_i$ is the $i$-th component of
$p$.

 It follows (see  Section~2.5)
  that for every unbounded set of tuples
$P \subset {\mathbb N}^{2K}$  the set of homomorphisms
 $$\Xi_P =  \{\xi_p \mid p \in P\}$$
 $G$-discriminates $G(U,T)$ into $G$. Therefore, (since $\lambda_\beta$ is monic),
 the family of $G$-homomorphisms
 $$ \Xi_{P,\beta} = \{\lambda_\beta  \xi_p \mid \xi_p \in \Xi_P
 \}$$
$G$-discriminates $G_S$ into $G$.

 One can give another description of the set $\Xi_{P,\beta}$ in terms of the basic
 automorphisms from the  basic sequence $\Gamma$. Observe first that
 $$ \lambda_\beta  \xi_p = \phi_{2K,p} \beta ,$$
 therefore
 $$\Xi_{P,\beta} = \{\phi_{2K,p} \beta \mid \ p \in P \}.$$

  We summarize the discussion above as follows.

\begin{theorem} \label{cy2} Let $G$ be a finitely generated
fully residually free group, $S = 1$  a regular standard quadratic
orientable equation,  and   $\Gamma$   its basic sequence of
automorphisms. Then for  any  solution $\beta:G_S \rightarrow G $ in
general position, any positive integer $J \geq 2$,  and  any
unbounded set $P \subset {\mathbb N}^{JK}$  the set of
$G$-homomorphisms $\Xi_{P,\beta}$ $G$-discriminates  $G_{R(S)}$ into
$G$. Moreover, for any fixed tuple $p'\in {\mathbb N}^{tK}$ the
family
 $$\Xi_{P,\beta,p'} = \{\phi _{tK,p'}\theta \mid \theta\in\Xi_{P,\beta}\}$$
$G$-discriminates  $G_{R(S)}$ into $G$.
\end{theorem}

For tuples $f = (f_1, \ldots,f_k)$ and $g = (g_1, \ldots,g_m)$
 denote the tuple
  $$fg = (f_1, \ldots,f_k,g_1, \ldots,g_m).$$
Similarly, for a set of tuples $P$ put
 $$f Pg =\{fpg \mid p \in P\}.$$

\begin{cy} \label{cy:cy2}
Let $G$ be a finitely generated fully residually free group, $S = 1$
a regular standard quadratic orientable equation,    $\Gamma$ the
basic sequence of automorphisms of $S$, and $\beta:G_S \rightarrow G
$ a solution of $S=1$ in general position. Suppose $P \subseteq
\mathbb{N}^{2K}$ is unbounded set, and $f \in \mathbb{N}^{Ks}$, $g
\in \mathbb{N}^{Kr}$ for some $r,s \in \mathbb{N}$. Then there
exists a number $N$ such that if $f$ is $N$-large and $s\geq 2$ then
the family
 $$\Phi_{P,\beta,f,g} = \{\phi _{K(r+s+2),q} \beta  \mid q \in fPg\}$$
$G$-discriminates  $G_{R(S)}$ into $G$.
\end{cy}
 {\it Proof.} By Theorem \ref{cy2} it suffices to show that if
 $f$ is $N$-large for some $N$ then
 $\beta_{f} = \phi_{2K,f}\beta$ is a solution of $S = 1$ in general
 position, i.e., the images of some particular finitely many  non-commuting
 elements from $G_{R(S)}$ do not commute in $G$. It has been shown above  that the set
 of solutions $\{\phi_{2K,h}\beta \mid h \in \mathbb{N}^{2K}\}$ is
 a discriminating set for $G_{R(S)}$. Moreover, for any finite set
 $M$ of non-trivial elements from $G_{R(S)}$ there exists a number $N$
 such that for any $N$-large tuple $h \in \mathbb{N}^{2K}$ the
 solution $\phi_{2K,h}\beta$ discriminates all elements from $M$ into
 $G$. Hence the result.
  \hfill $\Box$

\subsection{Small cancellation solutions of standard orientable
equations} \label{subsec:7.2-b}

Let $S(X) = 1$ be  a standard regular orientable quadratic
equation over $F$ written in the form (\ref{6-general}):
 $$\prod_{i=1}^{m}z_i^{-1}c_iz_i\prod_{i=1}^{n}[x_i,y_i]  =
\prod_{i=1}^{m}e_i^{-1}c_ie_i \prod_{i=1}^{n}[a_i,b_i].$$ In this
section we construct solutions in $F$ of $S(X) = 1$ which satisfy
some small cancellation conditions.

\begin{df}
\label{de:smallcan} Let $S = 1$ be a standard regular orientable
quadratic equation written in the form (\ref{6-general}). We say
that a solution $\beta: F_S \rightarrow F$ of $S= 1$  satisfies
the small cancellation condition $(1/\lambda)$ with respect to the
set $\bar {\mathcal W}_{\Gamma }$ (resp. $\bar {\mathcal
W}_{\Gamma ,L}$) if the following conditions are satisfied:
  \begin{enumerate}
  \item [1)] $\beta$ is in general position;
  \item [2)] for any 2-letter word $uv \in \bar {\mathcal
W_{\Gamma }}$ (resp. $uv \in {\mathcal W}_{\Gamma ,L}$)(in the
alphabet $Y$) cancellation in the word $u^\beta v^\beta$ does not
exceed $(1/\lambda) \min\{|u^\beta|, |v^\beta|\}$ (we assume here
and below that $u^\beta, v^\beta$ are given by their reduced forms
in $F$);
 \item [3)]   cancellation in a word $u^{\beta} v^{\beta}$
 does not exceed $(1/\lambda) \min\{|u^{\beta}|,|v^\beta|\}$
 provided $u, v$ satisfy one of the conditions below:
  \begin{itemize}
   \item [a)] $u=z_i, v=(z_{i-1}^{-1}c_{i-1}^{-1}z_{i-1})$,
    \item [b)] $u=c_i, v=z_i$,
    \item [c)]  $u=v=c_i$,
\end{itemize}
 \end{enumerate}
 (we assume here that $u^\beta, v^\beta$ are given by their
reduced forms in $F$).
\end{df}

{ \bf Notation:}  For a homomorphism  $\beta : F[X] \rightarrow F$
by $C_\beta$ we denote the set of all  elements that cancel in
$u^\beta v^\beta$  where $u, v$ are as in 2), 3) from Definition
\ref{de:smallcan} and the word that cancels in the product
$(c_2^{z_2})^\beta \cdot (dc_{m-1}^{-z_{m-1}})^\beta$.

\begin{lm}
\label{lm:32} Let $u,v$ be cyclically reduced elements of $G\ast
H$ such that $|u|,|v|\geq 2$. If for some $m, n > 1$ elements
$u^{m}$ and $v^{n}$ have a common initial segment of length
$|u|+|v|$, then $u$ and $v$ are both powers of the same element
$w\in G \ast H$. In particular, if both $u$ and $v$ are not proper
powers then $u = v$.
\end{lm}
 {\it Proof.}  The same argument as in the case of free groups.

\begin{cy} If $u,v\in F,\ [u,v]\neq 1,$ then for any $\lambda \geq
0$ there exist $m_0, n_0$ such that for any $m \geq m_0, n \geq
n_0$ cancellation between $u^m$ and $v^n$ is less than
$\frac{1}{\lambda} max\{|u^m|,|v^n|\}.$\end{cy}

\begin{lm}
\label{le:7.1.beta} Let $S(X) = 1$ be  a standard regular orientable
quadratic equation written in the form (\ref{6-general}):

 $$\prod_{i=1}^{m}z_i^{-1}c_iz_i\prod_{i=1}^{n}[x_i,y_i]  =
\prod_{i=1}^{m}e_i^{-1}c_ie_i \prod_{i=1}^{n}[a_i,b_i], \ \ n\geq 1,
$$ where  all  $c_i$ are cyclically reduced. Then there exists a
solution $\beta$ of this equation that satisfies the small
cancellation condition with respect to $\bar {\mathcal W}_{\Gamma ,L
}$.  Moreover,  for any word $w\in \bar {\mathcal W}_{\Gamma ,L}$
that does not contain elementary squares, the word $w^{\beta}$ does
not contain a cyclically reduced part of $A_i^{2\beta}$ for any
elementary period $A_i$.

\end{lm}
{\it Proof.} We will begin with a solution
 $$\beta_1: x_i \rightarrow a_i, y_i \rightarrow b_i, z_i \rightarrow e_i $$
  of $S = 1$ in $F$ in
 general position. We will show that for any $\lambda \in \mathbb{N}$ there are positive
integers  $m_i, n_i, k_i, q_j$ and a tuple $p=(p_1,\ldots p_{m})$
such that  the map $\beta: F[X] \rightarrow F$
 defined by

$$x_1^{\beta}=(\tilde b_1^{n_1}\tilde a_1)^{[\tilde a_1,\tilde b_1]^{m_1}}, \ \ y_1^{\beta}=
((\tilde b_1^{n_1}\tilde a_1)^{k_1}\tilde b_1)^{[\tilde a_1,\tilde
b_1]^{m_1}}, \ \ where \ \tilde a_1=x_1^{\phi _{m}\beta_1},\ \tilde
b_1=y_1^{\phi
 _{m}\beta_1}$$

  $$x_i^{\beta}=(b_i^{n_i}a_i)^{[a_i,b_i]^{m_i}}, \ \ y_i^{\beta}=
((b_i^{n_i}a_i)^{k_i}b_i)^{[a_i,b_i]^{m_i}}, \ \ i=2,\ldots n,$$

 $$z_i^{\beta}=c_i^{q_i}z_i^{\phi _{m}\beta _1}, \ \ i=1,\ldots m,$$
 is a  solution of $S= 1$ satisfying the small cancellation condition
 $(1/\lambda)$ with respect to $\bar {\mathcal W}_{\Gamma }$.
Moreover, we will show that one can choose the solution $\beta_1$
such that $\beta$ satisfies the small cancellation condition with
respect to $\bar {\mathcal W}_{\Gamma ,L }$.

 The solution $\beta _1$ is in general position,
therefore the neighboring items in the sequence
  $$c_1^{e_1}, \ldots, c_m^{e_m}, [a_1,b_1],\dots ,[a_n,b_n]$$
   do not commute. We have $[c_i^{e_i},c_{i+1}^{e_{i+1}}]\neq 1.$

There is a homomorphism $\theta _{\beta _1}:F_S\rightarrow \bar
F=F(\bar U,\bar T)$ into the group $\bar F$ obtained from $F$ by a
series of extensions of centralizers, such that $\beta=\theta
_{\beta _1}\psi _p$, where $\psi _p:\bar F\rightarrow F$. This
homomorphism $\theta _{\beta _1}$ is a  monomorphism on $F\ast
F(z_1,\dots, z_m)$ (this follows from the proof of Theorem 4 in
\cite{KMNull}, where the same sequence of extensions of
centralizers is constructed).

 The set  of
solutions $\psi_p$ for different tuples $p$ and numbers
$m_i,n_i,k_i, q_j$ is a discriminating family for $\bar F.$ We
just have to show that the small cancellation condition for
$\beta$ is equivalent to a finite number of inequalities in the
group $\bar F$.

 We have  $z_i^\beta =c_i^{q_i}z_i^{\phi _{m}\beta
_1}$ such that $\beta _1(z_i)=e_i$, and $p=(p_1,\dots ,p_{m})$ is
a large tuple.  Denote $\bar A_j= A_j^{\beta _1},\ j=1,\dots ,m.$
Then it follows from Lemma \ref{le:7.1.zforms}
 that

\medskip
$z_i^{\beta}=c_i^{q_i+1}e_i\bar
A_{i-1}^{p_{i-1}}c_{i+1}^{e_{i+1}}\bar A_{i}^{p_{i}-1}$, where
$i=2,\dots ,m-1$

\medskip
$z_m^{\beta}=c_m^{q_m+1}e_m\bar A_{m-1}^{p_{m-1}}a_1^{-1}\bar
A_{m}^{p_{m}-1}$,

where

\medskip
$\bar A_1=c_1^{e_1}c_2^{e_2}$, $\bar A_2=\bar A_1(p_1)= \bar
A_1^{-p_1}c_2^{e_2}\bar A_1^{p_1}c_3^{e_3},\dots $,
\medskip
$\bar A_{i}=\bar A_{i-1}^{-p_{i-1}}c_i^{e_i}\bar
A_{i-1}^{p_{i-1}}c_{i+1}^{e_{i+1}},\ i=2,\dots ,m-1$,
\medskip
$\bar A_{m}=\bar A_{m-1}^{-p_{m-1}}c^{e_m}\bar
A_{m-1}^{p_{m-1}}a_1^{-1}.$

 One can choose $p$ such that
$[\bar A_{i},\bar A_{i+1}]\neq 1, [\bar
A_{i-1},c_{i+1}^{e_{i+1}}]\neq 1, [\bar A_{i-1},c_{i}^{e_{i}}]\neq
1$ and $[\bar A_{m},[a_1,b_1]]\neq 1$, because their pre-images do
not commute in $\bar F$. We need the second and third inequality
here to make sure that $\bar A_{i}$ does not end with a power of
$\bar A_{i-1}$. Alternatively, one can prove by induction on $i$
that $p$ can be chosen to satisfy these inequalities.

 Then $c_i^{z_i^{\beta}}$
and $c_{i+1}^{z_{i+1}^{\beta}}$ have small cancellation, and
$c_m^{z_m^{\beta}}$ has small cancellation with $x_1^{\pm\beta},
y_1^{\pm\beta}$.

    Let
   $$ x_i^{\beta}=(b_i^{n_i}a_i)^{[a_i,b_i]^{m_i}}, \ \
   y_i^{\beta}=((b_i^{n_i}a_i)^{k_i}b_i)^{[a_i,b_i]^{m_i}}, \ \ i=2,\dots ,n
$$

for some positive integers  $m_i, n_i, k_i, s_j$ which values we
will specify in a due course.  Let $uv \in\bar{\mathcal
W}_{\Gamma}$. There are several cases to consider.

1) $uv = x_ix_i$. Then
 $$u^\beta v^\beta = (b_i^{n_i}a_i)^{[a_i,b_i]^{m_i}}
 (b_i^{n_i}a_i)^{[a_i,b_i]^{m_i}}.$$
Observe that the cancellation  between  $(b_i^{n_i}a_i)$ and
$(b_i^{n_i}a_i)$ is not more then $|a_i|$. Hence the cancellation
in $u^\beta v^\beta$ is not more then $|[a_i,b_i]^{m_i}| + |a_i|$.
We chose $n_i >> m_i$ such that

$$|[a_i,b_i]^{m_i}| + |a_i| < (1/\lambda) |(b_i^{n_i}a_i)^{[a_i,b_i]^{m_i}}|$$
which is obviously possible. Similar arguments prove the cases $uv
= x_iy_i$ and $uv = y_ix_i.$

2) In all other cases the cancellation in
 $u^\beta v^\beta$ does not exceed the cancellation between $[a_i,b_i]^{m_i}$ and
 $[a_{i+1},b_{i+1}]^{m_{i+1}}$, hence by Lemma  \ref{lm:32} it is not greater than
 $|[a_i,b_i]| +  |[a_{i+1},b_{i+1}]|.$

Let $u=z_i^{\beta}, v=c_{i-1}^{-z_{i-1}^{\beta}}.$ The
cancellation is the same as between $\bar A_{2i}^{p_{2i}}$ and
$\bar A_{i-1}^{-p_{i-1}}$ and, therefore, small.

Since $c_i$ is cyclically reduced, there is no cancellation
between $c_i$ and $z_{i}^{\beta}$.

The first statement of the lemma is proved.

We now will prove the second statement of the lemma. We have to
show that if $u=c_i^{z_i}$ or $u=x_j^{-1}$ and $v=c_1^{z_1}$, then
the cancellation between $u^{\beta}$ and $v^{\beta}$ is less than
$(1/\lambda)min \{|u|,|v|\}.$
 We can  choose the  initial solution $e_1,\dots, e_m, a_1,b_1,\dots ,a_n,b_n$ so
   that $[c_1^{e_1}c_2^{e_2}, c_3^{e_3}\ldots c_i^{e_i}]\neq 1$ ( $i\geq
   3$),
   $[c_1^{e_1}c_2^{e_2},[a_i,b_i]]\neq 1, (i=2,\dots ,n)$ and
   $[c_1^{e_1}c_2^{e_2}, b_1^{-1}a_1^{-1}b_1]\neq 1.$
Indeed, the equations $[c_1^{z_1}c_2^{z_2}, c_3^{z_3}\ldots
   c_i^{z_i}]=1$, $[c_1^{z_1}c_2^{z_2},[x_i,y_i]]=1, (i=2,\dots ,n)$ and
   $[c_1^{z_1}c_2^{z_2}, y_1^{-1}x_1^{-1}y_1]=1$ are not  consequences of the equation $S=1$, and,
therefore, there is a solution of $S(X)=1$ which does not satisfy
 any of  these equations.

To show that $u=c_i^{z_i^{\beta}}$  and $v=c_1^{z_1^{\beta}}$,
have small cancellation, we have to show  that $p$ can be chosen
so that $[\bar A_1, \bar A_{i}]\neq 1$ (which is obvious, because
the pre-images in $\bar G$ do not commute), and that $\bar
A_{i}^{-1}$ does not begin with a power of $\bar A_1$. The period
$\bar A_{i}^{-1}$ has form ($c_{i+1}^{-z_{i+1}}\ldots
c_3^{-z_3}\bar A_1^{-p_2}\ldots ).$ It begins with a power of
$\bar A_1$ if and only if $[\bar A_1, c_3^{e_3}\ldots
c_i^{e_i}]=1$, but this equality does not hold.

Similarly one can show, that the cancellation between
 $u=x_j^{-\beta}$ and $v=c_1^{z_1^{\beta}}$ is small. \hfill $\Box$

\begin{lm} \label{sc2}

Let $S(X) = 1$ be  a standard regular orientable quadratic
equation of the type {\rm  (\ref{eq:st9})}
$$\prod _{i=1}^{m}z_i^{-1}c_iz_i=c_1^{e_1}\ldots c_m^{e_m}=d,$$
 where  all  $c_i$ are cyclically reduced,
 and
 $$\beta_1:  z_i \rightarrow e_i $$
  a solution of $S = 1$ in $F$ in
 general position. Then for any $\lambda \in \mathbb{N}$ there is a  positive
integer $s$   and a tuple $p=(p_1,\ldots p_{K})$ such that the map
$\beta: F[X] \rightarrow F$
 defined by
$$z_i^{\beta}=c_i^{q_i}z_i^{\phi _K\beta _1}d^s, $$
  is a  solution of $S= 1$ satisfying the small cancellation condition
 $(1/\lambda)$ with respect to $\bar {\mathcal W}_{\Gamma ,L}$ with
 one exception when $u = d$ and $v =c_{m-1}^{-z_{m-1}}$ (in this
 case $d$ cancels out in $v^\beta$). Notice, however, that such
 word $uv$ occurs only in the product $wuv$ with  $w = c_2^{z_2}$,
 in which case cancellation
between $w^{\beta}$ and $dv^{\beta}$ is less than
$\min\{|w^{\beta}|,|dv^{\beta}|\}.$ Moreover,  for any word $w\in
\bar {\mathcal W}_{\Gamma ,L}$ that does not contain elementary
squares, the word $w^{\beta}$ does not contain a cyclically reduced
part of $A_i^{2\beta}$ for any elementary period $A_i$.
\end{lm}
{\it Proof.}  Solution $\beta$ is chosen   the same way as in the
previous lemma (except for the multiplication by  $d^s$) on the
elements $z_i,\ i\neq m$. We do not take $s$ very large, we just
need it to avoid cancellation between $z_2^{\beta}$ and $d$.
Therefore the cancellation between $c_i^{z_i^{\beta}}$ and
$c_{i+1}^{\pm z_{i+1}^{\beta}}$ is small for $i<m-1$. Similarly,
for $u=c_2^{z_2},\  v=d,\ w=c_{m-1}^{-z_{m-1}},$ we can make the
cancellation between $u^{\beta}$ and $dw^{\beta}$  less than
$min\{|u^{\beta}|,|dw^{\beta}|\}.$ \hfill $\Box$

\begin{lm}
\label{le:small-cancellation-beta-double} Let $U, V \in {\mathcal
W}_{\Gamma ,L}$ such that $UV = U \circ V$ and $UV \in {\mathcal
W}_{\Gamma ,L}$.

1. Let $n\neq 0.$ If $u$ is the last letter of $U$ and $v$ is the
first letter of $V$ then cancellation between  $U^\beta$ and
$V^\beta$ is equal to the cancellation between  $u^\beta$ and
$v^\beta$.

2. Let $n=0.$ If $u_1u_2$ are the last two letters of $U$ and
$v_1,v_2$ are the first two letters of $V$ then cancellation
between $U^\beta$ and $V^\beta$ is equal to the cancellation
between $(u_1u_2)^\beta$ and $(v_1v_2)^\beta$.

  \end{lm}
Since  $\beta$ has the small
  cancellation property with respect to
  $\bar {\mathcal W}_{\Gamma ,L},$
  this implies that the cancellation in $U^\beta V^{\beta}$ is equal to the cancellation in
  $u^\beta v^{\beta}$, which is equal to some element in
  $C_\beta$.
This proves the lemma.
 \hfill$\Box$

Let $w\in \bar{\mathcal W}_{\Gamma ,L}, \phi _j = \phi _{j,p},
W=w^{\phi _j}, $ and $A = A_j$.
  \begin{equation}
\label{eq:Adecomp3} W = B_1 \circ A^{q_1} \circ \ldots \circ B_k
\circ A^{q_k} \circ B_{k+1}
\end{equation}
 the canonical $N$-large $A$-representation of $W$ for some
positive integer $N$.

Since the occurrences $A^{q_i}$ above are stable we have $$B_1 =
{\bar B}_1 \circ A^{sgn(q_1)}, \ \ B_i = A^{sgn(q_{i-1})} \circ
{\bar B}_i \circ A^{sgn(q_i)} \ (2 \leq i \leq k), \ \ B_{k+1} =
A^{sgn(q_k)} \circ {\bar B}_{k+1}.$$ Denote
$A^{\beta}=c^{-1}A^{\prime} c,$ where $A^{\prime}$ is cyclically
reduced, and $c\in C_{\beta}.$ Then
$$B_1^\beta = {\bar B}_1^\beta c^{-1} (A^\prime)^{sgn(q_1)}c, \ \
B_i^\beta = c^{-1}(A^\prime)^{sgn(q_{i-1})} c {\bar B}_i^\beta
c^{-1} (A^\prime)^{sgn(q_i)}c, \ \ B_{k+1}^\beta  =c^{-1}
(A^\prime)^{sgn(q_k)} c{\bar B}_{k+1}^\beta.$$ By Lemma
\ref{le:small-cancellation-beta-double} we can assume that the
cancellation in the words above is small, i.e., it does not exceed
a fixed number $\sigma$ which is the maximum length of words from
$C_\beta$. To get an $N$-large canonical $A^\prime$-decomposition
of $W^\beta$ one has to take into account stable occurrences of
$A^\prime$. To this end, put $\varepsilon_i = 0$  if
${A^\prime}^{sgn(q_i)}$ occurs in the reduced form of ${\bar
B}_i^\beta c^{-1}(A^\prime)^{sgn(q_i)}$ as written (the
cancellation does not touch it), and put $\varepsilon_i =
sgn(q_i)$
 otherwise. Similarly, put $\delta_i = 0$ if
${A^\prime}^{sgn(q_i)}$ occurs in the reduced form of
$(A^\prime)^{sgn(q_{i})} c {\bar B}_{i+1}^\beta$ as written, and
put $\delta_i = sgn(q_{i})$ otherwise.

Now one can rewrite $W^\beta$ in the following form

\begin{equation}
\label{eq:7.3.can}
 W^\beta  = E_1 \circ
   (A^\prime)^{q_1 - \varepsilon_1 -\delta_1}
   \circ
  E_2
  \circ (A^\prime)^{q_2 - \varepsilon_2 -\delta_2} \circ \ldots
 \circ (A^{\prime})^{q_k -\varepsilon_k -\delta_k}  \circ E_{k+1},
 \end{equation}
 where $E_1=(B_1^\beta c^{-1}(A^\prime)^{\varepsilon_1}),\
 E_2=((A^\prime)^{\delta_1}cB_2^\beta
 c^{-1}(A^\prime)^{\varepsilon_2}),\ E_{k+1}=((A^\prime)^{\delta_{k}}
   cB_{k+1}^\beta).$

   Observe, that $d_i$ and $\varepsilon_i, \delta_i$ can be effectively computed from
 $W$ and $\beta$. It follows that one can effectively rewrite  $W^\beta$
 in the form  (\ref{eq:7.3.can}) and the form is unique.

The decomposition  (\ref{eq:7.3.can}) of  $W^\beta$ induces the
corresponding  $A^*$-decomposition of  $W$. This can be shown by an
argument similar to the one in Lemmas \ref{elper} and
\ref{le:A-A^*}, where it has been proven  that $A_{r+L}^*$-
decomposition induces the corresponding $A_r$-decomposition. To see
that  the argument works we need the last statement in Lemmas
\ref{le:7.1.beta} ($n>0$)  and \ref{sc2} ($n = 0$)  which ensure
that the "illegal" elementary squares do not occur because of the
choice of the solution $\beta$.

If the canonical $N$-large $A^*$-decomposition of $W$ has the form:
$$
D_1 (A^*)^{q_1} D_2 \ldots D_k(A^*)^{q_k} D_{k+1}$$ then the induced
one has the form:
\begin{equation}
\label{eq:induced} W = (D_1A^{* \varepsilon_1}) A^{* q_1 -
\varepsilon_1 - \delta_1} (A^{*  \delta_1} D_2 A^{*\varepsilon_2})
\ldots (A^{*  \delta_{k-1}} D_k A^{*\varepsilon_k}) A^{* q_k -
\varepsilon_k - \delta_k} (A^{*\delta_k}D_{k+1}).\end{equation}

We call this decomposition the {\it induced} $A^*$-decomposition
of $W$ with respect to $\beta$ and write it in the form:
\begin{equation}
\label{eq:induced*} W = D_1^* {(A^*)}^{q_1^*} D_2^*
 \ldots D_k^* {(A^*)}^{q_k^*}D_{k+1}^*,\end{equation}
where $D_i^* = {(A^*)}^{\delta_{i-1}} D_i
{(A^*)}^{\varepsilon_i}$, $q_i^* = q_i - \varepsilon_i -
\delta_i$, and, for uniformity, $\delta_1 = 0$ and
$\varepsilon_{k+1} = 0$.

\begin{lm}
\label{le:D-star-beta} For given positive integers $j$, $N$ and a
real number $\varepsilon > 0$ there is a constant  $C =
C(j,\varepsilon,N)
>0$ such that if $p_{t+1}-p_t
> C$ for every $t = 1, \ldots, j-1,$ and a word $W\in\bar{\mathcal
W}_{\Gamma ,L}$ has a canonical
 $N$-large $A^*$-decomposition (\ref{eq:induced*}), then this decomposition satisfies the following conditions:
\begin{equation} \label{eq:E-D} (D_1^*)^\beta = E_1
\circ_{\theta} (c R^\beta), \ \ (D_i^*)^\beta = (R^{-\beta}c^{-1})
\circ_{\theta} E_i \circ_{\theta} (cR^\beta), \ \
(D_{k+1}^*)^\beta = (R^{-\beta}c^{-1}) \circ_{\theta} E_{k+1},
\end{equation}
where $\theta  < \varepsilon |A'|$. Moreover, this constant $C$ can
be found effectively.
\end{lm}

{\it Proof.} Applying homomorphism $\beta$ to the reduced
 $A^*$-decomposition of $W$ (\ref{eq:induced*}) we can see
 that
 $$
 W^\beta = \left ((D_1^*)^\beta R^\beta c \right )
 {(A^\prime)}^{q_1^*} \left (c R^{\beta}(D_2^*)^\beta
 R^{-\beta}c^{-1}\right ){(A^\prime)}^{q_2^*} \ldots \left (cR^\beta
 (D_k^*)^\beta R^{-\beta}c^{-1}\right ){(A^\prime)}^{q_k^*} \left
 (cR^\beta (D_{k+1}^*)^\beta\right ) $$
 Observe that this decomposition has the same powers of $A^\prime$
 as the canonical $N$-large $A^\prime$-decomposition
 (\ref{eq:7.3.can}). From the uniqueness of such decompositions we
 deduce that
$$E_1 = (D_1^*)^\beta c^{-1} R^{-\beta}, \ \ E_i = c R^{\beta}(D_i^*)^\beta
 R^{-\beta}c^{-1}, \ \ E_{k+1} = cR^\beta (D_{k+1}^*)^\beta$$
 Put $\theta = |c| +|R^\beta|$. Rewriting the equalities above one can get
 $$
(D_1^*)^\beta = E_1 \circ_{\theta} (c R^\beta), \ \ (D_i^*)^\beta =
(R^{-\beta}c^{-1}) \circ_{\theta} E_i \circ_{\theta} (cR^\beta), \ \
(D_{k+1}^*)^\beta = (R^{-\beta}c^{-1}) \circ_{\theta} E_{k+1}.$$
Indeed,  in the decomposition (\ref{eq:7.3.can}) every occurrence
$(A^\prime)^{q_i - \varepsilon_i -\delta_i}$ is stable hence $E_i$
starts (ends) on $A^\prime$. The rank of $R$ is at most $rank(A)
-K+2$, and $\beta$ has small cancellation. Taking $p_{i+1}>>p_i$ we
obtain  $\varepsilon |A^\prime| > |c|+ |R^\beta|$. \hfill$\Box$

 Notice, that one can
effectively write down the induced $A^*$-decomposition of $W$ with
respect to $\beta$.

We summarize the discussion above in the following statement.

\begin{lm}\label{Claim1.b} For given positive integers $j$,  $N$
there is a constant  $C =  C(j,N)$  such that if $p_{t+1}-p_t >
C$, for  every $t = 1, \ldots, j-1,$ then for any   $W\in
\bar{\mathcal W}_{\Gamma , L}$    the following conditions are
equivalent:

\begin{enumerate}
\item  Decomposition {\rm (\ref{eq:Adecomp3})} is the canonical
(the canonical $N$-large) $A$-decomposition of $W$,

\item Decomposition {\rm (\ref{eq:7.3.can})} is the canonical (the
canonical $N$-large) $A^\prime$-decomposition of $W^\beta$,

\item  Decomposition {\rm (\ref{eq:induced})} is the canonical
(the canonical $N$-large) $A^*$-decomposition of $W.$

\end{enumerate}
\end{lm}

\subsection{Implicit Function Theorem for Quadratic
Equations}\label{se:7.3} In this section we prove Theorem
\ref{1,2,3,4} for orientable quadratic equations over a free group
$F = F(A)$. Namely, we prove the following statement.

\medskip
 {\it Let $S(X,A)=1$ be a regular standard orientable
quadratic equation over $F$. Then every equation $T(X,Y,A) = 1$
compatible with $S(X,A) = 1$ admits an effective complete
$S$-lift.}

 \medskip
{\bf A special discriminating set of solutions
 ${\mathcal L}$ and the corresponding   cut equation $\Pi$.}

 \medskip
Below we continue to use notations from the previous sections. Fix
a  solution $\beta$ of $S(X,A) = 1$ which satisfies the
cancellation condition $(1/\lambda)$ (with $\lambda
> 10$) with respect to $\bar {\mathcal W}_{\Gamma}$.

 Put
$$x_i^{\beta}=\tilde a_i,y_i^{\beta}=\tilde b_i,z_i^{\beta}=\tilde c_i. $$
 Recall that $$\phi_{j,p} = \gamma_j^{p_j} \ldots \gamma_1^{p_1} = \stackrel{\leftarrow}{\Gamma}_j^p$$
 where $j \in {\mathbb N}$,  $\Gamma_j = (\gamma_1,
 \ldots, \gamma_j)$ is the initial subsequence of length $j$ of the
sequence
 $\Gamma^{(\infty)}$, and $p = (p_1, \ldots,p_j) \in {\mathbb N}^j$.
  Denote by $\psi_{j,p}$ the following solution of $S(X) = 1$:
  $$\psi_{j,p}=\phi_{j,p}\beta .$$ Sometimes we
 omit $p$ in $\phi_{j,p}, \psi_{j,p}$ and simply write
$\phi_j, \psi_j$.

Below we continue to use notation:
  $$A = A_j,  \ A^* =
A^*_j = A^*(\phi_j) = R_j^{-1} \circ A_j \circ R_j, \ d = d_j =
|R_j|.$$
 Recall that $R_j$ has rank $\leq j-K +2$ (Lemma
\ref{le:cyclicAA*}).  By $A^\prime$ we denote the cyclically
reduced form of $A^\beta$ (hence of $(A^*)^\beta$).
  Recall that the set $C_\beta$ was defined right after Definition \ref{de:smallcan}.

Let
 $$ \Phi = \{\phi_{j,p} \mid j
\in {\mathbb N}, p \in {\mathbb N}^j \}.$$
  For    an arbitrary  subset ${\mathcal L}$ of $\Phi$ denote
  $${\mathcal L}^\beta = \{\phi \beta \mid \phi \in {\mathcal L}\}.$$

Specifying step by step various subsets of $\Phi$ we will
eventually ensure a very particular choice of a set of solutions
of $S(X) = 1$ in $F$.

Let $K = K(m,n)$ and $J \in \mathbb{N}, J \geq 3,$  a sufficiently
large positive integer which will be specified precisely in  due
course.  Put $L = JK$ and define ${\cal P}_1 = \mathbb{N}^L$,
 $${\mathcal L_1} = \{\phi_{L,p} \mid  p \in {\cal P}_1 \}.$$
By Theorem \ref{cy2} the set ${\mathcal L_1}^\beta$  is a
discriminating set of solutions of $S(X) = 1$ in $F$. In fact, one
can replace the set ${\cal P}_1$ in the definition of ${\mathcal
L_1}$ by any unbounded subset ${\cal P}_2 \subseteq {\cal P}_1$, so
that the new set is still discriminating. Now we construct by
induction a very particular unbounded subset ${\cal P}_2 \subseteq
\mathbb{N}^L$.
 Let $a \in  \mathbb{N}$ be a natural number and $h: \mathbb{N} \times \mathbb{N} \rightarrow \mathbb{N}$ a  function. Define a tuple
  $$p^{(0)} = (p^{(0)}_1, \ldots, p^{(0)}_L)$$
where
 $$ p^{(0)}_1 = a, \ \ p^{(0)}_{j+1} = p^{(0)}_j + h(0,j).$$
Similarly, if a tuple $p^{(i)} = (p^{(i)}_1, \ldots, p^{(i)}_L)$
is defined then put $p^{(i+1)} = (p^{(i+1)}_1, \ldots,
p^{(i+1)}_L)$, where
  $$ p^{(i+1)}_1 = p^{(i)}_1 + h(i+1,0), \ \ p^{(i+1)}_{j+1} = p^{(i+1)}_j + h(i+1,j).$$
This defines by induction an infinite set
 $${\cal P}_{a,h} = \{p^{(i)} \mid i \in \mathbb{N}\} \subseteq
\mathbb{N}^L$$
 such that  any infinite subset of ${\cal P}_{a,h}$ is also unbounded.

From now on fix a recursive monotonically increasing with respect
to both variables function $h$ (which will be specified in due
course) and put
$${\cal P}_2 =  {\cal P}_{a,h}, \ \ \ {\mathcal L_2} = \{\phi_{L,p} \mid  p \in {\cal P}_2 \}.$$

\begin{prop} \label{pr:disc-p}
Let $r \geq 2$ and $K(r+2) \leq L$ then there exists a number $a_0$
such that if $a \geq a_0$ and the function $h$ satisfies the
condition
 \begin{equation}
 \label{eq:h}
 h(i+1,j) > h(i,j) \ \ \text{for any} \  j = Kr+1,
 \ldots, K(r+2), i = 1,2, \ldots;
\end{equation}
then for any infinite subset ${\cal P} \subseteq {\cal P}_2$ the set
of solutions
 $${\mathcal L_{\cal P}}^\beta = \{\phi_{L,p}\beta \mid  p \in {\cal P} \}$$
is a discriminating set of solutions of $S(X,A) = 1$.
\end{prop}
 {\it Proof.} The result  follows from Corollary \ref{cy:cy2}.
 \hfill $\Box$

   Let  $\psi  \in {\mathcal L}_2^\beta$. Denote by $U_\psi$ the solution
  $X^{\psi}$ of the equation $S(X) = 1$ in $F$.  Since  $T(X,Y) = 1$
   is compatible with $S(X) = 1$ in $F$ the equation
$T(U_\psi,Y) = 1$ (in variables $Y$)  has a solution in $F$, say
$Y = V_\psi$.  Set
  $$\Lambda = \{(U_\psi, V_\psi) \mid \psi \in {\mathcal
  L}_2^\beta\}.$$
  It follows that every pair
$(U_\psi,V_\psi) \in \Lambda $ gives a solution of  the system
    $$ R(X,Y) = (S(X) = 1 \  \wedge \ T(X,Y) = 1).$$
      By Theorem
 \ref{th:cut} there exists a finite set ${\mathcal CE}(R)$ of cut equations
  which describes all solutions of $R(X,Y) = 1$ in $F$, therefore there exists a cut equation
$ \Pi_{\mathcal L_3, \Lambda}  \in {\mathcal CE}(R)$ and an
infinite subset ${\mathcal L}_3 \subseteq {\mathcal L_2}$ such
that $\Pi_{\mathcal L_3, \Lambda}$ describes all solutions of the
type $(U_\psi, V_\psi)$, where  $\psi \in  {\mathcal L}_3$.
 We state the precise formulation of this result  in the following
 proposition which, as we have mentioned already,  follows from Theorem \ref{th:cut}.

\begin{prop}
\label{prop:cut-L3}
 Let ${\mathcal L_2}$ and $\Lambda$ be as above. Then
 there exists  an infinite subset ${\cal P}_3 \subseteq {\cal P}_2$ and the corresponding set
 ${\mathcal L}_3  = \{\phi_{L,p} \mid p \in {\cal P}_3\} \subseteq
{\mathcal L_2}$,   a cut equation $\Pi_{\mathcal L_3, \Lambda} =
({\mathcal E}, f_X, f_M) \in {\mathcal CE}(R)$, and a
  tuple of words $Q(M)$ such that the following  conditions hold:
 \begin{enumerate}
 \item [1)] $f_X({\mathcal E}) \subset X^{\pm 1}$;
 \item [2)]  for every $\psi \in {\mathcal L}_3^\beta$  there exists a
tuple of words $P_\psi = P_\psi(M)$ and a solution $\alpha_\psi: M
\rightarrow F$ of  $\Pi_{\mathcal L_3, \Lambda}$ with respect to
$\psi: F[X] \rightarrow F$ such that:
 \begin{itemize}
  \item the solution  $U_\psi = X^\psi$ of $S(X) =1$ can be
presented as $U_\psi = Q(M^{\alpha_\psi})$ and the word
$Q(M^{\alpha_\psi})$ is reduced as written,
 \item $V_\psi = P_\psi(M^{\alpha_\psi})$.
 \end{itemize}
 \item [3)] there exists a tuple of words
$P$ such that for any solution (any group solution) $(\beta,
\alpha)$ of $\Pi_{\mathcal L_3, \Lambda}$  the pair $(U,V),$ where
$U = Q(M^\alpha)$ and $V = P(M^\alpha),$ is a solution of $R(X,Y)
= 1$ in $F$.
 \end{enumerate}
\end{prop}

Put
 $${\cal P} = {\cal P}_3, \ \ \ {\mathcal L} = {\mathcal L}_3, \ \ \  \Pi_{\mathcal L} = \Pi_{\mathcal L_3, \Lambda}.$$
   By Proposition
\ref{pr:disc-p}
 the set ${\mathcal L}^\beta$  is a discriminating
set of solutions of $S(X) = 1$ in $F$.

 \medskip
 {\bf The initial cut equation $\Pi_{\phi}$.}

 \medskip
 Now fix a tuple $p \in {\cal P}$ and the  automorphism $\phi= \phi_{L,p}  \in {\mathcal L}$.
Recall,  that for every $j \leq L$ the automorphism $\phi_j$ is
defined by $\phi_j = \stackrel{\leftarrow}{\Gamma}_j^{p_j}$, where
$p_j$ is the initial subsequence of $p$ of length $j$. Sometimes
we use  notation $\psi = \phi \beta, \psi_j = \phi_j  \beta$.

Starting with the cut equation $\Pi_{\mathcal L}$ we construct a
cut equation $\Pi_{\phi} = ( {\mathcal E}, f_{\phi,X}, f_M) $
which  is obtained from $\Pi_{\mathcal L}$  by replacing the
function $f_X: {\mathcal E} \rightarrow F[X] $ by a new function
$f_{\phi,X}: {\mathcal E} \rightarrow F[X]$, where $f_{\phi,X}$ is
the composition of $f_X$ and the automorphism $\phi$. In other
words, if an interval $e \in {\mathcal E}$ in $\Pi_{\mathcal L}$
has a label $x \in X^{\pm 1}$ then its label in $\Pi_{\phi}$ is
$x^\phi$.

Notice, that $\Pi_{\mathcal L}$ and $\Pi_{\phi}$ satisfy  the
following conditions:
 \begin{enumerate}
  \item [a)]  $\sigma ^{f_X \phi \beta} = \sigma ^{f_{\phi,X}
  \beta}$ for every $ \sigma \in {\mathcal E}$;
  \item [b)] the solution  of $\Pi_{\mathcal L}$ with respect to
  $\phi  \beta$ is also a solution of $\Pi_{\phi}$ with respect to
  $\beta$;
  \item [c)] any solution (any group solution) of $\Pi_\phi$ with
  respect to $\beta$ is a solution (a group solution) of $\Pi_{\mathcal L}$
  with respect to $\phi \beta$.
  \end{enumerate}

 The cut equation $\Pi_\phi$ has a very particular type. To deal with such cut
 equations we need the following definitions.

\begin{df}
Let $\Pi = ({\mathcal E}, f_X, f_M)$ be a cut equation. Then the
number $$length(\Pi) = max \{|f_M(\sigma)| \mid \sigma \in
{\mathcal E}\}$$ is called  the length of $\Pi$. We denote it by
$length(\Pi)$  or simply by $N_{\Pi}$.
\end{df}

Notice,  by construction,  $length(\Pi_\phi) =
length(\Pi_{\phi^\prime})$ for every $\phi, \phi^\prime \in
{\mathcal L}$. Denote
 $$N_{\mathcal L} = length(\Pi_\phi).$$

\begin{df}
\label{df:7.3.cutgamma} A  cut equation $\Pi = ( {\mathcal E},
f_{X}, f_M)$ is called a $\Gamma$-cut equation in {\em  rank} $j$
($rank(\Pi) = j$) and size $l$ if it satisfies the following
conditions:
\begin{enumerate}
\item [1)] let $W_\sigma = f_X(\sigma)$ for $\sigma \in {\mathcal
E}$ and  $N = (l+2)N_{\Pi}$.
  Then for every $\sigma \in {\mathcal E}$ $W_\sigma \in \bar {\mathcal W}_{\Gamma,L}$
  and one of the following conditions
  holds:
\begin{enumerate}
  \item [1.1)] $W_\sigma$ has $N$-large rank $j$ and its canonical
 $N$-large  $A_j$-decomposition has  size
$(N,2)$ i.e., $W_\sigma$ has the canonical $N$-large
$A_j$-decomposition
\begin{equation}
 \label{eq:Adecomp-1}
 W_{\sigma} = B_1 \circ A^{q_1}_j \circ \ldots B_k \circ A^{q_k}_j \circ B_{k+1},
 \end{equation}
with   $max_j(B_i) \leq 2$ and $q_i \geq N$;
 \item  [1.2)] $W_\sigma$ has rank $j$ and $\max_j(W_\sigma) \leq 2$;
\item [1.3)]  $W_\sigma$ has rank $< j$.
\end{enumerate}
Moreover,  there exists  at least one interval $\sigma \in
{\mathcal E}$ satisfying the condition 1.1).
 \item [2)] there exists a solution $\alpha : F[M]
\rightarrow F$ of the cut equation $\Pi$ with respect to the
homomorphism   $\beta:F[X] \rightarrow F$.
\end{enumerate}
\end{df}

\begin{lm} Let $l\geq 3$.
The cut equation $\Pi_\phi$ is a $\Gamma$-cut equation in rank $L$
and size $l$, provided  $$ p_L \geq (l+2)N_{\Pi_\phi}+3.$$
\end{lm}
{\it Proof.} By construction the labels of intervals from
$\Pi_\phi$ are precisely the words of the type $x^{\phi_L}$ and
every such word appears as a label. Observe, that
$rank(x_i^{\phi_L}) < L$ for every $i, 1 \leq i \leq n$ (Lemma
\ref{le:7.1.gammawords}, 1a). Similarly, $rank(x_i^{\phi_L}) < L$
for every $i < n$ and $rank(y_n^{\phi_L}) = L$ (Lemma
\ref{le:7.1.gammawords} 1b). Also, $rank(z_i^{\phi_L}) < L$ unless
$n = 0$ and  $i = m$, in the latter case $z_m^{\phi_L}) = L$ (Lemma
\ref{le:7.1.gammawords} 1c and 1d).
  Now consider the labels $y_n^{\phi_L}$ and $z_m^{\phi_L})$ (in the case $n =
  0$) of  rank $L$.  Again, it has been shown in Lemma
  \ref{le:7.1.gammawords} 1) that these labels have
 $N$-large  $A_L$-decompositions of size $(N,2)$, as required in 1.1)
 of the definition of a $\Gamma$-cut equation of rank $L$ and size
 $l$.

\hfill$\Box$

\medskip
 {\bf Agreement 1 on ${\cal P}.$} Fix an arbitrary integer $l$, $l \geq 5$.
 We may assume, choosing the
 constant $a$ to satisfy the condition
   $$a \geq (l+2)N_{\Pi_\phi}+3,$$
   that all tuples in the set ${\cal P}$ are
$[(l+2)N_{\Pi_\phi}+3]$-large. Denote $N = (l+2)N_{\Pi_\phi}.$

Now we introduce one  technical restriction on the set ${\cal P}$,
its real meaning will be clarified later.

\medskip
 {\bf Agreement 2 on ${\cal P}.$}
Let $r$ be an arbitrary fixed positive integer with $Kr \leq L$ and
$q$ be a fixed tuple of length $Kr$ which is an initial segment of
some tuple from
 ${\cal P}$. The choice of $r$ and $q$ will be clarified later.
  We may assume (suitably choosing the  function
  $h$) that all tuples from ${\cal P}$
 have $q$ as their initial segment. Indeed, it suffices to define
 $h(i,0) = 0$ and $h(i,j) = h(i+1,j)$ for all $i \in \mathbb{N}$ and
 $j = 1, \ldots, Kr$.

\medskip
 {\bf Agreement 3 on ${\cal P}.$} Let $r$ be the number from  Agreement
 2. By Proposition \ref{pr:disc-p} there exists a number $a_0$ such
 that for every infinite subset of ${\cal P}$ the corresponding set
 of solutions is a discriminating set. We may assume that $a > a_0$.

\medskip
{\bf Transformation $T^*$ of $\Gamma$-cut equations.}

 \medskip
 Now we describe a  transformation $T^*$ defined on
  $\Gamma$-cut
 equations and their solutions,  namely, given a $\Gamma$-cut equation $\Pi$ and its solution
 $\alpha$ (relative to the  fixed map $\beta:F[X] \rightarrow F$ defined above) $T^*$
 transforms $\Pi$ into a new $\Gamma$-cut equation $\Pi^* = T^*(\Pi)$ and $\alpha$ into a solution
 $\alpha^* = T^*(\alpha)$ of $T^*(\Pi)$ relative to $\beta$.

 Let $\Pi = ( {\mathcal E}, f_{X}, f_M)$ be a $\Gamma$-cut equation in rank $j$ and size
 $l$. The cut equation
 $$T^*(\Pi) = ( {\mathcal E^*}, f^*_{X^*}, f^*_{M^*})$$
 is defined as follows.

\medskip
{\bf Definition of the set ${\mathcal E}^*$.}

\medskip
 For $\sigma \in {\mathcal E}$ we denote $W_{\sigma} = f_X(\sigma)$. Put
 $${\mathcal E}_{j,N} = \{\sigma \in {\mathcal E} \mid W_\sigma \ \mbox{satisfies} \ 1.1) \}.$$
 Then ${\mathcal E} = {\mathcal E}_{j,N} \cup {\mathcal E}_{< j,N}$ where ${\mathcal E}_{< j,N}$ is
 the complement
 of $ {\mathcal E}_{j,N}$ in ${\mathcal E}$.

 Now let $\sigma \in  {\mathcal E}_{j,N}$.
Write the  word $W_\sigma^\beta$ in its canonical $A^\prime$
decomposition:
\begin{equation}
\label{eq:can-W-sigma-beta}
 W_\sigma^\beta = E_1 \circ {A^\prime}^{q_1} \circ E_2 \circ \cdots \circ E_k \circ
 {A^\prime}^{q_k} \circ E_{k+1}
 \end{equation}
 where $|q_i| \geq 1$, $E_i \neq 1$ for $2 \leq i \leq k$.

 Consider the partition $$f_M(\sigma)  = \mu_1 \ldots \mu_n$$ of $\sigma$.
  By the condition 2) of the definition of $\Gamma$-cut equations
   for the solution  $\beta:F[X] \rightarrow F$ there exists a solution
  $\alpha : F[M] \rightarrow F$ of the cut equation $\Pi$ relative to $\beta$. Hence
 $W_\sigma^\beta = f_M(M^\alpha)$ and the element
 $$f_M(M^\alpha) = \mu_1^\alpha \ldots \mu_n^\alpha$$ is reduced as written.
   It  follows that
   \begin{equation}
   \label{eq:7.3.13}  W_\sigma^\beta = E_1 \circ {A^\prime}^{q_1} \circ E_2 \circ \cdots  \circ E_k \circ
 {A^\prime}^{q_k} \circ E_{k+1}  = \mu_1^{\alpha} \circ  \cdots \circ \mu_n^{\alpha}.
\end{equation}

   We say that a variable $\mu_i$ is {\em long} if ${A^\prime}^{\pm (l+2)}$  occurs in
   $\mu_i^\alpha$ (i.e., $\mu_i^\alpha$ contains a stable occurrence of
   ${A^\prime}^{l}$),   otherwise it is called {\em short}.
   Observe, that the definition of long  (short) variables $\mu
   \in M$ does not depend on a choice of $\sigma$,  it depends only
   on the given homomorphism $\alpha$.
The graphical  equalities  (\ref{eq:7.3.13}) (when $\sigma$ runs
over ${\mathcal E}_{j,N}$)  allow one to effectively recognize long
and short variables in $M$. Moreover, since for every $\sigma \in
{\mathcal E}$ the length of the word $f_M(\sigma)$ is bounded by
$length(\Pi) = N_{\Pi}$ and $N = (l+2)N_{\Pi}$, every word
$f_M(\sigma)$ ($\sigma \in {\mathcal E}_{j,N}$) contains long
variables. Denote by $M_{\rm short}$, $M_{long}$ the sets of short
and long variables in $M$. Thus, $M = M_{\rm short} \cup M_{long}$
is a non-trivial partition of $M$.

Now we define the following property $P = P_{long,l}$ of
occurrences of powers of $A^\prime$ in $W_\sigma^\beta$: a given
stable occurrence ${A^\prime}^q$ satisfies $P$ if it occurs in
$\mu^\alpha$ for some long variable $\mu \in M_{long}$ and $q \geq
l$. It is easy to see that $P$ preserves correct overlappings.
Consider the set of stable occurrences ${\mathcal O}_{P}$ which
are maximal with respect to $P$. As we have mentioned already in
Section \ref{se:7.2.5}, occurrences from ${\mathcal O}_{P}$ are
pair-wise disjoint and this set is uniquely defined. Moreover,
$W_\sigma^\beta$ admits the  unique $A^\prime$-decomposition
relative  to the set ${\mathcal O}_{P}$:
 \begin{equation}
  \label{eq:7.3.can*}
   W_\sigma ^\beta = D_1 \circ (A^\prime)^{q_1}
   \circ D_2 \circ \cdots \circ D_k
 \circ (A^{\prime})^{q_k}  \circ
   D_{k+1},
  \end{equation}
where $D_i \neq 1$ for $i = 2, \ldots, k$. See Figure 13.
\begin{figure}[here]
\centering{\mbox{\psfig{figure=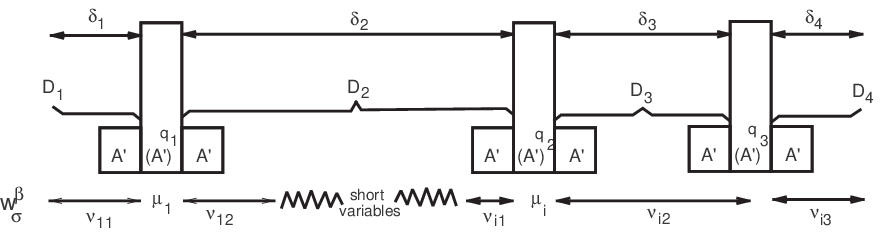,height=1.5in}}}
\caption{ Decomposition
(\ref{eq:7.3.can*})}\label{picture-decomposition-997}
\end{figure}

 Denote by $k(\sigma)$
the number of nontrivial elements among $D_1, \ldots, D_{k+1}$.

According to Lemma \ref{Claim1.b}  the $A^\prime$-decomposition
\ref{eq:7.3.can*} gives rise to the  unique associated
$A$-decomposition of $W_\sigma$ and hence the  unique associated
$A^*$-decomposition of $W_\sigma$.

Now with a given  $\sigma \in {\mathcal E}_{j,N}$ we
  associate a finite set of new intervals $E_\sigma$ (of the equation
  $T^*(\Pi)$):
  $$E_\sigma = \{\delta_1, \ldots , \delta_{k(\sigma)}\}$$
   and put
 $${\mathcal E}^* = {\mathcal E}_{<j,N} \cup \bigcup_{\sigma \in {\mathcal E}_{j,N}} E_\sigma.$$

\medskip
{\bf Definition of the set $M^*$}

 \medskip
 Let $\mu \in M_{long}$ and
 \begin{equation}
 \label{eq:7.3.14}
 \mu^\alpha = u_1\circ (A^\prime)^{s_1} \circ u_2 \circ \cdots \circ u_t\circ
 (A^{\prime})^{s_t}
 \circ u_{t+1}
 \end{equation}
  be the canonical $l$-large $A^\prime$-decomposition of $\mu^\alpha$.  Notice that
if $\mu$ occurs in  $f_M(\sigma)$ (hence $\mu^\alpha$ occurs in
$W_\sigma^\beta$) then this decomposition (\ref{eq:7.3.14})
   is precisely the $A^\prime$-decomposition of $\mu^\alpha$
  induced on $\mu^\alpha$ (as a subword of $W_\sigma^\beta$)
  from the $A^\prime$-decomposition (\ref{eq:7.3.can*}) of $W_\sigma^\beta$  relative to ${\mathcal O}_P$.

   Denote by  $t(\mu)$ the number of non-trivial elements among  $u_1, \ldots,  u_{t+1}$
 (clearly, $u_i \neq 1$ for $2 \leq i \leq t$).

We associate with each long variable $\mu$ a sequence of new
variables (in the equation $T^*(\Pi)$)
 $S_\mu = \{\nu_1, \ldots, \nu_{t(\mu)} \}$. Observe, since the decomposition
(\ref{eq:7.3.14}) of $\mu^\alpha$ is unique, the set $S_\mu$ is
well-defined (in particular, it does not depend on   intervals
$\sigma$).

 It is convenient to define here two
functions $\nu_{\rm left}$ and $\nu_{\rm right}$  on the set
$M_{long}$: if $\mu \in M_{long}$ then
 $$\nu_{\rm left}(\mu) =\nu_1, \ \ \ \nu_{\rm right}(\mu) =\nu_{t(\mu)}.$$

Now we define a new set of variable $M^*$ as follows:
$$ M^* = M_{\rm short}\cup \bigcup_{\mu \in M_{long}} S_\mu .$$

\medskip
{\bf Definition of the labelling function $f^*_{X^*}$}

\medskip
Put $X^* = X$. We define the labelling function $f^*_{X^*} :
{\mathcal E}^* \rightarrow F[X]$ as follows.

Let $\delta \in {\mathcal E}^*$. If $\delta \in {\mathcal
E}_{<j,N}$, then put $$f^*_{X^*}(\delta) = f_X(\delta).$$

 Let now  $\delta = \delta_i \in
E_\sigma$ for some $\sigma \in {\mathcal E}_{j,N}$. Then there are
three cases to consider.

a) $\delta$ corresponds to the consecutive occurrences of powers
${A^\prime}^{q_{j-1}}$ and ${A^\prime}^{q_{j}}$ in the
$A^\prime$-decomposition (\ref{eq:7.3.can*}) of $W_\sigma^\beta$
relative to ${\mathcal O}_P$. Here $j = i$ or $j = i-1$ with
respect to whether $D_1 = 1$ or $D_1 \neq 1$.

As we have mentioned before, according to Lemma \ref{Claim1.b} the
$A^\prime$-decomposition (\ref{eq:7.3.can*}) gives rise to the
unique associated $A^*$-decomposition of $W_\sigma$:
 \begin{equation}
 \label{eq:A-star}
 W_\sigma =
D_1^* \circ_d (A^*)^{q^*_1} \circ_d D^*_2 \circ \cdots \circ_d D^*_k
\circ_d (A^*)^{q^*_{k}} \circ_d D^*_{k+1}.
 \end{equation}

Now put
 $$f^*_X(\delta_i) = D_j^* \in F[X]$$
where $j = i$ if $D_1 = 1$ and  $j = i-1$ if  $D_1 \neq 1$. See
Figure 14.
\begin{figure}[here]
\centering{\mbox{\psfig{figure=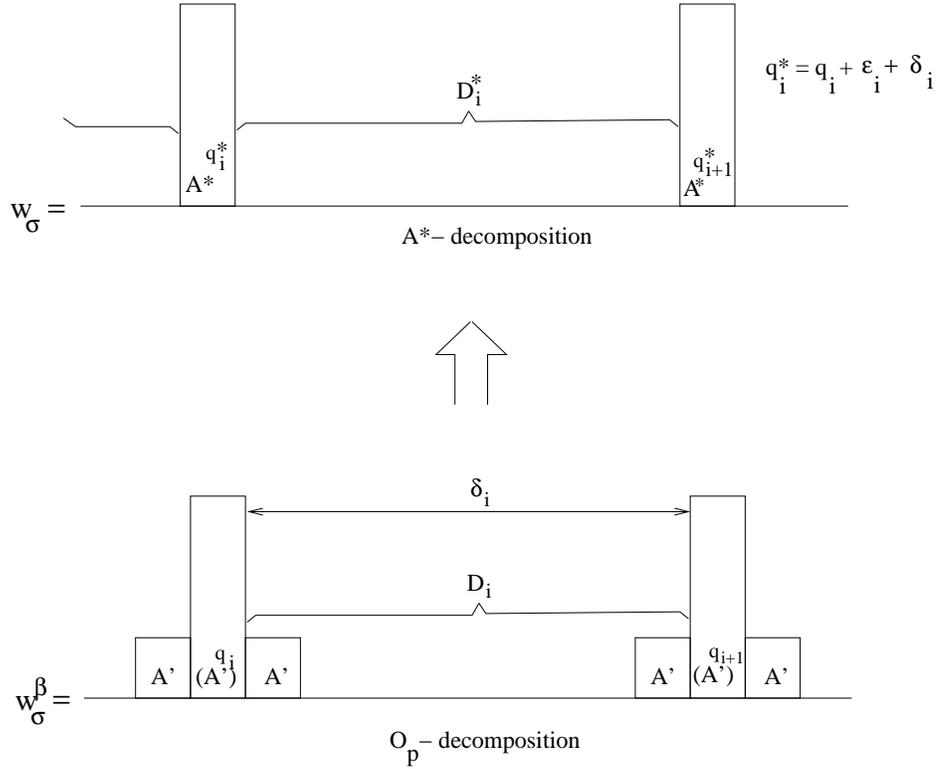,height=4in}}}
\caption{Defining $f^*_{X^*}$.} \label{fx-star}
\end{figure}

The other two cases are treated similarly to  case a).

 b) $\delta$ corresponds to the interval from the beginning of
$\sigma$ to the first $A^\prime$ power ${A^\prime}^{q_1}$ in the
decomposition (\ref{eq:7.3.can*}) of $W_\sigma^\beta$. Put
$$f^*_X(\delta )=D_1^*.$$

c) $\delta$ corresponds to the interval from the last occurrence
of  a power ${A^\prime}^{q_k}$ of $A^\prime$   in the
decomposition (\ref{eq:7.3.can*}) of $W_\sigma^\beta$ to the end
of the interval. Put
$$f^*_X(\delta )=D_{k+1}^*.$$

\medskip
{\bf   Definition of the function $f^*_{M^*}$.}

\medskip
Now we define the function $f^*: {\mathcal E}^* \rightarrow
F[M^*]$.

Let $\delta \in {\mathcal E}^*$. If $\delta \in {\mathcal
E}_{<j,N}$, then put
 $$f^*_{M^*}(\delta) = f_M(\delta)$$
 (observe that all variables in
$f_M(\delta)$ are short, hence they belong to $M^*$).

Let $\delta = \delta_i \in E_\sigma$ for some $\sigma \in {\mathcal
E}_{j,N}$. Again, there are three cases to consider.

a) $\delta$ corresponds to the consecutive occurrences of powers
${A^\prime}^{q_s}$ and ${A^\prime}^{q_{s+1}}$ in the
$A^\prime$-decomposition (\ref{eq:7.3.can*}) of $W_\sigma^\beta$
relative to ${\mathcal O}_P$.  Let the stable occurrence
${A^\prime}^{q_s}$ occur in $\mu_i^\alpha$ for a long variable
$\mu_i$, and the stable occurrence ${A^\prime}^{q_{s+1}}$ occur in
$\mu_j^\alpha$ for a long variable $\mu_j$.

Observe that
$$D_s = right(\mu_i)\circ  \mu_{i+1}^\alpha \circ \cdots \circ
\mu_{j-1}^\alpha  \circ left(\mu_j),$$ for some elements
$right(\mu_i), left(\mu_j) \in F$.

 Now put
$$f^*_{M^*}(\delta) = \nu_{i,right} \mu_{i+1} \ldots
\mu_{j-1} \nu_{j,left},$$
  See Figure 15.
 \begin{figure}[here]
\centering{\mbox{\psfig{figure=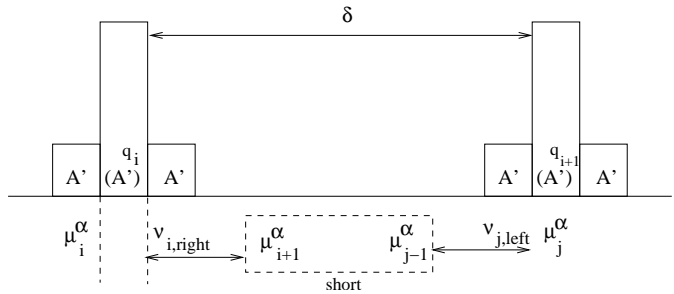,height=1.5in}}}
\caption{Defining $f^*_{M^*},$ case a)}\label{case-a}
\end{figure}

The other two cases are treated similarly to case a).

 b) $\delta$ corresponds to the interval from the beginning of
$\sigma$ to the first $A^\prime$ power ${A^\prime}^{q_1}$ in the
decomposition (\ref{eq:7.3.can*}) of $W_\sigma^\beta$. Put
 $$f^*_{M^*}(\delta) =  \mu_{1} \ldots
\mu_{j-1}\nu _{j,left}.$$

c) $\delta$ corresponds to the interval from the last occurrence
of  a power ${A^\prime}^{q_k}$ of $A^\prime$   in the
decomposition (\ref{eq:7.3.can*}) of $W_\sigma^\beta$ to the end
of the interval.

Denote $\Pi ^*=({\mathcal E}^*, f^*_{X^*}, f^*_{M^*}). $

Now we apply an auxiliary transformation $T'$ to the cut equation
$\Pi ^*$ as follows.    The resulting cut equation will be
 $$T'(\Pi ^*)= (\tilde{\mathcal E}, {\tilde f}_X, {\tilde f}_M),$$
with the same sets $X^*$ and $M^*$, and where ${\tilde f}_{X^*},
{\tilde f}_{M^*}$ are defined as follows. The transformation $T'$
can be applied only in the following two situations.

1) Suppose there are two intervals $\sigma, \gamma \in {\mathcal
E}^*$ such that
$$f^*_{M^*}(\sigma) = \mu \in M^{*\pm 1}, \ \ \ f^*_{M^*}(\gamma) = u \circ \mu \in F[{M^*}],$$
 for some  $u \in F[{M^*}]$ and $f^*_{X^*}(\sigma )=(A^*)^k$, $f^*_{X^*}(\gamma )=w\circ(A^*)^k$. Then put
 $${\tilde f}_{X^*}(\gamma) = w,  \ \ \
{\tilde f}_{M^*}(\gamma) = u,$$
$${\tilde f}_{X^*}(\delta) = f^*_{X^*}(\delta), \ \ \ {\tilde f}_{M^*}(\delta) = f^*_{M^*}(\delta) \ \ \ (\delta \neq
\gamma). $$

2)  Suppose there are two intervals $\sigma, \gamma \in {\mathcal
E}^*$ such that
$$f^*_{M^*}(\sigma) = \mu \in M^{*\pm 1}, \ \ \ f^*_{M^*}(\gamma) = \nu \circ \mu \in F[{M^*}],$$
and $f^*_{X^*}(\gamma )=(A^*)^k\circ f^*_{X^*}(\sigma ).$ Then put
 $${\tilde f}_{X^*}(\gamma ) = (A^*)^k, \ \ \
{\tilde f}_{M^*}(\gamma ) = \nu ,$$

$${\tilde f}_{X^*}(\delta) = f^*_{X^*}(\delta), \ \ \ {\tilde f}_{M^*}(\delta) = f^*_{M^*}(\delta) \ \ \ (\delta \neq
\gamma). $$

 We apply the transformation $T'$ consecutively to $\Pi^\ast$ until it is applicable.
 Notice, since $T'$ decreases the length of the element
$f^*_{M^*}(\gamma)$ it can only be applied a finite number of
times, say $s$, so $(T')^s(\Pi ^*)=(T')^{s+1}(\Pi ^*).$ Observe
also, that the resulting cut equation $(T')^s(\Pi ^*)$ does not
depend on a particular sequence of applications of the
transformation $T'$ to $\Pi^\ast$.  This implies that the  cut
equation $T^*(\Pi) = (T')^s(\Pi ^*) $ is well-defined.

\begin{Claim} \label{2.0.-1.}
 The homomorphism $\alpha^*: F[M^*]
\rightarrow F$  defined as (in the notations above):
$$\alpha^*(\mu) = \alpha(\mu) \ \ (\mu \in M_{\rm short}), $$
$$\alpha^*(\nu_{i,right}) = R^{-\beta}c^{-1} right(\mu_i) \ \ (\nu_i \in S_{\mu} \ for \ \mu \in M_{long})$$
$$\alpha^*(\nu_{i,left}) =  left(\mu_i)c R^{\beta}$$
is a solution of the cut equation $T^*(\Pi)$ with respect to
$\beta: F[X] \rightarrow F$.
\end{Claim}
{\em Proof.}   The statement follows directly from the construction.
  \hfill $\Box$

\medskip

\medskip
{\bf Agreement 4 on the set ${\mathcal P}$:} we assume (by
choosing the function $h$ properly, i.e., $h(i,j) > C(L,N+3)$, see
Lemma \ref{le:D-star-beta})
 that every tuple $p
\in {\mathcal P}$ satisfies the conditions of Lemma
\ref{le:D-star-beta}, so Claim \ref{2.0.-1.} holds for every $p
\in {\mathcal P}$.

\begin{df} \label{def-cut-of-rank}
Let $w \in \bar{\mathcal W}_{\Gamma,L}$. Let $1\le i \le K$. A
{\em cut of rank $i$} of $w$ is a decomposition $w = u \circ v$
where either $u$ ends with $A_i^{\pm2}$ or $v$ begins with
$A_i^{\pm2}$. In this event we say that $u$ and $v$ are obtained
by a cut (in rank $i$) from $w$.
\end{df}

\begin{df}
Given a 3-large tuple $p \in {\mathbb N}^L$, for any $0\le j\le L$
we define by induction (on $L-j$) a set of {\em patterns of rank
$j$} which are certain words in $F(X\cup C)$:
 \begin{enumerate}
 \item [(1)] Patterns of rank $L$ are precisely the letters from the alphabet
$X^{\pm1}$.
 \item [(2)] Now suppose $j=Ks+r$, where $0\le r <K$ and $Ks<L$.
  We represent $p$ as
\begin{equation} \label{p-repr}
p = p' q p'' \quad\text{where}\quad |p'| = Ks, \ |q| = K, \ |p''|
= L-Ks-K.
\end{equation}
\noindent Then a pattern of rank $j$ is either a word of the form
$u^{\phi_{K,q}}$ where $u$ is a pattern of rank $Ks+K$, or a subword
of $u^{\phi_{K,q}}$ formed by one or two cuts of ranks $>r$ (see
Definition \ref{def-cut-of-rank}).
 \end{enumerate}
\end{df}

\begin{rk} \label{re:8}
 $w \in \bar{\mathcal W}_{\Gamma,L}$ for any pattern $w$ of any
rank $j \leq  L$.
\end{rk}

 \begin{Claim}\label{2b}

Let $\phi_L=\phi_{L,p}$, where $p \in \mathbb{N}^L$ such that  $p_t
\geq (l+2)N_{\Pi}+3$ for $t = 1, \ldots,L$, and $l\geq 3$. Denote
$\Pi_L = \Pi_{\phi_L}$.

1) For $j\leq L$ the cut equation $\Pi _{L-j}=(T^*)^{j}(\Pi_L)$ is
well defined and it is a $\Gamma$-cut equation of rank $\leq L-j$
and size $l$. In particular, the sequence $\Sigma_{L,p}$ of
$\Gamma$-cut equations

\begin{equation}
\label{eq:7.3.15^*}
 \Sigma_{L,p}: \ \ \ \Pi_L \stackrel{T^*}{ \rightarrow} \Pi_{L-1}  \stackrel{T^*}{
\rightarrow} \ldots \stackrel{T^*}{ \rightarrow} \Pi_j \rightarrow
\ldots
\end{equation}
is well defined.

2) Let $j= Ks+r$, where $0\leq r<K$, $L=K(s+i)$, and $p'$ be from
the representation \eqref{p-repr}.  Denote $\phi_{Ks} =
\phi_{Ks,p'}$. Then the following are true:
\begin{enumerate}
 \item [(a)]
 for any interval $\sigma$ of $\Pi_j$ there is a pattern $w$ of rank $j$
 such that $f_X(\sigma) = w^{\phi_{Ks}}$;
  \item [(b)]  if $j = Ks$ ($r =0$) then for every interval
 $\sigma$ of the cut equation $\Pi _j$  the pattern
$w$, where $f_X(\sigma) = w^{\phi_{Ks}}$,  does not contain
$N$-large powers of elementary periods.
 \end{enumerate}\end{Claim}
 {\em Proof.}  Let $j= Ks+r$, $0\leq r<K$, $L=K(s+i)$. We  prove the claim by
  induction on $i$ and $m = K-r$ for $i  > 0$.

   Case $i=0$.  In this case $j = L$, so
 the labels of the intervals of $\Pi_L$ are of the  form $x^{\phi _L},\ x\in X$, and  the claim
 is obvious.

  Case  $i=1$.  We use induction on $m = 1, \ldots, K-1$ to prove that for every
  interval $\sigma$ from the cut equation
  $$ \Pi_{L-m} = ({\mathcal E}^{(L-m)}, f_{X}^{(L-m)},
 f_M^{(L-m)})$$
  the label $f_X^{(L-m)}(\sigma)$ is of the form $u^{\phi_{L-K}}$
  for some pattern  $u \in Sub(X^{\phi_K})$.

  Let $m = 1$. In this case $j = L-1$.
  For every $x \in X^{\pm 1}$ one can represent
 the  element $x^{\phi _L}$ as a product of elements of the type
 $y^{\phi_{L-K}}, y \in X^{\pm 1}$ (so the
 element $x^{\phi _L}$ is a word in the alphabet $X^{\phi _{L-K}}$).
 Indeed,
  $$x^{\phi _L}= (x^{\phi_K})^{\phi _{L-K}} = w^{\phi _{L-K}},$$
   where $w = x^{\phi_K}$ is a word in $X$. By Lemma \ref{le:A-A^*}
    there is a precise correspondence between stable $A_L^*$-decompositions of
 $$x^{\phi _L}= w^{\phi
_{L-K}}= D_1^{\phi _{L-K}}\circ _d A_L^{*q_1}\circ _dD_2^{\phi
_{L-K}}\ldots D_k^{\phi _{L-K}}\circ_d A_L^{*q_k}\circ
D_{k+1}^{\phi _{L-K}}$$ and stable $A_K$-decompositions of $w$
 $$w=D_1\circ
{A_K}^{q_1}\circ D_2\ldots D_k\circ A_K^{q_k}\circ D_{k+1}.$$  By
construction, application of the transformation $T^*$ to $\Pi _L$
removes  powers  $A_L^{*q_s}=A_K^{q_s\phi _{L-K}}$ which are
subwords of the word $w^{\phi _{L-K}}$ written in the alphabet
$X^{\phi _{L-K}}$. By construction the words $ D_s^{\phi _{L-K}}$
are the labels of the new intervals of the equation $\Pi _{L-1}$.
 Notice, that $D_s$ are subwords of $w= x^{\phi_K}$ which obtained
 from $w$ by one or two cuts in rank $L$. Hence $D_s$ are patterns
 in rank $L-1$, as required in 2a).

 Now we show that $\Pi_{L-1}$
  is a $\Gamma$-cut equation in rank $\leq L-1$ and size $l$. By 2a) and Remark \ref{re:8}
  $f_X(\sigma) \in \bar  {\mathcal W}_{\Gamma,L}$  for every
  interval $\sigma \in \Pi_{L-1}$. Thus the initial part of the first condition from the
  definition of $\Gamma$-cut equations is satisfied.
To show 1) it suffices to show that 1.1) in rank $L$ does not hold
for $\Pi _{L-1}$. Let $\delta \in {\mathcal E}^{L-1}$.  By the
construction $(A^\prime)^{l+2}$ does not occur in $\mu^\alpha$ for
any $\mu \in M^{L-1}$. Therefore  the maximal power of $A^\prime$
that can occur in $f_{M}(\delta)^\alpha$ is bounded from above by
$(l+1)|f_{M}(\delta)|$ which is less then $(l+2)length(\Pi
_{L-1})$. Hence there are no intervals in $\Pi_{L-1}$ which
satisfy the condition 1.1 from the definition of $\Gamma$-cut
equations. It follows that the rank of $\Pi_{L-1}$ is at most
$L-1$,
 as required.
Let $t$ be the rank of $\Pi_{L-1}$.  For an interval $\delta \in
\Pi_{L-1}$
 if the conditions 1.1) and
1.3) for $f_X(\delta)$ and the rank $t$ are not satisfied, then
the condition 1.2) is satisfied. Indeed, it is obvious from the
definition of patterns that either $f_X(\delta)$ has a non-trivial
$N$-large decomposition in rank $t$ or $\max_t(f_X(\delta)) \leq
2$. Finally, it has been shown in Claim 1 that $T^*(\Pi)$ has a
solution $\alpha^*$ relative to $\beta$.
 This proves  the condition 2) in the definition of the $\Gamma$-cut
 equation. Hence $\Pi_{L-1}$ is a $\Gamma$-cut equation of rank at
  $t \leq j-1$ and size $l$.

Suppose now by induction  on $m$ that for an interval $\sigma $ of
the cut equation $\Pi _j$ (for $m = L-j$) $f_X^{(j)}(\sigma)=u^{\phi
_{L-K}}$ for some  $u\in Sub(X^{\pm\phi _K}).$ Then either $\sigma$
does not change under $T^*$  or $f_X^{(j)}(\sigma )$ has
 a stable  $(l+2)$-large ${A_j}^*$-decomposition  in  rank
$j=r+(L-K)$ associated with long variables in $f_M^{(j)}(\sigma)$:
$$u^{\phi _{L-K}}=\bar D_1^{\phi _{L-K}}\circ _d A_j^{*q_1}\circ
_d\bar D_2^{\phi _{L-K}}\ldots \bar D_k^{\phi _{L-K}}\circ_d
A_j^{*q_k}\circ \bar D_{k+1}^{\phi _{L-K}},$$ and $\sigma$ is an
interval in $\Pi_j$. By Lemma \ref{le:A-A^*}, in this case there
is a stable $A_r$-decomposition of $u$:
$$u=\bar D_1\circ
A_r^{q_1}\circ \bar D_2\ldots \bar D_k\circ A_r^{q_k}\circ \bar
D_{k+1}.$$
 The application of the transformation
$T^*$ to $\Pi _j$ removes powers $A_j^{*q_s} = A_r^{q_s\phi
_{L-K}}$
 (since ${A_j}^*=A_r^{\phi _{L-K}}$)  which are subwords of the word $u^{\phi
_{L-K}}$ written in the alphabet $X^{\phi _{L-K}}$. By
construction the words $\bar D_s^{\phi _{L-K}}$ are the labels of
the new intervals of the equation $\Pi _{j-1}$, so they have  the
required form. This proves statement 2a) for $m+1$. Statement 1)
now follows from 2a) (the argument is the same as in rank $L-1$).
By induction the Claim holds for $m = K$, so the label
$f_X^{(L-K)}(\sigma)$ of an interval $\sigma$ in $\Pi _{L-K}$ is
of the form $u^{\phi _{L-K}},$ for some pattern $u$, where $u \in
Sub(X^{\pm\phi _K})$.
 Notice that  $ Sub(X^{\pm\phi _K})\subseteq {\mathcal
W}_{\Gamma,L}$ which proves  statement 2) (and, therefore,
statement 1) of the Claim for $i=1$.

Suppose, by induction, that labels of intervals in the cut
equation $\Pi _{L-Ki}$ have form $w^{\phi _{L-Ki}},$ $w$ is a
pattern in $\bar {\mathcal W}_{\Gamma ,L}.$ We can rewrite each
label in the form $v^{\phi _{L-K(i+1)}},$ where $v=w^{\phi _K}\in
\bar{\mathcal W}_{\Gamma ,L}$. Similarly to case $i=1$  we can
construct  the $T^*$-sequence
$$\Pi_{L-Ki}\rightarrow \ldots\rightarrow \Pi_{L-K(i+1)}$$
where each application of the transformation $T^*$ removes
subwords in the alphabet \newline $ X^{\phi _{L-K(i+1)}}$. The
argument above shows that the  labels of the  new intervals in all
cut equations $\Pi _{L-Ki-1)},\ldots ,\Pi _{L-K(i+1)}$ are of the
required  form $v^{\phi _{L-K(i+1)}},$ for  patterns $v$ where
$v\in\bar{\mathcal W}_{\Gamma ,L}.$ Following the proof it is easy
to see
 that in labels of intervals in $\Pi_{L-K(i+1)}$ the word $v$ does not contain $N$-large powers of $e^{\phi
_{L-K(i+1)}}$ for an elementary period $e$.

\hfill$\Box$

\begin{Claim}
\label{2.0.1}
 Let ${\cal P} \subseteq \mathbb{N}^L$ be an infinite set of $L$-tuples
 and for $p \in {\cal P}$ let
 $$\Sigma_{L,p}: \ \ \ \Pi_L^{(p)} \stackrel{T^*}{ \rightarrow} \Pi_{L-1}^{(p)}  \stackrel{T^*}{
\rightarrow} \ldots \stackrel{T^*}{ \rightarrow} \Pi_j^{(p)}
\rightarrow \ldots
$$
be the sequence  (\ref{eq:7.3.15^*}) of cut equations $ \Pi_j^{(p)}
= ({\cal E}^{j,p}, f_X^{j,p}, f_M^{j,p})$. Suppose that for a given
$j
>2K$ the following {\em ${\cal P}$-uniformity} property  $U({\cal P},j)$
(consisting of three conditions) holds:

\begin{enumerate}
 \item [(1)]  ${\cal E}^{j,p} = {\cal E}^{j,q}$ for every $p,q \in
 {\cal P}$, we denote this set by ${\cal E}^{j}$;
  \item [(2)] $f_M^{j,p} = f_M^{j,q}$ for every $p,q \in
 {\cal P}$;
  \item [(3)] for any $\sigma \in {\cal E}^{j}$ there exists a
  pattern $w_\sigma$ of rank $j$ such that for any $p \in {\cal P}$
  $f_X^{j,p}(\sigma) = w_\sigma^{\phi_{Kl,p'}}$ where $p'$  is the
  initial segment of $p$ of length  $Kl$, where $j = Kl+r$ and $0 < r
  \leq K$.
\end{enumerate}
Then there exists an infinite subset ${\cal P}'$ of ${\cal P}$ such
that the ${\cal P}'$-uniformity condition $U({\cal P}',j-1)$ holds
for $j-1$.
\end{Claim}
{\em Proof.} Follows from the construction.

 \hfill$\Box$

\medskip
{\bf Agreement 5 on the set ${\mathcal P}$:} we assume, in addition
to all the agreements above,  that  for the set ${\cal P}$ the
uniformity condition $U({\cal P},j)$ holds for every $j >2K$.
Indeed, by Claim \ref{2.0.1} we can adjust ${\cal P}$ consecutively
for each $j >2K$.

\begin{Claim} \label{2.0.0}  Let $\Pi =({\mathcal E},f_X,f_M)$ be a
$\Gamma$-cut equation in rank $j \geq 1$ from the sequence
(\ref{eq:7.3.15^*}). Then for every variable $\mu \in M$ there
exists a word ${\mathcal M}_\mu (M_{T(\Pi)}, X^{\phi_{j-1}},F)$
such that the following equality holds in the group $F$
$$\mu^\alpha = {\mathcal M}_\mu (M_{T(\Pi)}^{\alpha^*}, X^{\phi_{j-1}})^\beta.$$
Moreover, there exists an infinite subset $P' \subseteq P$ such
that
 the words ${\mathcal M}_\mu
(M_{T(\Pi)}, X)$ depend only on exponents $s_1, \ldots, s_t$ of
the canonical $l$-large decomposition (\ref{eq:7.3.14}) of the
words $\mu^\alpha$.
\end{Claim}
{\em Proof.}   The claim follows from the construction. Indeed, in
constructing
 $T^*(\Pi)$ we cut out leading periods of the type $(A_j^\prime)^s$
 from $\mu^\alpha$ (see (\ref{eq:7.3.14})). It follows that to get
 $\mu^\alpha$ back from $M_{T(\Pi)}^{\alpha^*}$ one needs to put
 the exponents $(A_j^\prime)^s$ back. Notice, that
 $$ A_j = A(\gamma_j)^{\phi_{j-1}}$$
  Therefore,  $$ (A_j)^s = A(\gamma_j)^{\phi_{j-1}\beta}$$
 Recall that $A_j^\prime$ is the cyclic reduced form of
 $A_j^\beta$, so
 $$ (A_j^\prime)^s = uA(\gamma_j)^{\phi_{j-1}\beta}v$$
  for some  constants $u, v \in C_{\beta}\subseteq F$. To see
  existence of the subset $P' \subseteq P$ observe that the length
  of the words $f_M(\sigma)$ does not depend on $p$, so there are
  only finitely many ways to cut out the leading periods  $(A_j^\prime)^s$
 from $\mu^\alpha$.  This proves the claim.
\hfill$\Box$

\medskip
{\bf Agreement 6 on the set ${\mathcal P}$:} we assume (replacing
$P$ with a suitable infinite subset) that every tuple $p \in
{\mathcal P}$ satisfies the conditions of  Claim \ref{2.0.0}.
Thus, for every $\Pi = \Pi_i$ from the sequence
(\ref{eq:7.3.15^*}) with a solution $\alpha$ (relative to $\beta$)
the solution $\alpha^*$ satisfies the conclusion of Claim
\ref{2.0.0}.

\begin{df}
We define a new transformation $T$ which is a modified version of
$T^\ast$. Namely, $T$ transforms cut equations and their solutions
$\alpha$ precisely as the transformation $T^\ast$, but it also
transforms the set of tuples ${\mathcal P}$ producing an infinite
subset ${\mathcal P}^\ast \subseteq {\mathcal P}$ which satisfies
the Agreements 1-6.
\end{df}
Now we define a sequence
\begin{equation}
\label{eq:7.3.15}
 \Pi_L \stackrel{T}{ \rightarrow} \Pi_{L-1}  \stackrel{T}{
\rightarrow} \ldots \stackrel{T}{ \rightarrow} \Pi_1
\end{equation} of $N$-large $\Gamma$-cut equations,
 where  $\Pi_L = \Pi_\phi$, and $\Pi_{i-1} = T(\Pi_i)$.
 From now on we fix the sequence (\ref{eq:7.3.15}) and refer to it as the
 {\em $T$-sequence.}

\begin{df} Let $\Pi = ({\mathcal E}, f_X,f_M)$ be a cut equation.
For a positive  integer $n$ by  $k_n(\Pi)$ we denote the number of
intervals $\sigma \in {\mathcal E}$ such that $|f_M(\sigma)| = n$.
The following finite sequence of integers
 $$Comp(\Pi) = (k_2(\Pi), k_3(\Pi), \ldots, k_{length(\Pi)}(\Pi))$$
is called the {\em complexity} of $\Pi$.
 \end{df}

 We well-order complexities of cut equations in the (right) shortlex order:
 if $\Pi$ and $\Pi^\prime$ are two cut equations then
  $Comp(\Pi) <  Comp(\Pi^\prime)$ if and only if $length(\Pi) <
 length(\Pi^\prime)$ or $length(\Pi) =  length(\Pi^\prime)$ and there exists $1 \leq i\leq length(\Pi)$ such
 that
  $k_j(\Pi) = k_j(\Pi^\prime)$
 for all $j > i$ but $k_i(\Pi) < k_i(\Pi^\prime)$.

Observe that intervals $\sigma \in {\mathcal E}$ with
$|f_M(\sigma)| = 1$ have no input into the complexity of a  cut
equation $\Pi = ({\mathcal E}, f_X,f_M)$. In particular, equations
with $|f_M(\sigma)| = 1$ for every $\sigma \in {\mathcal E}$ have
the minimal possible complexity among equations of a given length.
We will write $Comp(\Pi) = {\bf 0}$ in the case when $k_i(\Pi) =
0$ for every $i = 2, \ldots, length(\Pi)$.

\begin{Claim}\label{3.} Let $\Pi = ({\mathcal E}, f_X,f_M)$. Then the following holds:
 \begin{enumerate}
 \item $length(T(\Pi)) \leq length(\Pi)$;
  \item $Comp(T(\Pi)) \leq Comp(\Pi)$;

  \end{enumerate}
\end{Claim}

{\em Proof.}  By straightforward verification. Indeed, if $\sigma
\in {\mathcal E}_{<j}$
  then  $f_M(\sigma) = f_{M^*}^*(\sigma)$. If $\sigma \in {\mathcal E}_{j}$ and
   $\delta_i \in E_\sigma$ then   $$ f_{M^*}^*(\delta_i) = \mu_{i_1}^* \mu_{i_1+1} \ldots
\mu_{i_1+r(i)}^*,$$ where $\mu_{i_1} \mu_{i_1+1} \ldots
\mu_{i_1+r(i)}$ is a subword  of $ \mu_1 \ldots \mu_n$ and hence
$|f_{M^*}^*(\delta_i)| \leq |f_M(\sigma)|$, as required. $\Box$

  We need a few definitions related to the sequence (\ref{eq:7.3.15}). Denote by $M_j$ the set
of variables in the equation $\Pi_j$. Variables from $\Pi_L$ are
called {\em initial } variables. A variable $\mu$ from $M_j$ is
called {\em essential} if it occurs in some $f_{M_j}(\sigma)$ with
$|f_{M_j}(\sigma)| \geq 2$, such occurrence of $\mu$ is called
{\em essential}. By $n_{\mu,j}$ we denote the total number of all
essential occurrences of $\mu$ in $\Pi_j$. Then
 $$S({\Pi_j})=\sum_{i=2}^{N_{\Pi_j}} ik_i(\Pi_j) = \sum_{\mu \in M_j} n_{\mu,j} $$
is the total number of all essential occurrences of variables from
$M_j$ in $\Pi_j$.

\begin{Claim} \label{claim:S-Pi}
If  $1 \leq j \leq L$ then $S(\Pi_j) \leq 2S(\Pi_L)$.
\end{Claim}
{\em Proof.}
 Recall, that every variable
$\mu$ in $M_j$ either belongs to $M_{j+1}$ or it is replaced in
$M_{j+1}$ by the  set $S_\mu$ of new variables (see definition of
the function $f^*_{M^*}$ above). We refer to  variables from
$S_\mu$ as to {\em children} of $\mu$. A given occurrence of $\mu$
in some $f_{M_{j+1}}(\sigma)$, $\sigma \in {\mathcal E}_{j+1}$, is
called a {\em side occurrence} if it is either the first variable
or the last variable (or both) in  $f_{M_{j+1}}(\sigma)$. Now we
formulate several properties of variables from the sequence
(\ref{eq:7.3.15})  which come directly from the construction.  Let
$\mu \in M_j$. Then the following conditions hold:
 \begin{enumerate}
  \item every child of $\mu$ occurs only as a side variable in
  $\Pi_{j+1}$;
   \item every side variable $\mu$ has at most one essential
   child, say $\mu^*$. Moreover, in this event $n_{\mu^*,j+1}
   \leq n_{\mu,j}$;
   \item every initial variable $\mu$ has at most two essential
   children, say $\mu_{left}$ and $\mu_{right}$. Moreover, in this
   case $n_{\mu_{left},j+1} + n_{\mu_{right},j+1} \leq 2n_\mu$.
\end{enumerate}
Now the claim follows from the properties listed above. Indeed,
every initial variable from $\Pi_j$ doubles, at most, the number
of essential occurrences of its children in the next equation
$\Pi_{j+1}$, but all other variables (not the initial ones) do not
increase this number. $\Box$

 Denote by $width(\Pi)$ the {\em width} of $\Pi$ which is
   defined as
    $$width(\Pi) = \max_i {k_i(\Pi)}.$$

\begin{Claim} \label{claim:width}
For every $1 \leq j \leq L$ $width(\Pi_j) \leq 2S(\Pi_L)$
\end{Claim}
 {\em Proof.}  It follows directly from Claim \ref{claim:S-Pi}.
 $\Box$

 Denote by $\kappa(\Pi)$ the number of all $(length(\Pi)-1)$-tuples
of non-negative integers which are bounded by $2S(\Pi_L)$.

\begin{Claim}\label {4.} $Comp(\Pi_L) = Comp(\Pi_{\mathcal L})$.
\end{Claim}

{\em Proof.}   The complexity $Comp(\Pi_L)$ depends only on the
function $f_M$ in $\Pi_L$. Recall that $\Pi_L = \Pi_{\phi}$ is
obtained from the cut equation $\Pi_{\mathcal L}$ by changing only
the labelling function $f_X$, so $\Pi_{\mathcal L}$ and $\Pi_L$
have the same functions $f_M$, hence the same complexities. $\Box$

We say that a sequence
 $$ \Pi_L \stackrel{T}{ \rightarrow} \Pi_{L-1}  \stackrel{T}{
\rightarrow} \ldots  \rightarrow $$
 has {\em $3K$-stabilization} at $K(r+2)$ , where  $2 \leq r \leq L/K$,
  if
$$Comp(\Pi_{K(r+2)})= \ldots = Comp(\Pi_{K(r-1)}).$$
 In this event we denote
 $$K_0 = K(r+2), \ \ \ K_1 = K(r+1), \ \ \ K_2 = Kr, \ \ \ K_3 = K(r-1).$$
For the cut equation $\Pi_{K_1}$ by $M_{\rm veryshort}$ we denote
the subset of variables from $M(\Pi_{K_1})$  which occur unchanged
in $\Pi_{K_2}$ and are short in $\Pi_{K_2}$.

\begin{Claim}\label{5.} For a given $\Gamma$-cut equation $\Pi$ and a
positive integer $r_0 \geq 2$  if $L \geq Kr_0 + \kappa(\Pi)4K$
then for some $r \geq r_0$ either the sequence (\ref{eq:7.3.15})
has $3K$-stabilization at $K(r+2)$ or  $Comp(\Pi_{K(r+1)}) = {0}$.
\end{Claim}
 {\em Proof.}  Indeed, the claim follows by the "pigeon  hole" principle  from
Claims \ref{3.} and \ref{claim:width}  and the fact that there are
not more than $\kappa(\Pi)$ distinct complexities which are less
or equal to $Comp(\Pi)$.
  $\Box$

  Now we define a special set of solutions of the  equation $S(X) =
1$. Let $L = 4K + \kappa(\Pi)4K$,  $p$ be  a fixed $N$-large tuple
from ${\mathbb N}^{L-4K}$, $q$ be an arbitrary fixed $N$-large
tuple from ${\mathbb N}^{2K}$, and $p^*$ be an arbitrary $N$-large
tuple from ${\mathbb N}^{2K}$. In fact, we need $N$-largeness of
$p^*$ and $q$ only to formally satisfy the conditions of the
claims above. Put
$${\mathcal B}_{p,q,\beta} = \{\phi_{L-4K,p}  \phi_{2K,p^*}  \phi_{2K,q}  \beta
\mid p^* \in  {\mathbb N}^{2K} \}.$$

It follows from Theorem \ref{cy2} that ${\mathcal B}_{p,q,\beta}$
is a discriminating family of solutions of $S(X) = 1$.
 Denote $\beta_q =\phi_{2K,q} \circ
  \beta$. Then $\beta_q$ is a solution of $S(X) = 1$ in  general
  position and
   $${\mathcal B}_{q,\beta} = \{ \phi_{2K,p^*}  \beta_q
\mid p^* \in  {\mathbb N}^{2K} \}$$ is also a discriminating
family
 by Theorem \ref{cy2}.

Let $${\mathcal B} = \{\psi _{K_1}=\phi_{K(r-2),p'}  \phi_{2K,p^*}
 \phi_{2K,q}  \beta \mid p^* \in  {\mathbb N}^{2K} \},$$ where $p'$ is a
beginning of $p$.
\begin{prop}
\label{Or}

Let $L = 2K + \kappa(\Pi)4K$ and $\phi_L \in {\mathcal
B}_{p,q,\beta}$. Suppose  the sequence
$$ \Pi_L \stackrel{T}{ \rightarrow} \Pi_{L-1}  \stackrel{T}{ \rightarrow} \ldots
\rightarrow $$  of cut equations {\rm (\ref{eq:7.3.15})} has
$3K$-stabilization at $K(r+2), r\geq 2$.
 Then the set of  variables $M$ of the cut equation $\Pi_{K(r+1)}$ can be partitioned
 into three disjoint subsets
 $$M = M_{\rm veryshort} \cup M_{\rm free} \cup M_{\rm useless}$$
 for which the following holds:
\begin{enumerate}
  \item   there exists a finite  system of  equations
$\Delta(M_{\rm veryshort}) = 1$ over $F$ which has a solution in
$F$;
 \item   for every $\mu \in M_{\rm useless}$ there exists a word
$V_\mu \in F[X \cup M_{\rm free}\cup M_{\rm veryshort}]$ which
does not depend on tuples $p^*$ and $q$;
 \item  for every solution $\delta \in {\mathcal B}$, for every map
$\alpha_{\rm free} : M_{\rm free}
 \rightarrow F$, and every solution $\alpha_s: F[M_{\rm veryshort}] \rightarrow F$
  of the system $\Delta(M_{\rm veryshort}) = 1$  the map $\alpha: F[M] \rightarrow F$
  defined by
   \[ \mu^\alpha = \left\{\begin{array}{ll}
   \mu^{\alpha_{\rm free}} &\mbox{ if  $\mu \in M_{\rm free}$;}\\
   \mu^{\alpha_{s}} & \mbox{ if $\mu \in M_{\rm veryshort}$;}\\
     V_\mu(X^\delta, M_{\rm free}^{\alpha_{\rm free}},
   M_{\rm veryshort}^{\alpha_s}) & \mbox{ if $\mu \in M_{\rm useless}$.}
   \end{array}
 \right. \]
  is a group solution of $\Pi_{K(r+1)}$ with respect to $\beta$.
   \end{enumerate}
 \end{prop}
{\em Proof.}     Below we describe (in a series of claims
\ref{6.}-\ref{21}) some properties of partitions of intervals of
cut equations from the sequence (\ref{eq:7.3.15}):
  $$ \Pi_{K_1} \stackrel{T}{ \rightarrow} \Pi_{K_1-1}
\stackrel{T}{ \rightarrow} \ldots \stackrel{T}{
\rightarrow}\Pi_{K_2}.$$

Fix an arbitrary integer $s$  such that $K_1 \geq s \geq K_2$.

\begin{Claim}\label{6.}  Let $f_M(\sigma) = \mu_1 \ldots \mu_k$ be  a
partition of an interval $\sigma$ of rank $s$ in $\Pi_s$. Then:
 \begin{enumerate}
 \item
the variables  $\mu_2, \ldots, \mu_{k-1}$ are very short;
 \item  either $\mu_1$ or $\mu_k$, or both, are long variables.
 \end{enumerate}
\end{Claim}
 {\em Proof.}  Indeed,  if any of the variables $\mu_2, \ldots, \mu_{k-1}$
 is long then the interval $\sigma$ of $\Pi_s$ is replaced in $T(\Pi_s)$ by a set of intervals
 $E_\sigma$ such that $|f_M(\delta)| < |f_M(\sigma)|$ for every $\delta \in
 E_\sigma$.  This implies that complexity of $T(\Pi_s)$ is
smaller than of $\Pi_s$ - contradiction. On the other hand, since
$\sigma$ is a partition of rank $s$ some variables must be long -
hence the result. \hfill$\Box$

Let $f_M(\sigma) = \mu_1 \ldots \mu_k$ be  a partition of an
interval $\sigma$ of rank $s$ in $\Pi_s$. Then the variables
$\mu_1$ and $\mu_k$ are called {\em side variables}.

\begin{Claim} \label{7.} Let $f_M(\sigma) = \mu_1 \ldots \mu_k$ be  a
partition of an interval $\sigma$ of rank $s$ in $\Pi_s$.  Then
this partition will induce a partition of the form
$\mu_1'\mu_2\ldots \mu_{k-1}\mu_k'$ of some interval in rank $s-1$
in $\Pi_{s-1}$ such that  if $\mu_1$ is short in rank $s$ then
$\mu_1' = \mu_1$, if $\mu_1$ is long in $\Pi_s$ then  $\mu_1'$ is
a new variable which does not appear in the previous ranks.
Similar conditions hold for $\mu_k$.
\end{Claim}
 {\em Proof.}  Indeed, this follows from the
construction of the transformation $T$.\hfill$\Box$
\begin{Claim}\label{8.} Let $\sigma_1$ and
$\sigma_2$  be two intervals of ranks $s$ in $\Pi_s$ such that
$f_X(\sigma_1) = f_X(\sigma_2)$
 and
 $$f_M(\sigma_1) = \mu_1\nu_2\ldots \nu_k,\ \  f_M(\sigma_2) = \mu_1\lambda_2\ldots \lambda_l.$$
  Then   for any solution $\alpha$ of
$\Pi_s$ one has
$$\nu_k^\alpha = \nu_{k-1}^{-\alpha} \ldots \nu_2^{-\alpha}\lambda_2^{-\alpha}
\ldots \lambda_{l-1}^{-\alpha}\lambda_l^{-\alpha}$$ i.e,
$\nu_k^\alpha$ can be expressed via $\lambda_l^{\alpha}$ and a
product of images of short variables.\end{Claim}

\begin{Claim}\label{9.} Let $f_M(\sigma) = \mu_1 \ldots \mu_k$ be  a
partition of an interval $\sigma$ of rank $s$ in $\Pi_s$. Then for
any $u \in X \cup E(m,n)$ the word  $\mu_2^\alpha \ldots
\mu_{k-1}^\alpha $ does not contain a subword of the type
$c_1(M_u^{\phi_{K_1}})^{\beta}c_2,$ where $c_1,c_2\in C_\beta$,
and $M_u^{\phi_{K_1}}$ is the middle of $u$ with respect to
$\phi_{K_1}$.
\end{Claim}
 {\em Proof.}  By Corollary \ref{cy:middles} every word $M_u^{\phi_{K_1}}$
contains a big power (greater than $(l+2)N_{\Pi_s}$) of a period
in rank strictly  greater than $K_2$. Therefore, if
$(M_u^{\phi_{K_1}})^\beta$ occurs in the word $\mu_2^\alpha \ldots
\mu_{k-1}^\alpha $ then some of the variables $\mu_2, \ldots,
\mu_{k-1}$ are not short in some rank greater than $K_2$ -
contradiction. $\Box$

\begin{Claim}\label{15.} Let $\sigma$ be an interval   in $\Pi_{K_1}.$
 Then $f_X(\sigma) = W_{\sigma}$ can be
 written in the form
  $$W_\sigma = w^{\phi_{K_1}},$$
 and the following
holds:
 \begin{enumerate}
\item [(1)] the word $w$  can be uniquely written as $w=v_1\ldots
v_e,$ where  $v_1,\ldots v_e\in X^{\pm 1}\cup E(m,n)^{\pm 1}$, and
$v_iv_{i+1}\not \in E(m,n)^{\pm 1}$.

    \item [(2)]  $w$ is either a subword of a word from the list in Lemma
\ref{le:xyu} or  there exists $i$ such that $v_i=x_2x_1\prod
_{s=m}^{1}c_s^{-z_s}$ and $v_1\ldots v_i$, $v_{i+1}\ldots v_e$ are
 subwords of words from the list in Lemma \ref{le:xyu}. In addition,
$(v_iv_{i+1})^{\phi _K}=v_i^{\phi _K}\circ v_{i+1}^{\phi _K}.$

 \item  [(3)] if $w$ is  a subword of a word from the list in Lemma
\ref{le:xyu}, then  at most for two indices $i,j$ elements $v_i,
v_j$ belong to $ E(m,n)^{\pm 1},$ and, in this case $j=i+1.$
 \end{enumerate}
\end{Claim}

{\em Proof.}  The fact that $W_{\sigma}$ can be written in such a
form follows from Claim \ref{2b}.  Indeed, by Claim \ref{2b},
$W_{\sigma}=w^{\phi _{K_1}},$ where $w\in {\mathcal W}_{\Gamma
,L},$
 therefore it is either a subword of a word
from the list in Lemma \ref{le:xyu} or contains a subword from the
set $Exc$ from statement (3) of Lemma \ref{main}. It can contain
only one such subword, because two such subwords of a word from
$X^{\pm \phi _L}$ are separated by big (unbounded) powers of
elementary periods.

 The uniqueness
of $w$ in the first statement follows from the fact that
$\phi_{K_1}$ is an automorphism. Obviously, $w$ does not depend on
$p$.

 Property (3) follows from the comparison of
the set $E(m,n)$
 with the list from Lemma \ref{le:xyu}.\hfill$\Box$

\begin{Claim}
\label{claim:16} Let $ \Pi_{K_1}=({\mathcal E},f_X, f_M)$ and $\mu
\in M$ be  a long variable (in rank $K_1$) such that $f_M(\delta)
\neq \mu$  for any $\delta \in  {\mathcal E}$. If $\mu$ occurs as
the left variable in $f_M(\sigma)$ for some  $\delta \in {\mathcal
E}$ then it does not occur as the right variable in $f_M(\delta)$
for any $\delta \in {\mathcal E}$ (however,  $\mu^{-1}$ can occur
as the right variable). Similarly, If $\mu$ occurs as the right
variable in $f_M(\sigma)$  then it does not occur as the right
variable in any $f_M(\delta)$.

\end{Claim}
 {\em Proof.}    Suppose $\mu$ is a long variable such that
 $f_M(\sigma) = \mu \mu_2\ldots$ and  $f_M(\delta) = \ldots \mu_s\mu$ for some
intervals  $\sigma, \delta$  from $ \Pi_{K_1}$.  By Claim \ref{15.}
$W_\sigma = w^{\phi_{K_1}}$ for some $w=v_1\ldots v_e,$ where
$v_1,\ldots v_e\in X^{\pm 1}\cup E(m,n)^{\pm 1}$, and
$v_iv_{i+1}\not \in E(m,n)^{\pm 1}$. We divide the proof into three
cases.

{\em Case 1)}. Let  $v_1\neq z_i, y_n^{-1}$. Then  $W_{\sigma}$
begins with a big power of some period $A_j^*, \ j>K_2$ (see Lemmas
\ref{le:7.1.zforms} - \ref{le:7.1.xiforms}), therefore $\mu _1$
begins with a big power of $A_j^{*\beta}$. It follows that in the
rank $j$ the transformation $T$ decreases the  complexity of the
current cut equation. Indeed, when $T$ transforms $\mu$ and $\sigma$
it produces a new set of variables $S_\mu = \{\nu_1, \ldots,
\nu_{t(\mu)} \}$ and a new set of intervals $E_\sigma = \{\sigma_1,
\ldots , \sigma_{k(\sigma)}\}$ such that $f^*_X(\sigma_1)=A_j^{*k}$
for some $k \geq 1$ and $f^*_M(\sigma_1) = \nu_1$. Simultaneously,
when $T$ transforms $\delta$ it produces (among other things) a new
set of intervals $E_\delta = \{\delta_1, \ldots ,
\delta_{k(\delta)}\}$ such that $f^*_X(\delta_{k(\delta) -1})$ ends
on $A_j^{*k}$ and $f^*_M(\delta_{k(\delta) -1})$ ends on $\nu_1$.
Now the transformation $T'$ (part 1) applies to $\sigma_1$ and
$\delta_{k(\delta) -1}$ and decreases the complexity of the cut
equation - contradiction.

 {\em Case 2)}. Let  $v_{left}=z_i$. Then $\mu^\alpha$ begins with
 $z_i^\beta = c_i^{q_i}z_i^{\phi_m\beta_1}$ (see Lemma
 \ref{le:7.1.beta}) for some sufficiently large $q_i$. This implies
 that $c_i^{q_i}$ occurs in $f_M(\delta)^\alpha = f_X(\delta)^\beta$
 somewhere inside (since $f_M(\delta) \neq \mu$). On the other hand,  $f_X(\delta)
 \in \bar{\mathcal{W}}_{\Gamma,L}$, so $c_i^{q_i}$ can occur only at
 the beginning of $f_X(\delta)^\beta$
 (see Lemmas \ref{lm:W-Gamma-L} and \ref{le:xyu}) - contradiction.

  {\em Case 3)}. Let $v_{left}=y_n^{-1}$. Then
$W_{\delta}=\ldots x_n^{-1}\circ y_n^{-1}$. In this case, similar to
the case 1), after application of $T^*$ to the current cut equation
in the rank $K_2+m+4n-4$ one can apply the transformation $T'$ (part
2) which decreases  the complexity - contradiction.

 This proves the claim.  \hfill$\Box$

Our next goal is to  transform further the cut equation $\Pi
_{K_1}$ to the form where all  intervals are labeled by elements
$x^{\phi _{K_1}},\ x\in (X\cup E(m,n))^{\pm 1}$. To this end we
introduce several  new  transformations of $\Gamma$-cut equations.

Let $\Pi = ({\mathcal E},f_X, f_M)$ be a $\Gamma$-cut equation in
rank $K_1$ and size $l$ with a solution $\alpha:F[M] \rightarrow
F$ relative to $\beta : F[X] \rightarrow F$. Let $\sigma \in
{\mathcal E}$ and
   $$W_\sigma=(v_1\ldots v_e)^{\phi
_{K_1}}, \ \ \ e \geq 2,$$  be the canonical decomposition of
$W_\sigma $. For $i, 1 \leq i < e,$ put
 $$v_{\sigma,i,left }=v_1\ldots v_i,\ v_{\sigma,i,right}  = v_{i+1}\ldots v_e.$$
 Let, as usual,
 $$f_M(\sigma) = \mu_1 \ldots \mu_k.$$

We start with a {\bf transformation $T_{1,left}$}.  For $\sigma
\in {\mathcal E}$ and $1 \leq i < e$ denote by $\theta$ the
boundary between $v_{\sigma,i,left}^{\phi _{K_1}\beta}$ and
$v_{\sigma,i,right}^{\phi _{K_1}\beta}$  in the reduced form of
the product $v_{\sigma,i,left}^{\phi _{K_1}\beta}
v_{\sigma,i,right}^{\phi _{K_1}\beta}$. Suppose now that there
exist $\sigma$ and $i$ such that the following two conditions
hold:

\begin{enumerate}
 \item [C1)]
$\mu _1^\alpha$ almost   contains the beginning of the word
$v_{\sigma,i,left}^{\phi _{K_1}\beta}$ till the boundary $\theta$
(up to a very short end of it), i.e.,  there are elements $u_1,
u_2, u_3, u_4 \in F$ such that $v_{\sigma,i,left}^{\phi
_{K_1}\beta} = u_1\circ u_2 \circ u_3$, $ v_{i+1}^{\phi
_{K_1}\beta} = u_3^{-1}\circ u_4$, $u_1u_2u_4 = u_1\circ u_2 \circ
u_4$, and $\mu_1^\alpha $ begins with $u_1$, and $u_2$ is very
short (does not contain $A^{\pm l} _{K_2}$) or trivial.

\item [ C2)] the boundary $\theta$  does not lie inside
$\mu_1^{\alpha}$.
\end{enumerate}

In this event the transformation $T_{1,left}$ is applicable to
$\Pi$ as  described below. We consider three cases with respect to
the location of $\theta$ on $f_M(\sigma)$.

\begin{figure}[here]
\centering{\mbox{\psfig{figure=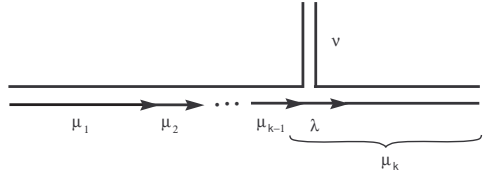,height=2.5in}}}\caption{ T2,
Case 1)} \label{T2}
\end{figure}
\begin{enumerate}

\item [Case 1)] $\theta$ is inside $\mu_k^\alpha$ (see Fig.
\ref{T2}). In this case we perform the following:

 a) Replace  the interval
$\sigma$ by two new intervals $\sigma_1, \sigma _2$ with the
labels $v_{\sigma,i,left}^{\phi _{K_1}},\ v_{\sigma,i,right}^{\phi
_{K_1}}$;

b) Put $f_M(\sigma _1)=\mu _1\ldots \mu _{k-1}\lambda\nu,$
$f_M(\sigma _2)=\nu ^{-1}\mu' _k,$  where $\lambda$ is a new very
short variable, $\nu$ is a new variable.

c) Replace everywhere $\mu _k$ by $\lambda\mu' _{k}$. This
finishes the description of the cut equation  $T_{1,left}(\Pi)$.

d) Define a solution $\alpha^\ast$ (with respect to $\beta$)
  of  $T_{1,left}(\Pi)$ in the natural way.
  Namely, $\alpha^\ast(\mu) = \alpha(\mu)$ for all variables $\mu$ which
  came unchanged from $\Pi$. The values
  $\lambda^{\alpha^\ast},  {\mu'}_{k}^{\alpha^\ast}$, $\nu ^{\alpha^\ast}$ are defined
  in the natural way, that is
${\mu'}_{k}^{\alpha^\ast}$ is the whole end part of $\mu_k^\alpha$
after the boundary $\theta$,
 $(\nu ^{-1}\mu'
_k)^{\alpha^\ast}=v_{\sigma,i,right}^{\phi_{K_1}\beta}$,
$\lambda^{\alpha^\ast} = \mu_k^\alpha ({\mu'}_{k}^\alpha)^{-1}$.

\item [Case 2)] $\theta$ is on the  boundary between  $\mu
_j^{\alpha}$ and $\mu _{j+1}^\alpha$ for some $j$. In this case we
perform the following:

a) We split the interval $\sigma$ into two new intervals $\sigma
_1$ and $\sigma _2$ with labels $v_{\sigma,i, left}^{\phi_{K_1}}$
and $v_{\sigma,i,right}^{\phi _{K_1}}$.

b) We introduce a new variable $\lambda$ and put $f_M(\sigma_1) =
\mu_1\ldots \mu_j\lambda$, $f_M(\sigma _2) =
\lambda^{-1}\mu_{j+1}\ldots \mu_k$.

 c) Define $\lambda^{\alpha^\ast}$ naturally.

\item [Case 3)] The boundary $\theta$  is contained inside
$\mu_i^\alpha$ for some $i (2 \leq i \leq r-1)$. In this case we
do the following:

 a) We split the interval $\sigma$ into two intervals $\sigma _1$ and
$\sigma _2$ with labels $v_{left}^{\phi_{K_1}}$ and
$v_{\sigma,i,right}^{\phi _{K_1}}$, respectively.

 b) Then we introduce three new variables $\mu _j', \mu _j'',
\lambda$, where  $\mu _j', \mu _j''$ are  ``very short", and add
equation $\mu _j = \mu _j'\mu _j''$ to the system $\Delta_{\rm
veryshort}$.

 c) We define
$f_M(\sigma _1) = \mu _1\ldots \mu_j'\lambda $, $f_M(\sigma _2) =
\lambda ^{-1}\mu _j''\mu _{i+1}\ldots \mu _k$.

 d) Define values of
$\alpha^\ast$ on the new variables naturally. Namely, put
$\lambda^{\alpha^\ast} $ to be equal to the terminal segment of
$v_{left}^{\phi_{K_1}\beta}$ that cancels in the product
$v_{left}^{\phi_{K_1}\beta}v_{\sigma,i,right}^{\phi _{K_1}\beta}$.
Now the values $\mu _j'^{\alpha^\ast}$ and $\mu
_j''^{\alpha^\ast}$ are defined to satisfy the equalities
 $$f_X(\sigma_1)^\beta = f_M(\sigma_1)^{\alpha^\ast}, f_X(\sigma_2)^\beta = f_M(\sigma_2)^{\alpha^\ast}.$$

\end{enumerate}

We described the transformation $T_{1,left}$. The transformation
$T_{1, right}$ is defined similarly. We denote both of them by
$T_1$.

Now we describe a {\bf transformation $T_{2,left}$}.

Suppose again that a cut equation $\Pi$ satisfies C1). Assume in
addition that for  these $\sigma$ and $i$ the following condition
holds:
\begin{enumerate}
  \item [C3)] the boundary $\theta$ lies inside $\mu_1^\alpha$.
  \end{enumerate}
Assume also that  one of the following three conditions holds:
 \begin{enumerate}
  \item [C4)] there are no intervals $\delta \neq \sigma$ in $\Pi$ such that
  $f_M(\delta)$ begins with $\mu_1$ or ends on $\mu_1^{-1}$;
  \item  [C5)] $v_{\sigma,i,left} \neq x_n$ (i.e., either  $i > 1$ or $i = 1$ but $v_1 \neq
 x_n$) and for every $\delta \in \mathcal{E}$ in $\Pi$ if $f_M(\delta)$ begins
with $\mu_1$ (or ends on $\mu_1^{-1}$) then  the canonical
decomposition of $f_X(\delta)$ begins with
$v_{\sigma,i,left}^{\phi_{K_1}}$ (ends with $v_{\sigma,i,left}^{-
\phi_{K_1}}$);
 \item  [C6)] $v_{\sigma,i,left} =  x_n$ ($i = 1$ and $v_1 = x_n$)
 and for every $\delta \in \mathcal{E}$ if $f_M(\delta)$ begins
with $\mu_1$ (ends with $\mu_i^{-1}$) then  the canonical
decomposition of $f_X(\delta)$ begins with $x_n^{\phi_{K_1}}$ or
with $y_n^{\phi_{K_1}}$ (ends with $x_n^{-\phi_{K_1}}$ or
$y_n^{-\phi_{K_1}}$).
 \end{enumerate}

In this event the transformation $T_{2,left}$ is applicable to
$\Pi$ as described below.

\begin{enumerate}
 \item [Case C4)] Suppose the condition C4) holds. In this case we
 do the following.

a)  Replace $\sigma$ by two new intervals $\sigma_1, \sigma _2$
with the labels $v_{\sigma,i,left}^{\phi _{K_1}},\
v_{\sigma,i,right}^{\phi _{K_1}}$;

b) Replace $\mu_1$ with  two new variables  $\mu' _1, \mu_1''$ and
put $f_M(\sigma_1) = \mu_1'$, $f_M(\sigma_2) = \mu_1''\mu_2 \ldots
\mu_k.$

c) Define $(\mu_1')^{\alpha^\ast}$ and $(\mu_1'')^{\alpha^\ast}$
such that $f_M(\sigma_1)^{\alpha^\ast} = v_{\sigma,i,left}^{\phi
_{K_1}\beta}$ and $f_M(\sigma_2)^{\alpha^\ast} =
v_{\sigma,i,right}^{\phi _{K_1}\beta}$.

 \item [Case C5)]  Suppose $v_{\sigma,i,left}\neq x_n$. Then do the
 following.

 a) Transform $\sigma$ as described in C4).

 b) If for some interval $\delta \neq \sigma$ the word
$f_M(\delta)$ begins with $\mu_1$ then replace $\mu _1$ in
$f_M(\delta)$ by the variable $\mu''_1$  and replace $f_X(\delta
)$ by $v_{\sigma,i,left}^{-\phi_{K_1}}f_X(\delta).$ Similarly
transform intervals $\delta$ that end with $\mu_1^{-1}$.

\item [Case C6)] Suppose $v_{\sigma,i,left}=x_n$. Then do the following.

 a)  Transform $\sigma$  as described in C4).

 b) If for some $\delta$ the word $f_M(\delta)$ begins with $\mu_1$  and $f_X(\delta)$ does not begin
with $y_n$ then  transform $\delta$ as described in Case C5).

 c)  Leave all other intervals  unchanged.

\end{enumerate}
We described the transformation $T_{2,left}$. The transformation
$T_{2, right}$ is defined similarly. We denote both of them by
$T_2$.

Suppose now that $\Pi = \Pi_{K_1}$.
 Observe that  the transformations
$T_1$ and $T_2$ preserve the properties described in Claims
\ref{3.}
 - \ref{4.} above.
 Moreover, for the
homomorphism $\beta : F[X] \rightarrow F$ we have constructed  a
solution $\alpha^* :F[M] \rightarrow F$ of $T_n(\Pi _{K_1})$  ($n
= 2,3$) such that the initial solution $\alpha$ can be
reconstructed from $\alpha^\ast$ and the equations $\Pi$ and
$T_n(\Pi)$.
 Notice also that the
length of the elements $W_{\sigma'}$ corresponding to new
intervals $\sigma$ are shorter than the length of the words
$W_\sigma$ of the original intervals $\sigma$ from which $\sigma'$
were obtained. Notice also that the transformations $T_1, T_2$
preserves the property of intervals formulated in the Claim
\ref{6.}.

\begin{Claim} \label{cl:T2-T3}

Let $\Pi$ be a cut equation which satisfies the conclusion of the
Claim \ref{6.}. Suppose  $\sigma$ is an interval in $\Pi$ such
that $W_{\sigma}$ satisfies the conclusion of Claim \ref{15.}. If
for some $i$
 $$(v_1 \ldots v_e)^{\phi_{K}} = (v_1 \ldots v_i)^{\phi_{K}} \circ
(v_{i+1} \ldots
 v_e)^{\phi_{K}}$$
  then either $T_1$ or $T_2$ is applicable to given
  $\sigma$ and $i$.
\end{Claim}

{\it Proof.} By Corollary \ref{R} the automorphism $\phi_{K_1}$
satisfies the
 Nielsen property with respect to $\bar {\mathcal W}_{\Gamma}$ with
 exceptions $E(m,n)$. By Corollary 12, equality
 $$(v_1 \ldots v_e)^{\phi_{K}} = (v_1 \ldots v_i)^{\phi_{K}} \circ
(v_{i+1} \ldots
 v_e)^{\phi_{K}}$$
implies that the element that is cancelled between $(v_1 \ldots
v_i)^{\phi_{K}\beta}$ and   $(v_{i+1} \ldots
 v_e)^{\phi_{K}\beta}$ is short in rank $K_2$. Therefore either $\mu
_1^{\alpha}$ almost
 contains
$(v_1 \ldots v_i)^{\phi_{K}\beta}$
 or $\mu _k^{\alpha}$ almost contains $(v_{i+1} \ldots
 v_e)^{\phi_{K}\beta}$. Suppose $\mu _1^{\alpha}$ almost contains
$(v_1 \ldots v_i)^{\phi_{K}\beta}$. Either we can apply
$T_{1,left}$, or the boundary $\theta$ belongs to $\mu
_1^{\alpha}$. One can verify using formulas from Lemmas
\ref{le:7.1.zforms}-\ref{le:7.1.xiforms} and \ref{main} directly
that in this case one of the conditions $C4)-C6)$ is satisfied,
and, therefore $T_{2,left}$ can be applied.
 \hfill $\Box$

\begin{lm} \label{n12} Given a cut equation $\Pi_{K_1}$ one can effectively
find a finite sequence of transformations $Q_1, \ldots, Q_s$ where
$Q_i \in \{T_1, T_2\}$ such that for every interval $\sigma$  of
the cut equation $\Pi_{K_1}^\prime = Q_s \ldots Q_1 (\Pi_{K_1})$
the label $f_X(\sigma)$ is of
 the form $u^{\phi _{K_1}}$, where $u \in X^{\pm 1}\cup E(m,n)$.

 Moreover, there exists an infinite subset $P'$ of the solution set $P$ of
$\Pi _{K_1}$ such that this sequence
 is the same for any solution in $P'$.
\end{lm}
{\it Proof.} Let $\sigma$ be an interval of the equation
$\Pi_{K_1}$. By Claim \ref{15.} the word $W_\sigma$ can be
uniquely written in the canonical decomposition form
 $$W_\sigma = w^{\phi_{K_1}}=(v_1 \ldots v_e)^{\phi_{K_1}},$$
 so that the conditions 1), 2), 3) of Claim \ref{15.} are satisfied.

It follows from the construction of $\Pi _{K_1}$ that either $w$
is a subword of a word between two elementary squares $x\neq c_i$
or begins and (or) ends with some power $\geq 2$ of an elementary
period. If $u$ is an elementary period,
 $u^{2\phi _K}=u^{\phi _K}\circ u^{\phi _K}$, except $u=x_n$,
 when the middle is exhibited in the proof of Lemma \ref{main}.
 Therefore, by Claim \ref{cl:T2-T3}, we can apply $T_1$ and $T_2$ and cut
$\sigma$ into
 subintervals $\sigma _i$
such that for any $i$ $f_X(\sigma _i)$
 does not contain powers $\geq 2$ of elementary periods.
All possible values of $u^{\phi _K}$ for
 $u\in E(m,n)^{\pm 1}$ are shown in the proof of Lemma
 \ref{main}. Applying $T_1$ and $T_2$ as in Claim \ref{cl:T2-T3}
we can split intervals (and their
 labels) into parts with labels of the form $x^{\phi _{K_1}},\
 x\in (X\cup E(m,n)),$
except for the following cases:

1. $w=uv$, where $u$ is $x_i^2, i<n,\ v\in E_{m,n},$ and $v$ has
at
 least three letters,

2.
$w=x_{n-2}^2y_{n-2}x_{n-1}^{-1}x_nx_{n-1}y_{n-2}^{-1}x_{n-2}^2,$

3. $w=x_{n-1}^2y_{n-1}x_n^{-1}x_{n-1}y_{n-2}^{-1}x_{n-2}^{-2},$

4. $y_{r-1}x_r^{-1}y_r^{-1},\ r<n,$

5. $w=uv$, where $u=(c_1^{z_1}c_2^{z_2})^2,\ v\in E(m,n),$ and $v$
is one of the following: $v=\prod _{t=1}^mc_t^{z_t}x_1^{\pm 1},$
$v=\prod _{t=1}^mc_t^{z_t}x_1^{\pm 1}\prod _{t=m}^1c_t^{-z_t}$, $v
= \prod_{t=1}^mc_t^{z_t}x_1\prod
_{t=m}^1c_t^{-z_t}(c_1^{z_1}c_2^{z_2})^{-2},$

6. $w=uv$, where $u=(c_1^{z_1}c_2^{z_2})^2,\ v\in E(m,n),$ and $v$
is one of the following: $v=\prod _{t=1}^mc_t^{z_t}x_1^{-1
}x_2^{-1}$ or $v=\prod _{t=1}^mc_t^{z_t}x_1^{-1}y_1^{-1}.$

7. $w=z_iv.$

Consider the first case. If $f_M(\sigma)=\mu _1\ldots\mu _k,$ and
$\mu _1^{\alpha}$ almost contains
 $$x_i^{\phi
_{K_1}}(A_{m+4i+K_2}^*)^{-p_{m+4i+K_2}+1}x_{i+1}^{\phi
_{K_2}\beta}$$ (which is a non-cancelled initial peace of
$x_i^{2\phi _{K_1}\beta}$ up to a very short part of it), then
either $T_{1,left}$ or $T_{2,left}$ is applicable and we split
$\sigma$ into two intervals $\sigma _1$ and $\sigma _2$ with
labels $x_i^{2\phi _{K_1}}$ and $v^{\phi _{K_1}}$.

Suppose $\mu _1^{\alpha}$ does not contain $x_i^{\phi
_{K_1}}(A_{m+4i+K_2}^*)^{-p_{m+4i+K_2}+1}x_{i+1}^{\phi
_{K_2}\beta}$ up to a very short part. Then $\mu _k^{\alpha}$
contains the non-cancelled left end $E$ of $v^{\phi _{K+1}\beta},$
and $\mu _k^{\alpha}E^{-1}$ is not very short. In this case
$T_{2,right}$ is applicable.

We can similarly consider all Cases 2-6.

Case 7. Letter $z_i$  can  appear only
 in the beginning of $w$ (if  $z_i^{-1}$ appears at the
 end of $w$, we can replace $w$ by $w^{-1}$)
 If $w=z_it_1\ldots t_s$ is the canonical decomposition, then
 $t_k=c_j^{\pm z_j}$ for each $k$. If $\mu _1^{\alpha}$ is longer
 than
the non-cancelled part of  $(c_i^{p}z_i)^{\beta}$, or
 the difference between $\mu _1^{\alpha}$
and $(c_i^{p}z_i)^{\beta}$ is very short, we can split $\sigma$
into two parts, $\sigma _1$ with label $f_X(\sigma _1)=z^{\phi
_{K_1}}$ and $\sigma _2$ with label $f_X(\sigma _2)=(t_1\ldots
t_s)^{\phi _{K_1}}.$

If the difference between $\mu _1^{\alpha}$ and
$(c_i^{p}z_i)^{\beta}$ is not very short, and $\mu _1^{\alpha}$ is
shorter than the non-cancelled part of  $(c_i^{p}z_i)^{\beta}$,
then there is no interval $\delta$ with $f(\delta )\neq f(\sigma
)$ such that $f_M(\delta )$ and $f_M(\sigma )$ end with $\mu _k,$
and we can split $\sigma$ into two parts using $T_1$, $T_2$  and
splitting $\mu _k.$

We have considered all possible cases. $\Box$

Denote the resulting cut equation by $\Pi '_{K_1}.$

\begin{cy} The intervals of $\Pi '_{K_1}$ are labelled by elements
$u^{\phi _{K_1}},$ where

 for $n=1$
$$u\in\{z_i,\ x_i,\ y_i,\  \prod c_s^{z_s},\  x_1\prod
_{t=m}^{1}c_t^{-z_t},\}
$$

for $n=2$ \begin{multline*}u\in\{z_i,\ x_i,\ y_i,\  \prod
c_s^{z_s},\ y_1x_1\prod _{t=m}^{1}c_t^{-z_t},\  y_1x_1,\ \prod
_{t=1}^{m}c_t^{z_t}x_1\prod _{t=m}^{1}c_t^{-z_t}\ , \prod
_{t=1}^{m}c_t^{z_t}x_1^{-1}x_2^{\pm 1},\\
\prod _{t=1}^{m}c_t^{z_t}x_1^{-1}x_2x_1,\  \prod
_{t=1}^{m}c_t^{z_t}x_1^{-1}x_2x_1\prod _{t=m}^{1}c_t^{-z_t},\
x_1^{-1}x_2x_1\prod _{t=m}^{1}c_t^{-z_t},\  x_2x_1\prod
_{t=m}^{1}c_t^{-z_t},\\ x_1^{-1}x_2,\ x_2x_1\},\end{multline*} and
for $n\geq 3$, \begin{multline*}u\in \{z_i,\ x_i,\ y_i,\
c_s^{z_s},\ y_1x_1\prod _{t=m}^{3}c_t^{-z_t},\ \prod
_{t=1}^{m}c_t^{z_t}x_1^{-1}x_2^{-1},\ y_rx_r, \ x_1\prod
_{t=m}^{1}c_t^{-z_t},\\ y_{r-2}x_{r-1}^{-1}x_r^{-1},\
y_{r-2}x_{r-1}^{-1},\ x_{r-1}^{-1}x_r^{-1},\ y_{r-1}x_r^{-1},\
r<n,\ x_{n-1}^{-1}x_nx_{n-1},\\
y_{n-2}x_{n-1}^{-1}x_nx_{n-1}y_{n-2}^{-1},\
y_{n-2}x_{n-1}^{-1}x_n^{\pm 1},\  x_{n-1}^{-1}x_n,\ x_nx_{n-1},\\
y_{n-1}x_n^{-1}x_{n-1}y_{n-2}^{-1},\ y_{n-1}x_n^{-1},
y_{r-1}x_r^{-1}y_r^{-1}\}.\end{multline*}
\end{cy}

{\it Proof.} Direct inspection from Lemma \ref{n12}.
 \hfill $\Box$

Below we suppose $n>0$.

We still want to reduce the variety of possible labels of
intervals in $\Pi '_{K_1}$. We cannot apply $T_1$, $T_2$ to some
of the intervals labelled by $x^{\phi _{K_1}}$, $x\in X\cup
E(m,n)$, because there are some cases when $x^{\phi _{K_1}}$ is
completely cancelled in  $y^{\phi _{K_1}}$, $x,y\in (X\cup
E(m,n))^{\pm 1}.$

We will change the basis of $F(X\cup C_S)$, and then apply
transformations $T_1$, $T_2$ to the labels written in the new basis.
Replace, first, the basis $(X\cup C_S)$ by a new basis $\bar X\cup
C_S$ obtained by replacing each variable $x_s$ by
$u_s=x_sy_{s-1}^{-1}$ for $s>1$, and replacing $x_1$ by
$u_1=x_1c_m^{-z_m}$:
  $$ \bar X = \{u_1 = x_1c_m^{-z_m}, u_2 = x_2y_{1}^{-1}, \ldots, u_n = x_ny_{n-1}^{-1},
  y_1, \ldots,y_n, z_1, \ldots, z_m\}.$$

Consider  the case $n\geq 3.$ Then the labels of the intervals
will be rewritten as $u^{\phi _{K_1}}$, where
$$u\in\{z_i,\ u_iy_{i-1},\ y_i,\ \prod _sc_s^{z_s},\
 y_1u_1\prod _{j=n-1}^1c_j^{-z_j}, \  u_1^{-1}y_1^{-1}u_2^{-1},\ y_ru_ry_{r-1}, u_r,\
 u_{r-1}^{-1}y_{r-1}^{-1}u_r^{-1},$$ $$ u_ry_{r-1}u_{r-1}y_{r-2},\
 u_2y_1u_1\prod _{j=n-1}^1c_j^{-z_j},\ r<n;\ \
y_{n-2}^{-1}u_{n-1}^{-1}u_ny_{n-1}u_{n-1}y_{n-2},$$ $$
 u_{n-1}^{-1}u_ny_{n-1}u_{n-1},\ u_{n-1}^{-1}u_ny_{n-1}, \
 u_{n-1}^{-1}y_{n-1}^{-1}u_n^{-1},$$ $$
 y_{n-2}^{-1}u_{n-1}^{-1}u_ny_{n-1},  \ u_ny_{n-1}u_{n-1}y_{n-2},\
 u_n^{-1}u_{n-1},\ u_n\}.$$

In the cases $n=1,2$ some of the labels above do not appear, some
coincide. Notice, that $x_n^{\phi _K}=u_n^{\phi _K}\circ
y_{n-1}^{\phi _K},$ and that the first letter of $y_{n-1}^{\phi
_K}$ is not cancelled in the products
$(y_{n-1}x_{n-1}y_{n-2}^{-1})^{\phi _K},$ $(y_{n-1}x_{n-1})^{\phi
_K}$ (see Lemma \ref{le:7.1.x1formsmneq0}). Therefore,  applying
transformations similar to $T_1$ and $T_2$ to the cut equation
$\Pi'_{K_1}$ with labels written in the basis $\bar X$, we can
split all the intervals with labels containing $(u_ny_{n-1})^{\phi
_{K_1}}$ into two parts and obtain a cut equation with the same
properties and intervals labelled by $u^{\phi _{K_1}},$ where

$$u\in\{z_i,\ u_iy_{i-1},\ y_i,\ \prod _sc_s^{z_s},\
 y_1u_1\prod _{j=n-1}^1c_j^{-z_j}, \  u_1^{-1}y_1^{-1}u_2^{-1},\ y_ru_ry_{r-1}, u_r,\
 u_{r-1}^{-1}y_{r-1}^{-1}u_r^{-1},$$ $$ u_ry_{r-1}u_{r-1}y_{r-2},\
 u_2y_1u_1\prod _{j=n-1}^1c_j^{-z_j},\ r<n;$$
$$y_{n-2}^{-1}u_{n-1}^{-1}u_n,\ y_{n-1}u_{n-1}y_{n-2},\
 u_{n-1}^{-1}u_n,\ y_{n-1}u_{n-1},
\ u_n\}.$$

Consider for $i<n$ the expression for $(y_iu_i)^{\phi
_K}=A_{m+4i}^{-p_{m+4i}+1} $
$$\circ x_{i+1}\circ A
_{m+4i-4}^{-p_{m+4i-4}}\circ x^{p_{m+4i-3}}\circ y_i\circ A
_{m+4i-2}^{p_{m+4i-2}-1}\circ x_i\circ \tilde y_{i-1}^{-1}.$$

Formula 3.a) from Lemma \ref{main} shows that $u_i^{\phi _K}$ is
completely cancelled in  the product $y_i^{\phi _K}u_i^{\phi _K}$.
This implies that $y_i^{\phi _{K}}=v_i^{\phi _{K}}\circ u_i^{-\phi
_{K}}$.

Consider also the product $ y_{i-1}^{-\phi _K}u_i^{-\phi _K}=$

\vspace{5mm} \noindent
 $ \left({\bf A_{m+4i-4}^{-p_{m+4i-4}+1}\circ x_i\circ
\tilde y_{i-1}} \circ x_i^{-1}A_{m+4i-4}^{p_{m+4i-4}-1}\right ) $

\vspace{5mm} \noindent
 $ \left ( A_{m+4i-4}^{-p_{m+4i-4}+1}x_i\circ {\bf
(x_{i}^{p_{m+4i-3}}y_{i-1}\ldots
*)^{p_{m+4i-1}-1}x_i^{p_{m+4i-3}}y_ix_{i+1}^{-1}A_{m+4i}^{p_{m+4i}-1}}\right
),$

\noindent \vspace{ 3mm}
 where the non-cancelled part is made bold.

Notice, that $(y_{r-1}u_{r-1})^{\phi _{K}}y_{r-2}^{\phi
_{K}}=(y_{r-1}u_{r-1})^{\phi _{K}}\circ y_{r-2}^{\phi _{K}},$
because $u_{r-1}^{\phi _{K}}$ is completely cancelled in the
product $y_i^{\phi _{K}}u_i^{\phi _{K}}$.

Therefore, we can again apply the transformations similar to $T_1$
and $T_2$ and split the intervals into the ones with labels
$u^{\phi _{K_1}}$, where \begin{multline*}u\in\{z_s,\  y_i,\ u_i,\
\prod _s c_s^{z_s},\ y_ru_r,\  y_1u_1\prod _{j=m-1}^1c_j^{-z_j}, \
u_{n-1}^{-1}u_n=\bar u_n,\\ 1\leq i\leq n, \ 1\leq j\leq m,\ 1\leq
r <n\}.\end{multline*}

Now we change the basis $\bar X$ with a new basis $\hat X$ replacing
$y_r, 1<r<n,$ by a new variable $v_r=y_ru_r$,   $y_1u_1\prod
_{j=m-1}^1c_j^{-z_j}$ by $v_1$, and  $u_{n-1}u_n$ by $\bar u_n$:

$$\hat X=\{z_1, \ldots, z_m,\ u_1, \ldots,u_{n-1}, \bar u_n=u_{n-1}u_n ,\ v_1 = y_1u_1\prod
_{j=m-1}^1c_j^{-z_j}, v_2 =y_2u_2, \ldots, v_n = y_nu_n, \ y_n \}.$$

Then $y_r^{\phi _{K}}=v_r^{\phi _{K}}\circ u_r^{-\phi _{K}},$ and
$y_1^{\phi _K}=v_1^{\phi _K}\circ c_1^{z_1^{\phi _K}}\circ
c_{m-1}^{z_{m-1}^{\phi _K}}\circ u_1^{-\phi _K}$ (if $n\neq 1$).
Formula 2.c) shows that $u_n^{\phi _{K}}=u_{n-1}^{\phi _{K}}\circ
(u_{n-1}^{-1}u_n)^{\phi _{K}}.$

Apply transformations similar to $T_1$ and $T_2$ to the intervals
with labels written in the new basis $\hat X$  and obtain intervals
with labels $u^{\phi _{K_1}},$ where
$$u\in \hat X\cup \{c_m^{z_m}\}.$$
Denote the resulting cut equation by ${\bar \Pi} _{K_1} = (\bar{
\mathcal E}, f_{\bar X},f_{\bar M})$. Let  $\alpha$ be the
corresponding  solution of ${\bar \Pi} _{K_1}$ with respect to
$\beta .$

Denote by $\bar M_{side}$ the set of long variables in $\bar \Pi
_{K_1}$, then $ \bar M =  \bar M_{\rm veryshort}  \cup {\bar
M}_{side}$.

Define a binary relation $\sim_{left}$ on $\bar M_{side}^{\pm 1}$
as follows.  For $\mu_1, \mu'_1 \in \bar M_{side}^{\pm 1}$ put
$\mu_1 \sim_{left} \mu_1^\prime$ if and only if there exist two
intervals $\sigma, \sigma'  \in \bar{E}$ with $f_{\bar X}(\sigma
)=f_{\bar X}(\sigma')$ such that
 $$f_{\bar M}(\sigma)=\mu _1\mu_2\ldots \mu_{r}, \ \ \ f_{\bar M}(\sigma ')=\mu _1'\mu _2'\ldots
 \mu_{r'}'$$
and either $\mu_r = \mu'_{r'}$ or $\mu_r,  \mu'_{r'} \in  M_{\rm
veryshort}.$ Observe that if $\mu_1 \sim_{left} \mu_1^\prime$ then
 $$\mu_1 = \mu_1^\prime \lambda_1 \ldots \lambda_t$$
for some $\lambda_1, \ldots, \lambda_t \in M_{\rm veryshort}^{\pm
1}.$ Notice, that $\mu \sim_{left} \mu$.

Similarly, we define a binary relation $\sim_{right}$ on $\bar
M_{side}^{\pm 1}$.  For $\mu_r, \mu'_{r'} \in \bar M_{side}^{\pm
1}$ put $\mu_r \sim_{right} \mu_{r'}^\prime$ if and only if there
exist two intervals $\sigma, \sigma'  \in \bar{E}$ with $f_{\bar
X}(\sigma )=f_{\bar X}(\sigma')$ such that
 $$f_{\bar M}(\sigma)=\mu _1\mu_2\ldots \mu_{r}, \ \ \ f_{\bar M}(\sigma ')=\mu _1'\mu _2'\ldots
 \mu_{r'}'$$
and either $\mu_1 = \mu'_1$ or $\mu_1,  \mu'_1 \in  M_{\rm
veryshort}.$  Again, if  $\mu_r \sim_{right} \mu_{r'}^\prime$ then
 $$\mu_r =  \lambda_1 \ldots \lambda_t\mu_{r'}^\prime$$
for some $\lambda_1, \ldots, \lambda_t \in M_{\rm veryshort}^{\pm
1}.$

  Denote by $\sim$ the transitive closure of
   $$\{(\mu,\mu') \mid \mu \sim_{left} \mu'\} \cup \{(\mu,\mu') \mid \mu
   \sim_{right} \mu'\} \cup \{(\mu,\mu^{-1}) \mid \mu \in \bar M_{side}^{\pm 1} \}.$$
   Clearly, $\sim$ is an equivalence relation on $\bar
  M_{side}^{\pm 1}$.
  Moreover, $\mu \sim \mu'$ if and only if there exists a sequence
  of variables
   \begin{equation}
   \label{eq:seq-mu}
   \mu = \mu_0, \mu_1, \ldots, \mu_k = \mu'
    \end{equation}
 from $\bar M_{side}^{\pm 1}$  such that either $\mu_{i-1} = \mu_i$, or $\mu_{i-1} = \mu_i^{-1}$,
   or  $\mu_{i-1} \sim_{left} \mu_i$, or $\mu_{i-1} \sim_{right} \mu_i$
    for $i= 1, \ldots, k.$ Observe that if $\mu_{i-1}$ and $\mu_i$ from (\ref{eq:seq-mu})
    are side variables of "different sides" (one is on  the left, and the other is on the right)
     then $\mu_i = \mu_{i-1}^{-1}$. This implies that replacing in
     the sequence (\ref{eq:seq-mu}) some elements $\mu_i$ with
     their inverses one can get a new sequence
   \begin{equation}
   \label{eq:seq-mu-2}
   \mu = \nu_0, \nu_1, \ldots, \nu_k = (\mu')^\varepsilon
    \end{equation}
 for some $\varepsilon \in \{1,-1\}$ where  $\nu_{i-1} \sim \nu_i$ and all the variables $\nu_i$ are
 of the same side.   It follows  that if $\mu$ is a left-side variable and $\mu \sim \mu'$ then
     \begin{equation}
     \label{eq:sim-mu}
     (\mu')^\varepsilon =  \mu \lambda_1 \ldots \lambda_t
      \end{equation}
for some $\lambda_j \in M_{\rm veryshort}^{\pm 1}.$

It follows from (\ref{eq:sim-mu}) that for a variable $\nu \in
\bar M_{side}^{\pm 1}$ all variables from the equivalence class
$[\nu]$ of $\nu$ can be expressed via $\nu$ and very short
variables from $M_{\rm veryshort}$. So if we fix a system of
representatives $R$ of $\bar M_{side}^{\pm 1}$ relative to $\sim$
then all other variables from $\bar M_{side}$ can be expressed as
in (\ref{eq:sim-mu}) via variables from $R$ and very short
variables.

This allows one to introduce a new transformation $T_3$ of cut
equations. Namely, we fix a set of representatives $R$ such that for
every $\nu \in R$ the element $\nu^\alpha$ has minimal length among
all the variables in this class. Now, using (\ref{eq:sim-mu})
replace every variable $\nu$ in every word
 $f_M(\sigma)$ of a cut equation $\Pi$ by its expression via the
 corresponding representative variable from $R$ and a product of
 very short variables.

Now we repeatedly apply the transformation $T_3$  till the
equivalence relations $\sim_{left}$ and $\sim_{right}$ become
trivial.  This process stops in finitely many steps since the
non-trivial relations decrease the number of side variables.

Denote the resulting equation again by ${\bar \Pi} _{K_1}$.

Now we introduce an equivalence relation on partitions of  ${\bar
\Pi} _{K_1}$. Two partitions $f_M(\sigma )$ and $f_M(\delta )$ are
equivalent ($f_M(\sigma ) \sim f_M(\delta )$) if $f_X(\sigma
)=f_X(\delta )$ and either the left side variables or the right
side variables of $f_M(\sigma )$ and $f_M(\delta )$ are equal.
Observe, that $f_X(\sigma )=f_X(\delta )$ implies $f_M(\sigma
)^\alpha=f_M(\delta )^\alpha$, so in this case the partitions
$f_M(\sigma )$ and $f_M(\delta )$ cannot begin with $\mu$ and
$\mu^{-1}$ correspondingly. It follows that if $f_M(\sigma ) \sim
f_M(\delta )$ then the left side variables and, correspondingly,
the right side variables of $f_M(\sigma )$ and $f_M(\delta )$ (if
they exist) are equal. Therefore, the relation $\sim$ is, indeed,
an equivalence relation on the set of partitions of ${\bar \Pi}
_{K_1}$.

 If an equivalence class of partitions
contains two distinct elements $f_M(\sigma )$ and $f_M(\delta )$
then the equality
 $$ f_M(\sigma )^\alpha = f_M(\delta )^\alpha$$
implies the corresponding  equation on the variables $\bar M_{\rm
veryshort},$ which is obtained by deleting all side variables
(which are equal) from $f_M(\sigma )$ and $f_M(\delta )$ and
equalizing the resulting words in very short variables.

Denote by $\Delta (\bar M _{\rm veryshort})=1$ this system.

Now we describe a transformation $T_4$. Fix a set of
representatives $R_p$ of partitions of ${\bar \Pi} _{K_1}$ with
respect to the equivalence relation $\sim$. For a given class of
equivalent partitions we take as a representative an interval
$\sigma$ with $f_M(\sigma )=\mu _{\rm left}\ldots \mu _{\rm
right}.$

 Below we say that:  a word $w \in F[X]$ is {\em very short} if the reduced form of $w^\beta$ does not
contain $(A_j^\prime)^3$  for any $j \geq K_2$; a word $v \in F$
is {\em very short} if it does not contain $(A_j^\prime)^3$ for
any $j \geq K_2$; we also say that $\mu ^{\alpha}$ almost contains
$u^{\beta}$ for some word $u$ in the alphabet $X$ if $\mu
^{\alpha}$ contains a subword which is the reduced form of
$f_1u^{\beta}f_2$ for some  $f_1,f_2\in C_{\beta}$.

{\bf Principal variables.}  A long variable $\mu _{\rm left}$ or
$\mu _{\rm right}$ for an interval $\sigma$ of ${\bar \Pi} _{K_1}$
  with  $f_M(\sigma )=\mu _{\rm
left}\ldots \mu _{\rm right}$ is called {\em principal }  in
$\sigma$ in the following cases.

1) Let $f_X(\sigma )=u_i\ (i\neq n)$, where $u_i=x_iy_{i-1}^{-1}$
for $i>1$ and $u_1=x_1c_m^{-z_m}$ for $m\neq 0$. Then (see Lemma
\ref{main})
\begin{multline*}u_i^{\phi
_{K_1}}=A_{K_2+m+4i}^{*-q_4+1}x_{i+1}^{\phi _{K_2}}y_i^{-\phi
_{K_2}}x_i^{-q_1\phi _{K_2}}\\ \left (x_i^{-\phi
_{K_2}}A_{K_2+m+4i-4}^{*q_0}A_{K_2+m+4i-2}^{*(-q_2+1)}y_i^{\phi
_{K_2}}x_i^{-q_1\phi _{K_2}}\right )^{q_3-1}
A_{K_2+m+4i-4}^{*q_0}.\end{multline*}

The variable $\mu _{\rm right}$  is {\em principal} in $\sigma$ if
and only if  $\mu _{\rm right}^{\alpha}$  almost contains a
cyclically reduced part of
$$\left (x_i^{-\psi
_{K_2}}A_{K_2+m+4i-4}^{*q_0\beta}A_{m+4i-2}^{*(-q_2+1)\beta
}y_i^{\psi _{K_2}}x_i^{-q_1\psi _{K_2}}\right )^{q} =
(x_i^{q_1}y_i)^{\psi _{K_2}}
(A^{*\beta}_{K_2+m+4i-1})^{-q}(y_i^{-1}x_i^{-q_1})^{\psi _{K_2}},$$
for some $q > 2$. Now, the variable  $\mu _{\rm left}$  is {\em
principal} in $\sigma$ if and only if  $\mu _{\rm right}$ is not
principal in $\sigma$.

2) Let $f_X(\sigma )=v_i$, where $v_i=y_iu_i \ (i\neq 1,n)$ and
$v_1=y_1u_1\prod _{j=m-1}^1c_j^{-z_j}$. Then (see formula 3.a)
from Lemma \ref{main})
$$v_i^{\phi _{K_1}}=A_{K_2+m+4i}^{*(-q_4+1)} x_{i+1}^{\phi _{K_2}} A_{K_2+m+4i-4}^{*(-q_0)}
x_i^{q_1\phi _{K_2}}y_i^{\phi _{K_2}} A_{K_2+m+4i-2}^{*(q_2-1)}
A_{K_2+m+4i-4}^{*-1}, \ \ i\neq 1,$$  and

$$v_1^{\phi _{K_1}}=A_{K_2+m+4}^{*(-q_4+1)} x_{2}^{\phi _{K_2}} A_{K_2+2m}^{*(-q_0)}
x_1^{q_1\phi _{K_2}}y_1^{\phi _{K_2}} A_{K_2+m+1}^{*(q_2-1)}x_1\Pi
_{j=n}^1c_j^{-z_j}.$$

\noindent
 A side variable $\mu_{\rm right}$ ($\mu_{\rm left}$) is
{\em principal} in $\sigma$ if and only if  $\mu_{\rm
right}^{\alpha}$ (correspondingly, $\mu_{\rm left}^{\alpha}$)
almost contains $(A^{ \beta}_{K_2+m+4i})^{-q},$ for some $q>2$. In
this case both variables $\mu_{\rm left}^{\alpha}$, $\mu_{\rm
right}^{\alpha}$ could be simultaneously  principal.

 3) Let  $f_X(\sigma )=\bar u_n=u_{n-1}u_n$. Then (by formula 3.c) from Lemma \ref{main})
 $$ \bar u_n^{\phi
_{K_1}}=A_{K_2+m+4n-8}^{*}
A_{K_2+m+4n-6}^{-q_2+1}(y_{n-1}^{-1}x_n^{-q_1})^{\phi
_{K_1}}A_{K_2+m+4n-8}^{*q_0}(x_n^{q_5}y_n)^{\phi
_{K_1}}A_{K_2+m+4n-2}^{*q_6-1}A_{K_2+m+4n-4}^{*-1}.$$

 A side variable $\mu_{\rm right}$ ($\mu_{\rm left}$) is {\em principal}
  in $\sigma$ if $\mu _{\rm right}^{\alpha}$
(correspondingly, $\mu _{\rm left} ^{\alpha}$) almost contains
$(A^{\beta}_{K_2+m+4n-2})^q,$  for some $q>2$. In this case both
variables $\mu_{\rm left}^{\alpha}$, $\mu_{\rm right}^{\alpha}$
could be simultaneously  principal.

\vspace{3mm}
 4)  Let $f_X(\sigma )= y_n$.
Then (by Lemma \ref{le:7.1.xiforms})
 $$y_n^{\phi
_{K_1}}=A_{K_2+m+4n-4}^{*q_0\beta}A_{K_1}^{*q_3\beta}x_n^{q_1\psi_{K_2}}y_1^{\psi_{K_2}}.$$

 The  variable $\mu _{\rm right}$
 ($\mu_{\rm left}$) is {\em principal} in $\sigma$ if $\mu_{\rm
right}^{\alpha}$ (correspondingly, $\mu_{\rm left}^{\alpha}$)
almost contains
$$(A^{\beta}_{K_1})^q,$$ for some $q$ such that $2q>p_{K_1}-2$.
In this case both variables $\mu_{\rm left}^{\alpha}$, $\mu_{\rm
right}^{\alpha}$ could be simultaneously  principal.

 \vspace{3mm}
 5) Let $f_X(\sigma )=z_j$, $j=1,\dots ,m-1$.
Then (by Lemma \ref{le:7.1.zforms}) $$z_j^{\phi _{K_1}}=c_jz_j^{\phi
_{K_2}}A_{K_2+j-1}^{*\beta p_{j-1}}c_{j+1}^{z_{j+1}^{\phi
_{K_2}}}A_{K_2+j}^{*\beta p_j-1}.$$

 A variable $\mu _{\rm left}$
($\mu _{\rm right}$) is {\em principal} if $\mu _{\rm
right}^{\alpha}$ (correspondingly, $\mu _{\rm left} ^{\alpha}$)
almost contains $(A^{\beta}_{K_2+j})^q,$ for some  $|q|>2.$ Both
left and right side variables can be simultaneously principal.

 \vspace{3mm}
 6)
Let $f_X(\sigma )=z_m$. Then (by Lemma \ref{le:7.1.zforms})
$$z_m^{\phi_{K_1}} = c_m^{K_2} z_m^{\phi_{K_2}} A_{K_2+m-1}^{\ast
p_{m-1}} x_1^{-\phi_{K_2}} A_{K_2 +m}^{\ast p_m-1}.$$

 In this case  $\mu _{\rm left}$ is {\em principal} in
$\sigma$ if and only if $\mu_{\rm left}$ is long (i.e.,  it is not
very short), and we define $\mu_{\rm right}$ to be always
non-principal. Observe that if  $\mu_{\rm left}$ is very short then
$\mu _{\rm right}^{\alpha}= f z_m^{\phi _{K_1}\beta}$ for a very
short $f \in F$.

 Let $f_X(\sigma )=z_m^{-1}c_mz_m$. Then (by Lemma \ref{le:7.1.zforms})
$$f_X(\sigma )^{\phi _{K_1}}=A_{K_2+m}^{*-p_{m}+1}x_1^{\phi
_{K_2}}A_{K_2+m}^{*p_{m}}.$$

The variable $\mu _{\rm left}$  is {\em principal} in $\sigma$
 if and only if the following two
conditions hold:  $\mu _{\rm left}^{\alpha}$ almost contains $(A^{
\beta}_{K_2+m})^q,$ for some $q$ with $ |q|>2$; $\mu_{\rm
left}^{-1}\neq fz_m^{\phi _{K_1}\beta}$ for a very short $f\in F$.

  Similarly, the variable $\mu _{\rm right}$  is
{\em principal} in $\sigma$  if and only if the following two
conditions hold:  $\mu _{\rm right}^{\alpha}$ almost contains
$(A^{\beta}_{K_2+m})^q,$ for some $q$ with $ |q|>2$;  $\mu _{\rm
right}^{\alpha}\neq fz_m^{\phi _{K_1}\beta}$ for a very short $f\in
F$.

 Observe, that in this case
the variables $\mu _{\rm left}$ and $\mu _{\rm right}$ can be
simultaneously principal in $\sigma$ and non-principal in $\sigma$.
The latter happens if and only if  $\mu _{\rm
right}^{\alpha}=f_1z_m^{\phi _{K_1}\beta}$ and
 $\mu _{\rm left}^{\alpha}=z_m^{-\phi _{K_1}\beta}f_2$ for
 some very short elements $f_1, f_2 \in F$. Therefore, if both
 $\mu _{\rm left}$ and $\mu _{\rm right}$ are non-principal then
 they can be expressed in terms of $z_m^{\phi _{K_1}}$ and very short variables.

\begin{Claim} \label{claim:one-p} For every interval $\sigma$
of $\bar\Pi _{K_1}$ its partition $f_M(\sigma)$
has at least one principal variable, unless this interval $\sigma$
and its partition $f_M(\sigma)$ are of those two  particular types
described in Case 6).
\end{Claim}
 {\em Proof.}  Let  $f_M(\sigma )=\mu _{\rm
left}\nu _1\ldots\nu _k \mu _{\rm right},$ where $\nu _1\ldots\nu
_k$ are very short variables. Suppose  $A_{r+K_2}$ is the oldest
period such that $f_X(\sigma)$ has $N$-large
$A_{r+K_2}$-decomposition.

If $r \neq 1$ then (see Lemmas \ref{le:7.1.zforms} -
\ref{le:7.1.xiforms})   $A_{r+K_2}$ contains some $N$-large exponent
of $A_{r-1+K_2}$. Therefore $\nu _1^{\alpha}\ldots\nu _k^{\alpha}$
does not contain $A_{r+K_2}'$, hence  either $\mu _{\rm left}$ or
$\mu _{\rm right}$ almost contains $A_{r+K_2}^{\beta q},$ where
$|q|>2$. This finishes all the cases except for the Case 1). In Case
1) a similar argument shows that $\nu _1^{\alpha}\ldots\nu
_k^{\alpha}$ does not contain $A_{r-1+K_2}'$, so one of the side
variables is principal.

If $r=1$, then $A_{1+K_2}$ contains some $N$-large exponent of
$A_{2+K_3}$. Again,  $\nu _1^{\alpha}\ldots\nu _k^{\alpha}$ does not
 contain $A_{1+K_2}'$, because the complexity of the cut equation
$\Pi _{K_1}$ does not changed in ranks from $K_0$ to $K_3$. Now, an
argument similar to the one above finishes the proof.   \hfill
$\Box$

\begin{Claim}
\label{claim:N-N} If both side variables of a partition
$f_M(\sigma)$  of an interval $\sigma$ from $\bar\Pi _{K_1}$  are
non-principal, then they are non-principal in every partition of an
interval from  $\bar\Pi _{K_1}$.
\end{Claim}
 {\em Proof.}  It follows directly from the description of the side
 variables $\mu _{\rm left}$ and $\mu _{\rm right}$ in the Case 6)
 of the definition of principal variables. Indeed, if
 $\mu _{\rm left}$ and $\mu _{\rm right}$  are both non-principal, then (see Case 6))
  each of them is either very short, or it is equal to $f_1 z_m^{\phi _{K_1}\beta}f_2$
  for some  very short $f_1, f_2  \in F$.  Clearly, neither of them
  could be principal in other partitions.
 \hfill $\Box$

\begin{Claim}\label{claim:p-n} Let $n\neq 0$. Then  a side
variable can be principal only in one class of equivalent partitions
of intervals from  $\bar\Pi _{K_1}$.
\end{Claim}
{\em Proof.}  Let  $f_M(\sigma )=\mu _{\rm left}\nu _1\ldots\nu _k
\mu _{\rm right},$ where $\nu _1\ldots\nu _k$ are very short
variables. Suppose  $A_{r+K_2}$ is the oldest period such that
$f_X(\sigma)$ has $N$-large $A_{r+K_2}$-decomposition.

In every case from the definition of principal variables (except for
Case 1))  a principal variable in $\sigma$ almost contains a cube
$(A_{r+K_2}^\prime)^3$.  In Case 1) the principal variable  almost
contains $(A_{r-1+K_2}^\prime)^3$, moreover, if $\mu _{\rm left}$ is
the principal variable then $\mu _{\rm left}^\alpha$  contains an
$N$-large exponent of $(A_{r+K_2}^\prime)$.

We consider only the situation when the partition $f_M(\sigma)$
satisfies Case 1), all other cases can be done similarly.

Clearly,  if $f_X(\sigma) = u_i$ then  a principal variable in
$\sigma$ does not
 appear as a principal variable in the  partition of any other
interval $\delta$ with $f_X(\delta) \neq u_i$, $f_X(\delta) \neq
v_i$. Suppose that a principal variable in $\sigma$ appears as a
principal variable of the partition of $\delta$ with $f_X(\delta) =
u_i$. Then partitions $f_M(\sigma)$ and $f_M(\delta)$ are
equivalent, as required. Suppose now that  a principal variable
$\mu$ in $\sigma$ appears as a principal variable of the partition
of $\delta$ with $f_X(\delta) = v_i$. If $\mu = \mu _{\rm right}$
then it cannot appear as the right principal variable, say
$\lambda_{\rm right}$, of $f_M(\delta)$. Indeed, $\mu _{\rm
right}^\alpha$ ends   (see Case 1) above) with almost all of the
word $(A_{K_2+m+4i-4}^{*q_0})^\beta$ (except, perhaps,  for a short
initial segment of it). But the write principal variable
$\lambda_{\rm right}$ should end (see Case 2) above ) with almost
all of the word $A_{K_2+m+4i-4}^{*-1}$ (except, perhaps,  for a
short initial segment of it), so $\mu _{\rm right} \neq \lambda_{\rm
right}$. Similarly, if the left side variable $\lambda_{\rm left}$
of $f_M(\delta)$ is principal in $\delta$ then $\mu _{\rm right}
\neq \lambda_{\rm left}$.
 Suppose now that  $\mu = \mu _{left}$, then $\mu _{\rm right}$ is not principle
 in $\sigma$, so it is not true that $\mu _{\rm right}$ almost
 contains the  cube of the cyclically reduced part of
$$ x_i^{-\psi
_{K_2}}A_{K_2+m+4i-4}^{*q_0\beta}A_{m+4i-2}^{*(-q_2+1)\beta
}y_i^{\psi _{K_2}}x_i^{-q_1\psi _{K_2}} .$$
 Then $\mu _{\rm left}$
is very long, so  it is easy to see that it does not  appear in the
partition of $\delta$ as a principle variable. This finishes the
Case 1).
   \hfill $\Box$

For the cut equation $\bar\Pi _{K_1}$ we construct a finite graph
$\Gamma = (V,E).$ Every vertex from $V$ is marked by variables
from $\bar M_{side}^{\pm 1}$ and letters
  from the alphabet $\{P,N\}$. Every edge from $E$ is colored either as red or blue.
  The graph $\Gamma$ is constructed as follows. Every partition $f_M(\sigma) = \mu_1 \ldots \mu_k$
   of $\bar\Pi _{K_1}$ gives two
vertices $v_{\sigma,left}$ and $v_{\sigma,right}$ into $\Gamma$,
so
  $$V = \cup_{\sigma} \{v_{\sigma,left},  v_{\sigma,right}\}.$$
We mark $v_{\sigma,left}$ by $\mu_1$ and  $v_{\sigma,right}$ by
$\mu_k$. Now we mark the  vertex $v_{\sigma,left}$
 by a letter $P$ or letter $N$  if $\mu_1$
   is correspondingly principal or non-principal in
 $\sigma$. Similarly, we mark $v_{\sigma,right}$  by $P$ or $N$ if
 $\mu_k$ is principal or non-principal in $\sigma$.

 For every $\sigma$ the
vertices $v_{\sigma,left}$ and $v_{\sigma,right}$ are connected by
a {\em red} edge. Also, we connect by a {\em blue} edge every pair
of
 vertices which are marked by variables $\mu, \nu$ provided $\mu =
 \nu$ or $\mu = \nu^{-1}$. This describes the graph $\Gamma$.

 Below we construct a new graph $\Delta$
  which is obtained from $\Gamma$ by deleting some  blue edges
  according to the following procedure.
Let $B$ be a maximal connected blue component of $\Gamma$, i.e., a
connected component of the graph obtained from $\Gamma$ by
deleting all red edges. Notice, that $B$ is a complete graph, so
every two vertices in $B$ are connected by a blue edge. Fix a
vertex $v$ in $B$ and consider the star-subgraph $Star_B$ of $B$
generated by all edges adjacent to $v$. If $B$ contains a vertex
marked by $P$ then we choose $v$ with label $P$, otherwise $v$ is
an arbitrary vertex of $B$. Now, replace $B$ in $\Gamma$ by the
graph $Star_B$, i.e., delete all edges in $B$ which are not
adjacent to $v$. Repeat this procedure for every maximal blue
component $B$ of $\Gamma$. Suppose that the blue component
corresponds to long bases of case 6) that are non-principal and
equal to $f_1z_m^{\phi _{K_1}}f_2$ for very short $f_1,f_2$. In
this case, we remove all the blue edges that produce  cycles if
the red edge from $\Gamma$ connecting non-principal $\mu _{\rm
left}$ and $\mu _{\rm right}$ is added to the component (if such a
red edge exists).

Denote the resulting graph by $\Delta$.

In the next claim we describe connected components of the graph
$\Delta$.

\begin{Claim}\label{20.}
Let $C$ be a connected component of $\Delta$. Then one of the
following holds:
 \begin{enumerate}
   \item  [(1)]\label{item:4} there is a vertex  in $C$ marked by a variable
   which does not occur as a principal variable in any partition
    of $\bar\Pi _{K_1}$. In particular, any component which
    satisfies one of the following conditions has such a vertex:
  \begin{itemize}
  \item  [a)] there is a vertex  in $C$ marked by a variable which is a short variable in some partition of
   $\bar\Pi _{K_1}$.
   \item [b)] there is a red  edge   in $C$ with both endpoints marked by $N$ (it corresponds to a partition described
   in Case 6 above);
    \end{itemize}
   \item [(2)]\label{item:3} both endpoints of every red edge in $C$ are marked by $P$.
  In this case $C$ is an isolated vertex;
  \item [(3)] \label{item:5} there is a vertex  in $C$  marked by a  variable
  $\mu$ and $N$ and if $\mu$ occurs as a label of an endpoint  of some red edge in $C$
   then the other endpoint of this edge is marked by $P$.
  \end{enumerate}
\end{Claim}

{\em Proof.}   Let $C$ be a connected component of $\Delta$.
Observe first, that if $\mu$ is a short variable  in
   $\bar\Pi _{K_1}$ then $\mu$ is not principle in $\sigma$ for
   any interval $\sigma$ from $\bar\Pi _{K_1}$, so there is no vertex in $C$
   marked by both $\mu$ and $P$.  Also, it follows from  Claim \ref{claim:N-N} that
if there is a red  edge  $e$ in $C$ with both endpoints marked by
$N$, then the variables assigned to endpoints of $e$ are
non-principle in any interval $\sigma$ of  $\bar\Pi _{K_1}$. This
proves the part "in particular" of 1).

  Now assume that  the component $C$  does not satisfy any
of the conditions (1), (2). We need to show that $C$ has type (3).
It follows that every variable which occurs as a label of a vertex
in $C$ is long and it labels, at least, one vertex in $C$ with
label $P$. Moreover,  there are non-principle occurrences of
variables in $C$.

 We summarize some properties of $C$ below:
 \begin{itemize}
  \item  There are no blue edges in $\Delta$ between
vertices with labels $N$ and $N$ (by construction).
 \item  There are no blue edges between
 vertices labelled by $P$ and $P$ (Claim  \ref{claim:p-n}).
  \item   There are no  red edges in
  $C$ between vertices labelled by $N$ and $N$ (otherwise 1) would hold).
   \item  Any reduced path in $\Delta$ consists of
   edges of alternating color (by construction).
 \end{itemize}

We claim that $C$ is a tree. Let $p =  e_1 \ldots e_k$ be a simple
loop in $C$ (every vertex in $p$ has degree 2 and the terminal
vertex of $e_k$ is equal to the starting point of $e_1$).

We show first that $p$ does not have red edges with endpoints
labelled by $P$ and $P$. Indeed, suppose there exists such an edge
in $p$. Taking cyclic permutation of $p$ we may assume that $e_1$
is a red edge with labels $P$ and $P$. Then $e_2$ goes from a
vertex with label $P$ to a vertex with label $N$. Hence the next
red edge $e_3$ goes from $N$ to $P$, etc. This shows that every
blue edge along $p$ goes from $P$ to $N$. Hence the last edge
$e_k$ which must be blue goes from $P$ to $N$ -contradiction,
since all the labels of $e_1$ are $P$.

It follows that both colors of edges and labels of vertices in $p$
alternate. We may assume now that $p$ starts with a vertex with
label $N$ and the first edge $e_1$ is red. It follows that the end
point of $e_1$ is labelled by $N$ and all blue edges go from $N$ to
$P$. Let $e_i$ be a blue edge from $v_i$ to $v_{i+1}$. Then the
variable $\mu_i$ assign to the vertex $v_i$ is principal in the
partition associated with the red edge $e_{i-1}$ , and the variable
$\mu_{i+1} = \mu_i^{\pm 1}$ associated with $v_{i+1}$ is a
non-principal side variable in the partition $f_M(\sigma)$
associated with the red edge $e_{i+1}$. Therefore, the the side
variable $\mu_{i+2}$ associated with the end vertex $v_{i+2}$ is a
principal side variable in the partition $f_M(\sigma)$ associated
with $e_{i+1}$. It follows from the definition of principal
variables that the length of $\mu_{i+2}^\alpha$ is much longer than
the length of $\mu_{i+1}^\alpha$, unless the variable $\mu_i$ is
described in the Case 1). However, in the latter case the variable
$\mu_{i+2}$ cannot occur in any other partition $f_M(\delta)$  for
$\delta \neq \sigma$. This shows that there no blue edges in
$\Delta$ with endpoints labelled by such $\mu_{i+2}$. This implies
that $v_{i+2}$ has degree one in $\Delta$ - contradiction wit the
choice of $p$. This shows that there are no vertices labelled by
such variables described in Case 1).  Notice also, that the length
of variables (under $\alpha$) is preserved along blue edges:
$|\mu_{i+1}^\alpha| = |(\mu_i^{\pm 1})^\alpha| = |\mu_i^\alpha|$.
Therefore,
 $$  |\mu_i^\alpha| = |\mu_{i+1}^\alpha| < |\mu_{i+2}^\alpha|$$
 for every $i$.

It follows that going along $p$ the length of $ |\mu_i^\alpha|$
increases, so $p$ cannot be a loop. This implies that $C$ is a
tree.

Now we are ready to show that the component $C$ has type 3) from
Claim  \ref{20.}. Let $\mu_1$ be a variable assigned to some vertex
$v_1$ in $C$ with label $N$.   If $\mu_1$ satisfies the condition 3)
 from Claim  \ref{20.}  then we are done. Otherwise, $\mu_1$ occurs as a
label of one of $P$-endpoints, say $v_2$  of a red edge $e_2$ in $C$
such that the other endpoint of $e_2$, say $v_3$   is non-principal.
Let $\mu_3$ be the label of $v_3$. Thus $v_1$ is connected to $v_2$
by a blue edge and $v_2$ is connected to $v_3$ by a red edge.   If
$\mu_3$ does not satisfy the condition 3)  from Claim \ref{20.} then
we can repeat the process (with $\mu_3$ in place of $\mu_1$). The
graph $C$ is finite, so  in finitely many steps either we will find
a variable that satisfies 3)  or we will construct a closed reduced
path in $C$.  Since $C$ is a tree the latter does not happen,
therefore $C$ satisfies 3), as required.

 $\Box$

\begin{Claim} \label{21} The graph $\Delta$ is a forest, i.e., it is union of trees.
\end{Claim}
 {\em Proof.}
  Let $C$ be a connected component of $\Delta$. If $C$ has type (3)
  then it is a tree, as has been shown in Claim \ref{20.}
    If $C$ of the type (2) then by Claim \ref{20.}  $C$ is an
    isolated vertex - hence a tree.
    If $C$ is of the type (1) then $C$ is a tree because each interval corresponding to this component has exactly one
    principal variable (except some particular intervals of type 6)
    that do not have principal variables at all and do not produce cycles), and the same long variable cannot be principal in two different
    intervals. Although the same argument as in (3) also
    works here.

 $\Box$

 Now  we define the sets $\bar M_{useless}, \bar M_{free}$ and assign values
 to variables from $\bar M = \bar M_{useless} \cup  \bar M_{free} \cup \bar M_{\rm
 veryshort}$. To do this we use the structure of connected components of
 $\Delta$. Observe first, that all occurrences of a given variable
 from ${\bar M_{sides}}$ are located in the same connected
 component.

 Denote by $\bar M_{free}$ subset of $\bar M$ which consists of variables
of the following types:
 \begin{enumerate}
  \item variables which do not occur as principal in any partition
  of $(\bar\Pi _{K_1})$;
   \item \label{item: 2-2}
   one (but not the other) of the variables  $\mu$ and $\nu$ if they
   are both principal side variables of a partition of the
   type (\ref{item:3}) and such that $\nu \neq \mu^{-1}$.
 \end{enumerate}

Denote by $\bar M_{useless} = \bar M_{side}- \bar M_{free}.$

 \begin{Claim}
For  every $\mu \in \bar M_{\rm useless}$ there exists a word
$V_\mu \in F[X \cup \bar M_{\rm free}\cup \bar M_{\rm veryshort}]$
such that for every map $\alpha_{\rm free} : \bar M_{\rm free}
 \rightarrow F$, and every solution $\alpha_s: F[\bar M_{\rm veryshort}] \rightarrow F$
  of the system $\Delta(\bar M_{\rm veryshort}) = 1$  the map $\alpha: F[\bar M] \rightarrow F$
  defined by
   \[ \mu^\alpha = \left\{\begin{array}{ll}
   \mu^{\alpha_{\rm free}} &\mbox{ if  $\mu \in \bar M_{\rm free}$;}\\
   \mu^{\alpha_{s}} & \mbox{ if $\mu \in \bar M_{\rm veryshort}$;}\\
     \bar V_\mu(X^\delta, \bar M_{\rm free}^{\alpha_{\rm free}},
   \bar M_{\rm veryshort}^{\alpha_s}) & \mbox{ if $\mu \in \bar M_{\rm useless}$.}
   \end{array}
 \right. \]
  is a group solution of $\bar \Pi_{K_1}$ with respect to $\beta$.

 \end{Claim}
 {\em Proof.}   The claim follows from Claims \ref{20.} and
 \ref{21}. Indeed,
take as values of short variables an arbitrary solution $\alpha
_s$ of the system $\Delta (\bar M_{\rm veryshort})=1$. This system
is obviously consistent, and we fix its solution. Consider
connected components of type (1) in Claim \ref{20.}.  If $\mu$ is
a principal variable for some $\sigma$ in such a component, we
express $\mu ^{\alpha}$ in terms of values of very short variables
$\bar M_{\rm veryshort}$ and elements $t^{\psi _{K_1}},$ $t\in X$
that correspond to labels of the intervals. This expression does
not depend on $\alpha _s, \beta $ and tuples $q, p^*.$ For
connected components of $\Delta $ of types (2) and (3) we express
values $\mu ^{\alpha}$  for $\mu\in M_{useless}$ in terms of
values $\nu ^{\alpha}$, $\nu\in M_{\rm free}$  and $t^{\psi
_{K_1}}$ corresponding to the labels of the intervals. $\Box$

We can now finish the proof of Proposition \ref{Or}. Observe, that
$M_{\rm veryshort}\subseteq \bar M_{\rm veryshort}.$ If $\lambda$
is an additional very short variable from $M^*_{\rm veryshort}$
that appears when transformation $T_1$ or $T_2$ is performed,
$\lambda ^{\alpha}$ can be expressed in terms $M_{\rm
veryshort}^{\alpha}$. Also, if a variable $\lambda$ belongs to
$\bar M_{\rm free}$ and does not belong to $M$, then there exists
a variable $\mu\in M$, such that $\mu ^{\alpha}=u^{\psi
_{K_1}}\lambda ^{\alpha},$ where $u\in F(X, C_S)$, and we can
place $\mu$ into $M_{\rm free}.$

Observe, that the argument above is based only on the tuple $p$,
it does not depend on the tuples $p^*$ and $q$. Hence the words
$V_\mu$ do not depend on $p^*$ and $q$.

 The Proposition is proved for $n\neq 0.$
If $n=0$, partitions of the intervals with  labels $z_{n-1}^{\phi
_{K_1}}$ and $z_{n}^{\phi _{K_1}}$ can have equivalent principal
right variables, but in this case the left variables will be
different and do not appear in other non-equivalent partitions.
The connected component of $\Delta$ containing these partitions
will have only four vertices  one blue edge.

In the case $n=0$ we transform equation $\Pi _{K_1}$ applying
transformation $T_1$ to the form when the intervals are labelled
by $u^{\phi _{K_1}},$ where $$ u\in\{z_1,\ldots ,z_m,
c_{m-1}^{z_{m-1}}, z_{m}c_{m-1}^{-z_{m-1}}\}.$$

If $\mu _{\rm left}$ is very short for the interval $\delta$
labelled by $(z_mc_{m-1}^{-z_{m-1}})^{\phi _{K_1}},$ we can apply
$T_2$ to
 $\delta$, and split it into intervals with labels $z_m^{\phi
_{K_1}}$ and $c_{m-1}^{-z_{m-1}^{\phi _{K_1}}}.$ Indeed, even if
we had to replace $\mu_{\rm right}$ by the product of two
variables, the first of them would be very short.

If $\mu _{\rm left}$ is not very short for the interval $\delta$
labelled by \[(z_mc_{m-1}^{-z_{m-1}})^{\phi _{K_1}}= c_mz_m^{\phi
_{K_2}}A_{m-1}^{*p_{m-1}-1},\] we do not split the interval, and
$\mu _{\rm left}$ will be considered as the principal variable for
it. If $\mu _{\rm left}$ is not very short for the interval
$\delta$ labelled by $z_m^{\phi _{K_1}}=z_m^{\phi
_{K_2}}A_{m-1}^{*p_{m-1}}$, it is a principal variable, otherwise
$\mu _{\rm right}$ is principal.

If an interval $\delta$ is labelled by $(c_{m-1}^{z_{m-1}})^{\phi
_{K_1}}=A_{m-1}^{*-p_{m-1}+1}c_m^{-z_m^{\phi
_{K_2}}}A_{m-1}^{*p_{m-1}},$ we consider $\mu _{\rm right}$
principal if $\mu _{\rm right}^{\alpha}$ ends with
$(c_m^{-z_m^{\phi _{K_2}}}A_{m-1}^{*p_{m-2}})^{\beta},$ and the
difference is not very short. If $\mu _{\rm left}^{\alpha}$ is
almost $z_m^{-\phi _k\beta}$ and $\mu _{\rm right}^{\alpha}$ is
almost $z_m^{\phi _k\beta}$, we do not call any of the side
variables principal. In all other cases $\mu _{\rm left}$ is
principal.

Definition of the principal variable in the interval with label
$z_i^{\phi _{K_1}}$,  $i=1,\dots m-2$ is the same as in 5) for
$n\neq 0.$

A variable can be principal only in one  class of equivalent
partitions. All the rest of the proof is the same as for $n>0.$

 \hfill $\Box$

 Now we continue the proof of Theorem \ref{1,2,3,4}.
 Let  $L = 2K + \kappa(\Pi)4K$ and
 $$\Pi_\phi= \Pi_L\rightarrow \Pi_{L-1}\rightarrow\ldots  \rightarrow \ldots  $$
be the sequence of $\Gamma$-cut equations (\ref{eq:7.3.15}). For a
$\Gamma$-cut equation $\Pi_j$ from (\ref{eq:7.3.15}) by $M_j$ and
$\alpha_j$ we denote the corresponding set of variables and the
solution relative to $\beta$.

By Claim \ref{5.} in the sequence (\ref{eq:7.3.15}) either there
is $3K$-stabilization at $K(r+2)$ or $Comp(\Pi_{K(r+1)}) = 0$.

Case 1. Suppose there is $3K$-stabilization at $K(r+2)$ in the
sequence (\ref{eq:7.3.15}).

By Proposition \ref{Or} the set of variables $M_{K(r+1)}$ of the
cut equation $\Pi_{K(r+1)}$ can be partitioned into three subsets
 $$M_{K(r+1)} = M_{\rm veryshort} \cup M_{\rm free} \cup M_{\rm
 useless}$$
  such that there exists a finite consistent system of equations $\Delta(M_{\rm veryshort}) = 1$
   over $F$ and words $V_\mu \in F[X, M_{\rm free}, M_{\rm
   veryshort}]$, where $\mu \in M_{\rm useless}$, such that
for every solution $\delta \in {\mathcal B}$, for every map
$\alpha_{\rm free} : M_{\rm free}
 \rightarrow F$, and every solution $\alpha_{short}: F[M_{\rm veryshort}] \rightarrow F$
  of the system $\Delta(M_{\rm veryshort}) = 1$  the map $\alpha_{K(r+1)}: F[M] \rightarrow F$
  defined by
   \[ \mu^{\alpha_{K(r+1)}} = \left\{\begin{array}{ll}
   \mu^{\alpha_{\rm free}} &\mbox{ if  $\mu \in M_{\rm free}$;}\\
   \mu^{\alpha_{short}} & \mbox{ if $\mu \in M_{\rm veryshort}$;}\\
     V_\mu(X^\delta, M_{\rm free}^{\alpha_{\rm free}},
   M_{\rm veryshort}^{\alpha_s}) & \mbox{ if $\mu \in M_{\rm useless}$}
   \end{array}
 \right. \]
  is a group solution of $\Pi_{K(r+1)}$ with respect to $\beta$.
Moreover, the words $V_\mu$  do not depend on tuples $p^*$ and
$q$.

By Claim \ref{2.0.0} if $\Pi = ({\mathcal E}, f_X, f_M)$ is
 a $\Gamma$-cut equation and  $\mu \in M$ then
there exists a word ${\mathcal M}_\mu (M_{T(\Pi)}, X)$ in the free
group $F[M_{T(\Pi)} \cup X]$ such that
$$\mu^{\alpha_{\Pi}} = {\mathcal M}_\mu (M_{T(\Pi)}^{{\alpha_{T(\Pi)}}},
X^{\phi_{K(r+1)}})^\beta, $$  where $\alpha_\Pi$ and
$\alpha_{T(\Pi)}$ are the corresponding solutions of $\Pi$ and
$T(\Pi)$ relative to $\beta$.

Now, going along the sequence (\ref{eq:7.3.15})  from
$\Pi_{K(r+1)}$ back to the cut equation $\Pi_L$ and using
repeatedly the remark above for each $\mu \in M_L$ we obtain a
word
$${\mathcal M^\prime}_{\mu,L} (M_{K(r+1)}, X^{\phi_{K(r+1)}}) = {\mathcal M^\prime}_{\mu,L}
(M_{\rm useless},M_{\rm free}, M_{\rm veryshort},
X^{\phi_{K(r+1)}})$$
 such that
$$\mu^{\alpha_L} = {\mathcal M^\prime}_{\mu,L} (M_{K(r+1)}^{\alpha_{K(r+1)}}, X^{\phi_{K(r+1)}})^{\beta}.$$
Let $\delta = \phi_{K(r+1)} \in {\mathcal B}$ and put
$${\mathcal M}_{\mu,L}(X^{\phi_{K(r+1)}}) = {\mathcal M^\prime}_{\mu,L} (V_\mu(X^{\phi_{K(r+1)}},M_{\rm free}^{\alpha_{\rm free}},
M_{\rm veryshort}^{\alpha_{short}}),M_{\rm free}^{\alpha_{\rm
free}}, M_{\rm veryshort}^{\alpha_{short}}, X^{\phi_{K(r+1)}}).$$
 Then for every $\mu \in M_L$
$$\mu^{\alpha_L}  = {\mathcal M}_{\mu,L}(X^{\phi_{K(r+1)}})^\beta$$
If we denote by ${\mathcal M}_{L}(X)$ a tuple of words
$${\mathcal M}_{L}(X) = ({\mathcal M}_{\mu_1,L}(X), \ldots,{\mathcal
M}_{\mu_{|M_L|},L}(X)),$$
  where $\mu_1, \ldots, \mu_{|M_L|}$ is some fixed ordering of $M_L$
  then
    $$M_L^{\alpha_L} ={\mathcal M}_{L}(X^{\phi_{K(r+1)}})^\beta.$$
Observe, that  the words  ${\mathcal M}_{\mu,L}(X)$, hence
${\mathcal M}_{L}(X)$ ( where $ X^{\phi_{K(r+1)}}$  is replaced by
$X$) are the same for  every $\phi_L \in {\mathcal B}_{p,q}$.

It follows from property c) of the cut equation $\Pi_\phi$ that
the solution ${\alpha_L}$ of $\Pi_\phi$ with respect to $\beta$
gives rise to a group solution of the original cut equation
$\Pi_{\mathcal L}$ with respect to $\phi_L\circ \beta$.

Now, property c) of the initial cut equation $\Pi_{\mathcal L} =
(\mathcal E,f_X,f_{M_L})$
 insures that for every $\phi_L \in {\mathcal B}_{p,q}$ the pair
$(U_{\phi_L\beta},V_{\phi_L\beta})$ defined by
  $$U_{\phi_L\beta} =  Q(M_L^{\alpha_L}) =  Q({\mathcal M}_L(X^{\phi_{K(r+1)}}))^\beta, $$
$$V_{\phi_L\beta}= P(M_L^{\alpha_L}) = P({\mathcal
M}_L(X^{\phi_{K(r+1)}}))^\beta.$$
 is a solution  of the  system$S(X) = 1 \wedge T(X,Y) = 1.$

  We claim that
 $$Y(X) = P({\mathcal M}_L(X))$$
  is a solution of the equation $T(X,Y) = 1$ in
$F_{R(S)}$.  By Theorem \ref{cy2} ${\mathcal B}_{p,q,\beta}$ is a
discriminating family of solutions for the group $F_{R(S)}$. Since
$$T(X,Y(X))^{\phi\beta} = T(X^{\phi\beta},Y(X^{\phi\beta})) = T(X^{\phi\beta}, {\mathcal
M}_L(X^{\phi\beta})) = T(U_{\phi_L\beta}, V_{\phi_L\beta}) = 1$$
for any ${\phi\beta} \in {\mathcal B}_{p,q,\beta}$ we deduce that
$T(X,Y_{p,q}(X)) = 1$ in $F_{R(S)}$.

Now we need to  show that $T(X,Y) = 1$ admits a complete $S$-lift.
 Let $W(X,Y) \neq 1$ be an inequality such that  $T(X,Y) = 1 \wedge W(X,Y) \neq 1$
  is compatible with $S(X) = 1$. In this event, one may assume
  (repeating the argument from the beginning of this section) that the set
   $$
   \Lambda = \{(U_\psi,V_\psi) \mid \psi \in {\mathcal L}_2
   \}$$
   is such that every pair
 $(U_\psi,V_\psi) \in \Lambda$ satisfies the formula $T(X,Y) = 1 \wedge W(X,Y) \neq
 1.$  In this case, $W(X,Y_{p,q}(X)) \neq
 1$ in $F_{R(S)}$, because its image in $F$ is non-trivial:
 $$W(X,Y_{p,q}(X))^{\phi\beta} = W(U_\psi,V_\psi) \neq 1.$$
 Hence $T(X,Y) =1$ admits a complete lift into  generic
point of $S(X) = 1$.

Case 2. A similar argument applies when $Comp(\Pi_{K(r+2)}) = 0$.
Indeed, in this case for every $\sigma \in {\mathcal E}_{K(r+2)}$
the word $f_{M_{K(r+1)}}(\sigma)$ has length one, so
$f_{M_{K(r+1)}}(\sigma) =\mu$ for some $\mu \in M_{K(r+2)}.$ Now
one can replace the word $V_\mu \in F[X \cup M_{\rm free} \cup
M_{\rm veryshort}]$ by the label $f_{X_{K(r+1)}}(\sigma)$ where
$f_{M_{K(r+1)}}(\sigma) =\mu$ and then repeat the argument.

 \hfill$\Box$

\subsection{Non-orientable quadratic equations}
 Consider now the equation
\begin{equation}\label{8}
\prod_{i=1}^{m}z_i^{-1}c_iz_i\prod_{i=1}^{n}x_i^2  = c_1\ldots
c_m\prod_{i=1}^{n}a_i^2,
\end{equation}
where $a_i, c_j$ give a solution in general position (in all the
cases when it exists). We will now prove Theorem \ref{1,2,3,4} for
a  regular standard non-orientable quadratic equation over $F$.

{\em Let  $S(X,A)=1$ be  a regular standard non-orientable
quadratic equation over $F$. Then every equation $T(X,Y,A) = 1$
compatible with $S(X,A)=1$ admits an complete $S$-lift.}

The proof of the theorem is similar to the proof in the orientable
case, but the basic sequence of automorphisms is different. We
will give a sketch of the proof in this section.

It is more convenient to consider a non-orientable equation in the
form

\begin{equation}\label{8'}
S=\prod_{i=1}^{m}z_i^{-1}c_iz_i\prod_{i=1}^{n}[x_i,y_i]x_{n+1}^2 =
c_1\ldots c_m\prod_{i=1}^{n}[a_i,b_i]a_{n+1}^2,
\end{equation}

or
\begin{equation}\label{8''}
S=\prod_{i=1}^{m}z_i^{-1}c_iz_i\prod_{i=1}^{n}[x_i,y_i]x_{n+1}^2x_{n+2}^2
= c_1\ldots c_m\prod_{i=1}^{n}[a_i,b_i]a_{n+1}^2a_{n+2}^2.
\end{equation}

Without loss of generality we consider equation (\ref{8''}).
 We define a basic sequence
$$\Gamma =(\gamma_1, \gamma_2, \ldots, \gamma_{K(m,n)})$$
 of $G$-automorphisms of the free $G$-group $G[X]$  fixing the left side of the equation (\ref{8''}).

We assume that
 each $\gamma \in \Gamma$    acts identically on all the generators
from $X$ that are not mentioned in the description of $\gamma$.

Automorphisms $\gamma _i,\ i=1,\ldots ,m+4n-1$ are the same as in
the orientable case.

\medskip \noindent
Let $ n = 0$. In this case $K=K(m,0) = m+2.$ Put

\medskip
 $\gamma _{m}: \
z_m\rightarrow  z_m(c_m^{z_m}x_1^2),\ \ \
 x_1\rightarrow x_1^{(c_m^{z_m}x_1^2)};$

\smallskip
$\gamma _{m+1}:\ x_1\rightarrow x_1(x_1x_{2}),\ \ \
x_{2}\rightarrow (x_1x_{2})^{-1}x_{2};$

\smallskip
$\gamma _{m+2}:\ x_1\rightarrow x_1^{(x_1^2x_{2}^2)},\ \ \
x_{2}\rightarrow x_{2}^{(x_1^2x_{2}^2)}.$

\medskip \noindent
Let $ n\geq 1$. In this case $K=K(m,n) = m+4n+2.$ Put

\medskip
 $\gamma _{m+4n}: \
x_n\rightarrow (y_nx_{n+1}^2)^{-1}x_n,\ \ \
 y_n\rightarrow y_n^{(y_nx_{n+1}^2)},\ \ \
 x_{n+1}\rightarrow x_{n+1}^{(y_nx_{n+1}^2)};$

\smallskip
$\gamma _{m+4n+1}:\ x_{n+1}\rightarrow x_{n+1}(x_{n+1}x_{n+2}),\ \
\ x_{n+2}\rightarrow (x_{n+1}x_{n+2})^{-1}x_{n+2};$

\smallskip
$\gamma _{m+4n+2}:\ x_{n+1}\rightarrow
x_{n+1}^{(x_{n+1}^2x_{n+2}^2)},\ \ \ x_{n+2}\rightarrow
x_{n+2}^{(x_{n+1}^2x_{n+2}^2)}.$

These automorphisms induce automorphisms on $G_S$ which we denote
by the same letters.

  Let $\Gamma = (\gamma_1, \ldots, \gamma_K)$ be the
 basic sequence of automorphisms for $S = 1$. Denote by
$\Gamma_{\infty}$ the
 infinite periodic sequence with period $\Gamma$, i.e.,
 $\Gamma_{\infty} = \{\gamma_i\}_{i \geq 1}$ with $\gamma_{i+K} =
\gamma_i$.
 For $j \in {\bf N}$ denote by $\Gamma_j$ the initial segment of
 $\Gamma_{\infty}$ of length $j$. Then for a given $j$ and  $p \in
{\bf N}^j$ put
 $$ \phi_{j,p} =
\stackrel{\leftarrow}{\Gamma}_j^{\stackrel{\leftarrow}{p}}.$$

Let $$\Gamma _P=\{\phi _{j,p}|p\in P\}.$$ We can prove the
analogue of Theorem \ref{cy2}, namely, that  a family of
homomorphisms $\Gamma _P\beta$ from $G_S=G_{R(S)}$ onto $G$, where
$\beta$ is a solution in general position, and $P$ is unbounded,
is a discriminating family.

The rest of the proof is the same as in the orientable case.

\subsection{Implicit function theorem: NTQ systems} \label{se:7.4}

\begin{df}
Let $G$ be a group with a generating set $A$. A system of
equations $S = 1$  is called {\em triangular quasiquadratic}
(shortly, TQ) if it can be partitioned into the following
subsystems

\medskip
$S_1(X_1, X_2, \ldots, X_n,A) = 1,$

\medskip
$\ \ \ \ \ S_2(X_2, \ldots, X_n,A) = 1,$

$\ \ \ \ \ \ \ \ \ \  \ldots$

\medskip
$\ \ \ \ \ \ \ \ \ \ \ \ \ \ \ \ S_n(X_n,A) = 1$

\medskip \noindent
 where for each
$i$ one of the following holds:
\begin{enumerate}
\item [1)] $S_i$ is quadratic  in variables $X_i$;
 \item [2)] $S_i= \{[y,z]=1, [y,u]=1 \mid y, z \in X_i\}$ where $u$ is a
group word in $X_{i+1} \cup  \ldots \cup X_n \cup A$ such that its
canonical image  in $G_{i+1}$ is not a proper power. In this case
we say that $S_i=1$ corresponds to an extension of a centralizer;
 \item [3)] $S_i= \{[y,z]=1 \mid y, z \in X_i\}$;
 \item [4)] $S_i$ is the empty equation.
  \end{enumerate}

Define $G_{i}=G_{R(S_{i}, \ldots, S_n)}$ for $i = 1, \ldots, n$
and put $G_{n+1}=G.$ The  TQ system $S = 1$ is called {\em
non-degenerate} (shortly, NTQ) if each system  $S_i=1$, where
$X_{i+1}, \ldots, X_n$ are viewed as the corresponding constants
from $G_{i+1}$ (under the canonical maps $X_j \rightarrow
G_{i+1}$, $j = i+1, \ldots, n$, has a solution in $G_{i+1}$. The
coordinate group of an NTQ system is called an {\em NTQ group}.

 An NTQ system $S = 1$ is called {\em regular} if for each $i$ the system $S_i=1$ is either of the type 1) or
 4), and in the former case the quadratic equation $S_i=1$ is in
 standard form and regular (see Definition \ref{regular}).
\end{df}

 One of the results to be proved in this section is the following.

\begin{theorem} \label{reg} Let $U(X,A)=1$ be a regular NTQ-system. Every equation $V(X,Y,A)=1$ compatible with
$U=1$ admits a complete $U$-lift.\end{theorem}

{\em Proof.}   We use induction on the number $n$ of levels in the
system $U=1$. We construct a solution tree $T_{sol}(V(X,Y,A)\wedge
U(X,Y))$ with parameters $X=X_1\cup\ldots\cup X_n.$ In the
terminal vertices of the tree there are generalized equations
$\Omega _{v_1},\ldots ,\Omega _{v_k}$ which are equivalent to cut
equations $\Pi _{v_1},\ldots ,\Pi _{v_k}$.

If $S_1(X_1,\ldots ,X_n)=1$ is an empty equation, we can take
Merzljakov's words (see Theorem \ref{merzl}) as values of
variables from $X_1$, express $Y$ as functions in $X_1$ and a
solution of some $W(Y_1,X_2,\ldots ,X_n)=1$ such that for any
solution of the system

$\ \ \ \ \ S_2(X_2, \ldots, X_n,A) = 1$

$\ \ \ \ \ \ \ \ \ \  \ldots$

\medskip
$\ \ \ \ \ \ \ \ \ \ \ \ \ \ \ \ S_n(X_n,A) = 1$

\noindent equation $W=1$ has a solution.

Suppose, now that $S_1(X_1,\ldots ,X_n)=1$ is a regular quadratic
equation. Let $\Gamma $ be a basic sequence of automorphisms for
the equation $S_1(X_1,\ldots ,X_n,A)=1.$
 Recall that $$\phi_{j,p} = \gamma_j^{p_j} \ldots \gamma_1^{p_1} =
 \stackrel{\leftarrow}{\Gamma}_j^p,$$
 where $j \in {\mathbb N}$,  $\Gamma_j = (\gamma_1,
 \ldots, \gamma_j)$ is the initial subsequence of length $j$ of the
sequence
 $\Gamma^{(\infty)}$, and $p = (p_1, \ldots,p_j) \in {\mathbb N}^j$.
  Denote by $\psi_{j,p}$ the following solution of \newline $S_1(X_1) = 1$:
  $$\psi_{j,p}=\phi_{j,p}\alpha,$$
where $\alpha$   is a composition of a solution of $S_1=1$ in
$G_2$ and a solution from a generic family of solutions of the
system

$\ \ \ \ \ S_2(X_2, \ldots, X_n,A) = 1$

$\ \ \ \ \ \ \ \ \ \  \ldots$

\medskip
$\ \ \ \ \ \ \ \ \ \ \ \ \ \ \ \ S_n(X_n,A) = 1$

\noindent in $F(A).$ We can always suppose that $\alpha$ satisfies
a small cancellation condition with respect to $\Gamma .$

  Set $$ \Phi = \{\phi_{j,p} \mid j \in {\mathbb N}, p \in {\mathbb N}^j \}$$
  and let ${\mathcal L}^{\alpha }$ be an  infinite subset of $\Phi ^{\alpha}$ satisfying one of the cut equations above. Without loss of generality we can suppose it satisfies $\Pi _1$.
By Proposition \ref{Or} we can express variables from $Y$ as
functions of the set of $\Gamma $-words in $X_1$, coefficients,
variables $M_{free}$ and variables $M_{veryshort}$, satisfying the
system of equations $\Delta (M_{veryshort})$   The system $\Delta
(M_{veryshort})$ can be turned into a generalized equation with
parameters $X_2\cup\ldots \cup X_n,$ such that for any solution of
the system

$\ \ \ \ \ S_2(X_2, \ldots, X_n,A) = 1$

$\ \ \ \ \ \ \ \ \ \  \ldots$

\medskip
$\ \ \ \ \ \ \ \ \ \ \ \ \ \ \ \ S_n(X_n,A) = 1$

the system $\Delta (M_{veryshort})$ has a solution. Therefore, by
induction, variables $(M_{veryshort})$ can be found as elements of
$G_2$, and variables $Y$ as elements of $G_1$. $\Box$

\begin{lm}
All stabilizing automorphisms (see \cite{Gr})) of the left side of
the equation
\begin{equation}\label{10}
c_1^{z_1}c_2^{z_2}(c_1c_2)^{-1}=1
\end{equation}
have the form $z_1^{\phi}=c_1^kz_1(c_1^{z_1}c_2^{z_2})^n,
z_2^{\phi}=c_2^mz_2(c_1^{z_1}c_2^{z_2})^n.$ All stabilizing
automorphisms of the left side of the equation
\begin{equation}\label{11}
x^2c^{z}(a^2c)^{-1}=1\end{equation} have the form
$x^{\phi}=x^{(x^2c^z)^n}, z^{\phi}=c^kz(x^2c^z)^n$. All
stabilizing automorphisms  of the left side of the equation
\begin{equation}\label{12}
x_1^2x_2^2(a_1^2a_2^2)^{-1}=1\end{equation} have the form
$x_1^{\phi}=(x_1(x_1x_2)^m)^{(x_1^2x_2^2)^n},
x_2^{\phi}=((x_1x_2)^{-m}x_2)^{(x_1^2x_2^2)^n}.$
\end{lm}
{\em Proof.}    The computation of the automorphisms can be done
by utilizing the Magnus software system.
$\Box$


If a quadratic equation $S(X) = 1$ has only commutative solutions
then the radical $R(S)$ of $S(X)$ can be described (up to a linear
change of variables) as follows (see \cite{KMNull}):
 $$Rad(S)=ncl\{[x_i,x_j], [x_i,b],  \mid \ i,j =1,\ldots ,k\},$$
 where $b$ is an element (perhaps, trivial) from $F$. Observe, that if $b$ is not
 trivial then $b$ is not a proper power in $F$. This shows that $S(X) = 1$ is
 equivalent to the system
 \begin{equation}
\label{eq:U}
 U_{com}(X) =  \{[x_i,x_j]=1, [x_i,b]=1,  \mid \ i,j =1,\ldots ,k\}.
\end{equation}
The system $U_{com}(X) = 1$ is equivalent to a single equation,
which we also denote by $U_{com}(X)=1.$
 The coordinate group $H = F_{R(U_{com})}$ of the system $U_{com}= 1$, as well
 as of the corresponding  equation, is $F$-isomorphic to the free extension of the
 centralizer $C_F(b)$ of rank $n$. We need the following notation to deal with $H$.
 For a set $X$ and  $b \in F$    by  $A(X)$ and $A(X,b)$ we denote
  free abelian groups with  basis $X$ and $X \cup \{b\}$,  correspondingly. Now,
  $H \simeq F \ast_{b = b} A(X,b)$. In particular, in the case when $b = 1$ we
  have $H = F \ast A(X)$.

 \begin{lm}
\label{simple} Let $F = F(A)$ be a non-abelian free group and
$V(X,Y,A) = 1$, $W(X,Y,A) = 1$ be equations over $F$.  If a
formula
 $$\Phi=\forall X (U_{com}(X) = 1\rightarrow \exists Y (V(X,Y,A)=1 \wedge W(X,Y,A)\not
=1))$$
 is true in $F$ then there exists a finite number of $<b>$-embeddings
 $\phi_k: A(X,b) \rightarrow A(X,b) \ ( k \in K)$   such that:

 1)  every formula
 $$\Phi _k = \exists Y (V(X^{\phi _k}, Y,A)=1 \wedge W(X^{\phi _k},Y,A)\not =1)$$
 holds in the coordinate group $H = F\ast_{b = b} A(X,b)$;

 2)  for any solution  $\lambda: H \rightarrow F$  there exists a solution
 $\lambda^{\ast} : H \rightarrow F$ such that $\lambda = \phi_k
 \lambda^{\ast}$ for some $k \in K$.
 \end{lm}
{\em Proof.}   We  construct a set of initial parameterized
generalized equations ${\mathcal GE}(S) = \{\Omega_1, \ldots,
\Omega_r\}$ for $V(X,Y,A) = 1$ with respect to the set of
parameters $X$.
 For each $\Omega \in {\mathcal GE}(S)$ in  Section \ref{se:5.5} we
constructed  the finite tree $T_{sol}(\Omega)$ with respect to
parameters $X$. Observe, that non-active part
$[j_{v_0},\rho_{v_0}]$ in the root equation $\Omega =
\Omega_{v_0}$ of the tree $T_{sol}(\Omega)$ is  partitioned into a
disjoint union of closed sections  corresponding to $X$-bases and
constant bases (this follows from the construction of the initial
equations in the set  ${\mathcal GE}(S)$). We label every closed
section $\sigma$ corresponding to a variable  $x \in X^{\pm 1}$ by
$x$, and every constant section corresponding to a constant $a$ by
$a$.  Due to our construction of the tree $T_{sol}(\Omega)$ moving
along a brunch $B$ from the initial vertex $v_0$ to a terminal
vertex $v$ we transfer all the bases from the non-parametric part
into parametric part until, eventually,  in  $\Omega_v$ the whole
interval consists of the parametric part. For a terminal vertex
$v$ in $T_{sol}(\Omega )$ equation $\Omega _v$ is periodized (see
Section 5.4). We can consider the correspondent periodic structure
$\mathcal P$ and the subgroup $\tilde Z_2$. Denote the cycles
generating this subgroup by $z_1,\ldots ,z_m$.  Let $x_i=b^{k_i}$
and $z_i=b^{s_i}$. All $x_i$'s are cycles, therefore the
corresponding system of  equations can be written as a system of
linear equations with integer coefficients in variables
$\{k_1,\ldots ,k_n\}$ and variables $\{s_1,\ldots ,s_m\}:$
\begin{equation}\label{abel}
k_i=\sum _{j=1}^m \alpha _{ij}s_j +\beta _i,\ i=1,\ldots ,n.
\end{equation}

We can always suppose $m\leq n$ and at least for one equation
$\Omega _v$ $m=n$, because otherwise the solution set of the
irreducible system $U_{com}=1$ would be represented as a union of
a finite number of proper subvarieties.

We will show now that all the tuples $(k_1,\ldots ,k_n)$ that
correspond to some system (\ref{abel}) with $m<n$ (the dimension
of the subgroup $H_v$ generated by $\bar k-\bar\beta=k_1-\beta
_1,\ldots ,k_n-\beta _n$ in this case is less than $n$), appear
also in the union of  systems (\ref{abel}) with $m=n$. Such
systems have form $\bar k-\bar\beta _{q} \in H_q$, $q$ runs
through some finite set $Q$, and where $H_q$ is a subgroup of
finite index in $Z^n=<s_1>\times\ldots\times <s_n>$. We use
induction on $n$.
 If for some terminal vertex $v$, the system (\ref{abel})
has $m<n$, we can suppose without loss of generality  that the set
of tuples $H$ satisfying this system  is defined by the equations
$k_r=\ldots ,k_n=0$. Consider just the case $k_n=0$. We will show
that all the tuples $\bar k_0=(k_1,\ldots ,k_{n-1},0)$ appear in
the systems (\ref{abel}) constructed for the other terminal
vertices with $n=m$. First, if $N_q$ is the  index of the subgroup
$H_q,$ $N_q\bar k\in H_q$ for each tuple $\bar k$. Let $N$ be the
least common multiple of $N_1,\ldots ,N_Q$. If a tuple
$(k_1,\ldots ,k_{n-1},tN)$ for some $t$ belongs to $\bar\beta
_q+H_q$ for some $q$, then $(k_1,\ldots ,k_{n-1},0)\in
 \bar\beta _q+H_q$, because $(0,\ldots ,0,tN)\in H_q$. Consider the set $K$ of all tuples $(k_1,\ldots ,k_{n-1},0)$
such that $(k_1,\ldots ,k_{n-1},tN)\not\in  \bar\beta _q+H_q$ for
any $q=1,\ldots ,Q$ and $t\in {\mathbb Z}$ . The set
$\{(k_1,\ldots ,k_{n-1},tN)|(k_1,\ldots ,k_{n-1},0)\in K, t\in
{\mathbb Z}\}$ cannot be a discriminating set for $U_{comm}=1$.
Therefore it satisfies some proper equation. Changing variables
$k_1,\ldots ,k_{n-1}$ we can suppose that for an irreducible
component the equation has form
 $k_{n-1}=0$. The contradiction arises from the fact that we cannot obtain a discriminating set for $U_{comm}=1$ which does not belong to
$\bar\beta _q+H_q$ for any $q=1,\ldots ,Q.$

Embeddings $\phi _k$ are given by the systems (\ref{abel}) with
$n=m$ for generalized equations $\Omega _v$ for all terminal
vertices $v$. $\Box$

There are two more important generalizations of the implicit
function theorem, one -- for arbitrary NTQ-systems, and another --
 for arbitrary systems. We need a few more definitions  to explain
 this. Let $U(X_1, \ldots, X_n,A) = 1$ be an NTQ-system:
$$
\begin{array}{rr}
S_1(X_1, X_2, \ldots, X_n,A) & = 1 \\
S_2(X_2, \ldots, X_n,A) &  = 1\\
\vdots  &  \\
S_n(X_n,A) &  = 1
\end{array}
$$
and $G_{i}=G_{R(S_{i}, \ldots, S_n)}$, $G_{n+1} = F(A)$.

A $G_{i+1}$-automorphism $\sigma$ of $G_i$ is called a {\em
canonical automorphism }
 if the following holds:
 \begin{enumerate}
  \item [1)] if $S_i$ is quadratic  in variables $X_i$ then $\sigma$ is
  induced by a $G_{i+1}$-automorphism of the  group
$G_{i+1}[X_i]$ which fixes $S_i$;
 \item [2)] if $S_i= \{[y,z]=1, [y,u]=1 \mid y, z \in X_i\}$ where $u$ is a
group word in $X_{i+1} \cup  \cdots \cup X_n \cup A$,  then $G_i=
G_{i+1} \ast_{u = u} Ab(X_i\cup \{u\})$, where $Ab(X_i\cup \{u\})$
is a free abelian group with basis $X_i\cup \{u\}$, and in this
event $\sigma$ extends an automorphism of $Ab(X_i\cup \{u\})$
(which fixes $u$);
 \item [3)] If $S_i= \{[y,z]=1 \mid y, z \in X_i\}$ then $G_i =G_{i+1} \ast  Ab(X_i)$,  and in this event $\sigma$ extends an
automorphism of $Ab(X_i)$;
 \item [4)] If $S_i$ is the empty equation then $G_i = G_{i+1}[X_i]$, and in this
 case $\sigma$ is just the identity automorphism of $G_i$.
\end{enumerate}

Let $\pi_i$   be a fixed $G_{i+1}[Y_{i}]$-homomorphism
$$\pi_i : G_i[Y_i]  \rightarrow G_{i+1}[Y_{i+1}],$$
where $\emptyset = Y_1 \subseteq Y_2 \subseteq \ldots \subseteq
Y_n \subseteq Y_{n+1}$ is an ascending chain of finite sets of
parameters, and $G_{n+1} = F(A)$.  Since the system $U = 1$ is
non-degenerate such homomorphisms $\pi_i$ exist. We assume also
that if $S_i(X_i) = 1$ is a standard quadratic equation (the case
1) above) which has a non-commutative solution in $G_{i+1}$, then
$X^{\pi_i}$ is a non-commutative  solution  of $S_i(X_i) = 1$ in
$G_{i+1}[Y_{i+1}].$

A {\em fundamental sequence} (or a {\em fundamental set}) of
solutions of the system $U(X_1,\dots ,X_n,A)=1$ in $F(A)$ with
respect to the  fixed homomorphisms $\pi_1, \ldots, \pi_n$ is a
set of all solutions of $U = 1$ in $F(A)$ of the form
$$\sigma _1\pi _1\cdots \sigma _n\pi _n\tau ,$$
where  $\sigma_i$ is $Y_i$-automorphism of $G_i[Y_i]$ induced by a
canonical automorphism of $G_i$, and  $\tau$ is an
$F(A)$-homomorphism $\tau: F(A \cup Y_{n+1}) \rightarrow F(A)$.
 Solutions
from a given fundamental set of $U$ are called {\em fundamental }
solutions.

Below we describe two  useful  constructions.  The first one is a
{\it normalization} construction which allows one to rewrite
effectively an NTQ-system $U(X) = 1$ into a  normalized NTQ-system
$U^*=1$. Suppose we have an NTQ-system $U(X)=1$ together with a
fundamental sequence of solutions which we denote $\bar V(U)$.

 Starting from the bottom we
replace each non-regular quadratic equation $S_i = 1$ which has a
non-commutative solution by a system of equations effectively
constructed as follows.

1)  If $S_i = 1$  is in  the form
 $$c_1^{x_{i1}}c_2^{x_{i2}}=c_1c_2,$$
  where $[c_1,c_2]\not =1$,  then we replace it  by a system
  $$\{x_{i1}=z_1 c_1z_3,
x_{i2}=z_2 c_2z_3, [z_1,c_1]=1, [z_2,c_2]=1, [z_3,c_1c_2]=1\}.$$

 2) If $S_i = 1$  is in  the form
 $$x_{i1}^2c^{x_{i2}}= a^2c,$$
  where $[a,c]\not = 1$, we replace it by  a system
   $$\{x_{i1}= a^{z_1}, x_{i2}=z_2 cz_1, [z_2, c]=1,[z_1,a^2c]=1\}.$$

 3)  If $S_i = 1$  is in  the form
 $$x_{i1}^2x_{i2}^2=a_1^2a_2^2$$
 then we replace it  by the system
 $$\{x_{i1}=(
a_1z_1)^{z_2}, x_{i2}=(z_1^{-1} a_2)^{z_2}, [z_1, a_1 a_2]=1,
[z_2,a_1^2a_2^2]=1\}.$$

The normalization  construction effectively provides an NTQ-system
 $U^* = 1$ such that each fundamental  can be obtained
 from a solution of $U^*=1$.
We refer to  this system as to the normalized system of $U$
corresponding to the fundamental sequence. Similarly, the
coordinate group of the normalized system is called the {\it
normalized} coordinate group of $U = 1$.

\begin{lm}
\label{le:7.5embed} Let $U(X) = 1$ be an NTQ-system, and $U^* = 1
$ be  the normalized system corresponding to the fundamental
sequence $\bar V(U)$. Then the following holds:

1) The coordinate group $F_{R(U)}$ canonically  embeds into
$F_{R(U^*)}$;

2) The system $U^* = 1$ is an NTQ-system of the type

\medskip
$S_1(X_1, X_2, \ldots, X_n,A) = 1$

\medskip
$\ \ \ \ \ S_2(X_2, \ldots, X_n,A) = 1$

\medskip
$\ \ \ \ \ \ \ \ \ \  \ldots$

\medskip
$\ \ \ \ \ \ \ \ \ \ \ \ \ \ \ \ S_n(X_n,A) = 1$

\medskip\noindent
in which every $S_i = 1$ is either a regular quadratic equation or
an empty equation or a system of the type $$ U_{com}(X,b) =
\{[x_i,x_j]=1, [x_i,b]=1,  \mid \ i,j =1,\ldots ,k\}$$ where $b
\in G_{i+1}$.

3) Every solution $X_0$ of $U(X)=1$ that belongs to the
fundamental sequence $\bar V(U)$ can be obtained from a solution
of the system $U^*=1$.
\end{lm}
{\em Proof.}   Statement 1) follows from the normal forms of
elements in free constructions or from the fact that applying
standard automorphisms $\phi_L$ to  a non-commuting solution  (in
particular, to a basic one) one obtains a discriminating set of
solutions (see Section 7.2). Statements 2) and 3) are obvious from
the normalization construction.$\Box$

\begin{df} A family of solutions $\Psi$ of a regular NTQ-system $U(X,A)=1$ is called {\em generic}
if for any equation $V(X,Y,A)=1$ the following is true: if for any
solution from $\Psi$ there exists a solution of
$V(X^{\psi},Y,A)=1$, then $V=1$ admits a complete $U$-lift.

A family of solutions $\Theta$ of a regular  quadratic equation
$S(X)=1$ over a group $G$  is called {\em generic} if for any
equation $V(X,Y,A)=1$ with coefficients in $G$ the following is
true: if for any solution $\theta\in\Theta$ there exists a
solution of $V(X^{\theta},Y,A)=1$ in $G$, then $V=1$ admits a
complete $S$-lift.

A family of solutions $\Psi$ of an NTQ-system $U(X,A)=1$ is called
{\em generic} if $\Psi=\Psi _1\ldots \Psi _n,$ where $\Psi _i$ is
a generic family of solutions of $S_i=1$ over $G_{i+1}$ if $S_i=1$
is a regular quadratic system, and $\Psi _i$ is a discriminating
family for $S_i=1$ if it is  a system of the type $U_{com}$.

\end{df}

The second construction is a {\it correcting extension of
centralizers} of a normalized NTQ-system $U(X) = 1$ relative to an
equation $W(X,Y,A) = 1$, where $Y$ is a tuple of new variables.
Let $U(X) = 1$ be an NTQ-system in the normalized form:

\medskip
$S_1(X_1, X_2, \ldots, X_n,A) = 1$

\medskip
$\ \ \ \ \ S_2(X_2, \ldots, X_n,A) = 1$

\medskip
$\ \ \ \ \ \ \ \ \ \ \dots$

\medskip
$\ \ \ \ \ \ \ \ \ \ \ \ \ \ \ \ S_n(X_n,A) = 1$

\medskip\noindent
So every $S_i = 1$ is either a regular quadratic equation or an
empty equation or a system of the type
$$ U_{com}(X,b) =  \{[x_i,x_j]=1, [x_i,b]=1,  \mid \ i,j =1,\ldots ,k\}$$ where
$b \in G_{i+1}$.  Again, starting from the bottom we find the
first equation $S_i(X_i) = 1$ which is in the form $U_{com}(X) =
1$ and replace it with a new centralizer extending  system ${\bar
U}_{com}(X) = 1$ as follows.

We construct $T_{sol}$ for the system $W(X,Y) = 1 \wedge U(X)=1$
with parameters $X_i,\ldots ,X_n$. We obtain generalized equations
corresponding to final vertices. Each of them consists of a
periodic structure on $X_i$ and generalized equation on
$X_{i+1}\ldots X_n$. We can suppose that for the periodic
structure the set of cycles $C^{(2)}$ is empty. Some of the
generalized equations have a solution over the extension of the
group $G_i$. This extension is  given by the relations $\bar
U_{com}(X_i)=1, S_{i+1}(X_{i+1},\ldots , X_n)=1,\ldots ,
S_n(X_n)=1$, so that there is an embedding $\phi_k : A(X,b)
\rightarrow A(X,b)$.   The others provide a proper (abelian)
equation $E_j(X_i) = 1$  on $X_i$.  The argument above shows that
replacing each  centralizer extending system $S_i(X_i) = 1$ which
is in the form $U_{com}(X_i) = 1$ by a new  system of the type
${\bar U}_{com}(X _i) = 1$ we eventually rewrite the system $U(X)
= 1$ into  finitely many new ones ${\bar U}_1(X) = 1, \ldots,
{\bar U}_m(X) = 1 $. We denote this set of NTQ-systems by
${\mathcal C}_W (U)$. For every NTQ-system ${\bar U}_m(X) = 1  \in
{\mathcal C}_W(U)$ the embeddings $\phi_k$ described above give
rise to embeddings ${\bar \phi}:F_{R(U)} \rightarrow F_{R({\bar
U})}$. Finally, combining normalization and correcting extension
of centralizers (relative to $W = 1$)  starting with an NTQ-system
$U = 1$ and a fundamental sequence of its solutions $\bar V(U)$ we
can obtain a finite set
 $${\mathcal NC}_W(U) =   {\mathcal C}_W(U^*) $$
 which comes equipped with a  finite set of embeddings
 $\theta_{i}: F_{R(U)}  \rightarrow
 F_{R({\bar U}_i)}$ for each ${\bar U}_i \in {\mathcal NC}_W(U)$.
These embeddings are called {\em correcting normalizing
embeddings}. The construction implies the following result.

\begin{theorem}\label{tqe}(Parametrization theorem)
Let $U(X,A) =1$ be an NTQ-system with a fundamental sequence of
solutions $\bar V(U)$. Suppose a formula
 $$\Phi=\forall X (U(X) = 1 \rightarrow \exists Y (W(X, Y,A)=1  \wedge W_1(X, Y,A)\not =1)$$
 is true in $F$. Then for every ${\bar U}_i \in {\mathcal NC}_W(U)$
 the formula
 $$\exists Y (W(X^{\theta_{i}}, Y,A)=1 \wedge W_1(X^{\theta _{i}}, Y,A)\not =1)$$
is true in the group $G_{R({\bar U}_{i})}$ for every correcting
normalizing embedding $\theta_i:  F_{R(U)}  \rightarrow F_{R({\bar
U}_i)}$. This formula can be effectively verified and solution $Y$
can be effectively found.

Futhermore, for every fundamental solution $\phi:
F_{R(U)}\rightarrow F$ there exists a fundamental solution $\psi$
of one of the systems $\bar U_i=1,$  where ${\bar U}_i \in
{\mathcal NC}_W(U)$ such that $\phi =\theta _i\psi .$
\end{theorem}
 As a corollary of this theorem and results from Section 5 we
obtain the following theorems.

\begin{theorem}  Let $U(X,A) =1$
be an NTQ-system and $\bar V(U)$  a fundamental set of solutions
of $U = 1$ in $F = F(A)$. If a formula
 $$\Phi=\forall X (U(X) = 1 \rightarrow \exists Y (W(X, Y,A)=1  \wedge W_1(X, Y,A)\not =1)$$
 is true in $F$ then one can effectively find finitely many NTQ
 systems $U_1 = 1, \ldots, U_k = 1$ and embeddings $\theta_i:  F_{R(U)}  \rightarrow F_{R(U_i)}$
 such that the formula
 $$\exists Y (W(X^{\theta_{i}}, Y,A)=1 \wedge W_1(X^{\theta _{i}}, Y,A)\not =1)$$
is true in each group $F_{R({ U}_{i})}$.  Furthermore,  for every
solution $\phi: F_{R(U)}  \rightarrow F$ of $U = 1$ from $\bar
V(U)$  there exists $i \in \{1, \ldots, k\}$   and a fundamental
solution $\psi: F_{R({ U}_{i})} \rightarrow F$ such that $\phi =
\theta_i \psi$. \end{theorem}

\begin{theorem} Let $S(X)=1$ be an arbitrary  system of equations over $F$.
If a formula
$$ \Phi = \forall X \exists Y (S(X) = 1 \ \ \rightarrow \ \ (W(X, Y,A)=1  \wedge W_1(X, Y,A)\not =1))$$
is  true in $F$  then one can effectively find finitely many NTQ
 systems $U_1 = 1, \ldots, U_k = 1$ and $F$-homomorphisms $\theta_i:  F_{R(S)}  \rightarrow F_{R(U_i)}$
 such that the formula
 $$\exists Y (W(X^{\theta_{i}}, Y,A)=1 \wedge W_1(X^{\theta _{i}}, Y,A)\not =1)$$
is true in each group $F_{R({ U}_{i})}$.  Furthermore,  for every
solution $\phi: F_{R(S)}  \rightarrow F$ of $S = 1$    there
exists $i \in \{1, \ldots, k\}$   and a fundamental solution
$\psi: F_{R({ U}_{i})} \rightarrow F$ such that $\phi = \theta_i
\psi$. \end{theorem}

\end{document}